\input amstex
\input xy
\xyoption{all}

\documentstyle{amsppt}
\nologo
\NoBlackBoxes

\font\ninemi=cmmi9

\font\twelvemi=cmmi12

\newcount\refCount
\def\newref#1 {\advance\refCount by 1
\expandafter\edef\csname#1\endcsname{\the\refCount}}

\newref BER      
\newref BOR      
\newref BOL      
\newref BGR      
\newref TVS      
\newref BRTHREE  
\newref COM      
\newref EMERJ    
\newref EMERI    
\newref FEDI     
\newref FETH     
\newref ING      
\newref LAZ      
\newref MEH      
\newref MOR      
\newref PRA      
\newref ROB      
\newref SCHNA    
\newref SCHTBD   
\newref SCHTAN   
\newref SCHTUF   
\newref SCHTIW   
\newref SCHTPF   
\newref SCHTNEW  
\newref SERLG     

\def\plim#1{\underset{\longleftarrow\atop #1}\to{\lim} \, }
\def\ilim#1{\underset{\longrightarrow\atop #1}\to{\lim} \, }

\def\id{\roman{id}}
\def\iso{\buildrel \sim \over \longrightarrow}
\def\leftiso{\buildrel \sim \over \longleftarrow}

\def\cotimes{\mathop{\hat{\otimes}}\nolimits}

\def\Hom{\mathop{\roman{Hom}}\nolimits}
\def\Lin{\mathop{\Cal{L}}\nolimits}
\def\End{\mathop{\roman{End}}\nolimits}

\def\Res{\mathop{\roman{Res}}\nolimits}
\def\MaxSp{\mathop{\roman{Sp}}\nolimits}

\def\SL{\mathop{\roman{SL}}\nolimits}
\def\GL{\mathop{\roman{GL}}\nolimits}

\def\Fun{\Cal F}
\def\Aff{\Cal C^{\alg}}
\def\An{\Cal C^{\an}}
\def\La{\Cal C^{\la}}
\def\Con{\Cal C}

\def\DCon{\Cal D}
\def\DAn{\Cal D^{\an}}
\def\DLa{\Cal D^{\la}}

\def\core#1#2{\roman{co}_{#1}(#2)}

\def\lie#1{\goth #1}

\def\open#1{#1^{\circ}} 

\def\Z{{\Bbb Z}}
\def\R{{\Bbb R}}
\def\Q{{\Bbb Q}}

\def\F{{\Bbb F}}

\def\G{{\Bbb G}}
\def\N{{\Bbb N}}
\def\A{{\Bbb A}}
\def\B{{\Bbb B}}
\def\H{{\Bbb H}}
\def\J{{\Bbb J}}

\def\X{{\Bbb X}}
\def\Y{{\Bbb Y}}

\def\sup{\roman{sup}}
\def\ord{\roman{ord}}

\def\la{\roman{la}}
\def\lalg{\roman{lalg}}

\def\an{\roman{an}}
\def\lf{\roman{lf}}
\def\sm{\roman{sm}}
\def\con{\roman{con}}
\def\op{\roman{op}}
\def\ev{\roman{ev}}

\def\alg{\roman{alg}}

\def\Ad{\roman{Ad}}

\def\Gr{\roman{Gr}}
\def\orbit{o}
\def\strict{\roman{st.sm}}

\def\d{\, d\,}

\def\section #1{\vskip 2mm {(\bf #1)}}

\topmatter

\title
\textfont1 = \twelvemi
locally analytic vectors in 
representations of locally
$p$-adic analytic groups
\endtitle

\rightheadtext{\textfont1 = \ninemi
locally analytic vectors in representations of $p$-adic groups}

\author
Matthew Emerton
\endauthor

\affil
Northwestern University 
\endaffil

\thanks
The author would
like to acknowledge the support 
of the National Science Foundation (award number DMS-0070711)
\endthanks

\address
Northwestern University \newline\indent
Department of Mathematics\newline\indent
2033 Sheridan Rd.\newline\indent
Evanston, IL 60208-2730, USA
\endaddress

\email
emerton\@math.northwestern.edu
\endemail

Draft: May 4, 2004

\medskip

\toc
\head 0. Terminology, notation and conventions \page{7}\endhead
\head 1. Non-archimedean functional analysis \page{10}\endhead
\head 2. Non-archimedean function theory \page{27} \endhead
\head 3. Continuous, analytic, and locally analytic vectors \page {43} \endhead
\head 4. Smooth, locally finite, and locally algebraic vectors \page{71}\endhead
\head 5. Rings of distributions \page{81} \endhead
\head 6. Admissible locally analytic representations \page{100} \endhead
\head 7. Representations of certain product groups \page{124}\endhead
\head {} References \page{135} \endhead
\endtoc
\endtopmatter

\document
Recent years have seen the emergence of a new branch of
representation theory:
the theory of representations of locally $p$-adic analytic
groups on locally convex $p$-adic topological vector spaces
(or ``locally analytic representation theory'', for short).
Examples of such representations are provided by finite dimensional
algebraic representations of $p$-adic reductive groups, and also
by smooth representations of such groups (on $p$-adic vector spaces).
One might call these the ``classical'' examples of such representations.
One of the main interests of the theory (from the point of view of
number theory) is that it provides a setting
in which one can study $p$-adic completions of the classical representations
\cite{\BRTHREE},
or construct ``$p$-adic interpolations'' of them
(for example, by defining locally analytic analogues of the
principal series, as in \cite{\SCHTAN}, or by 
constructing representations via
the cohomology of arithmetic quotients of symmetric spaces, as in
\cite{\EMERI}).

Locally analytic representation theory also plays an important
role in the analysis
of $p$-adic symmetric spaces; indeed, this analysis
provided the original
motivation for its development.
The first ``non-classical'' examples in the theory were found by Morita,
in his analysis of the $p$-adic upper half-plane (the $p$-adic symmetric
space attached to $\GL_2(\Q_p)$) \cite{\MOR},
and further examples were found by Schneider and Teitelbaum
in their analytic investigations of the $p$-adic 
symmetric spaces of $\GL_n(\Q_p)$ (for arbitrary $n$) \cite{\SCHTBD}.
Motivated in part by the desire to understand these examples,
Schneider and Teitelbaum have recently initiated a systematic study
of locally analytic representation theory
\cite{\SCHTBD, \SCHTAN, \SCHTUF, \SCHTIW, \SCHTNEW}.
In particular, they have introduced the important notions
of admissible and strongly admissible
locally analytic representations of a 
locally $p$-adic analytic group
(as well as the related notion of admissible continuous 
representations of such a group).

The goal of this paper is to provide the foundations for the 
locally analytic representation theory that is required in
the papers \cite{\EMERJ, \EMERI}.   In the course of
writing those papers we have
found it useful to
adopt a particular point of view on locally analytic representation
theory: namely,  we regard a locally analytic representation as 
being the inductive limit of its 
subspaces of analytic vectors (of
various ``radii of analyticity''), 
and we use the analysis of these subspaces as one of the basic tools in our
study of such representations. 
Thus in this paper we present a development of locally
analytic representation theory built around this point of view.
Some of the material that we present is entirely new
(for example, the notion of essentially admissible
representation, which plays a key role in \cite{\EMERJ},
and the results of section~7, which are used in \cite{\EMERI});
other parts of it can be found
(perhaps up to minor variations) in the papers of Schneider and
Teitelbaum cited above, or in the thesis of
Feaux de Lacroix \cite{\FETH}.  
We have made a deliberate effort to keep the paper reasonably self-contained,
and we hope that this will be of some benefit to the reader.

We will now give a more precise description of the view-point 
on locally analytic representation theory that this paper adopts,
and that we summarised above.

Let $L$ be a finite extension of $\Q_p,$ and
let $K$ be an extension of $L$, complete with respect
to a discrete valuation extending the discrete valuation
of $L$.  We let $G$ denote a locally $L$-analytic
group, and consider representations of $G$ by continuous
operators on Hausdorff locally convex topological $K$-vector spaces.  
Our first goal is to define, for any
such representation $V$,
the topological $K$-vector
space $V_{\la}$ of locally analytic vectors
in $V$.  As a vector space, $V_{\la}$ is defined to be the
subspace of $V$ consisting of those vectors $v$ for which
the orbit map $g \mapsto gv$ is a locally analytic $V$-valued
function on $G$.  The non-trivial point in the definition
is to equip $V_{\la}$ with an appropriate topology.
In the references \cite{\SCHTBD,\SCHTNEW}, the authors endow
this space (which they denote by $V_{\an}$, rather than $V_{\la}$) with
the topology that it inherits as a closed subspace of the space
$\La(G,V)$ of locally analytic $V$-valued functions on $G$
(each vector $v \in V_{\la}$ being identified with the corresponding
orbit map $\orbit_v: G \rightarrow V$).
We have found it advantageous to endow $V_{\la}$ with a finer topology,
with respect to which it is exhibited
as a locally convex inductive limit of Banach spaces.
(In some important situations --
for example, when $V$ is of compact type,
or is a Banach space equipped with an admissible continuous
$G$-action -- we prove that
the topology that we consider coincides with that
considered by Schneider and Teitelbaum.)

Suppose first that $\G$ is an affinoid
rigid analytic group defined over $L$, and that $G$ is the
group of $L$-valued points of $\G$. 
If $W$ is a Banach
space equipped with a representation of $G$, then we say that this
representation is
$\G$-analytic if for each $w \in W$ the orbit map $\orbit_w: G \rightarrow W$
given by $w$ extends to a $W$-valued rigid analytic function on $\G$.
For any $G$-representation $V$, 
we define $V_{\G-\an}$ to be the locally convex inductive limit
over the inductive system of $G$-equivariant maps $W \rightarrow V$,
where $W$ is a Banach space equipped with a $\G$-analytic action of $G$.

We now consider the case of an arbitrary locally $L$-analytic group $G$.
Recall that a chart $(\phi,H,\H)$
of $G$ consists of an open subset $H$ of $G$, an affinoid space
$\H$ isomorphic to a closed ball, and a locally
analytic isomorphism $\phi: H \iso \H(L)$.  If $H$ is furthermore
a subgroup of $G$, then the fact that $\H(L)$ is Zariski dense in $\H$
implies that there is at most one rigid analytic group structure
on $\H$ inducing the given group structure on $H$.  If such a group
structure exists, we refer to the chart $(\phi, H,\H)$ as an analytic
open subgroup of $G$.  We will typically suppress reference to the
isomorphism $\phi$, and so will speak of an analytic open subgroup $H$
of $G$, letting $\H$ denote the corresponding rigid analytic group,
and identifying $H$ with the group of points $\H(L)$.

For any $G$-representation on a Hausdorff convex $K$-vector space $V$,
and any analytic open subgroup $H$ of $\G$,
we can define as above the space $V_{\H-\an}$ of $\H$-analytic vectors in $V$.
(If we ignore questions of topology, then $V_{\H-\an}$ consists of those locally
analytic vectors with ``radius of analyticity'' bounded below by $H$.)
We define $V_{\la}$ to the be the locally convex inductive limit over all
locally analytic open subgroups $H$ of $G$ of the spaces $V_{\H-\an}$.

The representation $V$ of $G$ is said to be locally analytic
if $V$ is barrelled, and if the natural map $V_{\la} \rightarrow V$
is a bijection.  If $V$ is an $LF$-space (we recall the meaning of
this, and some related, functional analytic terminology in subsection~1.1
below), then we can show that
if this map is a bijection, it is in fact a topological isomorphism.  
Thus given a locally analytic representation of $G$ on
an $LF$-space $V$, we may write $V \iso \ilim{n} V_{\H_n-\an},$
where $H_n$ runs over a cofinal sequence of analytic open subgroups
of $G$.

The category of admissible locally analytic $G$-representations,
introduced in \cite{\SCHTNEW}, admits a useful description
from this point of view.  We show that a locally analytic
$G$-representation on a Hausdorff convex $K$-vector space
$V$ is admissible if and only if $V$ is an
$LB$-space, such that for each analytic open subgroup $H$ of $G$,
the space $V_{\H-\an}$ admits a closed $H$-equivariant embedding into 
a finite direct sum of copies of the space $\An(\H,K)$ of rigid analytic
functions on~$\H$.

Recall that in \cite{\SCHTNEW}, a locally analytic
$G$-representation $V$ is defined
to be admissible if and only if $V$ is of compact type, and if
the dual space $V'$ is a coadmissible module under the action
of the ring $\DAn(H,K)$ of locally analytic distributions on $H$,
for some (or equivalently, every) compact open subgroup $H$ of $G$.
For this definition to make sense (that is, for the notion of
a coadmissible $\DAn(H,K)$-module to be defined), the authors
must prove that the ring $\DAn(H,K)$ is a Fr\'echet-Stein algebra,
in the sense of \cite{\SCHTNEW, def., p.~8}.  This result
\cite{\SCHTNEW, thm.~5.1} is the main theorem of that reference.

In order to establish our characterisation of admissible locally
analytic representations, we are led to give an alternative
proof of this theorem, and an alternative description of
the Fr\'echet-Stein structure on $\DLa(H,K)$,
which is more in keeping with our point of
view on locally analytic representations.  While the proof of
\cite{\SCHTNEW} relies on the methods of \cite{\LAZ}, we rely
instead on the methods used in \cite{\BER} to prove
the coherence of the sheaf of rings $\Cal D^{\dagger}$.

We also introduce the category of essentially admissible
locally analytic $G$-representations.  To define this category,
we must assume that the centre $Z$ of $G$ is topologically
finitely generated.
(This is a rather mild condition, which is satisfied, for example, if
$G$ is the group of $L$-valued points of a linear algebraic group over $L$.)
Supposing that this is so,
we let $\hat{Z}$ denote the rigid analytic space parameterising
the locally analytic characters of $Z$,
and let $\An(\hat{Z},K)$ denote the Fr\'echet-Stein algebra of $K$-valued
rigid analytic functions on $\hat{Z}$. 

Let $V$ be a convex $K$-vector space of compact type equipped with
a locally analytic $G$-representation, and suppose that 
$V$ may be written as a union
$V = \ilim{n} V_n,$ where each $V_n$ is a $Z$-invariant $BH$-subspace
of $V$.
The $\DLa(H,K)$-action
on the dual space $V'$ then extends naturally to an
action of the completed tensor product algebra
$\An(\hat{Z},K) \cotimes_{K} \DLa(H,K)$.
Our proof of the fact that $\DLa(H,K)$ is Fr\'echet-Stein generalises
to show that this completed tensor product is also Fr\'echet-Stein.
We say that $V$ is an essentially admissible locally analytic
representation of $G$ if, furthermore, $V'$ is a coadmissible
module with respect to this Fr\'echet-Stein algebra, for some
(or equivalently, any) compact open subgroup $H$ of $G$.

It is easy to show, using the characterisation of admissible locally 
analytic representations described above,
that any such locally analytic representation of $G$ is 
essentially admissible.  Conversely, if
$V$ is any essentially admissible locally analytic
representation of $G$, and if $\chi$ is a $K$-valued point of $\hat{Z}$,
then the closed subspace $V^{\chi}$ of $V$ on which
$Z$ acts through $\chi$ is an admissible locally analytic representation
of $G$.  
The general
theory of Fr\'echet-Stein algebras \cite{\SCHTNEW, \S 3} implies
that the category
of essentially admissible locally analytic $G$-representations
is abelian. 

Since the category of coadmissible $\An(\hat{Z},K)$-modules
is equivalent to the category of coherent rigid analytic sheaves
on the rigid analytic space $\hat{Z}$, one may think of
a coadmissible $\An(\hat{Z},K)\cotimes_K \DLa(H,K)$-module
as being a ``coherent sheaf'' of coadmissible $\DLa(H,K)$-modules
on $\hat{Z}$.  Thus, roughly speaking,
one may think of an essentially admissible
locally analytic $G$-representation $V$ as being a family
of admissible locally analytic $G$-representations parameterised
by the space $\hat{Z},$ whose fibre over a point $\chi \in \hat{Z}(K)$
is equal to $V^{\chi}$.

The category of essentially admissible locally
analytic representations provides the setting for the
Jacquet module construction for locally analytic representations
that is the subject of the paper \cite{\EMERJ}.
These functors are in turn applied in {\cite \EMERI} to construct
``eigenvarieties'' (generalising the eigencurve of
\cite{\COM}) that $p$-adically interpolate
systems of eigenvalues attached to automorphic
Hecke eigenforms on reductive groups over number fields.

Let us point out that functional analysis currently
provides the most important technical tool in the theory of locally
analytic representations. 
Indeed, since continuous $p$-adic valued functions on a $p$-adic group 
are typically not locally integrable for Haar measure (unless they
are locally constant), there is so far no real analogue in this 
theory of the harmonic analysis which plays such
an important role in the theory of smooth representations
(although one can see some shades of harmonic analysis in the theory:
the irreducibility result of \cite{\SCHTAN, thm.~6.1}
depends for its proof on Fourier analysis in the non-compact
picture, and Hecke operators make an appearance in the construction
of the Jacquet module functor of \cite{\EMERJ}).
Thus one relies on softer functional analytic methods
to make progress.   This paper is no exception;
it relies almost entirely on such methods.

A more detailed summary of the paper now follows,
preceding section by section.

In section~1 we develop the non-archimedean functional analysis
that we will require in the rest of the paper.
Subsection~1.1 is devoted to recalling various pieces of terminology that
we will need, and to proving some results for which we could not
find references in the literature.  None of the results are difficult,
and most or all are presumably well-known to experts.
In subsection~1.2 we recall the theory of Fr\'echet-Stein
algebras developed in \cite{\SCHTAN, \S 3}, and prove some
additional results that we will require in our applications of
this theory. 

In section~2 we recall the basics of non-archimedean function theory.
Subsection~2.1 recalls the basic definitions regarding spaces
of continuous,
rigid analytic, and locally analytic functions with values in locally
convex $K$-vector spaces, and establishes some basic properties of
these spaces that we will require.	
Subsection~2.2 introduces the corresponding spaces of distributions.
In subsection~2.3 we recall the definition and basic properties
of the restriction of scalars functor, in both the rigid analytic 
and the locally analytic setting.  

In section~3 we present our construction of 
the space of locally analytic vectors attached to 
a representation of a non-archimedean locally $L$-analytic group $G$.

After some preliminaries in subsections~3.1 and~3.2, 
in subsection~3.3 we suppose that $G$ is the group of points
of an affinoid
rigid analytic group defined over $L$, and define the space
of analytic vectors.
In subsection~3.4 we extend this construction
to certain non-affinoid rigid analytic groups.  

In subsection~3.5, we return to the situation in which $G$ is
a locally $L$-analytic group, and construct the space of locally
analytic vectors attached to any $G$-representation.
In subsection~3.6 we recall the notion of locally analytic
representation, and also introduce the related notion of
analytic representation, and establish some basic properties
of such representations.

Section~4 begins by recalling, in subsection~4.1, the notion
of smooth and locally finite vectors in a $G$-representation.
The main point of this subsection is to prove some
simple facts about representations in which every
vector is smooth or locally finite.

In subsection~4.2 we assume that
$G$ is the group of $L$-valued points of a reductive linear
algebraic group $\G$ over $L$.  For any finite dimensional algebraic
representation $W$ of $\G$ over $K$, and for any $G$-representation
$V$, we define the space $V_{W-\alg}$ of
locally $W$-algebraic vectors in $V$, and study some of its
basic properties.
As in \cite{\PRA},
the representation $V$ is said to be locally algebraic
if every vector of $V$ is $W$-locally algebraic vector for some
representation $W$ of $\G$.   
We prove that any irreducible locally
algebraic representation is isomorphic to the tensor product of a
smooth representation of $G$ and a finite dimensional
algebraic representation of $\G$  (first proved
in \cite{\PRA}).

One approach to analysing representations $V$ of $G$, the importance
of which has been emphasised by Schneider and Teitelbaum,
is to pass to the dual space $V'$, and to regard $V'$ as a module
over an appropriate ring of distributions on $G$. 
The goal of section~5 is to recall this approach, and 
to relate it to the view-point of section~3.

In subsection~5.1 we prove some simple forms of Frobenius
reciprocity, and apply these to obtain a uniform development
of the dual point of view for continuous, analytic, and
locally analytic representations.  
In subsection~5.2 we recall the description of algebras of
analytic distributions
via appropriate completions of universal enveloping algebras.
In subsection~5.3 we use this description, together with the
methods of \cite{\BER, \S 3}, to present a new construction of
the Fr\'echet-Stein structure on the ring $\DLa(H,K)$ of locally
analytic distributions on any compact open subgroup $H$ of $G$.
In fact, we prove a slightly more general result, which implies
not only that $\DLa(H,K)$ is a Fr\'echet-Stein algebra, but also
that the completed tensor product $A \cotimes_K \DLa(H,K)$ is
Fr\'echet-Stein, for a fairly general class
of Fr\'echet-Stein $K$-algebras $A$, namely, those $K$-Fr\'echet
algebras $A$ that admit an integral Fr\'echet-Stein structure
(a notion introduced in definition~5.3.21).

In section~6 we study the various admissibility conditions
that have arisen so far in locally analytic representation theory.

In subsection~6.1 we present our alternative definition
of the category of admissible locally analytic $G$-representations,
and prove that it is equivalent to the definition presented
in \cite{\SCHTNEW}. 

In subsection~6.2, we recall the notion
of strongly admissible locally analytic $G$-representation
(introduced in \cite{\SCHTAN}) and also the notion of admissible
continuous $G$-representation (introduced in \cite{\SCHTIW}).
We prove that any strongly admissible locally analytic $G$-representation
is an admissible $G$-representation, and also that if $V$
is an admissible continuous $G$-representation, then $V_{\la}$
is a strongly admissible locally analytic $G$-representation.
(Both these results are established independently in~\cite{\SCHTNEW}.)

In subsection~6.3 we study admissible smooth and admissible
locally algebraic representations.  The main result
concerning smooth representations is that a locally
analytic smooth representation
is admissible as a locally analytic representation
if and only if it is admissible as a smooth representation,
and that it is then necessarily equipped with its finest
convex topology.
(This has also been proved independently in \cite{\SCHTNEW}.)

In subsection~6.4 we suppose that the centre $Z$ of $G$
is topologically finitely generated.  We are then able to
define the category of essentially
admissible locally analytic $G$-representations.  
As already indicated, these are locally analytic $G$-representations
on a space $V$ of compact type whose dual space $V'$ is a coadmissible
$\An(\hat{Z},K)\cotimes_K \DLa(H,K)$-module, for one (or equivalently
any) compact open subgroup $H$ of $G$.  (The $K$-algebra
$\An(\hat{Z},K)$ admits an integral Fr\'echet-Stein structure,
and so the results of subsection~5.3
imply that $\An(\hat{Z},K)\cotimes_K \DLa(H,K)$ is a Fr\'echet-Stein
algebra.)

In subsection~6.5 we introduce some simple notions related to
invariant open lattices in $G$-representations.  In particular,
we say that a $G$-invariant open lattice is admissible if 
the (necessarily smooth) action of $G$ induced on $L / \pi$ (where
$\pi$ is a uniformiser of the ring of integers of $K$) is admissible.
We characterise strongly admissible locally analytic representations
as being those essentially admissible locally analytic representations
that contain a separated open $H$-invariant admissible lattice, for
some (or equivalently, every) compact open subgroup $H$ of $G$.

In section~7, we discuss the representations of groups of the form
$G \times \Gamma,$ where $G$ is a locally $L$-analytic group and
$\Gamma$ is a Hausdorff locally compact topological group that admits a
countable basis of neighbourhoods of the identity consisting of
open subgroups.  (The motivating example is a group of the form
$\G(\Bbb A_f),$ where $\G$ is a reductive group defined over some number
field $F$, $\Bbb A_f$ denotes the ring of finite ad\`eles of $F$, and $L$ is
the completion of $F$ at some finite prime.)  In subsection~7.1
we introduce the notion of a strictly smooth representation of $\Gamma$
on a Hausdorff convex $K$-vector space $V$.  This is a useful strengthening of
the notion of a smooth representation in the context of topological
group actions, defined by requiring that the natural map
$\ilim{H} V^H \rightarrow V$ be a topological isomorphism. 
(Here the locally convex
inductive limit is indexed by the open subgroups $H$ of $\Gamma$.)
In subsection~7.2 we extend the various notions of
admissibility introduced in section~6 to the context of $G\times \Gamma$
representations.  Briefly, a topological representation of $G\times \Gamma$
on a Hausdorff convex $K$-vector space $V$ is said to be admissible locally
analytic, essentially admissible locally analytic, or admissible continuous,
if the $\Gamma$-action on $V$ is strictly smooth, and if the closed subspace
of invariants $V^H$ has the corresponding property as a $G$-representation,
for each open subgroup $H$ of $\Gamma$.  The subsection closes by considering
locally algebraic representations in this context.

{\it Acknowledgments.}  
I would like to thank David Ben-Zvi, Kevin Buzzard, Brian Conrad,
Robert Kottwitz, Peter Schneider, Jeremy Teitelbaum,
and Bertrand Toen for helpful discussions on various aspects
of this paper.  

\head 0. Terminology, notation and conventions \endhead

We now describe our terminological and notational conventions
for the paper.

{\bf (0.1)}
Throughout the paper we fix a finite extension $L$
of the field of $p$-adic numbers $\Q_p$ (for some fixed prime
$p$),
and let $\ord_L: L^{\times} \rightarrow \Q$
denote the discrete valuation on $L$,
normalised so that $\ord_L(p) = 1$.
We also fix an extension $K$ of $L$,
spherically complete with respect to a a non-archimedean valuation
$\ord_K: K^{\times} \rightarrow \R$ extending $\ord_L$.
We let $\Cal O_L$ and $\Cal O_K$ denote the rings of integers of $L$
and $K$ respectively.  We also let $|\,\text{--}\,|$ denote the
absolute value on $K$ induced by the valuation $\ord_K$.

{\bf (0.2)}  A topological $K$-vector space $V$ is called locally
convex if there is a neighbourhood basis of the origin consisting of
$\Cal O_K$-submodules.  
Let us remark that the theory of locally convex vector spaces
over the non-archimedean field $K$ can be developed in a fashion
quite similar
to that of locally convex real or complex topological vector
spaces.  In particular, the hypothesis that $K$ is spherically
complete ensures that the Hahn-Banach theorem holds for convex
$K$-vector spaces \cite{\ING}.  Thus we will on occasion cite a
result in the literature dealing with real or complex convex spaces
and then to apply the corresponding non-archimedean analogue,
leaving it to the reader to check that the archimedean proof
carries over the the non-archimedean setting.
The lecture notes of Schneider \cite{\SCHNA} present
a concise development of the foundations of non-archimedean
functional analysis.

We will often abbreviate the phrase ``locally convex
topological $K$-vector space'' to ``convex $K$-vector space''.

{\bf (0.3)}
The category of convex $K$-vector spaces and continuous $K$-linear maps
is additive, and all morphisms have kernels and cokernels.  We say
that a continuous $K$-linear map $\phi: V \rightarrow W$ in this
category is strict if its image is equal to its coimage.  
Concretely, this means that the quotient of $V$ modulo the kernel
of $\phi$ embeds as a topological subspace of $W$.  In particular,
a surjection of convex $K$-vector spaces
is strict if and only if it is an open mapping.

We will say that a sequence of continuous $K$-linear maps of convex $K$-vector
spaces is exact if it is exact as a sequence of $K$-linear maps,
and if furthermore all the maps are strict.

{\bf (0.4)} Unless otherwise stated, any finite dimensional
$K$-vector space is regarded as a convex $K$-vector space,
by equipping it with its unique Hausdorff topology.

{\bf (0.5)} If $V$ and $W$ are two $K$-vector space we let
$\Hom(V,W)$ denote the $K$-vector space of $K$-linear maps from
$V$ to $W$.  The formation of $\Hom(V,W)$ is contravariantly
functorial in $V$ and covariantly functorial in $W$.

{\bf (0.6)} An important case of the discussion of~(0.5) occurs 
when $W = K$. The space $\Hom(V,K)$ is the dual space to $V$, which
we denote by $\check{V}$.  The formation of $\check{V}$ is
contravariantly functorial in $V$.

Passing to the transpose yields a natural injection $\Hom(V,W)
\rightarrow \Hom(\check{W},\check{V})$.

{\bf (0.7)}
If $V$ and $W$ are two convex $K$-vector spaces,
then we let $\Lin(V,W)$ denote the $K$-vector subspace $\Hom(V,W)$
consisting of continuous $K$-linear maps.
As with $\Hom(V,W),$ the formation of
$\Lin(V,W)$ is contravariantly functorial in $V$ and covariantly
functorial in $W$.

The space $\Lin(V,W)$ admits various convex topologies.  
The two most common are the strong topology (that is,
the topology of uniform convergence on bounded subsets of $V$),
and the weak topology (that is, the topology of pointwise
convergence).  We let $\Lin_b(V,W)$, respectively $\Lin_s(V,W)$,
denote the space $\Lin(V,W)$ equipped with the strong,
respectively weak, topology. The formation of $\Lin_b(V,W)$
and $\Lin_s(V,W)$ are again both functorial in $V$ and $W$.

{\bf (0.8)} An important case of the discussion of~(0.7) occurs 
when $W = K$.  The space $\Lin(V,K)$ is the topological dual space
to $V$, which we will denote by $V'$.
The Hahn-Banach theorem shows that if $V$ is Hausdorff, then
the functionals in $V'$ separate the elements of $V$.
The formation of $V'$ is contravariantly functorial in $V$.
We let $V'_b$ and $V'_s$ denote $V$ equipped with the strong,
respectively weak, topology.  The formation of either of these
convex spaces is contravariantly functorial in $V$.

Passing to the transpose yields the two injective
maps $\Lin(V,W) \rightarrow \Lin(W_b',V_b')$ and 
$\Lin(V,W) \rightarrow \Lin(W_s',V_s')$.
(The injectivity follows from the fact that 
the elements of $W'$ separate the elements of $W$.)

There are natural $K$-linear injections
$V \rightarrow (V'_s)'$ and $V\rightarrow
(V_b')'$.  The first is always a bijection, while the second
is typically not.  If it is, the space $V$ is called
semi-reflexive.

There is a refinement of the notion of semi-reflexivity.
Namely, we say that the Hausdorff convex $K$-vector space $V$ is reflexive
if the map
$V \rightarrow (V_b')'$ is a topological isomorphism, when the
target is endowed with its strong topology.

{\bf (0.9)} If $V$ and $W$ are two convex $K$-vector spaces then
there are (at least) two natural locally convex topologies that can be placed
on the tensor product $V \otimes_K W.$  The inductive tensor product
topology is universal for separately continuous bilinear maps of
convex $K$-vector spaces $V \times W \rightarrow U,$ while the projective
tensor product topology is universal for jointly continuous bilinear
maps $V \times W \rightarrow U.$
We let $V \otimes_{K,\iota} W$ denote $V \otimes_K W$ equipped with
its inductive tensor product topology, and $V \otimes_{K,\pi} W$
denote $V\otimes_K W$ equipped with its projective tensor product
topology.  (This is the notation of \cite{\SCHNA, \S 17}.) 
We denote the completions of these two
convex $K$-vector spaces by $V\cotimes_{K,\iota} W$ and $V \cotimes_{K,\pi} W$
respectively.  The identity map on $V\otimes_K W$ induces
a natural continuous bijection
$V \otimes_{K,\iota} W \rightarrow V \otimes_{K,\pi} W,$ 
which extends to a continuous linear map
$V \cotimes_{K,\iota} W \rightarrow V \cotimes_{K,\pi} W.$ 

In many contexts, the bijection $V\otimes_{K,\iota} W \rightarrow
V \otimes_{K,\pi} W$ is a topological isomorphism.  This is the case,
for example, if $W$ is finite dimensional, or if $V$ and $W$ are
Fr\'echet spaces.  (See proposition~1.1.31 below
for an additional such case.)  In these situations, we write simply
$V\otimes_K W$ and $V\cotimes_K W$ to denote the topological tensor
product and its completion.

{\bf (0.10)} If $G$ is a group then we let $e$ denote the
identity element of $G$.

{\bf (0.11)} If $V$ is a $K$-vector space and $G$ is a group,
then we will refer to a left action of $G$ on $V$ by $K$-linear
 endomorphisms of $V$ simply as a $G$-action on $V$, or
as a $G$-representation on $V$.  (Thus by convention our
group actions are always left actions.)

If $V$ is furthermore a topological $K$-vector space,
and $G$ acts on $V$ through continuous endomorphisms of $V$,
we will refer to the $G$-action as a topological $G$-action.

Finally, if $V$ is a topological $K$-vector space
and $G$ is a topological group, then we have the notion
of a continuous $G$-action on $V$: this is an action for which
the action map $G \times V \rightarrow V$ is continuous.
(Such an action is thus in particular a topological $G$-action
on $V$.)

One can also consider an intermediate notion of continuity
for the action of a topological group $G$ on a topological $K$-vector
space $V$,
in which one asks merely that the action map $G\times V \rightarrow V$
be separately continuous.
We remind the reader that if $V$ is a barrelled convex $K$-vector space
and $G$ is locally compact then this {\it a priori}
weaker condition in fact implies the joint continuity of
the $G$-action.  (Indeed, since $G$ is locally compact, we see that
the action of a neighbourhood of the identity in
$G$ is pointwise bounded, and thus equicontinuous,
since $V$ is barrelled. The claim is then a consequence of
lemma~3.1.1 below.)

{\bf (0.12)} If $V$ and $W$ are two $K$-vector spaces each equipped
with a $G$-action, then the $K$-vector space
$\Hom(V,W)$ is equipped with a $G$-action,
defined by the condition that
$g(\phi(v)) = (g \phi)(g v),$ for $g\in G,$ $\phi
\in \Hom(V,W)$ and $v\in V$.  An element of $\Hom(V,W)$ is fixed
under this $G$-action precisely if it is $G$-equivariant.
We let $\Hom_G(V,W)$ denote the subspace of $\Hom(V,W)$ consisting
of $G$-fixed elements; it is thus the space of $G$-equivariant
$K$-linear maps from $V$ to $W$.

{\bf (0.13)} A particular case of the discussion of~(0.12) occurs
when $V$ is a $K$-vector space equipped with a $G$-action,
and $W = K$, equipped with the trivial $G$-action.
The resulting $G$-action on $\check{V}$ is called
the contragredient $G$-action, and is characterised by the condition
$\langle g \check{v}, g v\rangle = \langle \check{v} , v \rangle$ for any
elements $\check{v}\in \check{V}$ and $v\in V$.  In other words, it is the
action obtained on $\check{V}$ after converting the transpose action
of $G$ on $\check{V}$ (which is a right action) into a left action.

If $W$ is a second $K$-vector space equipped with a $G$-action,
then the natural map $\Hom(V,W) \rightarrow \Hom(\check{W},
\check{V})$ is $G$-equivariant, if we endow the target with
the $G$-action of~(0.12) induced by the contragredient $G$-action
on each of $\check{W}$ and $\check{V}$.

{\bf (0.14)}
If $V$ and $W$ are two convex $K$-vector spaces each equipped with
a topological $G$-action then $\Lin(V,W)$ is a $G$-invariant subspace
of $\Hom(V,W).$  We let $\Lin_G(V,W)$ denote the subspace of $G$-invariant
elements; this is then the space of continuous $G$-equivariant $K$-linear
maps from $V$ to $W$.

Since the formation of $\Lin_b(V,W)$ and $\Lin_s(V,W)$ is functorial
in $V$ and $W$, we see that the $G$-action on either of these convex
spaces is topological.

{\bf (0.15)}
A particular case of the discussion of~(0.14) occurs when $V$ is
a convex $K$-vector space equipped with a topological $G$-action, and $W = K$
equipped with the trivial $G$-action.
Thus we see that $V'$ is $G$-invariant with respect to the contragredient action
of $G$ on $\check{V}$, and
that the contragredient action of $G$ on either
$V'_b$ or $V'_s$ is again topological.

Furthermore, if $W$ is a second convex $K$-vector space equipped with
a topological $G$-action, then the natural maps
$\Lin(V,W) \rightarrow \Lin(W'_b,V'_b)$ and $\Lin(V,W)
\rightarrow \Lin(W'_s,V'_s)$ are $G$-equivariant.  (Here the $G$-action
on the target of either of these maps is defined as in~(0.14) using
the contragredient $G$-action on $W'$ and $V'$.)

{\bf (0.16)} We will frequently be considering groups or
$K$-vector spaces equipped with additional topological structures,
and consequently we will sometimes use the adjective ``abstract'' to indicate
that a group or vector space is not equipped with any additional
such structure.  Similarly, we might write that a $K$-linear
map between convex $K$-vector spaces is an isomorphism of abstract
$K$-vector spaces if it is a $K$-linear isomorphism that is not
necessarily a topological isomorphism.

{\bf (0.17)}
A crucial point in non-archimedean analysis is the distinction
between locally analytic and rigid analytic
spaces (and hence between locally analytic and rigid analytic functions).
Naturally enough, in the theory of locally analytic representations
it is the notion of locally analytic function that comes to the fore.
However, any detailed investigation of locally analytic phenomena
ultimately reduces
to the consideration of functions given by power series,
and hence to the consideration of algebras of rigid analytic functions.
We refer to \cite{\SERLG} and \cite{\BGR} for 
the foundations of the locally and rigid analytic theories respectively.

\head 1. Non-archimedean functional analysis \endhead

\section{1.1} In this subsection we recall some functional analytic
terminology, and establish some simple results that we will require.
We begin with a discussion of $FH$-spaces, 
$BH$-spaces, spaces of $LF$-type, spaces of $LB$-type, $LF$-spaces, 
$LB$-spaces, and spaces of compact type.
Our discussion of these notions owes quite a lot to that of \cite{\FETH, \S 1}

\proclaim{Definition 1.1.1} Let $V$ be a Hausdorff locally convex
topological $K$-vector space. 
We say that $V$ is an $FH$-space if it admits
a complete metric that induces a locally convex topology
on $V$ finer than than its given topology.  We refer to
the topological vector space structure on $V$
induced by such a norm as a latent Fr\'echet space structure on $V$.

If $V$ admits a latent Fr\'echet space structure that can
be defined by a norm (that is, a latent Banach space structure),
then we say that $V$ is a $BH$-space.
\endproclaim

\proclaim{Proposition 1.1.2} (i) An $FH$-space $V$
admits a unique latent Fr\'echet space structure.  We let
$\overline{V}$ denote $V$ equipped with its latent Fr\'echet
space structure.

(ii) If $f:V_1 \rightarrow V_2$ is a continuous morphism
of $FH$-spaces then $f:\overline{V}_1 \rightarrow \overline{V}_2$
is a continuous morphism of Fr\'echet spaces.
\endproclaim
\demo{Proof} Suppose that $V_1$ and $V_2$ are $FH$-spaces
and that $f:V_1\rightarrow V_2$ is a continuous map between them,
and let $\overline{V}_1$ and $\overline{V}_2$ denote {\it some}
latent Fr\'echet space structure on each of these spaces.
Then the identity map induces a continuous bijection
$$\overline{V}_1 \times\overline{V}_2 \rightarrow V_1\times V_2.$$
Since $V_2$ is Hausdorff, the graph $\Gamma_f$ of $f$ is
closed in $V_1\times V_2,$ and so also in $\overline{V}_1\times
\overline{V}_2.$  The closed graph theorem now shows that
$f:\overline{V}_1 \rightarrow \overline{V}_2$ is continuous.
This proves~(ii).  To prove~(i), take $V_1=V_2=V,$
let $\overline{V}_1$ and $\overline{V}_2$ be two latent Fr\'echet
space structures on $V$, and let $f$ be the identity map
from $V$ to itself.
\qed\enddemo

\proclaim{Proposition 1.1.3} (i) If $V$ is an $FH$-space
(respectively a $BH$-space) and $W$ is a
closed subspace of $V$ then $W$ is also an $FH$-space
(respectively a $BH$-space), and
$\overline{W}$ is a closed subspace of $\overline{V}$.

(ii) If $V_1 \rightarrow V_2$ is a surjective map of Hausdorff topological
$K$-vector spaces and $V_1$ is an $FH$-space
(respectively a $BH$-space) then $V_2$ is an $FH$-space
(respectively a $BH$-space).
\endproclaim
\demo{Proof} (i) The identity map from $V$ to itself is
a continuous map $\overline{V} \rightarrow V$.  Thus the
preimage of $W$ under this map is a closed subspace
of $\overline{V}$, and thus equips $W$ with a latent Fr\'echet
space structure, which is in fact a Banach space structure
if $\overline{V}$ is a Banach space.

(ii) Let $W$ denote the kernel of the
composite $\overline{V}_1 \rightarrow V_1 \rightarrow V_2$, which
is a continuous map.
Since $V_2$ is Hausdorff we see that $W$ is closed in $\overline{V}_1,$
and so the map $\overline{V}_1/W \rightarrow V_2$ is a continuous
bijection whose source is either a Fr\'echet space
or a Banach space, depending on our hypothesis on $V_1$.
This equips $V_2$ with either
a latent Fr\'echet space structure or 
a latent Banach space structure.
\qed\enddemo

\proclaim{Definition 1.1.4} If $V$ is a Hausdorff topological
$K$-vector space, an $FH$-subspace (respectively a $BH$-subspace)
$W$ of $V$ is a subspace $W$ of $V$ that becomes an $FH$-space
(respectively a $BH$-space) when
equipped with its subspace topology.
\endproclaim

\proclaim{Proposition 1.1.5} If $W_1$ and $W_2$ are
$FH$-subspaces (respectively $BH$-subspaces)
of the Hausdorff topological $K$-vector
space $V$ then so are $W_1\bigcap W_2$ and $W_1 + W_2$.
\endproclaim
\demo{Proof} Let $W$ denote the kernel of the morphism
$$\overline{W}_1\oplus \overline{W}_2 \buildrel (w_1,w_2) \mapsto w_1 + w_2
\over \longrightarrow V.\tag 1.1.6$$
Since this is a continuous morphism into the Hausdorff space $V$,
we see that $W$ is a closed subspace of $\overline{W}_1\oplus \overline{W}_2,$
and so is either a Fr\'echet space or a Banach space, depending
on our hypothesis on $W_1$ and $W_2$.  

The continuous
morphism $$\overline{W}_1\oplus \overline{W}_2 \buildrel (w_1,w_2) \mapsto
w_1 \over \longrightarrow V$$ induces a continuous bijection
from $W$ to $W_1\bigcap W_2,$ and thus equips $W_1\bigcap W_2$ with either
a latent Fr\'echet space structure or a latent 
Banach space structure.  The morphism~(1.1.6) induces a continuous
bijection from $(\overline{W}_1\oplus \overline{W}_2)/W$ to $W_1+W_2,$ and
thus equips $W_1+W_2$ with either a latent Fr\'echet space structure
or a latent Banach space structure.  This
proves the proposition.
\qed\enddemo

\proclaim{Proposition 1.1.7} If $f:V_1\rightarrow V_2$ is
a continuous map of Hausdorff topological $K$-vector spaces
and $W$ is an $FH$-subspace (respectively a $BH$-subspace) of $V_1$,
then $f(W)$ is an $FH$-subspace (respectively a $BH$-subspace) of $V_2$.
\endproclaim
\demo{Proof} This follows from proposition~1.1.3~(ii).
\qed\enddemo

\proclaim{Proposition 1.1.8} Let $V$ be a Hausdorff topological
$K$-vector space, and let $W$ be a finite dimensional $K$-vector
space.  If $U$ is an $FH$-subspace (respectively a $BH$-subspace)
of $V$ then $U\otimes_K W$ is
an $FH$-subspace (respectively a $BH$-subspace)
of $V\otimes_K W$. Conversely, if $U'$ is an $FH$-subspace (respectively
a $BH$-subspace) of $V\otimes_K W,$ then there is an $FH$-subspace
(respectively a $BH$-subspace) $U$
of $V$ such that $U' \subset U\otimes_K W$.
\endproclaim
\demo{Proof}
Let $U$ be an $FH$-subspace or a $BH$-subspace of $V$.
The latent Fr\'echet or Banach space structure $\overline{U}$ on $U$ induces
a latent Fr\'echet or Banach space structure $\overline{U}\otimes_K W$ on
$U\otimes_K W,$ proving the first assertion.

Suppose now that $U'$ is an arbitrary $FH$-subspace or $BH$-subspace
of $V\otimes_K W.$
What we have just proved shows that $U'\otimes_K \check{W}$ is either an
$FH$-subspace or a $BH$-subspace
of $V\otimes_K W \otimes_K \check{W}$.  Let $U$ be the image of
$U'\otimes_K \check{W}$ under the natural map $V\otimes_K W \otimes_K
\check{W} \rightarrow V$, given by contracting $W$ against $\check{W}$.
Proposition~1.1.7 shows that $U$ is either an $FH$-subspace or a
$BH$-subspace of $V,$ and it is
immediately seen that there is an inclusion $U' \subset U\otimes_K W$.
\qed\enddemo

\proclaim{Definition 1.1.9}
We say that a Hausdorff convex $K$-vector space $V$ is of $LF$-type
(respectively of $LB$-type) if we may write $V =
\bigcup_{n=1}^{\infty} V_n,$ where $V_1 \subset V_2 \subset \cdots
\subset V_n \subset \cdots $ is an increasing sequence of $FH$-subspaces
(respectively $BH$-subspaces) of $V$.
\endproclaim

\proclaim{Proposition 1.1.10} If $V$ is a Hausdorff convex $K$-vector
space of $LF$-type, say $V = \bigcup_{n=1}^{\infty} V_n,$
for some increasing sequence $V_1\subset V_2
\subset \cdots \subset V_n \subset \cdots$ 
of $FH$-subspaces of $V$,
then any $FH$-subspace of $V$ 
is contained in $V_n$ for some $n$.
\endproclaim
\demo{Proof} 
This follows from \cite{\TVS, prop.~1, p.~I.20}.
\qed\enddemo

If $V$ is a Hausdorff locally convex topological $K$-vector
space, if $A$ is a bounded $\Cal O_K$-submodule of $V$, and if $V_A$ denotes
the vector subspace of $V$ spanned by $A$, then the 
gauge of $A$ defines a norm on $V_A$, with respect to which
the map $V_A \rightarrow V$ becomes continuous.  If $A$
is semi-complete, then $V_A$ is complete with respect to
this norm \cite{\TVS, cor., p.~III.8},
and hence is a $BH$-subspace of $V$.

Conversely, if $W$ is a $BH$-subspace of $V$, then the
image of the unit ball of $\overline{W}$ is a bounded
$\Cal O_K$-submodule $A$ of $V$, and the map $\overline{W} \rightarrow V$
factors as $\overline{W} \rightarrow V_B \rightarrow V$
for any bounded $\Cal O_K$-submodule $B$ of $V$ containing $A$. (Here
$V_B$ is equipped with its norm topology.)  In particular,
if $V$ is semi-complete, then the spaces $V_A$, as $A$ ranges
over all closed bounded $\Cal O_K$-submodules of $V$,
are cofinal in directed set of all $BH$-subspaces of $V$.

\proclaim{Proposition 1.1.11}
If $V$ is a semi-complete locally convex Hausdorff $K$-vector space, and 
if is $A$ is a bounded subset of $V$, then there is a $BH$-subspace $W$
of $V$ containing $A$, such that $A$ is bounded when regarded as a subset of $\overline{W}$.
If furthermore $V$ is either Fr\'echet or of compact type, then
$W$ can be chosen so that the topologies induced on $A$ by $V$ and
$\overline{W}$ coincide. 
\endproclaim
\demo{Proof}
The first claim follows from the preceding discussion.
In the case where $V$ is a Fr\'echet space, the second claim
of the theorem follows from
\cite{\SCHNA, lem.~20.5}.  Suppose now that $V$ is of compact type,
and write $V = \ilim{n} V_n$, where the transition maps
$V_n \rightarrow V_{n+1}$ are injective and compact maps between
Hausdorff convex spaces.
We may find an integer $n$ 
and a bounded subset $A_{n-1}$ of $V_{n-1}$ such that 
$A$ is equal to the image of $A_{n-1}$ under the injection
$V_{n-1} \rightarrow V$ \cite{\SCHNA, lem.~16.9~(ii)}.  
Since the injection $V_{n-1} \rightarrow V_n$ is compact,
the image $A_n$ of $A_{n-1}$ in $V_n$ is c-compact
(being closed, since it is the preimage of the closed subset $A$
of $V$).  The map $V_n \rightarrow V$ is thus a continuous
bijection from the c-compact $\Cal O_K$-module $A_n$ to the
$\Cal O_K$-module $A$, and so is a topological
isomorphism.  Thus if we take $\overline{W} = V_n,$ we prove
the second claim in the case where $V$ is of compact type.
\qed\enddemo

If $V$ is a Hausdorff convex $K$-vector space, then there is
a continuous bijection $\ilim{W} \overline{W} \rightarrow V,$
where $W$ runs over all $BH$-subspaces of $V$.  
The spaces for which this map is a topological isomorphism
are said to be ultrabornological \cite{\TVS, ex.~20, p.~III.46}.
In particular, semi-complete bornological Hausdorff convex
spaces are ultrabornological
\cite{\TVS, cor., p.~III.12}.
Since the property of being either barrelled or bornological is preserved
under the passage to locally convex inductive limits, we see that
any ultrabornological
Hausdorff convex $K$ vector space is necessarily
both barrelled and bornological.

Suppose that $V$ is a convex $K$-vector space described as a
locally convex inductive limit $V = \ilim{i\in I} V_i,$
and that $U$ is a closed subspace of $V$.  If $U_i$ denotes
the preimage of $U$ in $V_i$, for each $i\in I$, then the
natural map
$\ilim{i\in I} U_i \rightarrow U$ is a continuous
bijection, but in general not a topological isomorphism.

\proclaim{Proposition 1.1.12} Suppose that $V$ is a
Hausdorff ultrabornological convex $K$-vector space, and write
$V = \ilim{W}\overline{W},$ where $W$ ranges over all
$BH$-subspaces of $V$.  If $U$ is a closed subspace of $V$
that is also ultrabornological, 
then the natural map $\ilim{W} \overline{U\cap W} \rightarrow U$
is a topological isomorphism.
\endproclaim
\demo{Proof}
Since $U$ is ultrabornological, by assumption, it is the
direct limit of its $BH$-subspaces.  Each of these subspaces
is a $BH$-subspace of $V$, and so $\overline{U\cap W}$
ranges over all $BH$-subspaces of $U$, as $W$ ranges
over all $BH$-subspaces of $V$.
\qed\enddemo

The following definition provides a relative version of the
notion of $BH$-space.

\proclaim{Definition 1.1.13} A continuous linear map between Hausdorff
topological $K$-vector spaces $U \rightarrow V$ is said to be
$BH$ if it admits a factorisation of the form $U \rightarrow \overline{W}
\rightarrow V,$ where $\overline{W}$ is a $K$-Banach space.

Equivalently, $U\rightarrow V$ is a $BH$-map if its image
lies in a $BH$-subspace $W$ of $V$, and if it lifts to a continuous
map $U \rightarrow \overline{W}$.
\endproclaim

Obviously, any map whose source or target is a $K$-Banach space
is a $BH$-map.

\proclaim{Lemma 1.1.14} Any compact map between
locally convex Hausdorff topological $K$-vector spaces
is a $BH$-map.
\endproclaim
\demo{Proof}  This is standard.  A proof is given in the course
of the argument of \cite{\SCHNA, p.~93}.
\qed\enddemo

\proclaim{Lemma 1.1.15} If $V$ is a semi-complete 
convex Hausdorff $K$-vector space, then a map $U \rightarrow 
V$ is a $BH$-map if and only if some open neighbourhood of
the origin in $U$ has bounded image in $V$.
\endproclaim
\demo{Proof}
Since Banach spaces contain bounded open neighbourhoods
of the origin, it is clear that any $BH$-map $U\rightarrow 
V$ has the stated property.  Conversely, suppose
that some open neighbourhood of $U$ maps into a bounded
subset $A$ of $V$.  We may as well suppose that $A$ is a closed
$\Cal O_K$-submodule of $V$.  Then $U \rightarrow V$ factors
through the natural map $V_A \rightarrow V$, and so is
a $BH$-map.  (Here we are using the notation and results
discussed prior to the statement of proposition~1.1.11.)
\qed\enddemo

\proclaim{Definition 1.1.16} (i) We say that a locally
convex topological
$K$-vector space is an $LF$-space (respectively, $LB$-space)
if it is isomorphic to a 
locally convex inductive limit of a sequence of $K$-Fr\'echet spaces
(respectively $K$-Banach spaces).

(ii) We say that a Hausdorff locally convex topological $K$-vector space is
of compact type if it isomorphic to the locally convex inductive
limit of a sequence of convex $K$-vector spaces in which the transition maps
are compact.
\endproclaim

Any convex $K$-vector space $V$ of compact type is certainly an $LB$-space.
Write $V \iso \ilim{n} V_n$,
where $\{V_n\}$ is an inductive sequence
of $K$-Banach spaces with compact transition maps.  If $W_n$
denotes the quotient of $V_n$ by the kernel of the natural map
$V_n \rightarrow V$, then $\{W_n\}$ again forms an inductive sequence
of $K$-Banach spaces (since $V$ is Hausdorff by assumption), whose transition
maps are of compact type and injective.  Furthermore, the natural
map $\ilim{n} V_n \rightarrow \ilim{n} W_n$ is obviously an isomorphism.
Thus any compact type convex $K$-vector space may in fact be written
as the locally convex inductive limit of an inductive
sequence of $K$-Banach spaces having compact injective transition maps.
Conversely, the locally convex inductive limit
of any inductive sequence of $K$-Banach spaces having compact injective
maps is of compact type (and in particular is Hausdorff). 

We refer to \cite{\SCHTAN, p.~445} for the proof of the
preceding claim, and for a fuller discussion
of spaces of compact type.   Of particular importance is the
fact that any such space is reflexive, and that passage 
to the strong dual induces an anti-equivalence of categories between
the category of compact convex $K$-vector spaces and the category
of nuclear $K$-Fr\'echet spaces \cite{\SCHTAN, thm.~1.3}.

Proposition~1.1.10 shows that the Hausdorff $LB$-spaces
are precisely the Hausdorff ultrabornological convex spaces
of $LB$-type.

We recall the following version of the open mapping theorem.

\proclaim{Theorem 1.1.17}  A continuous surjection $V \rightarrow W$
from a convex Hausdorff $LF$-space to a convex Hausdorff space
is strict if and only if $W$ is an $LF$-space.
\endproclaim
\demo{Proof}  The forward implication
is a special case of \cite{\TVS, cor., p.~II.34}.
To prove the converse implication, we must show that any Hausdorff
quotient of a convex Hausdorff $LF$-space is again an $LF$-space.
Write $V = \ilim{n} V_n$ as an inductive limit of a sequence
of Fr\'echet spaces, and let $V \rightarrow W$ be a strict surjection.
For each value of $n$,
let $W_n$ denote the quotient of the Fr\'echet space $V_n$ by the
kernel of the composite $V_n \rightarrow V \rightarrow W;$ then
$W_n$ is again a Fr\'echet space.  The surjection from $V$ to $W$
factors as
$$V \rightarrow \ilim{n} W_n \rightarrow W,$$
and the second morphism is a continuous bijection.  Since the composite
of both morphisms is a strict surjection, by assumption, so is this 
second morphism.  Thus it is a topological isomorphism, and hence $W$
is an $LF$-space, as claimed.
\qed\enddemo

In particular, a continuous bijection between $LF$-spaces
is necessarily a topological isomorphism.

\proclaim{Proposition 1.1.18} If $V$ is a Hausdorff convex $K$-vector space that
is both an $LB$-space and normable then $V$ is a Banach space.
\endproclaim
\demo{Proof} Write $V \iso \ilim{n} V_n$.  We can and do assume, without
loss of generality, that the transition morphisms $V_n \rightarrow
V_{n+1}$ are injective. 
Passing to duals yields
a continuous bijection $$V'_b \rightarrow \plim{n} (V_n')_b.\tag 1.1.19$$
(Recall that the subscript $b$ denotes that each dual is equipped with its
strong topology.) 
Since $V$ is normable, its dual $V'_b$ is a Banach space.
Similarly, each dual $(V_n')_b$ is a Banach space, and so $\plim{n} (V_n')_b$
is a Fr\'echet space.  Thus the open mapping theorem shows that~(1.1.19)
is an isomorphism.

Let $p_n: V'_b \rightarrow (V_n')_b$ denote the natural map. Let $U$
denote the unit ball in $V'_b$, and $U_n$ the unit ball in $(V_n')_b$.
Taking into account the fact that~(1.1.19) is an isomorphism, together with the
definition of the projective limit topology, we see that we may find some $n$
and some $\alpha \in K^{\times}$ such that $$p_n^{-1}(\alpha U_n) \subset U.
\tag 1.1.20$$
Since $U$ contains no non-trivial subspace of $V'_b,$ we see from~(1.1.20)
that $p_n$ has a trivial kernel, and so is injective.  Thus~(1.1.20) shows that
in fact $p_n$ is a topological embedding.  

Let $\hat{V}$ denote the Banach space obtained by completing $V$.
The natural map $\hat{V}'_b \rightarrow V'_b$ is an isomorphism,
and so the natural map $\hat{V}'_b \rightarrow (V_n')_b$ is a topological
embedding.  We deduce from \cite{\TVS, cor.~3, p.~IV.30} that the
injection $V_n \rightarrow \hat{V}$ is a topological embedding,
and hence has a closed image.  It also has a dense image
(since $V'\rightarrow V_n'$ is injective), and so is an isomorphism.
We conclude that the natural map $V\rightarrow \hat{V}$ is an isomorphism,
as required.
\qed\enddemo

\proclaim{Proposition 1.1.21} If $V$ is a Hausdorff $LF$-space that
admits a topological embedding into a convex $K$-vector space of compact type,
then $V$ is itself of compact type.
\endproclaim
\demo{Proof}
Suppose that $V$ admits an embedding into the space $W$ of compact type,
and write $W = \ilim{n} W_n,$
where for each natural number $n$ the map $W_n \rightarrow W_{n+1}$
is a compact, injective $K$-linear map of $K$-Banach spaces.
If we let $V_n$ denote the preimage in $V$ of $W_n$,
then there is a continuous bijection
$\ilim{n}V_n \rightarrow V.$  Now for each $n$ the map $V_n \rightarrow
V_{n+1}$ is again compact, and so $\ilim{n}V_n$ is of compact type.
Since $V$ is assumed to be a Hausdorff $LF$-space, theorem~1.1.17
shows that $\ilim{n}V_n \rightarrow V$ must be an isomorphism,
and so $V$ is of compact type, as required.
\qed\enddemo

\proclaim{Proposition 1.1.22}
If $V$ is a Hausdorff convex $K$-vector space that is both an $LF$-space
and semi-complete, then writing $V = \ilim{n} V_n$ as the inductive limit
of a sequence of Fr\'echet spaces, we have that
for any convex $K$-vector space $W$, the natural map
$\Lin_b(V,W) \rightarrow \plim{n} \Lin_b(V_n,W)$ is a topological
isomorphism.
\endproclaim
\demo{Proof}
Since $V = \ilim{n} V_n,$ we see that for any convex $K$-vector
space $W$ the natural map
$\Lin_b(V,W) \rightarrow \plim{n} \Lin_b(V_n,W)$ is certainly a continuous
bijection.  Propositions~1.1.10 and~1.1.11 show that it
is in fact a topological isomorphism, as claimed.
\qed\enddemo

\proclaim{Corollary 1.1.23}
If $V$ is a semi-complete locally convex Hausdorff $LB$-space,
and if $W$ is a Fr\'echet space, then $\Lin_b(V,W)$ is a Fr\'echet space.
\endproclaim
\demo{Proof}
If we write $V = \ilim{n} V_n$ as the locally convex inductive limit of
a sequence of Banach spaces, then proposition~1.1.22
shows that $\Lin_b(V,W) \iso \plim{n} \Lin_b(V_n,W)$.  Since
a projective limit of a sequence of Fr\'echet spaces is again
a Fr\'echet space, it suffices to prove the corollary in the
case when $V$ is a Banach space.  Since $W$ is a Fr\'echet space,
it is obvious that
$\Lin_b(V_n,W)$ has a countable basis of neighbourhoods of the
origin.
That it is complete follows from \cite{\SCHNA, prop.~7.16}.
\qed\enddemo

\proclaim{Proposition 1.1.24} If $V$ is a semi-complete locally
convex Hausdorff $LB$-space, and if $W$ is a Fr\'echet space,
then any morphism $V \rightarrow W$ is a $BH$-morphism.
\endproclaim
\demo{Proof}
Suppose given such a morphism.
By lemma~1.1.15, we must show that there is an open neighbourhood
of the origin of $V$ that has bounded image in $W$.
The transpose of the given map induces a continuous map
$$W'_b \rightarrow V'_b.\tag 1.1.25$$
Corollary~1.1.23 shows that $V'_b$
is a Fr\'echet space.  Thus \cite{\TVS, lem.~1, p.~IV.26}
implies that some neighbourhood of $W'_b$ maps into a
bounded subset of $V'_b$.  Since $V$ is bornological
(being an $LB$-space), any bounded subset of $V'_b$
is equicontinuous.  Thus~(1.1.25) maps the pseudo-polar
of some bounded $\Cal O_K$-submodule of $W$ 
into the pseudo-polar of an open lattice in $V$.
From \cite{\SCHNA, prop.~13.4} we conclude that
the original map $V \rightarrow W$ maps some
open sublattice of $V$ into a bounded $\Cal O_K$-submodule
of $W$, and we are done.
\qed\enddemo

We now establish some results pertaining to topological tensor products.

\proclaim{Proposition 1.1.26}
Suppose that $U$, $V$ and $W$ are Hausdorff locally convex $K$-vector spaces,
with $U$ bornological, and  $V$ and $W$ complete.  If we are given
an injective continuous map $V \rightarrow W$, then the induced map
$U \cotimes_{K,\pi} V \rightarrow U \cotimes_{K,\pi} W$ is again injective.
\endproclaim
\demo{Proof}
Since $U$ is bornological, the double duality map
$U\rightarrow (U'_b)'_b$ is a topological embedding \cite{\SCHNA, lem.~9.9},
which by \cite{\SCHNA, cor.~17.5~(ii)}
induces a commutative diagram
$$\xymatrix{ U\cotimes_{K,\pi} V \ar[r]\ar[d] & U \cotimes_{K,\pi} W \ar[d] \\
(U'_b)'_b \cotimes_{K,\pi} V \ar[r] & (U'_b)'_b \cotimes_{K,\pi} W}$$
in which the vertical arrows are again topological embeddings.
Thus it suffices to prove the proposition with $U$ replaced by $(U'_b)'_b$.
Propositions~7.16 and~18.2 of \cite{\SCHNA} now yield a commutative diagram
$$\xymatrix{ (U'_b)'_b \cotimes_{K,\pi} V \ar[r] \ar[d] & 
(U'_b)'_b \cotimes_{K,\pi} W \ar[d] \\
\Lin(U'_b,V) \ar[r] & \Lin(U'_b,W),}$$
in which the vertical arrows identify each of the topological tensor products
with the subspace of the corresponding spaces of linear maps consisting
of completely continuous maps.  These arrows are thus injective,
and so is the lower
horizontal arrow (being induced by the injection $V \rightarrow W$).
Consequently the upper horizontal arrow is also injective,
as claimed.
\qed\enddemo

\proclaim{Corollary 1.1.27} If $U_0 \rightarrow U_1$ and $V_0 \rightarrow V_1$
are injective continuous maps of $K$-Banach spaces, then the induced map
$U_0 \cotimes_K V_0 \rightarrow U_1 \cotimes_K V_1$ is injective.
\endproclaim
\demo{Proof}
Banach spaces are bornological and complete.  Thus this follows
by applying proposition~1.1.26 to the pair of maps
$U_0 \cotimes_K V_0 \rightarrow U_0 \cotimes_K V_1$
and $U_0 \cotimes_K V_1 \rightarrow U_1 \cotimes_K V_1.$
\qed\enddemo

\proclaim{Proposition 1.1.28} If $V$ and $W$ are two nuclear $K$-Fr\'echet
spaces then the completed topological tensor product
$V\cotimes_K W$ is again a nuclear $K$-Fr\'echet space.
\endproclaim
\demo{Proof} The completed tensor product of two Fr\'echet
spaces is again a Fr\'echet space.
Propositions~19.11 and~20.4 of \cite{\SCHNA} show that if
$V$ and $W$ are nuclear, then so is $V\cotimes_K W$.
\qed\enddemo

\proclaim{Proposition 1.1.29}
If $V$ and $W$ are two Fr\'echet spaces,
each described as a projective limit of a sequence of Fr\'echet spaces,
$V \iso \plim{n} V_n$ and $W \iso \plim{n} W_n$,
then there is a natural isomorphism
$V\cotimes_K W \iso \plim{n} V_n \cotimes_K W_n$.
\endproclaim
\demo{Proof}
The functoriality of the formation of projective tensor
products and of completions shows that there are natural maps 
$$V\otimes_K W \rightarrow V\cotimes_K W \rightarrow
\plim{n} V_n\cotimes_K W_n.$$
A consideration of the definition of the the projective topology on
$V\otimes_K W$, and hence of the topology on the Fr\'echet
space $V\cotimes_K W$, shows that the second of these maps is
a topological isomorphism.
\qed\enddemo

\proclaim{Lemma 1.1.30}
If $V$ and $W$ are two convex $K$-vector spaces expressed as
locally convex inductive limits,
say $V \iso \ilim{i \in I} V_i$ and $W \iso \ilim{j \in J} W_j$,
then there is a natural isomorphism
$\ilim{(i,j) \in I \times J} V_i \otimes_{K,\iota} W_j 
\iso V\otimes_{K,\iota} W.$
\endproclaim
\demo{Proof}
For each $(i,j) \in I \times J$, the functoriality of
the formation of inductive tensor products yields a continuous map
$V_i \otimes_{K,\iota} W_j \rightarrow V \otimes_{K,\iota} W.$
Passing to the inductive limit yields a continuous map
$\ilim{(i,j) \in I \times J} V_i \otimes_{K,\iota} W_j 
\rightarrow V\otimes_{K,\iota} W,$
which is in fact a bijection (since on the level of abstract
$K$-vector spaces, the formation of tensor products commutes with
passing to inductive limits).
A consideration of the universal property that defines the
inductive tensor product topology shows that this map is
in fact a topological isomorphism.
\qed\enddemo

\proclaim{Proposition 1.1.31}
If $V$ and $W$ are two semi-complete $LB$-spaces,
then the natural bijection $V\otimes_{K,\iota} W \rightarrow
V\otimes_{K,\pi} W$ is a topological isomorphism.  (Thus
from now on we write simply $V\otimes_K W$ to denote the
topological tensor product of $V$ and $W$, and $V \cotimes_K W$
to denote its completion.)
\endproclaim
\demo{Proof}
This follows from the argument used to prove
\cite{\TVS, thm.~2, p.~IV.26}.
Corollary~1.1.23 replaces the citation of ``IV, p.~23, prop.~3''
in that argument, while proposition~1.1.24 replaces the citation
of ``lemma 1''.
\qed\enddemo

Not only is the preceding proposition proved in the same way
as \cite{\TVS, thm.~2, p.~IV.26}, it generalises that result.
Indeed, the strong dual of a reflexive Fr\'echet space is 
a complete $LB$-space (as follows from 
propositions~2 and~4 of \cite{\TVS, pp.~IV.22--23}).

As was already mentioned in~(0.8),
when the hypotheses of the preceding proposition apply,
we write simply $V\otimes_K W$ to denote the tensor product
of $V$ and $W$,
equipped with its inductive (or equivalently, its projective)
tensor product topology.

\proclaim{Proposition 1.1.32}
Let $V$ and $W$ be convex $K$-vector spaces of compact type,
and let $\ilim{n} V_n \iso V$ and $\ilim{n} W_n \iso W$
be expressions of $V$ and $W$ respectively as the 
locally convex inductive limit of Banach spaces,
with injective and compact transition maps. 

(i) There is a natural isomorphism
$\ilim{n} V_n \cotimes_K W_n \iso 
V \cotimes_K W.$
In particular, $V \cotimes_K W$ is of compact type.

(ii) The compact space $V \cotimes_K W$ is dual to the
nuclear Fr\'echet space $V'_b \cotimes_K W'_b.$
\endproclaim
\demo{Proof} 
We begin with~(i).  
Consider the commutative diagram
$$\xymatrix{ \ilim{n} V_n \otimes_K W_n \ar[r] \ar[d]
 & V \otimes_K W \ar[d] \\
\ilim{n} V_n \cotimes_K W_n \ar[r] & V\cotimes_K W.}\tag 1.1.33$$
Lemma~1.1.30 implies that the top horizontal arrow is a topological
isomorphism.  By definition, the right-hand vertical arrow
is a topological embedding that identifies its target
with the completion of its source.  We will show that the
same is true of the left-hand vertical arrow.  This will
imply that the bottom horizontal is also a topological
isomorphism, as required.

The fact that the composite of the top horizontal arrow
and the right-hand vertical arrow in~(1.1.33) is a topological embedding
implies that the left-hand vertical arrow is also a topological
embedding.  It clearly has dense image.
From corollary~1.1.27 and \cite{\SCHNA, lem.~18.12} we see
that $\ilim{n} V_n \cotimes_K W_n$ is the inductive
limit of a sequence of Banach spaces with compact and injective
transition maps, hence is of compact type,
and so in particular is complete.
Thus the left-hand vertical arrow of~(1.1.33) does identify its target with
the completion of its source,  and the discussion of
the preceding paragraph shows that~(i) is proved.

In order to prove part~(ii),
note that
the strong duals $V'_b$ and $W'_b$ are nuclear Fr\'echet spaces
(and so in particular are reflexive).
It follows from \cite{\SCHNA, prop.~20.13} that there is a
topological isomorphism $V \cotimes_K W \iso
(V'_b \cotimes_K W'_b)'_b.$ 
Part~(ii) follows upon
appealing to the
reflexivity of the convex space $V\cotimes_K W$. 
(Alternatively, proposition~1.1.28
shows that
$V'_b \cotimes_K W'_b$ is again a nuclear Fr\'echet
space.  This also provides an alternative proof that $V\cotimes_K W$
is of compact type.)
\qed\enddemo

\proclaim{Lemma 1.1.34} If $U$ and $V$ are
convex $K$-vector spaces, with $U$ barrelled,
and if $A$ and $B$ are bounded $\Cal O_K$-submodules
of $U$ and $V$ respectively, then $A\otimes_{\Cal O_K} B$
is a bounded $\Cal O_K$-submodule of $U\otimes_{K,\imath} V.$
\endproclaim
\demo{Proof}
To prove the lemma,
it suffices to show that for any
separately continuous $K$-bilinear pairing
$U\times V \rightarrow W$,
the induced map $U\otimes_{K,\imath} V
\rightarrow W$ takes $A\otimes_{\Cal O_K} B$
to a bounded subset of $W$.
The given pairing realises $B$ as a pointwise
bounded set of continuous maps $U \rightarrow W$.
Since $U$ is barrelled, it follows that
in fact $B$ induces an equicontinuous set
of maps $U\rightarrow W$.  Thus $A \otimes_{\Cal O_K} B$
does indeed have bounded image in $W$,
as required.
\qed\enddemo

\proclaim{Proposition 1.1.35}
If $U$, $V$ and $W$ are convex $K$-vector spaces,
with $U$ barrelled, 
then the natural isomorphism of abstract $K$-vector spaces
$\Hom(U\otimes_K V, W) \iso \Hom(U,\Hom(V,W))$
(expressing the adjointness of $\otimes$ and $\Hom$)
induces a continuous bijection of convex $K$-vector spaces
$$\Lin_b(U\otimes_{K,\imath} V,W)
\rightarrow \Lin_b(U,\Lin_b(V,W)).$$
\endproclaim
\demo{Proof}
It follows from \cite{\TVS, prop.~6, p.~III.31} that 
the adjointness of $\otimes$ and $\Hom$ induces
an isomorphism of abstract $K$-vector spaces
$$\Lin(U\otimes_{K,\imath} V,W)
\iso \Lin(U,\Lin_b(V,W)).$$
It remains to be shown that this map is continuous, when
both the source and the target are equipped with
their strong topology.
This follows from lemma~1.1.34.
\qed\enddemo

We end this subsection with a presentation of some miscellaneous results and
definitions.

\proclaim{Proposition 1.1.36} Let $V$ and $W$ be Hausdorff convex
$K$-vector spaces,
and suppose that every bounded subset of $W'_b$ is equicontinuous
as a set of functionals on $W$.  (This holds if $W$ is barrelled
or bornological, for example.)  Then passing to the transpose
induces a topological embedding of convex spaces 
$\Lin_b(V,W) \rightarrow \Lin_b(W'_b,V'_b).$
\endproclaim
\demo{Proof}
As observed in~(0.8), passing to the transpose always induces
an injection $\Lin_b(V,W) \rightarrow \Lin_b(W'_b,V'_b).$
We must show that under the assumptions of the proposition this
is a topological embedding. 

Let $B'$ be bounded in $W'_b$ and let $U'$ be a neighbourhood of zero in
$V'_b$.  Then by definition of the strong topology on $V',$
there exists a bounded subset $B$ of $V$ such that
$U' \supset \{ v' \in V' \, | \, v'(B) \subset \Cal O_K\}.$
A neighbourhood basis of the origin in $\Lin_b(W'_b,V'_b)$
is thus given by the sets
$S_{B,B'}' =
\{\phi': W'_b \rightarrow V'_b \, | \, \phi'(b')(B) \subset \Cal O_K
\text{ for all } b' \in B' \},$
as $B$ and $B'$ range over the bounded subsets of $V$ and $W'_b$ 
respectively.
If $S_{B,B'}$ denotes the preimage of $S_{B,B'}'$ in $\Lin(V,W),$
then $S_{B,B'} = \{ \phi: V \rightarrow W \, | \, b'(\phi(B)) \subset
\Cal O_K \text{ for all } b' \in B' \}.$

If $B$ is a bounded subset of $V$ and if $U$ is a neighbourhood of
zero in $W$, then write
$T_{B,U} = \{\phi: V \rightarrow W \, | \, \phi(B) \subset U\}$.
As $B$ runs over all bounded subsets of $V$ and $U$ runs over all
neighbourhoods of zero in $W$ the sets $T_{B,U}$ form a neighbourhood
basis of zero in $\Lin_b(V,W)$.  
Thus we must compare the collections of subsets $T_{B,U}$ and $S_{B,B'}$
of $\Lin(V,W)$.  Since $W$ is locally convex, we may restrict our
attention to those $T_{B,U}$ for which $U$ is an open $\Cal O_K$-submodule
of $W$.

If $U$ is an open $\Cal O_K$-submodule of $V$ and we set
$B'_U = \{ w' \in W \, | \, w'(U)
\subset \Cal O_K\}$ then $B'_U$ is a bounded subset of $W'_s,$
and the Hahn-Banach theorem shows that
$U = \{ w \in W \, | \, b'(w) \in \Cal O_K \text{ for all } b' \in B'_U\}.$
We conclude that $T_{B,U} = S_{B,B'_U}.$

Conversely, if $B'$ is an arbitrary bounded subset of $W'_b$ then by assumption
$B'$ is an equicontinuous set of functionals on $W,$
and so there exists an open $\Cal O_K$-submodule $U$ in $W$ such
that $b'(U) \subset \Cal O_K$ for every $b' \in B'$.  Thus
$B' \subset B'_U,$ and so $S_{B,B'} \supset S_{B,B'_U} = T_{B,U}$.
This result, together with that of the preceding paragraph, shows that
$\Lin_b(V,W) \rightarrow \Lin_b(W'_b,V'_b)$ is a topological embedding,
as required.
\qed\enddemo

Proposition~1.1.36 is an extension of the well-known fact that for any convex
$K$-vector space $W$ satisfying the hypothesis of the proposition, the
double duality map $W \rightarrow (W'_b)'_b$ is a topological
embedding.  (See \cite{\TVS, p.~IV.15}, or \cite{\SCHNA, lem.~9.9}.)

\proclaim{Definition 1.1.37} If $V$ is a locally convex
topological $K$-vector space, then we define the bounded-weak
topology on $V'$ to be the finest locally convex topology on $V'$
which induces on each bounded subset of $V'_b$ a topology
coarser than that induced by $V'_s$.
We let $V'_{bs}$ denote $V'$ equipped with the bounded-weak topology.
\endproclaim

\proclaim{Proposition 1.1.38} Suppose given $V$ and $W$ and
locally convex spaces, with $W$ barrelled and complete,
and a continuous linear map $V \rightarrow W$.
If this map has 
the property that it maps any bounded subset of $V$
into a c-compact subset of $W$, then its transpose induces
a continuous map $W'_{bs} \rightarrow V'_b$.
\endproclaim
\demo{Proof} By \cite{\ROB, prop.~5, p.~153}, we see
that the transpose map $W'_b \rightarrow V'_b$ also takes
bounded subsets into c-compact subsets.

Let $B$ be a bounded subset of $W'_b$.
Since $W$ is barrelled, the classes of equicontinuous, strongly bounded,
and weakly bounded subsets of $W'$ coincide.
Also, any weakly closed weakly bounded subset of $W'$,
being also equicontinuous,
is in fact weakly c-compact.

We see in particular that the weak closure of a
strongly bounded set is again strongly bounded,
and thus that any strongly bounded subset of $W'$
is contained in a strongly bounded subset that is weakly
c-compact.  
Let $B$ be such a subset of $W'$.  The image of $B$
is then c-compact with respect to the topology on $V'_s$.
In particular it is closed in $V'_s$, and so also in
$V'_b$.  On the other hand this image is contained
in a c-compact subset of $V'_b$, and thus is itself
a c-compact subset of $V'_b$.  The topologies
on the image of $B$ induced by $V'_b$ and $V'_s$ thus
coincide, and so the transpose map $W'_s \rightarrow V'_b$
is continuous when restricted to $B$.  By definition,
we conclude that it yields a continuous map
$W'_{bs} \rightarrow V'_b,$ as claimed.
\qed\enddemo

\proclaim{Definition 1.1.39}  We say that a Hausdorff convex $K$-vector
space $V$ is hereditarily complete if any Hausdorff quotient of $V$
is complete.
\endproclaim

Fr\'echet spaces are hereditarily complete, as are spaces
of compact type.  (In fact, the strong dual
to any reflexive Fr\'echet space is hereditarily complete, by \cite{\ROB,
cor., p.~114} and the discussion on page~123 of this reference.)

\proclaim{Proposition 1.1.40} Any closed subspace of a hereditarily
complete Hausdorff convex $K$-vector space is again hereditarily
complete.
\endproclaim
\demo{Proof}
Let $V \rightarrow W$ be a closed embedding, with $W$ a
hereditarily complete convex $K$-vector space.
If $U$ is a closed subspace of $V$, then $V/U$ is a closed
subspace of $W/U$.  Since $W/U$ is complete, by assumption,
the same if true of $V/U$.  Thus $V$ is also hereditarily complete.
\qed\enddemo

\section{1.2} In this subsection we recall
some basic aspects of the theory
of Banach and Fr\'echet $K$-algebras, and especially the notion of
Fr\'echet-Stein structure, introduced in \cite{\SCHTNEW}.

We adopt the convention that all modules are left modules.

\proclaim{Definition 1.2.1} (i) 
A locally convex topological $K$-algebra is
a locally convex topological $K$-vector space $A$ equipped with
a $K$-algebra structure, such that the multiplication map
$A \times A \rightarrow A$ is jointly continuous.
(Equivalently, such that the multiplication map induces a continuous
morphism $A\otimes_{K,\pi} A \rightarrow A$.)

If $A$ is a $K$-Banach space (respectively, a $K$-Fr\'echet space),
then we will refer to $A$ as a $K$-Banach algebra (respectively,
a $K$-Fr\'echet algebra).

(ii) If $A$ is a locally convex topological $K$-algebra, then
a locally convex topological $A$-module is a locally
convex topological $K$-vector space $M$ equipped with an $A$-module
structure (compatible with its $K$-vector space structure), such
that the multiplication map $A \times M \rightarrow M$ is
jointly continuous.   (Equivalently, such that the multiplication
map induces a continuous morphism
$A \otimes_{K,\pi} M \rightarrow M$.)

If $M$ is a $K$-Banach space (respectively, a $K$-Fr\'echet space),
then we will refer to $M$ as an $A$-Banach module (respectively,
an $A$-Fr\'echet module).

(iii) If $A$ is a locally convex topological $K$-algebra, and
$M$ is a locally convex topological $A$-module, then we say that
$M$ is finitely generated over $A$ (as a topological module) if there exists a strict
surjection of $A$-modules $A^n \rightarrow M$ for some $n\geq 0$.
\endproclaim

Note that in part~(iii), when we say that a topological module is finitely
generated, we require that it be presentable topologically (not just
algebraically) as a quotient of some finite power of $A$.

We remark that it follows from \cite{\BGR, prop.~1.2.1/2} that
if $A$ is a Banach algebra in the sense defined above,
then $A$ admits a norm with respect to which it becomes a
Banach algebra in the usual sense \cite{\BGR, def.~3.7.1/1}.

\proclaim{Lemma 1.2.2} If $A$ is a locally convex topological
$K$-algebra, and $M$ is a locally convex topological $A$-module,
then there is a unique locally convex topological $A$-module structure
on the completion $\hat{M}$ of $M$, such that the natural morphism
$M \rightarrow \hat{M}$ becomes $A$-linear.
\endproclaim
\demo{Proof}
Clear.
\qed\enddemo

\proclaim{Lemma 1.2.3}
Let $A \rightarrow B$ be a continuous homomorphism of locally convex
topological $K$- algebras, and let $M$ be a locally convex topological
$A$-module.  If $B \otimes_A M$ is endowed with the quotient topology
obtained by regarding it as a quotient of
$B \otimes_{K,\pi} M$, then it becomes a locally convex topological $B$-module.
\endproclaim
\demo{Proof}
Clear.
\qed\enddemo

In the context of lemma~1.2.3, we will always regard $B\otimes_A M$ as
a topological $B$-module by endowing with it the quotient topology
induced from $B\otimes_{K,\pi} M.$  We let $B\cotimes_A M$
denote the completion of the locally convex $B$-module $B\otimes_A M$.

\proclaim{Proposition 1.2.4} Let $A$ be a Noetherian $K$-Banach
algebra.

(i) Every finitely generated $A$-module has a unique
structure of $K$-Banach space making it an $A$-Banach module.

(ii) If $f: M \rightarrow N$ is an $A$-linear morphism between
finitely generated $A$-modules, then $f$ is continuous and strict 
with respect to the Banach space structures on $M$ and $N$ given by~(i).
\endproclaim
\demo{Proof}
Part~(i) is a restatement of \cite{\BGR, props.~3.7.3/3}.
Part~(ii) is a restatement of \cite{\BGR, prop. 3.7.3/2, cor.~3.7.3/5}.
\qed\enddemo

This proposition shows in particular that if $A$ is Noetherian,
then the forgetful functor from
the category of finitely generated $A$-Banach modules to
the category of finitely generated abstract $A$-modules is an equivalence
of categories.  In particular, since the latter category is
abelian ($A$ being Noetherian), the same is true of the former.

\proclaim{Proposition 1.2.5} Let $A \rightarrow B$ be a continuous
homomorphism of locally convex topological $K$-algebras, and suppose
that $B$ is hereditarily complete (in the sense of definition~1.1.39).
If $M$ is a finitely generated locally convex
topological $A$-module, then the natural map
$B\otimes_A M \rightarrow B\cotimes_A M$ is surjective, and 
$B\cotimes_A M$ is a finitely generated topological $B$-module.
\endproclaim
\demo{Proof} By assumption there is a strict surjection $A^n \rightarrow M$
for some $n\geq 0,$ and hence a strict surjection
$B^n \rightarrow B\otimes_A M$.  Since $B$, and hence $B^n,$ is hereditarily
complete, the quotient of $B^n$ by the closure of the kernel of this
surjection is complete.  Thus $B\otimes_A M \rightarrow B\cotimes_A M$
is surjective, as claimed.
\qed\enddemo

Note that if $A$ and $B$ are Noetherian Banach algebras,
then \cite{\BGR, prop.~3.7.3/6} shows that the map $B\otimes_A M
\rightarrow B \cotimes_A M$ is in fact an isomorphism.

\proclaim{Definition 1.2.6} Let $A$ be a locally convex topological
$K$-algebra.  A weak Fr\'echet-Stein structure on $A$ consists of
the following data:

(i) A sequence of locally convex topological
$K$-algebras $\{A_n\}_{n \geq 1},$ such that each
$A_n$ is hereditarily complete (in the sense of definition~1.1.39).

(ii) For each $n\geq 1,$ a continuous $K$-algebra homomorphism 
$A_{n+1} \rightarrow A_n$, which is a $BH$-map of convex $K$-vector spaces
(in the sense of definition~1.1.13).

(iii) An isomorphism of locally convex topological $K$-algebras
$A \iso \plim{n} A_n$ (the projective limit being
taken with respect to the maps of~(ii)), such that each
each of the induced maps $A \rightarrow A_n$ has dense image.

We say that $A$ is a weak Fr\'echet-Stein $K$-algebra if $A$
admits a weak Fr\'echet-Stein structure.
\endproclaim

Note that if $A$ admits a weak Fr\'echet-Stein structure,
then it is the projective limit of a sequence of Hausdorff
convex $K$-vector spaces
with $BH$-transition maps, and so is a Fr\'echet space.
Thus a weak Fr\'echet-Stein algebra is, in particular,
a Fr\'echet algebra.

Suppose given a sequence of locally convex topological
$K$-algebras $\{A_n\}_{n\geq 1}$ satisfying conditions~(i)
and~(ii) of definition~1.2.6, and an isomorphism
of locally convex topological $K$-algebras $A \iso \plim{n} A_n$.
If we let $B_n$ denote the closure of the image of $A$ in $A_n,$
then the projective system $\{B_n\}_{n\geq 1}$
also satisfies conditions~(i) (by proposition~1.1.40) and~(ii)
(since clearly the $BH$-map $A_{n+1} \rightarrow A_n$ induces
a $BH$-map $B_{n+1} \rightarrow B_n$).
The natural map $A \iso \plim{n} B_n$ is also an isomorphism,
and by construction it satisfies condition~(iii) of definition~1.2.6. 
Thus $A$ is a weak Fr\'echet-Stein $K$-algebra,
and that the projective system $\{B_n\}_{n\geq 1}$
forms a weak Fr\'echet-Stein structure on $A$.

We say that two weak Fr\'echet-Stein structures $A \iso \plim{n} A_n$
and $A \iso \plim{n} B_n$ on $A$ are equivalent if the
projective systems $\{A_n\}_{n \geq 1}$ and $\{B_n\}_{n\geq 1}$ are
isomorphic in the category of projective systems of
topological $K$-algebras.

\proclaim{Proposition 1.2.7} If $A$ is a locally convex topological $K$-algebra,
then any two weak Fr\'echet-Stein structures on $A$ are equivalent.
\endproclaim
\demo{Proof}
Let $A \iso \plim{n} A_n$ and $A \iso \plim{n} B_n$ be two
weak Fr\'echet-Stein structures on $A$.  
Since each transition map $A_{n+1} \rightarrow A_n$ is
assumed to be $BH$-maps with dense image, we may factor 
each map as $A_{n+1} \rightarrow V_n \rightarrow A_n,$
where $V_n$ is a Banach space, and each map has dense
image.  Similarly, we may factor each map $B_{n+1} \rightarrow B_n$
as $B_{n+1} \rightarrow W_n \rightarrow B_n,$ where $W_n$ is 
a Banach space, and each map has dense image.

Thus we obtain a pair of
projective systems of Banach spaces $\{V_n\}_{n\geq 1}$
and $\{W_n\}_{n\geq 1},$ each with dense transition maps,
the first being equivalent to $\{A_n\}_{n \geq 1}$ and
the second being equivalent to $\{B_n\}_{n\geq 1}$.
Since the projective limit of each
of $\{V_n\}_{n\geq 1}$ and $\{W_n\}_{n\geq 1}$ is the same
Fr\'echet space $A$, they are necessarily equivalent. 
It follows that $\{A_n\}_{n\geq 1}$ and $\{B_n\}_{n\geq 1}$
are equivalent as projective systems of convex $K$-vector spaces.
Since $A$ has dense image in each of the topological
algebras $A_n$ and $B_n$, they are then necessarily
equivalent as projective systems of topological $K$-algebras.
\qed\enddemo

\proclaim{Definition 1.2.8} Let $A$ be a weak Fr\'echet-Stein algebra,
and let $A \iso \plim{n} A_n$ be a choice of
a weak Fr\'echet-Stein structure on $A$.

If $M$ is a locally convex topological $A$-module,
then
we say that $M$ is coadmissible (with respect to the given
weak Fr\'echet-Stein structure on $A$) if we may find the following data:

(i) A sequence
$\{M_n\}_{n \geq 1}$, such that $M_n$ is a finitely generated locally
convex topological $A_n$-module for each~$n\geq 1$. 

(ii)
An isomorphism of topological $A_n$-modules $A_n \cotimes_{A_{n+1}} M_{n+1} \iso M_n$
for each~$n\geq 1.$  

(iii) An isomorphism of topological $A$-modules $M \iso \plim{n} M_n.$  (Here the projective
limit is taken with respect to the transition maps $M_{n+1} \rightarrow M_n$
induced by the isomorphisms of~(ii).)

We will refer to such a projective sequence $\{M_n\}_{n\geq 1}$ as above
as an $\{A_n\}_{n\geq 1}$-sequence (for the coadmissible module $M$).
\endproclaim

If $A$ is equipped with a weak Fr\'echet-Stein structure $A =\plim{n} A_n$,
then since each $A_n$ is hereditarily complete, proposition~1.2.5
implies that for any finitely generated locally convex topological
$A_{n+1}$-module $M_n$, the
completed tensor product $A_n \cotimes_{A_{n+1}} M_n$ is finitely
generated.  Thus hypothesis~(ii) of definition~1.2.8 is reasonable.
Note that it also implies that $M_n$ is Hausdorff, for each $n\geq 1$.

If we are in the situation of definition~1.2.8,
then for any $n\geq 1,$ we may construct a commutative diagram
$$\xymatrix{ A_{n+1}^m \ar[r] \ar[d] & A_n^m \ar[d] \\ M_{n+1} \ar[r] & M_n,}$$
in which the vertical arrows are strict surjections.
Since the upper horizontal arrow is $BH$,
by assumption,
the same is true of the lower horizontal arrow. 
Thus $M \iso \plim{n} M_n$ is isomorphic to the projective limit of Hausdorff
locally convex spaces under $BH$-transition maps, and so is a
Fr\'echet space.

Let us remark that the isomorphism $A \iso \plim{n} A_n$ witnesses
$A$ as a coadmissible topological $A$-module.

\proclaim{Proposition 1.2.9} If $A$ is a weak Fr\'echet-Stein algebra,
and if $M$ is a locally convex topological $A$-module that is coadmissible
with respect to one choice of weak Fr\'echet-Stein structure on $A$,
then $M$ is coadmissible with respect to any such choice.
\endproclaim
\demo{Proof}
Let $A \iso \plim{n \geq 1} A_n$ and $A \iso \plim{n \geq 1} B_n$ be two weak
Fr\'echet-Stein structures on $A$. 
Proposition~1.2.7 shows that
we may find monotonically increasing functions $\phi$ and $\psi$
mapping $\N$ to itself, 
each of whose images is a cofinal (equivalently, infinite)
subset of $\N$,
and such that for 
each $n\geq 1$, the map $A \rightarrow A_n$ (respectively
$A \rightarrow B_n$) factors through the map $A \rightarrow B_{\phi(n)}$
(respectively $A \rightarrow A_{\psi(n)}$).

Suppose that $\{M_n\}_{n \geq 1}$
is an $\{A_n\}_{n\geq 1}$-sequence realising $M$ as a coadmissible
module for the weak Fr\'echet-Stein structure on $A$ induced by
the projective sequence $\{A_n\}_{n\geq 1}$.
For each $n \geq 1,$
define $N_n := B_n \cotimes_{A_{\psi(n)}} M_{\psi(n)}.$ 
Proposition~1.2.5 shows that $N_n$ is a finitely generated
$B_m$-module. 
Since $\psi(n+1) \geq \psi(n)$,
the transition map $M_{\psi(n+1)} \rightarrow M_{\psi(n)}$ induces
a continuous $B_{n+1}$-linear map $N_{n+1} \rightarrow N_n$.
This in turn induces 
a continuous $B_n$-linear map $B_n\cotimes_{B_{n+1}} N_{n+1} \rightarrow N_n$,
which sits in the commutative diagram
$$\xymatrix{B_n \cotimes_{B_{n+1}} N_{n+1} \ar[r] \ar@{=}[d] &
N_n \ar@{=}[ddd] \\
B_n\cotimes_{B_{n+1}}(B_{n+1}\cotimes_{A_{\psi(n+1)}} M_{\psi(n+1)} )
\ar[d]^-{\sim} & \\
B_n \cotimes_{A_{\psi(n+1)}} M_{\psi(n+1)} \ar[d]^-{\sim} & \\
B_n\cotimes_{A_{\psi(n)}} (A_{\psi(n)}\cotimes_{A_{\psi(n+1)}} M_{\psi(n+1)})
\ar[r] & B_n\cotimes_{A_{\psi(n)}} M_n.}$$
By assumption, the natural map
$A_{\psi(n)}\cotimes_{A_{\psi(n+1)}} M_{\psi(n+1)}\rightarrow M_{\psi(n)}$
is an isomorphism,
and thus the same is true of the lower, and so also of the upper,
horizontal arrow in this diagram.  Altogether, we see that the sequence
of topological modules $\{N_n\}_{n\geq 1}$ satisfies conditions~(i) and~(ii)
of definition~1.2.8.

The projection $M \rightarrow M_{\psi(n)}$ induces a continuous $A$-linear map
$M \rightarrow N_n$, for each $n \geq 1$ (compatible with the transition
map $N_{n+1} \rightarrow N_n$), and consequently there is a continuous
$A$-linear map
$M \rightarrow \plim{n} N_n$.  We will prove that this is a topological
isomorphism,
and hence show that $\{N_n\}_{n\geq 1}$ also satisfies condition~(iii)
of definition~1.2.8, and thus is a $\{B_n\}_{n\geq 1}$-sequence
realising $M$ as a coadmissible module with respect to the weak Fr\'echet-Stein
structure on $M$ induced by the projective sequence $\{B_n\}_{n\geq 1}$.

Note that there is the obvious natural map $M_{\psi(n)} \rightarrow N_n$.
On the other hand,  since $\psi(\phi(n)) \geq n$ for any $n \geq 1$,
we see that for any such value of $n$ there is a natural map 
$N_{\phi(n)} = B_{\phi(n)}\cotimes_{A_{\psi(\phi(n))}} M_{\psi(\phi(n))}
\rightarrow A_n \cotimes_{A_{\psi(\phi(n))}} M_{\psi(\phi(n))}
\rightarrow M_n$
(the second map being induced by the fact that $A \rightarrow A_n$
factors through $A \rightarrow B_{\phi(n)},$ and the third map by
the transition function $M_{\psi(\phi(n))} \rightarrow M_n$).
Thus the projective systems $\{M_n\}_{n\geq 1}$ and $\{N_n\}_{n\geq 1}$
are equivalent, and in particular the natural
map $M \rightarrow \plim{n} N_n$ is a topological isomorphism.   This
completes the proof of the proposition.
\qed\enddemo

\proclaim{Definition 1.2.10} Let $A$ be a locally convex
topological $K$-algebra.
A Fr\'echet-Stein structure on $A$ is a weak Fr\'echet-Stein
structure $A \iso \plim{n} A_n$ on $A$, such that for each $n\geq 1,$
the algebra $A_n$ is a left
Noetherian $K$-Banach algebra, and the transition map
$A_{n+1} \rightarrow A_n$ is right flat.
\endproclaim

The notion of Fr\'echet-Stein structure and Fr\'echet-Stein
algebra is introduced in \cite{\SCHTNEW, def., p.~7}.   
The notion of weak Fr\'echet-Stein structure is a useful
generalisation, which provides additional
flexibility when it comes to constructing
and analysing Fr\'echet-Stein algebras and their coadmissible
modules.

Suppose that $ A \iso \plim{n} A_n $ is a Fr\'echet-Stein
structure on the locally convex topological $K$-algebra $A$.
If $M \iso \plim{n} M_n$ is a coadmissible $A$-module,
then proposition~1.2.4 implies that each $M_n$ is an $A_n$-Banach module,
while the remark following the proof of proposition~1.2.5 implies
that the natural map $A_n \otimes_{A_{n+1}} M_{n+1} \rightarrow M_n$
is an isomorphism.  

The following theorem incorporates the most important results
of \cite{\SCHTNEW, \S 3}.  

\proclaim{Theorem 1.2.11}  Let $A$ be a Fr\'echet-Stein algebra,
and let $A \iso \plim{n} A_n$ be a choice of
weak Fr\'echet-Stein structure on $A$.

(i) If $M$ is a coadmissible $A$-module, and if $\{M_n\}_{n\geq 1}$
is an $\{A_n\}_{n\geq 1}$-sequence for which there is a topological isomorphism
$M \iso \plim{n} M_n$,
then for each value of~$n$, the 
natural map $A_n \cotimes_A M \rightarrow M_n$ is an isomorphism.
Consequently, the natural map
$A \rightarrow \plim{n} A_n\cotimes_A M$ is an isomorphism.

(ii) The category of coadmissible $A$-modules is closed under passing to
finite direct sums,
kernels, cokernels, closed submodules, and Hausdorff quotient modules.
In particular, it is abelian.
\endproclaim
\demo{Proof} Let $A \iso \plim{m} B_n$ be a Fr\'echet-Stein
structure on $A$; such a structure exists, by assumption.
Following the proof of proposition~1.2.9,
we use proposition~1.2.7 to choose
monotonically increasing functions $\phi$ and $\psi$
mapping $\N$ to itself, 
each of whose images is a cofinal (equivalently, infinite)
subset of $\N$,
such that for 
each $n\geq 1$, the map $A \rightarrow A_n$ (respectively
$A \rightarrow B_n$) factors through the map $A \rightarrow B_{\phi(n)}$
(respectively $A \rightarrow A_{\psi(n)}$).
For each $n \geq 1,$
we define $N_n := B_n \cotimes_{A_{\psi(n)}} M_{\psi(n)}.$ 
As in the proof of proposition~1.2.9, we find that
$\{N_n\}_{n\geq 1}$ is a $\{B_n\}_{n\geq 1}$-sequence
realising $M$ as a coadmissible module with respect to the
Fr\'echet-Stein structure on $A$ induced by the projective sequence
$\{B_n\}$.
It follows from \cite{\SCHTNEW, cor.~3.1} that for any 
$n\geq 1$, the natural map
$B_n \otimes_A M \rightarrow N_n$ is a topological isomorphism.
(Actually, that result as stated shows only
that this map is a continuous $B_n$-linear bijection between
finitely generated topological $B_n$-modules.  However,
such a map is necessarily a topological isomorphism, by
proposition~1.2.4.)
Since $N_n$ is in fact complete, the source of this isomorphism
is isomorphic to $B_n \cotimes_A M.$
Thus the natural map
$B_n\cotimes_A M \rightarrow N_n = B_n\cotimes_{A_{\psi(n)}} M_{\psi(n)}$
is an isomorphism for any $n\geq 1$.
Applying this with $n$ replaced by $\phi(n)$,
we obtain the required isomorphism
$$\multline A_n\cotimes_A M \iso
A_n\cotimes_{B_{\phi(n)}} (B_{\phi(n)}\cotimes_A M)  \\ \iso
A_n\cotimes_{B_{\phi(n)}} (B_{\phi(n)} \cotimes_{A_{\psi(\phi(n))}}
M_{\psi(\phi(n))}) 
\iso A_n \cotimes_{A_{\psi(\phi(n))}} M_{\psi(\phi(n))}
\iso M_n \endmultline $$
(the final isomorphism following from our assumption that
$\{M_n\}$ is an $\{A_n\}$-sequence).
This proves~(i).

Part~(ii) is a restatement of \cite{\SCHTNEW, cor.~3.4~(i),(ii)}.
\qed\enddemo

Note that~(i), combined with proposition~1.2.4, implies that
any $A$-linear map between coadmissible locally convex topological
$A$-modules is necessarily continuous and strict.  In particular,
the forgetful functor from the category of coadmissible locally
convex topological $A$-modules to the category of abstract $A$-modules
is fully faithful.  Part~(ii) then implies that any finitely presented
$A$-module admits a unique topology with respect to
which it becomes a coadmissible topological $A$-module.
(This is a restatement of \cite{\SCHTNEW, cor.~3.4~(v)}.)

\proclaim{Definition 1.2.12}
We say that a topological $K$-algebra $A$ is a nuclear 
Fr\'echet algebra
if there exists a projective system
$\{A_n\}_{n \geq 1}$
of $K$-Banach algebras
with compact transition maps, and an isomorphism
of topological $K$-algebras
$A \iso \plim{n} A_n$.
\endproclaim

It follows from \cite{\SCHNA, cor.~16.6} that any
topological $K$-algebra $A$ satisfying the conditions
of the preceding definition is in particular a nuclear
Fr\'echet space.  

If we let $B_n$ denote the closure of the image of $A$
under the natural map $A \rightarrow A_n$, for each $n \geq 1,$
the space $B_n$ is a Banach subalgebra,
and the natural map $A \rightarrow \plim{n} B_n$ is 
an isomorphism.  Thus, in definition~1.2.12, it is no
loss of generality to require that $A$ have dense image
in each of the topological $K$-algebras $A_n$.
Since compact type maps are $BH$-maps (by lemma~1.1.14),
and since Banach space and compact type space are both
hereditarily complete, we find that a nuclear Fr\'echet
algebra is in particular a weak Fr\'echet-Stein algebra.

\proclaim{Lemma 1.2.13} The completed tensor product of 
a pair of nuclear Fr\'echet algebras over $K$ 
is again a nuclear Fr\'echet algebra.
\endproclaim
\demo{Proof}
If $A \iso \plim{n} A_n$ and $B \iso \plim{n} B_n$
are two nuclear Fr\'echet algebras, written as projective
limits, with compact transition maps, of $K$-Banach algebras,
then propositions~1.1.29 and \cite{\SCHNA, lem.~18.12}
together yield a topological isomorphism
$A \cotimes_K B \iso \plim{n} A_n \cotimes_K B_n,$
where the transition maps are again compact.
This proves the lemma.
\qed\enddemo

\proclaim{Proposition 1.2.14} Let $A$ be a nuclear Fr\'echet
algebra, and write $A \iso \plim{n} A_n$
as in definition~1.2.12, with the further hypothesis
that $A$ has dense image in each $A_n$.  Let $A^{\op}$ denote
the opposite algebra to $A$.  If $V$ is a compact
type convex $K$-vector space, then the following structures
on $V$ are equivalent:

(i) An $A$-module structure for which
the multiplication map
$A\times V \rightarrow V$ is separately continuous.

(ii) A topological $A^{\op}$-module structure on the dual space $V'_b$.

(iii) An isomorphism $V \iso \ilim{n} V_n$, where $\{V_n\}$
is an inductive sequence of 
$K$-Banach spaces, in which each $V_n$ is endowed with the structure
of a topological
$A_n$-module, and such that the transition maps $V_n \rightarrow V_{n+1}$
are compatible with the map $A_{n+1} \rightarrow A_n$. 
\endproclaim
\demo{Proof}
Proposition~1.1.35 shows that giving a separately
continuous bilinear map
$$A\times V \rightarrow V \tag 1.2.15$$
is equivalent
to giving a continuous map $A \rightarrow \Lin_b(V,V),$
and (taking into account that $A$ and $V'_b$ are both
nuclear Fr\'echet spaces, so that a separately continuous bilinear
map
$$A\times V'_b \rightarrow V'_b \tag 1.2.16$$
is automatically jointly continuous)
that giving a jointly continuous bilinear map as in~(1.2.16)
is equivalent to giving a map
$A \rightarrow \Lin_b(V'_b,V'_b).$  Proposition~1.1.36,
and the fact that $V$ is reflexive, shows that the natural
map $\Lin_b(V,V) \rightarrow \Lin_b(V'_b,V'_b)$ is a topological
isomorphism.  Thus passing to duals shows that the existence
of the separately continuous bilinear map~(1.2.15)
is equivalent to the existence of the jointly continuous bilinear
map~(1.2.16).  It is straightforward to check that~(1.2.15)
induces an $A$-module structure on $V$ if and only if~(1.2.16)
induces an $A^{\op}$-module structure on $V'_b$.  Thus~(i) and~(ii)
are equivalent.

If~(iii) holds, then each $V_n$ is in particular
a topological $A$-module, and the transition maps $V_n \rightarrow V_{n+1}$
are compatible with the $A$-module structure on source and target.
Passing to the inductive limit
in $n$, we find that $V$ is an $A$-module, and that
the map $A\times V \rightarrow V$ describing this module
structure is separately continuous.  Thus~(iii) implies~(i).

We now turn to showing that~(i) implies~(iii).
Suppose to begin with that $V$ is any compact type convex $K$-vector space,
and consider the space $\Lin_b(A,V)$. 
By \cite{\SCHNA, cor.~18.8} there is an isomorphism
$\Lin_b(A,V) \iso A'_b\cotimes_K V$, and so by
proposition~1.1.32~(i) we see that $\Lin_b(A,V)$ is of compact type.
By proposition~1.1.10, any map $A \rightarrow V$ factors through a map
$A_n \rightarrow V_n$ for some $n \geq 1$,
and hence there is a continuous bijection
$$\ilim{n} \Lin_b(A_n, V_n) \rightarrow \Lin_b(A,V).\tag 1.2.19$$
(In fact this is a topological isomorphism, since its source
and target are both $LB$-spaces.)
The image of $\Lin_b(A_n,V_n)$ in $\Lin_b(A,V)$ is
an $A_n$-invariant $BH$-subspace of $\Lin_b(A,V)$.
The right regular representation of $A_n$ on $\Lin_b(A_n,V_n)$
makes the latter a topological $A_n$-Banach module,
and thus $\Lin_b(A,V)$ satisfies the condition of~(iii).

Now suppose that $V$ satisfies condition~(i).
Proposition~1.1.35 also yields a continuous injection
$$V \rightarrow \Lin_b(A,V).\tag 1.2.18$$
Evaluation at $1 \in A$ gives a splitting of the injection~(1.2.18),
and hence this map is in fact a closed embedding.
If we give $\Lin_b(A,V)$
the $A$-module structure induced by the right regular
representation of $A$, then~(1.2.18) is furthermore a map of $A$-modules.
The discussion of the preceding paragraph shows that
$V$ is thus a closed $A$-submodule of an $A$-module
satisfying condition~(iii).  Thus $V$ also satisfies condition~(iii).
More precisely, if we let $V_n$ denote the preimage of
$V$ in $\Lin_b(A_n,V_n)$, then $V_n$ is a closed $A_n$-invariant 
submodule of $\Lin_b(A_n,V_n)$ (because $A$ is dense in $A_n$),
and $V = \ilim{n} V_n$.  
\qed\enddemo

\head 2. Non-archimedean function theory \endhead

\section{2.1}  In this subsection we recall some notions
of non-archimedean function theory.

\proclaim{Definition 2.1.1} If $X$ is a set and $V$ is
a $K$-vector space, we let $\Fun(X,V)$ denote the $K$-vector
space of $V$-valued functions on $X$.
\endproclaim

The formation of $\Fun(X,V)$
is evidently covariantly functorial in $V$ and contravariantly
functorial in $X$. In particular, for each $x\in X$ we obtain
a natural map $\Fun(X,V) \rightarrow \Fun(x,V) \iso V,$
that we denote by $\ev_x$.  (Of course, this is simply
the evaluation map at $x$.)

\proclaim{Definition 2.1.2}
If $X$ is a topological space and $V$ is a Hausdorff
locally convex topological $K$-vector space, we
let $\Con(X,V)$ denote the $K$-vector space of continuous
$V$-valued functions on $G$, equipped with the
(Hausdorff locally convex)
topology of uniform convergence on compact sets.
\endproclaim

The vector space $\Con(X,V)$ is a subspace of $\Fun(X,V)$,
and the formation of $\Con(X,V)$
is covariantly functorial in $V$ and contravariantly
functorial in $X$. In particular, if $x\in X$ then
the restriction of $\ev_x$ to $\Con(X,V)$ is a continuous
map to $V$, which we denote by the same symbol $\ev_x$.

If $X$ is compact and $V$ is a Banach space then so is $\Con(X,V)$;
in particular, if $X$ is compact then $\Con(X,K)$ is a Banach space.

Multiplying $K$-valued functions by vectors in $V$
induces a jointly continuous
bilinear map $\Con(X,K) \times V \rightarrow \Con(X,V),$
and hence a continuous map $\Con(X,K)\otimes_{K,\pi} V
\rightarrow \Con(X,V).$
The example of \cite{\SCHNA, pp.~111--112} shows that if
$V$ is complete, then this map induces an isomorphism
$$\Con(X,K)\cotimes_{K,\pi} V \iso \Con(X,V).\tag 2.1.3$$

\proclaim{Proposition 2.1.4}
If $X$ is a compact topological space,
and if $V$ is a Fr\'echet space (respectively a Banach space),
then $\Con(X,V)$ is again a Fr\'echet space (respectively a Banach space).
\endproclaim
\demo{Proof}
Since $V$ is in particular assumed to be complete,
the isomorphism~(2.1.3) shows
that $\Con(X,V)$ is the completed tensor product of a Banach space
and a Fr\'echet space (respectively a Banach space),
and hence is itself a Fr\'echet space (respectively a Banach space).
\qed\enddemo

\proclaim{Corollary 2.1.5} If $X$ is a locally compact topological
space that can be written as a disjoint union of compact open
subsets and $V$ is a Fr\'echet space 
then $\Con(X,V)$ is barrelled.
\endproclaim
\demo{Proof}
Let $X = \coprod_i X_i$ be the hypothesised decomposition of $X$
into a disjoint union of compact open subsets.  Then
we obtain a natural isomorphism $\Con(X,V) \iso \prod_i \Con(X_i,V)$.
Proposition~2.1.4 implies in particular that each of the spaces
$\Con(X_i,V)$ is barrelled.
Thus $\Con(X,V)$ is isomorphic to a product of barrelled spaces,
hence is itself barrelled \cite{\TVS, cor., p.~IV.14}.
\qed\enddemo

Let $X$ be a compact topological space.
If $V$ is a Hausdorff convex $K$-vector space and $W$ is a
$BH$-subspace of $V$, then the image of the natural
map $\Con(X,\overline{W}) \rightarrow \Con(X,V)$ is
a $BH$-subspace of $V$.
(We refer to definition~1.1.1 for the
notion of a $BH$-subspace of $V$.  Note in particular that if 
$W$ is a $BH$-subspace of $V$ then
$\overline{W}$ denotes the latent Banach space structure on $W$.)

\proclaim{Proposition 2.1.6}
If $X$ is a compact topological space,
and if $V$ is either a Fr\'echet space or a space of compact type,
then any $BH$-subspace of $\Con(X,V)$ is contained in
the image of the natural map $\Con(X,\overline{W}) \rightarrow \Con(X,V)$,
for some $BH$-subspace $W$ of $V$.  In particular,
the natural map
$\ilim{W} \Con(X,\overline{W}) \rightarrow \Con(X,V)$
(where the locally convex inductive limit is
taken over all $BH$-subspaces $W$ of $V$)
is a continuous bijection.
If $V$ is a Fr\'echet space, then this map
is even a topological isomorphism.
\endproclaim
\demo{Proof}
Since $V$ is in particular assumed to be complete, 
the same is true of $\Con(X,V)$.  Thus any $BH$-subspace of
$\Con(X,V)$ is contained in the subspace generated by some 
closed bounded $\Cal O_K$-submodule of $\Con(X,V)$.  Any such set
is in turn contained in the set $\Con(X,A)$, for some closed bounded
$\Cal O_K$-submodule $A$ of $V$.   
Proposition~1.1.11 shows that
we may find a $BH$-subspace $W$ of $V$ containing $A$,
such that $A$ is bounded when regarded as a subset of $\overline{W}$
and such that the topology induced on $A$ by $\overline{W}$ agrees with
the topology induced on $A$ by $V$.  We may therefore regard
$\Con(X,A)$ as a bounded subset of $\Con(X,\overline{W}).$
The first claim of the proposition follows.

If $U$ ranges over the directed set of all $BH$-subspaces of $\Con(X,V)$,
then the natural map $\ilim{U} \overline{U} \rightarrow \Con(X,V)$
is a continuous bijection.   If $V$ is a Fr\'echet space, then
proposition~2.1.4 implies that $\Con(X,V)$ is also, and so in particular
is ultrabornological.  Thus in this case, the preceding natural map
is even a topological isomorphism.
The result of the preceding paragraph
shows that the natural map $\ilim{W} \Con(X,\overline{W})
\rightarrow \ilim{U}\overline{U}$ is a also a topological isomorphism.
The remaining claims of the proposition now follow.
\qed\enddemo

If $V$ is a convex space of compact type, then in particular
$V$ is of $LB$-type,
and the preceding proposition then implies
that the convex space $\Con(X,V)$ is also of $LB$-type.

\proclaim{Proposition 2.1.7}
If $V$ and $W$ are two Hausdorff convex $K$-vector spaces,
and if $X$ is a locally compact topological space,
then the natural map
$\Lin_b(V,W) \rightarrow \Lin_b(\Con(X,V),\Con(X,W))$ (induced
by the functoriality of the formation of $\Con(X,\text{--}\,)$)
is continuous.
If $X$ is non-empty then it is even a topological embedding.
\endproclaim
\demo{Proof}
If $C$ is a compact subset of $X$ and $U$ is a neighbourhood
of zero in $V$ (respectively $W$) then let $S_{C,U}$ denote
the subset of $\Con(X,V)$ (respectively $\Con(X,W)$) defined by
$S_{C,U} = \{f \, | \, f(C) \subset U\}$.  This set 
is a neighbourhood of zero in $\Con(X,V)$ (respectively $\Con(X,W)$),
and such sets form a neighbourhood basis of zero, if
$C$ is allowed to run over all compact subsets of
$X$ and $U$ to run over all neighbourhoods of zero in $V$
(respectively $W$).
In light of this it is straightforward to check that a subset $B$ of
$\Con(X,V)$ is bounded if and only if for every compact subset $C$ of
$X$, the subset $\{f(c) \, | \, f\in B, c\in C \}$ 
is a bounded subset of $V$.

If $B$ is bounded subset of $\Con(X,V)$ and $S_{C,U}$ is a neighbourhood
of zero in $\Con(X,W)$ of the type described in the preceding
paragraph, let $T_{B,C,U}$ denote the subset of $\Lin(\Con(X,V),
\Con(X,W))$ defined by $T_{B,C,U} = \{\phi  \, | \,
\phi(B) \subset S_{C,U}\}.$ This set is 
a neighbourhood of zero in $\Lin_b(\Con(X,V),\Con(X,W))$,
and such sets form a neighbourhood basis of zero
if $B$ is allowed to run over all  bounded subsets $B$ of $\Con(X,V)$,
$C$ over all compact subsets of $X$, and $U$ over all neighbourhoods
of zero in $W$.

Fixing such a $B$, $C$, and $U$, let $D$ denote the subset of $V$
defined by $D = \{f(c) \, | \, f\in B, c \in C\}$.  As observed
above, $D$ is a bounded subset of $V$.
The preimage of $T_{B,C,U}$ under the natural map
$\Lin(V,W) \rightarrow \Lin(\Con(X,V),\Con(X,W))$ is equal to
$\{\phi  \, | \, \phi(D) \subset U\}$.  Since $D$
is bounded in $V$ and $U$ is open in $W$ this a neighbourhood
of zero in $\Lin_b(V,W)$.  Thus we have proved that
the natural map $\Lin_b(V,W) \rightarrow \Lin_b(\Con(X,V),\Con(X,W))$
is continuous. 

If $X$ is non-empty then this map is also evidently injective,
and so it remains
to prove that any open subset of $\Lin_b(V,W)$ can be obtained
by pulling back an open subset of $\Lin_b(\Con(X,V),\Con(X,W))$.

Let $D$ be an arbitrary bounded subset of $V$, and define $B$
to be that subset of $\Con(X,V)$ consisting of the constant functions
corresponding to the elements of $D$.  In this case we see that
for any non-empty compact subset $C$ of $X$, the set $D$ is recovered
as the set $\{f(c) \, | \, f \in B, c \in C\}.$  Thus our preceding
calculations show that for any neighbourhood $U$ of zero in $W$,
the preimage of $T_{B,C,U}$ in $\Lin(V,W)$ is precisely equal
to $\{\phi \, | \, \phi(D) \subset U\}.$  Since sets of this type
form a neighbourhood basis of zero in $\Lin_b(V,W)$, we find that
the map $\Lin_b(V,W) \rightarrow \Lin_b(\Con(X,W),\Con(X,W))$ is
a topological embedding, as claimed.
\qed\enddemo

\proclaim{Definition 2.1.8}
If $V$ is a Hausdorff convex $K$-vector space, then let
$c_0(\N,V)$ denote the space of sequences in $V$ that converge
to zero, equipped with the (locally convex Hausdorff) topology
of uniform convergence on finite sets.
\endproclaim

If $\hat{\N} = \N \bigcup \{ \infty \}$
denotes the one-point compactification
of the discrete set $\N$, then we may regard $c_0(\N,V)$
as a closed subspace of 
$\Con(\hat{\N},V)$ consisting of those functions
that vanish at infinity.  

\proclaim{Definition 2.1.9} If $V$ is a $K$-Banach space
and $\X$ is an affinoid rigid analytic space defined over $L$, 
we define the $K$-Banach space $\An(\X,V)$
of $V$-valued rigid analytic functions
on $\X$ to be the completed tensor product
$\An(\X,K)\cotimes V.$  
\endproclaim

The formation of $\An(\X,V)$ is evidently covariantly functorial
in $V$ and contravariantly functorial in $\X$.

\proclaim{Proposition 2.1.10} If 
$\X$ is an affinoid rigid analytic space defined
over $L$, then there is an isomorphism $\An(\X,K) \iso
c_0(\N,K).$
\endproclaim
\demo{Proof}
This is standard.
\qed\enddemo

\proclaim{Definition 2.1.11}
If $V$ is a Hausdorff locally convex topological
$K$-vector space and $\X$ is an affinoid rigid analytic
space defined over $L$, we define the locally convex
space $\An(\X,V)$ of $V$-valued rigid
analytic functions on $\G$ to be the
locally convex inductive limit of Banach spaces
$$\An(\X,V) := \ilim{W} \An(\X,\overline{W}),$$
where $W$ runs over the directed
set of all $BH$-subspaces $W$ of $V$.
\endproclaim

The formation of
$\An(\X,V)$ is covariantly functorial in $V$ (use proposition~1.1.7)
and contravariantly functorial in $\X$.
Note that if $V$ is a $K$-Banach space then definitions~2.1.9
and~2.1.11 yield naturally isomorphic objects, since $V$ then forms
a final object in the directed set of all $BH$-subspaces of $V$.

If $W$ is a $BH$-subspace of $V$, then the continuous injection
$\overline{W} \rightarrow V$ induces a natural map
$\An(\X,\overline{W}) = \An(\X,K) \cotimes_K \overline{W}
\rightarrow \An(\X,K)\cotimes_{K,\pi} V.$

Passing to the inductive limit
in $W$ we obtain a natural map
$$\An(\X,V) \rightarrow \An(\X,K) \cotimes_{K,\pi} V.\tag 2.1.12 $$

\proclaim{Proposition 2.1.13}
Let $V$ be a Hausdorff convex space, and let $\X$ be
an affinoid rigid analytic space over $L$.

(i)
The map~(2.1.12) is injective.

(ii)
If $V$ is furthermore a Fr\'echet space,
then the map~(2.1.12) is a topological isomorphism.
\endproclaim
\demo{Proof}
Let $\hat{V}$ denote the completion of $V$.
If $W$ is a $BH$-subspace of $V$, then the injection
$\overline{W} \rightarrow V \rightarrow \hat{V}$ yields an
injection $c_0({\N}, \overline{W}) \rightarrow c_0({\N},\hat{V})$.
Passing to the locally convex inductive limit over all $BH$-subspaces
$W$ of $V$, we obtain an injection
$$\ilim{W} c_0({\N}, \overline{W}) \rightarrow c_0({\N},\hat{V}).
\tag 2.1.14$$

If we choose an isomorphism
$\An(\X,K) \iso c_0(\N,K)$, as in proposition~2.1.10,
then the isomorphism~(2.1.3) allows us to rewrite this injection
as an injection 
$$\ilim{W} \An(\X,K)\cotimes_K \overline{W} \rightarrow
\An(\X,K)\cotimes_{K,\pi} V.$$
This is the map~(2.1.12), and so we have proved part~(i)
of the proposition.

If $V$ is a Fr\'echet space,
then proposition~2.1.6 shows that the map~(2.1.14)
is a topological isomorphism,
and hence that~(2.1.12) is a topological isomorphism in this case.
This proves part~(ii) of the proposition.
\qed\enddemo

Note that the preceding result 
shows that $\An(\X,V)$ is Hausdorff, 
since it injects continuously into the
Hausdorff space $\An(\X,K) \cotimes_{K,\pi} V.$

Since the proof identifies $\An(\X,V)$ with $c_0(\N,V),$ when $V$
is either Fr\'echet or of compact type,
we also see that if $V$ is Fr\'echet 
then $\An(\X,V)$ is again Fr\'echet.  On the other hand,
if $V$ is of $LB$-type, then proposition~1.1.10, together with
the definition of $\An(\X,V)$ as an inductive limit,
shows that $\An(\X,V)$ is an $LB$-space.

We require the notion of a relatively compact morphism of rigid analytic
spaces.  Although this can be defined more generally \cite{\BGR, \S 9.6.2},
the following special case of this definition will suffice for our
purposes.

\proclaim{Definition 2.1.15}  We say that a morphism $\X \rightarrow
\Y$ of affinoid rigid analytic varieties over $L$ is relatively compact if 
it fits into a commutative diagram
$$\xymatrix{ \X \ar[r] \ar[d]  & \Y \ar[d] \\
\B_n(r_1) \ar[r] & \B_n(r_2),}$$
for some $r_1 < r_2$, in which the right-hand vertical arrow is a closed
embedding.
(Here for any $n$ and $r$ we let $\B_n(r):= \{ (x_1,\ldots,x_n) \,
| \, |x_i| \leq r \text{ for all } i = 1,\ldots, n\}$ denote
the $n$-dimensional closed ball of radius $r$ centred at the origin of
$\A^n$.)
\endproclaim

\proclaim{Proposition 2.1.16} If $\X \rightarrow \Y$ is a relatively compact
morphism of affinoid rigid analytic spaces over $L$,
then the induced map $\An(\Y, K) \rightarrow \An(\X,K)$ is compact.
\endproclaim
\demo{Proof}
The diagram whose existence is guaranteed by definition~2.1.15 yields
a diagram of morphisms of Banach spaces
$$\xymatrix{ \An(\B(r_2),K) \ar[r]\ar[d] & \An(\B(r_1),K) \ar[d] \\
\An(\Y,K) \ar[r] & \An(\X,K),}$$
in which the left-hand vertical arrow is surjective, and hence open.
Thus to prove the proposition, it suffices to show that the upper
morphism is compact.  This is standard, and is in any case
easily checked.
\qed\enddemo

We will need to consider topological vector-space valued analytic
functions on rigid analytic spaces that are not necessarily
affinoid.  The class of spaces encompassed by the following
definition is suitable for our purposes.

\proclaim{Definition 2.1.17} Let $\X$ be a rigid analytic space
defined over $L$.

(i) We say that $\X$ is
$\sigma$-affinoid if there is an increasing sequence
$\X_1 \subset \X_2 \subset \cdots \subset \X_n \subset \cdots$
of admissible affinoid open subsets of $\X$ such that
$\X = \bigcup_{n=1}^{\infty} \X_n$.

(ii) We say that $\X$ is strictly $\sigma$-affinoid if there
is an increasing sequence 
$\X_1 \subset \X_2 \subset \cdots \subset \X_n \subset \cdots$
of admissible affinoid open subsets of $\X$ such that
each of these inclusions is relatively compact, and such that
$\X = \bigcup_{n=1}^{\infty} \X_n$.
\endproclaim

The basic example of a $\sigma$-affinoid rigid analytic space
that we have in mind is that of an open polydisk, which may
be written as a union of an increasing sequence of closed polydisks.
Of course, any open polydisk is even strictly $\sigma$-affinoid.

If $\X = \bigcup_{n=1}^{\infty} \X_n$ is $\sigma$-affinoid, then
any admissible affinoid open subset $\Bbb Y$ of $\X$ is contained
in $\X_n$ for some $n$.  Thus the sequence
of admissible affinoid open subsets $\{\X_n\}$ of $\X$ is cofinal
in the directed set of all admissible affinoid open subsets of $\X$.

\proclaim{Definition 2.1.18}
Suppose that $\X$ is a $\sigma$-affinoid rigid analytic space defined
over $L$. If $V$ is a Hausdorff locally analytic convex $K$-vector space,
then we define the convex $K$-vector space $\An(\X,V)$
of analytic $V$-valued functions
on $\X$ to be the projective limit
$\plim{\Y} \An(\Y,V),$ where $\Y$ runs over all admissible affinoid
open subsets of $\X$.
\endproclaim

Since $\An(\X,V)$ is defined as the projective limit of Hausdorff
locally convex spaces, it is again Hausdorff locally convex.
The formation of $\An(\X,V)$ is covariantly functorial in $V$
and contravariantly functorial in $\X$.
The remarks preceding definition~2.1.18 show that the projective limit
that appears in this definition
may be taken over any increasing sequence $\{\X_n\}$ which witnesses the
fact that $\X$ is $\sigma$-affinoid.

If $V$ is a Fr\'echet space, then
each $\An(\X_n,V)$ is a Fr\'echet space
(by proposition~2.1.13~(ii)),
and thus $\An(\X,V)$ is a Fr\'echet space, being the projective
limit of a sequence of Fr\'echet spaces.  In particular,
$\An(\X,K)$ is a nuclear Fr\'echet space. 
If $\X$ is strictly $\sigma$-affinoid,
then proposition~2.1.16 shows that $\An(\X,K)$ is the projective
limit of a sequence of Banach spaces equipped with compact transition
maps.  Thus if $\X$ is strictly $\sigma$-affinoid, then
$\An(\X,K)$ is a nuclear Fr\'echet space.

\proclaim{Proposition 2.1.19} Let $\X$ be a $\sigma$-affinoid rigid
analytic space over $L$.  If $V$ is a $K$-Fr\'echet space,
then there is a natural isomorphism $\An(\X,V) \iso
\An(\X,K) \cotimes_K V.$ 
\endproclaim
\demo{Proof}
Write $X = \bigcup_{n\geq 1} \X_n,$ where $\X_1 \subset \cdots
\X_n \subset \X_{n+1} \cdots $ is an increasing sequence of 
admissible open affinoid subsets of $\X$.
Definition 2.1.18, together with proposition~2.1.13~(iii),
yields an isomorphism
$$\An(\X,V) = \plim{n} \An(\X_n,V) \iso \plim{n} (\An(\X_n,K)\cotimes_K V).$$
The claimed isomorphism is now provided by proposition~1.1.29.
\qed\enddemo

We now discuss the evaluation of rigid analytic functions on $\X$
at the $L$-valued points of $\X$.

\proclaim{Proposition 2.1.20} If $\X$ is a $\sigma$-affinoid rigid
analytic space over $L$, 
if $V$ is any Hausdorff convex $K$-vector space, and if we write $X:= \X(L)$,
then the evaluation at the points of $X$ yields
a natural continuous map $\An(\X,V) \rightarrow \Con(X,V)$.
If $X$ is Zariski dense in $\X$ then this map is injective.
\endproclaim
\demo{Proof}
We first suppose that $\X$ is affinoid, and
that $V$ is a Banach space.  Then by definition
$\An(\X,V) = \An(\X,K)\cotimes V.$  If we take into account the
natural isomorphism provided by~(2.1.3), then we that the continuous map
$$\An(\X,K) \rightarrow \Con(X,K) \tag 2.1.21$$
induced by evaluating rigid
analytic functions on points in $X$ yields the continuous
map $$\An(\X,V) \rightarrow \Con(X,V)\tag 2.1.22 $$
of the proposition.  The naturality of~(2.1.22) is evident.
Furthermore,
by definition $X$ is Zariski dense in $\X$ if and only if the map~(2.1.21)
is injective, in which case the map~(2.1.22) is also injective.  

Now consider the case where $\X$ is affinoid but $V$ is arbitrary.
The preceding paragraph gives a continuous map
$\An(\X,\overline{W}) \rightarrow \Con(X,\overline{W})$ for each $BH$-subspace
$W$ of $V$, which is injective if $X$ is Zariski dense in $\X$.
Composing this with the continuous injection $\Con(X,\overline{W})
\rightarrow \Con(X,V)$ yields a continuous map
$\An(\G,\overline{W}) \rightarrow \Con(G,V)$, which is injective
if $X$ is Zariski dense in $\X$.  Finally,
taking the inductive limit over all $BH$-subspaces $W$ yields
the required continuous map
$\An(\G,V) \rightarrow  \Con(G,V),$ which is injective
if $X$ is Zariski dense in $\X$.  Its naturality is clear.

Now suppose that $\X$ is $\sigma$-affinoid, and
write $\X = \bigcup_{n=1}^{\infty},$ with 
$\{\X_n\}$ an increasing sequence of admissible affinoid open subsets
of $\X$.  Set $X_n = \X_n(L),$ so that $X = \bigcup_{n=1}^{\infty}
X_n.$  Each $X_n$ is a compact open subset of $X$,
and the sequence $X_n$ is cofinal in the directed set of
compact subsets of $X$.
Thus $\Con(X,V) = \plim{n} \Con(X_n,V),$  and the map
$\An(\X,V) \rightarrow \Con(X,V)$ under consideration is the
projective limit of the maps $\An(\X_n,V) \rightarrow \Con(X_n,V)$.
The preceding paragraph shows that these maps
are continuous, and thus so is their projective limit.
If $X$ is Zariski dense in $\X$, then
$X_n = X \cap \X_n$ is Zariski dense in $\X_n$ for each sufficiently
large value of $n$,
and so for each such $n$, the
map $\An(\X_n,V) \rightarrow \Con(X_n,V)$ is injective.  
Passing to the projective limit over $n$, we find that the
map $\An(\X,V) \rightarrow \Con(X,V)$ is injective, as required.
\qed\enddemo

\proclaim{Proposition 2.1.23}
Let $V$ be a Hausdorff convex $K$-vector space and let $W$
be a closed subspace of $V$.  Let $\X$ be a $\sigma$-affinoid rigid
analytic space defined over $L$ for which the set $X:= \X(L)$ is
Zariski dense in $\X$.  Then the diagram
$$\xymatrix{\An(\X,W) \ar[r] \ar[d]^-{(2.1.20)} & \An(\X,V) \ar[d]^-{(2.1.20)}\\
\Con(X,W) \ar[r] & \Con(X,V) }$$ 
is Cartesian on the level of abstract vector spaces. 
If $V$ (and hence $W$) is a Fr\'echet space,
then the horizontal arrows of this diagram are closed immersions,
and the diagram
is Cartesian on the level of topological vector spaces.
\endproclaim
\demo{Proof} Let us suppose first that $V$ is a Fr\'echet space,
and that $\X$ is affinoid. 
We must show
that the natural map $\An(\X,W) \rightarrow \An(\X,V)$ identifies
$\An(\X,W)$ with the closed subspace of $\An(\X,V)$ consisting of
functions that are $W$-valued when evaluated at points of $X$.

Consider the exact sequence $0 \rightarrow W \rightarrow V \rightarrow V/W 
\rightarrow 0,$ in which the inclusion is a closed embedding
and the surjection is strict.
After taking the completed tensor product of this exact sequence
with $V$, and applying the isomorphisms provided
by proposition~2.1.13~(iii), we obtain an exact sequence
$$ 0 \rightarrow \An(\X,W) \rightarrow \An(\X,V) \rightarrow \An(\X,V/W)
\rightarrow 0$$ having the same properties.  Thus to prove the proposition,
it suffices to observe that (since $X$ is Zariski dense in $\X$)
an element of $\An(\X,V/W)$ that vanishes at every point of 
$X$ is necessarily the zero element.

Now suppose that $V$ is an arbitrary Hausdorff convex $K$-vector
space (and continue to suppose that $\X$ is affinoid).
If $U$ runs through the members of the directed set
of $BH$-subspaces of $V$ then $U \cap W$ runs through the
members of the directed set of $BH$-subspaces of $W$.  Applying
the result already proved, we find that the diagram
$$\xymatrix{\An(\X,\overline{U\cap W}) \ar[r] \ar[d] &
\An(\X,\overline{U}) \ar[d]\\
\Con(X,\overline{U\cap W}) \ar[r] & \Con(X,\overline{U}) }$$ 
is Cartesian.  The diagram
$$\xymatrix{\Con(X,\overline{U\cap W}) \ar[r]\ar[d] & \Con(X,\overline{U}) 
\ar[d] \\ \Con(X,W) \ar[r] & \Con(X,V)}$$
is also evidently Cartesian,
and so we conclude that the diagram
$$\xymatrix{\An(\X,\overline{U\cap W}) \ar[r] \ar[d] &
\An(\X,\overline{U}) \ar[d] \\
\Con(X,W) \ar[r] & \Con(X,V)}$$
is Cartesian.
Passing to the inductive limit over $U$ yields the proposition.
(Note that because we pass to a locally convex inductive limit,
we cannot conclude that the diagram of the proposition is Cartesian
as a diagram of topological vector spaces.)

Finally, if we assume that $\X$ is $\sigma$-affinoid,
and write $\X = \bigcup \X_n$, where $\{\X_n\}$ is an increasing
sequence of admissible affinoid open subspaces, then
the diagram
$$\xymatrix{\An(\X,W) \ar[r] \ar[d] & \An(\X,V) \ar[d]\\
\Con(X,W) \ar[r] & \Con(X,V) }$$ 
is obtained as a projective limit over $n$ of the diagrams
$$\xymatrix{\An(\X_n,W) \ar[r] \ar[d] & \An(\X_n,V) \ar[d]\\
\Con(X_n,W) \ar[r] & \Con(X_n,V) }.$$   (Here we write 
$X_n:= \X_n(L)$.)  Since each of these is Cartesian, by what
has already been proved, the same
is true of their projective limit.  Since furthermore the
projective limit of closed embeddings is a closed embedding,
we see that if $V$ and $W$ are Fr\'echet spaces,
then the horizontal
maps in this diagram are closed embeddings, and that in this
case the diagram is Cartesian even on the level of topological
$K$-vector spaces.
\qed\enddemo

\proclaim{Proposition 2.1.24} If $\X$ is a
$\sigma$-affinoid rigid analytic space
defined over $L$, and if each of $V$ and $W$ is a Fr\'echet space,
then the natural map
$$\Lin_b(V,W) \rightarrow \Lin_b(\An(\X,V),\An(\X,W))$$
(induced by the functoriality of the construction of $\An(\X,\text{--}\,)$)
is continuous.
\endproclaim
\demo{Proof}
We prove the proposition first in the case where $\X$ is affinoid.
If we fix an isomorphism $\An(\X,K) \iso c_0(\N,K)$, as we may do,
by proposition~2.1.12,
then~(2.1.3) and proposition~2.1.14~(ii) together yield
an isomorphism $\An(\X,K) \iso c_0(\N,V)$, functorial in $V$.
The proposition, in the case where $\X$ is affinoid,
now follows from proposition~2.1.9
(applied to $X = \hat{\N}$).

Suppose now that $X$ is $\sigma$-affinoid,
and write $\X = \bigcup_{n=1}^{\infty}$ as the union of
an increasing sequence of admissible open affinoid subspaces.
Since by definition $\An(\X,W) := \plim{n} \An(\X_n,W)$,
there is a natural isomorphism
$$\Lin_b(\An(\X,V),\An(\X,W)) \iso \plim{n} \Lin_b(\An(\X,W), \An(\X_n,W)).$$
Thus it suffices to show that, for each $n \geq 1$,
the natural map
$\Lin_b(V,W) \rightarrow \Lin_b(\An(\X,V),\An(\X_n,W))$
is continuous.  If we factor this map as
$\Lin_b(V,W) \rightarrow \Lin_b(\An(\X_n,V), \An(\X_n,W))
\rightarrow \Lin_b(\An(\X,V),\An(\X_n,W)),$
then its continuity is a consequence
of the result of the preceding paragraph.
\qed\enddemo

Suppose now that $X$ is a locally $L$-analytic manifold.  (We follow
\cite{\SCHTAN} in always assuming that such an $X$ is strictly paracompact --
see the discussion of \cite{\SCHTAN, p.~5}.)
Let $\{X_i\}_{i\in I}$ be a partition of $X$ into a disjoint
union of charts $X_i$; we refer to such a partition as an
analytic partition of $X$.  By definition each $X_i$ is equipped
with an identification $\phi_i:X_i \iso \X_i(L),$ where $\X_i$ is
a rigid analytic closed ball defined over $L$.  The collection
of all analytic partitions of $X$ form a directed set
(under the relation of refinement).

Let $\{X_i\}_{i \in I}$ be an analytic partition of $X$,
and let $\{X_j\}_{j \in J}$ be an analytic partition of $X$
that refines it.  Let $\sigma: J \rightarrow I$ be the map
describing the refinement, so that $\X_j \subset \X_{\sigma(j)}$
for each $j \in J$.  If each of these inclusions is a relatively compact
morphism, then we say that the refinement is relatively compact. 

The following definition is taken from \cite{\FETH, 2.1.10}.
Another account appears in \cite{\SCHTAN, p.~5}.

\proclaim{Definition 2.1.25} 
Let $V$ be a Hausdorff convex $K$-vector space and
$X$ be a locally $L$-analytic manifold.

We define the convex $K$-vector space
$\La(X,V)$ of locally analytic $V$-valued functions on
$X$ to be the locally convex inductive limit
$\ilim{\{X_i, W_i\}_{i\in I}}\prod_{i \in I} \An(\X_i,\overline{W}_i),$
taken over the directed set 
of collections of pairs $\{X_i, W_i\}_{i\in I}$, where
$\{X_i\}_{i\in I}$ is an analytic partition of $X$,
and each $W_i$ is a $BH$-subspace of $V$.
\endproclaim

The formation of $\La(X,V)$ is covariantly functorial
in $V$ (by proposition~1.1.7) and contravariantly functorial in $X$.

If $\{X_i\}_{i \in I}$ is any
partition of $X$ into disjoint open subsets, then the natural map
$\La(X,V) \rightarrow \prod_{i\in I} \La(X_i,V)$ is an isomorphism
\cite{\FETH, 2.2.4}.

If $X$ is compact, then any analytic partition of
$X$ is finite, and so we may replace the variable $BH$-subspaces
$W_i$ by a $BH$-subspace $W$ containing each of them.
Thus taking first the inductive limit in $W$,
we obtain a natural isomorphism
$\La(X,V) \iso \ilim{\{X_i\}_{i\in I}}\prod_i \An(\X_i,V),$
where the inductive limit is taken over the directed set of
analytic partitions $\{X_i\}_{i\in I}$ of $X$.
Alternatively, we may first take the inductive limit with respect
to the analytic partitions of $X$, and so obtain a natural isomorphism
$\ilim{W} \La(X,\overline{W}) \iso \La(X,V)$, where
the inductive limit is taken over all $BH$-subspaces $W$ of $V$.

Maintaining the assumption that $X$ is compact, 
choose a sequence
$(\{X_i\}_{i \in I_n})_{n\geq 1}$ of analytic partitions
of $X$ that is cofinal in the directed set of all such analytic partitions,
and such that for each $n\geq 1,$ the partition
$\{X_i\}_{i \in I_{n+1}}$ is a relatively compact refinement of the partition
$\{X_i\}_{i \in I_n}$.
Proposition~2.1.16 then shows that the natural map
$\prod_{i \in I_n} \An(\X_i,K) \rightarrow \prod_{i \in I_{n+1}}
\An(\X_i,K)$ is a compact map.  Passing to the locally convex
inductive limit, we find that $\La(\X,K)$ is of compact type.

If $X$ is arbitrary, then we see from the preceding paragraph
that $\La(\X,K)$ is a product of spaces of compact type.
Thus $\La(\X,K)$ is reflexive (hence barrelled) and complete,
being a product of reflexive and complete spaces.

\proclaim{Proposition 2.1.26} If $V$ is a Hausdorff convex $K$-vector space and
$X$ is a locally $L$-analytic manifold,
then evaluation at points of $X$
induces a continuous injection $\La(X,V) \rightarrow \Con(X,V)$,
that is natural, in the sense that it is compatible with the functorial
properties of its source and target.
Furthermore, this injection has dense image.
\endproclaim
\demo{Proof} If $\{X_i\}_{i\in I}$ is a analytic partition
of $X$ then proposition~2.1.20
yields a continuous injection $\An(\X_i,V) \rightarrow \Con(X_i,V),$
and so a continuous injection
$$\prod_{i\in I} \An(\X_i,V) \rightarrow \prod_{i\in I} \Con(X_i,V)
\iso \Con(X,V).$$
Passing to the inductive limit yields the required map of the proposition.
Its naturality is clear.

To see that the image is dense, note that locally constant functions
are contained in $\La(X,V)$, and are dense in $\Con(X,V)$.
\qed\enddemo

This result shows that $\La(X,V)$ is Hausdorff, since it admits a
continuous injection into the Hausdorff space $\Con(X,V)$.

\proclaim{Proposition 2.1.27} If $X$ is a locally analytic $L$-manifold,
if
$V$ is a Hausdorff convex $K$-vector space, and if
$W$ is a closed subspace of $V$,
then the diagram
$$\xymatrix{\La(X,W) \ar[r] \ar[d]^-{(2.1.26)} & \La(X,V) \ar[d]^-{(2.1.26)}\\
\Con(X,W) \ar[r] & \Con(X,V) }$$ 
is Cartesian on the level of abstract $K$-vector spaces. 
\endproclaim
\demo{Proof}
Let $\{X_i\}_{i\in I}$ be an analytic partition of
$X$.  For each $i \in I$ we see by proposition~2.1.23 that
the diagram
$$\xymatrix{\An(\X_i,W) \ar[r] \ar[d] & \An(\X_i,V) \ar[d] \\
\Con(X_i,W) \ar[r] & \Con(X_i,V)}$$
is Cartesian as a diagram of abstract $K$-vector spaces.
Taking the product over all $i\in I,$ taking into account the
isomorphism $\Con(X,W) \iso \prod_{i\in I} \Con(X_i,W),$ and
passing to the locally convex inductive limit over all analytic
partitions of $X$, we deduce the proposition.
\qed\enddemo

\proclaim{Proposition 2.1.28}
If $V$ is a Hausdorff locally convex topological $K$-vector space,
and if $X$ is a compact locally $L$-analytic manifold, then there
is a continuous injection 
$\La(X,V) \rightarrow \La(X,K) \cotimes_{K,\pi} V.$
If $V$ is of compact type,
then it is even a topological isomorphism.  In particular, in this case
$\La(X,V)$ is again a space of compact type.
\endproclaim
\demo{Proof}
If $\{X_i\}_{i \in I}$ is an analytic partition of $X$ (necessarily
finite, since $X$ is compact), then~(2.1.12) yields a natural map
$$\multline \prod_{i\in I} \An(\X_i,V) \rightarrow
\prod_{i\in I} (\An(\X_i,K) \cotimes_{K,\pi} V) \\
\iso (\prod_{i \in I} \An(\X_i,V)) \cotimes_{K,\pi} V \rightarrow
\La(X,K) \cotimes_{K,\pi} V. \endmultline $$
Passing to the locally convex inductive limit over all
analytic partitions, we obtain a map
$\La(X,V) \rightarrow \La(X,K) \cotimes_{K,\pi} V.$
To see that it is injective, we may compose it with the natural
map $\La(X,K) \cotimes_{K,\pi} V \rightarrow \Con(X,K) \cotimes_{K,\pi}
V \iso \Con(X,\hat{V})$ (where $\hat{V}$ denotes the completion of
$V$, and the isomorphism is provided by~(2.1.3)).  The map
that we obtain
is then seen to be equal to the composite $\La(X,V) \rightarrow \Con(X,V)
\rightarrow \Con(X,\hat{V})$, where the first arrow is the injection
of~2.1.26, and the second arrow is induced by the embedding $V \rightarrow
\hat{V}$.  Since this latter composite is injective, we include that
the map $\La(X,V) \rightarrow \La(X,K) \cotimes_{K,\pi} V$ that
we have constructed is also injective.

Suppose now that $V$ is of compact type, and write $V = \ilim{n} V_n$
as the inductive limit of a sequence of Banach spaces, with injective
and compact transition maps.  Choose a cofinal sequence
$(\{X_i\}_{i\in I_m})_{m\geq 1}$ of
analytic partitions of $X$, such
that each partition provides a relatively compact refinement of the partition
that precedes it.
We may compute $\La(X,V)$ as the locally convex inductive limit
$\ilim{\{X_i\}_{i\in I_m}, n} (\prod_{i\in I_m} \An(\X_i,K)) \cotimes_K V_n.$
Part~(i) of proposition~1.1.32 now shows that the injection
$\La(X,V) \rightarrow \La(X,K) \cotimes_K V$ is a topological isomorphism,
and that $\La(X,K) \cotimes_K V$ is of compact type, as claimed.
\qed\enddemo

The final statement of the preceding proposition is a restatement
of \cite{\SCHTAN, lem.~2.1}.
The discussion that follows proposition~2.1.7 below
shows that if $V$ is complete and of $LB$-type, then
the map of the preceding proposition is bijective.

\proclaim{Corollary 2.1.29} If $W \rightarrow V$ is a closed
embedding of Hausdorff convex $K$-vector space of compact type,
and if $X$ is a compact locally $L$-analytic manifold,
then the natural injection $\La(X,W) \rightarrow \La(X,V)$ is a closed
embedding.
\endproclaim
\demo{Proof}
Propositions~2.1.27 and~2.1.28 together imply that the source
and image of this injection are both spaces of compact type.
By theorem~1.1.17, it is necessarily a closed embedding.
\qed\enddemo

\proclaim{Proposition 2.1.30} If $V$ is a Hausdorff convex space
of $LB$-type,
and if $X$ is a compact locally $L$-analytic manifold,
then $\La(X,V)$ is an $LB$-space. 
\endproclaim
\demo{Proof}
Write $V = \ilim{n} V_n$
as the inductive limit of a sequence of Banach spaces, with injective
and compact transition maps.  Choose a cofinal sequence
$(\{X_i\}_{i\in I_m})_{m\geq 1}$ of
analytic partitions of $X$.
We may compute $\La(X,V)$ as the locally convex inductive limit
$\ilim{\{X_i\}_{i\in I_m}, n} (\prod_{i\in I_m} \An(\X_i,K)) \cotimes_K V_n.$
Thus it is an $LB$-space.
\qed\enddemo

\proclaim{Proposition 2.1.31} If $X$ is a compact locally analytic $L$-manifold
and if $V$ and $W$ are two Hausdorff convex $K$-vector spaces of compact
type, then the natural map
$\Lin_b(V,W) \rightarrow \Lin_b(\La(X,V),\La(X,W)),$ induced
by the functoriality of the construction of $\La(X,\text{--}\,)$,
is a continuous map of convex $K$-vector spaces.
\endproclaim
\demo{Proof}
Since $V$ and $W$ are of compact type, proposition~2.1.28
shows that the same is true of $\La(X,V)$ and $\La(X,W)$,
and that furthermore, the latter spaces are isomorphic to
$\La(X,K) \cotimes_{K} V$ and $\La(X,K)\cotimes_{K} V$ respectively.
From \cite{\SCHNA, prop.~20.9} we deduce isomorphisms
$V'_b\cotimes_{K,\pi} W \iso \Lin_b(V,W)$ and
$$(\La(X,K) \cotimes_K V)'_b \cotimes_{K,\pi} \La(X,K) \cotimes_{K,\pi} W
\iso \Lin_b(\La(X,V),\La(X,W)).$$
Taking into account the isomorphism
$(\La(X,K) \cotimes_K V)'_b \iso \La(X,K)'_b \cotimes_{K,\pi} V'_b$
yielded by proposition~1.1.32~(ii),
we find that the natural map 
$\Lin_b(V,W) \rightarrow \Lin_b(\La(X,V),\La(X,W))$
that is under consideration may be rewritten as a map
$$\multline V'_b\cotimes_{K,\pi} W \rightarrow
\La(X,K)'_b\cotimes_{K,\pi} V'_b \cotimes_{K,\pi} \La(X,K) \cotimes_{K,\pi} W
\\
\iso \La(X,K)'_b\cotimes_{K,\pi} \La(X,K) \cotimes_{K,\pi} V'_b
\cotimes_{K,\pi} W.\endmultline \tag 2.1.32$$
This is immediately seen to be obtained from the map
$$K \rightarrow \La(X,K)'_b\cotimes_{K,\pi} \La(X,K)\tag 2.1.33$$
(obtained by taking $V = W = K$ in~(2.1.32),
which sends $1 \in K$ to the element
of the tensor product corresponding to the identity map on $\La(X,K)$)
by tensoring through with $V'_b\cotimes_{K,\pi} W$, and then completing.
Since~(2.1.33) is certainly continuous, we conclude that
the same is true of~(2.1.32).
\qed\enddemo

\section{2.2}
In this subsection we recall the various types of 
distributions that are relevant to non-archimedean function theory.

\proclaim{Definition 2.2.1} If $X$ is a topological
space, then we let $\DCon(X,K)$ denote the dual space to the 
convex $K$-vector space $\Con(X,K)$.  This is the space of $K$-valued measures on
$X$.
\endproclaim

\proclaim{Definition 2.2.2} If $\X$ is a $\sigma$-affinoid rigid analytic 
space defined over $L$, then
we let $\DAn(\X,K)$ denote the dual space to the convex $K$-vector space
$\An(\X,K)$.  This is the space of $K$-valued analytic
distributions on $\X$.
\endproclaim

\proclaim{Definition 2.2.3} If $X$ is a locally $L$-analytic manifold,
then we let $\DLa(X,K)$ denote the dual space to the convex space
of compact type $\La(X,K)$.  This is the space of $K$-valued
locally analytic distributions on $X$.
\endproclaim

In these definitions we have not specified any particular topology
on the dual spaces under consideration.  As with the dual space
to any convex space, they admit various locally convex topologies.
Frequently we will endow these spaces with
their strong topologies, in which case we add a subscript `$b$' to
emphasise this.  If $X$ is
a compact topological space $X$ then $\DCon(X,K)_b$ is a Banach space,
if $\X$ is an affinoid rigid analytic space over $L$ then $\DAn(\X,K)_b$
is also a Banach space, if $\X$ is a strictly $\sigma$-affinoid rigid
analytic space then
$\DAn(\X,K)_b$ is of compact type, 
if $X$ is a locally $L$-analytic space then $\DLa(X,K)_b$ is 
reflexive, and if $X$ is furthermore compact
then $\DLa(X,K)_b$ is a nuclear Fr\'echet space.  (The third and
fifth claims follow from the fact that
the dual of a nuclear Fr\'echet space is a space of compact type,
and conversely. The fourth follows from the fact that $\DLa(X,K)$
is dual to the reflexive space $\La(X,K)$.)

If $\mu$ is an element of $\DCon(X,K)$ and $f$ is an element
of $\Con(X,K)$, then it is sometimes suggestive to write
$\int_X f(x) d \mu(x)$ to indicate the evaluation of the functional
$\mu$ on $f$. 
Similar notation can be used for the evaluation of elements of
$\DAn(\X,K)$ (respectively $\DLa(X,K)$)
on elements of $\An(\X,K)$
(respectively $\La(X,K)$) for $\sigma$-affinoid
rigid analytic spaces $\X$ (respectively
locally analytic spaces $X$) over $L$.

If $X$ is a topological space
and $x$ is an element of $X$ then the functional $\ev_x$
defines an element of $\DCon(X,K),$
which we denote (as usual) by $\delta_x$.
If we let $K[X]$ denote the vector space spanned on $X$ over $K$
then the map $x \mapsto \delta_x$ gives a map
$K[X] \rightarrow \DCon(X,K)$, which is easily seen to be an embedding
(since a continuous function on $X$ can assume arbitrary values
at a finite collection of points of $X$).  The image of this
map is weakly dense in $\DCon(X,K)$ (since a continuous function
that vanishes at every point of $X$ vanishes), but is typically
not strongly dense (since $\Con(X,K)$ is typically not reflexive).

Similarly, if $\X$ is a rigid analytic space over $L$,
if $X = \X(L),$ and if $x$ is an element of $X,$ then $\ev_x$
defines an element $\delta_x$ of $\DAn(\X,K)$.  If $X$ is Zariski
dense in $\X$ then the elements $\delta_x$ are weakly dense
in $\DAn(\X,K)$ (since by assumption a rigid analytic function
which vanishes at every point of $X$ vanishes on $\X$), but are
typically not strongly dense in $\DAn(\X,K)$ (since $\An(X,K)$ is
typically not reflexive).
However,
if $\X$ is as strictly $\sigma$-affinoid, then the nuclear Fr\'echet space
$\An(\X,K)$ is reflexive, and so the elements $\delta_x \in
\DAn(\X,K)$ corresponding to the functionals $\ev_x$
are both weakly and strongly dense in $\DAn(\X,K)$.
Similarly,
if $X$ is a locally $L$-analytic manifold then $\La(X,K)$ is
reflexive, and so the elements $\delta_x$
defined by the functionals $\ev_x$ are both
weakly and strongly dense in $\DLa(X,K)$.
(This is \cite{\SCHTAN, lem.~3.1} .)

\proclaim{Lemma 2.2.4} If $X$ is a topological space,
then the map $X \rightarrow \DCon(X,K)_s$, given
by $x \mapsto \delta_x$, is continuous.
\endproclaim
\demo{Proof}
Let $\xi$ denote the topology on $X$ induced by this
injection.  By definition of the weak topology,
for any $x \in X$, a basis of $\xi$-neighbourhoods of 
$x$ is provided by sets of the form $\{ y \in X \, | \,
|f(x) - f(y)| < \epsilon \}$, where $f$ is an element
of $\Con(X,K)$, and $\epsilon$ a positive real number.
Thus the topology $\xi$ is coarser than
the given topology on $X$, and the proposition follows.  
\qed\enddemo

If $X$ is a locally $L$-analytic manifold, then
the continuous injection $\La(X,K) \rightarrow \Con(X,K)$ of proposition~2.1.26
induces a continuous map
$$\DCon(X,K)_b \rightarrow \DLa(X,K)_b.\tag 2.2.5$$
Also, the association of $\delta_x$ to $x$ induces an injective map
$$X \rightarrow \DLa(X,K)_b.\tag 2.2.6$$

\proclaim{Proposition 2.2.7}  Let $X$ be a locally
$L$-analytic manifold.

(i) The map~(2.2.5) is injective,
with dense image.  

(ii) The map~(2.2.6) is a topological embedding.
\endproclaim
\demo{Proof}
Since the map of proposition~2.1.26
has dense image, the map~(2.2.5) is injective.
Since the map of proposition~2.1.26 is injective,
and since $\La(X,K)$ is reflexive,
the map~(2.2.5) induces an injection on dual spaces, and so has dense
image.  This proves~(i).

We turn to proving~(ii). 
Let $\xi$ denote the given topology
on $X$, and let $\xi'$ denote the topology induced on $X$ by
regarding it as a subspace of $\DLa(X,K)$ via the map~(2.2.6).
Any locally analytic function $f\in \La(X,K)$
induces a continuous functional on $\DLa(X,K)$, which thus restricts to a
$\xi'$-continuous function on $X$.  From the definition of~(2.2.6), this
is exactly the function $f$ again.  Thus we see in particular that
any $\xi$-locally constant function on $X$ is $\xi'$-continuous.
Since a basis of $\xi$-open subsets of $X$ can be cut out via
locally constant functions, we see that the topology $\xi'$ is
finer than $\xi.$

Since $\La(X,K)$ is reflexive, any
closed and bounded subset of $\La(X,K)$ is c-compact.  Also, $\Con(X,K)$
is complete and barrelled (the latter by corollary~2.1.5).
The hypotheses of proposition~1.1.38 are
thus satisfied by the map of proposition~2.1.26,
and we infer from that proposition that
the transpose map~(2.2.5) factors as $\DCon(X,K)_b \rightarrow \DCon(X,K)_{bs}
\rightarrow \DLa(X,K)_b$.  We will show that the map $x \mapsto \delta_x$
induces a continuous map $X \rightarrow \DCon(X,K)_{bs},$ and thus
that~(2.2.6) is continuous.  This will imply that $\xi$ is finer
than $\xi'$, and so prove~(ii).

Let $X = \coprod_i X_i$ be a partition of $X$ into a disjoint union
of compact open subsets.  
It suffices to show that the map
$X_i \rightarrow \DCon(X,K),$ given by $x \mapsto \delta_x$,
is continuous for each $i$.
The image of $X_i$ in
$\DCon(X,K)$ 
is bounded as a subset of $\DCon(X,K)_b$ (since $X_i$ is
compact), and so by definition of the bounded-weak topology,
the topology on this image induced by $\DCon(X,K)_{bs}$
is coarser than the weak topology.  Thus it suffices to show
that the map $X_i \rightarrow \DCon(X,K)_s$ is continuous.
This map factors as $X_i \rightarrow \DCon(X_i,K)_s
\rightarrow \DCon(X,K)_s$, and so by lemma~2.2.4 
we are done. 
\qed
\enddemo 

If $V$ is any Hausdorff convex $K$-vector space, then by composing
elements of $\Lin(\DLa(X,K)_b, V)$ with
the map~(2.2.6), we obtain 
a $K$-linear map
$$\Lin(\DLa(X,K)_b,V) \rightarrow \Con(X,V).\tag 2.2.8$$

\proclaim{Proposition~2.2.9} The map~(2.2.8) is continuous,
when the source is given its strong topology.
\endproclaim
\demo{Proof}
If we let $\{X_i\}_{i \in I}$ denote a partition of $X$
into compact open subsets, then the isomorphisms
$\bigoplus_{i \in I} \DLa(X_i,K)_b \iso \DLa(X,K)_b$
and $\prod_{i \in I} \Con(X_i,V) \iso \Con(X,V)$
allow us to write~(2.2.8) as the product of the maps
$\Lin_b(\DLa(X_i,K)_b,V) \rightarrow \Con(X_i,V).$
Hence we may restrict our attention to the case where $X$ is
compact.  If $\hat{V}$ denotes the completion of $V$,
we may embed~(2.2.8) into the commutative diagram
$$\xymatrix{\Lin_b(\DLa(X,K)_b,V) \ar[r] \ar[d] &  \Con(X,V) \ar[d] \\
\Lin_b(\DLa(X,K)_b,\hat{V}) \ar[r] &  \Con(X,\hat{V}),}$$
in which the vertical arrows are embeddings.  Thus we may
assume in addition that $V$ is complete.
Taking into account \cite{\SCHNA, cor.~18.8} and the isomorphism~(2.1.3),
we may rewrite~(2.2.8) (in the case where $X$ is compact and $V$
is complete) as the completed tensor product with $V$ of the natural
map $\La(X,K) \rightarrow \Con(X,K)$.  Thus it is continuous, as claimed.
\qed\enddemo

\proclaim{Proposition~2.2.10} There is a continuous
linear map $$\La(X,V) \rightarrow \Lin_b(\DLa(X,K)_b,V),$$
uniquely determined by the fact that its composite with~(2.2.8)
is the map of proposition~2.1.26.
If furthermore $V$ is of $LB$-type, then this map is a bijection.
\endproclaim
\demo{Proof}
The required map is constructed in
\cite{\SCHTAN, thm.~2.2} and the discussion that follows
the proof of that theorem, and by that theorem is shown to be
bijective if $V$ is of $LB$-type.  We will give another description
of this map, which shows that it is continuous.
Let $\{X_i\}_{i\in I}$
be an analytic partition of $X$, and let $\{W_i\}_{i\in I}$ be
a collection of $BH$-subspaces of $V$.
For each $i \in I$, the natural map
$\An(\X_i,K) \rightarrow \La(X_i,K)$ and the map
$\overline{W}_i \rightarrow V$ together induce a continuous map
$$\multline \An(\X_i,\overline{W}_i) = \An(\X_i,K)\cotimes_K \overline{W}_i
\rightarrow \La(X_i,K)\cotimes_K \overline{W}_i  \\
\iso \Lin_b(\DLa(X_i,K),\overline{W}_i)
\rightarrow \Lin_b(\DLa(X_i,K),V)\endmultline$$
(the isomorphism following from \cite{\SCHNA, cor.~18.8}).
Taking the product of these maps over all $i \in I$,
we obtain a continuous map
$$\multline \prod_i \An(\X_i,\overline{W}_i) \rightarrow 
\prod_i \Lin_b(\DLa(X_i,K),V)  \\ \iso \Lin_b(\bigoplus_{i\in I}
\DLa(X_i,K),V) \iso \Lin_b(\DLa(X,K),V).\endmultline$$
Passing to the direct limit over all such collections of pairs
$\{(X_i,W_i)\}_{i \in I}$, we obtain a continuous map
$$\La(X,K) \rightarrow \Lin_b(\DLa(X,K),V).\tag 2.2.12$$
If $f \in \La(X,K)$, then its image 
$\phi \in \Lin(\DLa(X,K)_b,V)$ is characterised by
the property $\phi(\delta_x) = f(x)$, and so we see that~(2.2.12)
coincides with the map constructed in \cite{\SCHTAN},
and that the composite of this map with~(2.2.8) is
the map of proposition~2.1.26, as claimed.
\qed\enddemo

The preceding result relates to proposition~2.1.28, in
the following way.
Suppose that $V$ is complete.  Then, taking into account the isomorphism
of \cite{\SCHNA, cor.~18.8}, the map of proposition~2.2.10 may be rewritten
as a map
$$\La(X,V) \rightarrow \La(X,K) \cotimes_{K,\pi} V.$$
This is easily checked to be equal to the map of proposition~2.1.28.
Thus if $V$ is complete and of $LB$-type, 
then we see that the map of proposition~2.1.28 is a bijection.

\section{2.3} In this brief subsection we discuss two
``change of field'' functors.

If $F$ is a finite extension of $K$, then the valuation on $K$
extends in a unique fashion to a spherically complete
non-archimedean valuation on $F$.  We can now define
an ``extension of scalars'' functor as follows.
If $V$ is any locally convex $K$-vector space $V$, then
the tensor product $F\otimes_K V$ is naturally a locally
convex $F$-vector space.  Furthermore, it is Hausdorff, normable,
metrisable, complete, barrelled, bornological, an $LF$-space,
and $LB$-space, or of compact type, if and only if the same is true
of $V$.   (If $F$ is an arbitrary extension of $K$,
spherically complete with respect to a non-archimedean
valuation extending that on $K$, then we can defined
an extension of scalars functor by sending $V$ to
$F\otimes_{K,\pi} V.$  However, this functor is not
so easily studied if $F$ is of infinite degree over $K$.)

If $V$ is any Hausdorff locally convex $K$-vector space,
then given any topological space $X$, any $\sigma$-affinoid
rigid analytic space $\X$ over $L$, or any locally $L$-analytic
space $X$, there are natural isomorphisms
$F\otimes_K \Con(X,V) \iso \Con(X, F\otimes_K V),$
$F\otimes_K \An(\X,V) \iso \An(\X, F\otimes_K V),$
and $F\otimes_K \La(X,V) \iso \La(X,F\otimes_K V)$.

Suppose now that $E$ is a closed subfield of $L$
(so that there are inclusions $\Q_p \subset E \subset L$).
If $\X$ is a rigid analytic space over $E$, then we let
$\X_{/L}$ denote the corresponding rigid analytic
space over $L$, obtained by extending scalars.
If $A$ is an affinoid $E$-Banach algebra,
then $(\MaxSp A)_{/L} = \MaxSp L\otimes_E A.$
Extension of scalars defines a functor from the
category of rigid analytic spaces over $E$
to the category of rigid analytic spaces
over $L$.  
This functor admits a right adjoint,
referred to as ``restriction of scalars''.
If $\X$ is a rigid analytic space over $L$, then
we let $\Res^{L}_{E} \X$ denote the rigid analytic space over $E$
obtained by restriction of scalars.

The right adjointness property of restriction of scalars implies
that for any rigid analytic space $\X$ over $L$,
there is a natural map
$$(\Res^{L}_{E} \X)_{/L} \rightarrow \X.\tag 2.3.1 $$
If $L$ is a a Galois extension of $E$,
then there is furthermore a natural isomorphism
$(\Res^{L}_{E} \X)_{/L} \iso \prod_{\tau} \X^{(\tau)},$
where $\tau$ ranges over the elements of the Galois 
group of $L$ over $E$, and $\X^{(\tau)}$ denotes
the rigid analytic space over $L$ obtained by twisting
$\X$ via the automorphism $\tau: L \iso L$.
(The map~(2.3.1) is obtained by projecting onto
the factor corresponding to the identity automorphism.)

We may similarly define an extension of scalars
functor from the category of locally $E$-analytic
spaces to the category of locally $L$-analytic
spaces, and a corresponding right adjoint
restriction of scalars functor.
Note that if $\X$ is a rigid analytic space over $L$,
then $\Res^L_E \X (E) = \X(L)$, and hence that
if $X$ is a locally $L$-analytic space, then
$\Res^L_E X$ has the same underlying topological space
as $X$.  Thus restriction of scalars may equally
well be regarded as the forgetful functor,
``regard a locally $L$-analytic space as a locally
$E$-analytic space''.

If $X$ is a locally $L$-analytic space, and we first restrict
and then extend scalars, then there is a natural map
$$(\Res^{L}_{E} X)_{/L} \rightarrow X.\tag 2.3.2 $$

\head 3. Continuous, analytic, and locally analytic vectors \endhead

\section{3.1}  In this subsection we study the left and
right regular actions of a group on the function spaces
associated to it.  We begin with a lemma.

\proclaim{Lemma 3.1.1} Let $V$ be a topological $K$-vector
space and $G$ a topological group, and consider a map
$$G \times V \rightarrow V\tag 3.1.2$$
such that each $g \in G$ acts as a linear transformation
of $V$.  Then the map~(3.1.2) is continuous if and only if it
satisfies the following conditions:

(i) for each $v\in V$ the map $g \mapsto g v$
from $G$ to $V$ induced
by~(3.1.2) is continuous at the identity $e \in G$;

(ii) for each $g \in G$ the
linear transformation $v \mapsto g v$ of $V$ induced by~(3.1.2)
is continuous (that is, the $G$-action on $V$ is topological);

(iii) the map~(3.1.2) is continuous at the point $(e,0)$ of $G\times V$.
\endproclaim
\demo{Proof}
It is clear that these conditions are necessary, and we
must show that these conditions are sufficient.
We begin by showing that if $v$ is an element of $V$,
then the map~(3.1.2)
is continuous at the point $(e,v) \in G \times V$. 

Let $M$ be a neighbourhood of zero in $V$, and choose a neighbourhood
$M'$ of zero such that $M'+M' \subset M$. By~(iii),
we may find a neighbourhood $H$ of $e$ in $G$ and a neighbourhood
$M''$ of zero in $V$ such that $g M'' \subset M'$ if
$g \in H$.
Since the
map $g \mapsto g v$ is continuous at the element~$e$,
by~(i), we may 
also find a neighbourhood $H'$ of $e$ contained in $H$
such that $H' v \subset v+M'$.  Thus we see that the image of
$H' \times (v+M'')$ is contained in $v+M$, proving that~(3.1.2)
is continuous at the point $(e,v) \in G \times V$.

Now let $(g,v)$ be an arbitrary element of $G \times V$, and
let $M$ be neighbourhood of zero in $V$.  Replacing $v$ by $g v$
in the discussion of the preceding paragraph, we may find an open
subset $H'$ of $e\in G$ and a neighbourhood $M''$ of zero in
$V$ such that the image of 
$H'\times (g v + M'')$ is contained in $g v + M.$ Since the 
action of $g$ on $V$ is continuous, by~(ii), we may
find a neighbourhood $M'''$ of zero such that
$g M''' \subset M''$.  
Now $H' g $ is a neighbourhood of $g$ and
the image of $H' g \times (v + M''')$ under~(3.1.2)
is contained in $g v + M$.
This proves that~(3.1.2) is continuous
at the point $(g,v) \in G \times V$.
\qed\enddemo

Note that condition~(iii) of the preceding lemma holds if a
neighbourhood of $e$ in $G$ acts equicontinuously on $V$.
Suppose, for example, that $G$ is a locally compact abelian
group, and $V$ is a barrelled convex $K$-vector space
equipped with a $G$-action for which the action map~(3.1.2) is
separately continuous.  Any relatively compact neighbourhood
of $e$ then acts in a pointwise bounded fashion on $V$,
and hence (since $V$ is barrelled) acts equicontinuously on $V$.
Thus all three conditions of lemma~3.1.1 are satisfied,
and~(3.1.2) is jointly continuous.

\proclaim{Corollary 3.1.3} If $V$ is a convex $K$-vector space equipped
with a topological action of the topological group $G$, and if
$H$ is an open subgroup of $G$ such that the resulting $H$-action on
$V$ is continuous, then the $G$-action on $V$ is continuous.
\endproclaim
\demo{Proof}
By assumption, condition~(ii) of lemma~3.1.1 is satisfied
for $G$, and conditions~(i) and~(iii) are satisfied with $H$ in
place of $G$.  But these latter two conditions depend only on
an arbitrarily small neighbourhood of the identity $e\in G$,
and so they hold for $G$ since they hold for $H$.  Lemma~3.1.1
now implies that the $G$-action on $V$ is continuous.
\qed\enddemo

\proclaim{Lemma 3.1.4} 
If $V$ is a topological $K$-vector space equipped with a
continuous action of a topological group $G$, then any compact
subset of $G$ acts equicontinuously on $V$.
\endproclaim
\demo{Proof}
Let $C$ be a compact subset of $G$,
and let $M$ be a neighbourhood of zero in $V$.
If $g \in C$, then there exists a neighbourhood $U_g$ of $g$ in $G$,
and a neighbourhood $M_g$ of zero in $V$, such that $U_g M_g \subset M$.
A finite number of the $U_g$ suffice to cover to compact set $C$.
If we let $M'$ denote the intersection of the corresponding finite
number of neighbourhoods $M_g$, we obtain a neighbourhood of zero
with the property that $C M' \subset M$.  This proves the lemma.
\qed\enddemo

Suppose now that $G$ is an abstract group and that $V$
is an abstract $K$-vector space.  The left and right regular
actions of $G$ on itself are the maps $G \times G \rightarrow G$
defined as $(g_1, g_2) \mapsto g_2^{-1} g_1$ and $(g_1,g_2) \mapsto
g_1 g_2$ respectively.  (We regard these as right actions of $G$
on itself.)  They commute one with the other.
Via functoriality of the formation of the $K$-vector
space of functions $\Fun(G,V)$,
we obtain (left) actions of $G$ on the space $\Fun(G,V)$, which
we again refer to as the left and right regular actions.
These actions also commute with one another.

If $G$ is a topological group (respectively
the $L$-valued points of a $\sigma$-affinoid rigid analytic group $\G$
defined over $L$, respectively a locally $L$-analytic group) and
if $V$ is a Hausdorff locally convex topological $K$-vector space, then
we similarly may consider the left and right regular representations
of $G$ on $\Con(G,V)$ (respectively on $\An(\G,V)$,
respectively on $\La(G,V)$).

\proclaim{Proposition 3.1.5} If $G$ is a locally compact topological
group and $V$ is a Hausdorff locally convex $K$-vector space,
then the left and right regular actions
of $G$ on $\Con(G,V)$ are both continuous.
\endproclaim
\demo{Proof} 
We will show that the right regular action is continuous.
The proof of the analogous result for the left regular action
proceeds along identical lines.

If $C$ is a compact subset of $G$ and
$M$ is a neighbourhood of zero in $V$,
let $U_{C,M}$ denote the subset of $\Con(G,V)$ consisting
of functions $\phi$ such that $\phi(C) \subset M$.
As $C$ runs over all compact subsets of $G$ and
$M$ runs over the neighbourhoods of zero in $V$,
the spaces $U_{C,M}$ form a basis of neighbourhoods of zero in $\Con(G,V)$.
If $H$ is a compact open subset of $G$ then
$C H$ is a compact subset of $G,$ and the right regular action
of $H$ on $\Con(G,V)$ takes $U_{ C H,M}$ to $U_{C,M}$.
Thus we see that $H$ acts equicontinuously on $\Con(G,V)$.
Choosing $H$ so that it contains a neighbourhood of $e \in G$,
we see that conditions~(ii) and~(iii) of lemma~3.1.1 hold.
It remains to be shown that condition~(i) of that lemma holds.

Let $\phi$ be an element of $\Con(G,V)$, $C$ a compact subset of $G$,
and $M$ a neighbourhood of zero in $V$.
We will find a neighbourhood
$U$ of the identity in $G$ such that $H (\phi + U_{C H,M}) \subset U_{C,M}$.
This will show that the right regular action is continuous at the point
$(e,\phi)$ of $G \times \Con(G,V)$, verifying condition~(i).

Since $\phi$ is continuous and $C$ is compact,
we may find an open neighbourhood $U$ of $e$ such that 
$\phi(g) - \phi(g h) \in M$ for any $g\in C$ and $h \in U$.
(This is a variation on the fact that
for a function on a compact space, continuity implies uniform continuity.)
Thus $h \phi \subset \phi + U_{C,M}$ for every $h\in H$, and so
$H (\phi + U_{C H,M}) \subset \phi + U_{C,M},$ as required.
\qed\enddemo

\proclaim{Proposition 3.1.6} If $V$ is a Hausdorff convex $K$-vector space,
$\G$ is a $\sigma$-affinoid rigid analytic group defined over $L$,
and $G = \G(L),$
then the left and right regular actions of $G$ on $\An(\G,V)$
are continuous.
\endproclaim
\demo{Proof} If we write $\G = \bigcup_{n=1}^{\infty},$
so that $G = \bigcup_{n=1}^{\infty} G_n$,
with $G_n := \G_n(L)$ open in $G$ for each $n \geq 1$,
then $\An(\G,V) = \plim{n} \An(\G_n,V)$, and so it suffices to
show that $G_n$ acts continuously on $\An(\G_n,V)$ for each
$n \geq 1$.  This reduces us to the case where $\G$ is affinoid.
We suppose this to be the case for the remainder of the proof.

Since $\An(\G,V)$ is barrelled, the discussion following lemma~3.1.1
shows that it suffices to show that the $G$-action on $\An(\G,V)$
is separately continuous.
The definition of $\An(\G,V)$ as an inductive limit then
allows us to reduce first to the case where $V$ is a Banach
space.   Finally, the definition of $\An(\G,V)$ for a Banach
space $V$ allows us to reduce to the case where $V = K$. 
The proof then proceeds along the same lines as that of
the preceding proposition.  (In place of the fact that
a continuous function on a compact space is uniformly continuous,
one should use the fact that an analytic function on an affinoid
space is Lipschitz \cite{\BGR, prop.~7.2.1/1}.)
\qed\enddemo

\section{3.2}
Now suppose that $V$ is an abstract $K$-vector space
equipped with an action of an abstract group $G$.  Each element $v$ of $V$ gives
rise to a function $\orbit_v \in \Fun(G,V),$ defined by
$\orbit_v(g) = g v$; we refer to $\orbit_v$ as the orbit map
of the vector $v$.  Formation of the orbit map yields
a $K$-linear embedding $\orbit: V \rightarrow \Fun(G,V)$, which
is a section to the map $\ev_e:\Fun(G,V) \rightarrow V$.

\proclaim{Lemma 3.2.1} The map $\orbit:V\rightarrow \Fun(G,V)$
is $G$-equivariant when $\Fun(G,V)$ is equipped with
the right regular $G$-action.  Thus
the restriction of $\ev_e$ to $\orbit(V)$
is a $G$-equivariant isomorphism when $\orbit(V)$ is endowed
with the right regular $G$-action.
\endproclaim
\demo{Proof} This lemma is an immediate consequence of the fact
that $\orbit$ is a section to $\ev_e,$ together with 
the formula $\orbit_v(g g') = \orbit_{g' v}(g)$.
\qed\enddemo 

There is a convenient characterisation of the image of $\orbit$.
First note that, in addition to the left and right regular actions of
$G$ on $\Fun(G,V)$, we obtain another action of $G$ induced
by its action on $V$ and the functoriality of $\Fun(G,V)$ in $V$, which
we refer to as the pointwise action of $G$. This commutes
with both the left and right regular actions.  Since the
left and right regular actions also commute one with the
other, we obtain an action of $G\times G \times G$ on
$\Fun(G,V),$ the first factor acting via the pointwise action,
the second via the left regular action, and the third via
the right regular action.  Explicitly, if $(g_1,g_2,g_3) \in G\times
G \times G$ and $\phi \in \Fun(V,G),$ then
$$((g_1,g_2,g_3) \cdot \phi) (g') = g_1 \cdot \phi(g_2^{-1} g' g_3).
\tag 3.2.2 $$

Let $\Delta_{1,2}: G \rightarrow G\times G \times G$ denote
the map $g \mapsto (g,g,1)$.  Then $\Delta_{1,2}$ induces
a $G$-action on $\Con(G,V)$ which commutes with the right regular
action of $G$.  Explicitly, this action is defined by
$$(g \cdot \phi)(g')= g\cdot \phi(g^{-1} g'). \tag 3.2.3 $$

We first explain how to untwist the $\Delta_{1,2}(G)$-action
on $\Fun(G,V)$.

\proclaim{Lemma 3.2.4} If we define a $K$-linear map
$\Fun(G,V) \rightarrow \Fun(G,V)$ by sending a
function $\phi$ to the function $\tilde{\phi}$ defined by
$\tilde{\phi}(g) = g^{-1}\cdot\phi(g),$ then this 
map is an isomorphism.  It is furthermore $G$-equivariant,
provided we equip its source with the $\Delta_{1,2}(G)$-action
and its target with the left regular $G$-action.
\endproclaim
\demo{Proof}
This is easily checked by the reader.
\qed\enddemo

We observe for later that we could also consider the $\Delta_{1,3}(G)$-action
on $\Fun(G,V),$ and that this action is made isomorphic to the right
regular $G$-action on $\Fun(G,V)$ via the isomorphism
$\phi \mapsto (\tilde{\phi}: g \mapsto g\phi(g)).$
For now, we require the following corollary.

\proclaim{Corollary 3.2.5} The image of $\orbit$ is equal to
the space $\Fun(G,V)^{\Delta_{1,2}(G)}$
of $\Delta_{1,2}(G)$-fixed vectors in $\Fun(G,V)$.
\endproclaim
\demo{Proof} An element $\phi$ of $\Fun(G,V)$ is fixed under
the action~(3.2.3) if and only if the element
$\tilde{\phi}$ that it is sent to under the isomorphism of lemma~3.2.4
is fixed under the left regular action of $G$; equivalently,
if and only if $\tilde{\phi}$ is constant.  The explicit
formula relating $\phi$ and $\tilde{\phi}$ shows that $\tilde{\phi}$
is constant if and only if $\phi(g) = g \phi(e) = \orbit_{\phi(e)}(g)$
for all $g\in G$.
\qed\enddemo

If $V$ is any $K$-vector space (not necessarily equipped with
a $G$-action), we
may consider the double function space
$\Fun(G,\Fun(G,V)).$ If we endow $\Fun(G,V)$ with its
right regular $G$-action, then the preceding constructions apply,
with $\Fun(G,V)$ in place of $V$.
There is a natural isomorphism
$\imath: \Fun(G\times G, V) \iso
\Fun(G,\Fun(G,V)),$ 
defined by $(\imath \phi)(g_1)(g_2) = \phi(g_1,g_2),$
and we will reinterpret the preceding constructions
in terms of this isomorphism.

Let $\alpha:G\times G \rightarrow G$ denote the map
$\alpha: (g_1,g_2) \mapsto g_2g_1$ and $\beta:G \rightarrow G\times G$
denote the map $\beta: g \mapsto (e,g)$.  These maps induce maps
$\alpha^*: \Fun(G,V) \rightarrow \Fun(G\times G,V)$ and 
$\beta^*:\Fun(G\times G,V) \rightarrow \Fun(G,V),$  and we obtain
the following commutative diagrams:
$$\xymatrix{\Fun(G,V) \ar[d]^{\alpha^*} \ar[dr]^{\orbit} & \\
\Fun(G\times G,V) \ar[r]^{\imath} & \Fun(G,\Fun(G,V))}\tag 3.2.6$$
and
$$\xymatrix{\Fun(G\times G,V) \ar[r]^{\imath} \ar[rd]^{\beta^*} &
\Fun(G,\Fun(G,V)) \ar[d]^{\ev_e} \\ & \Fun(G,V).}\tag 3.2.7$$
(The commutativity of each diagram is easily checked by the reader.)

Suppose now that $G$ is a locally compact topological group
and that $V$ is a Hausdorff convex $K$-vector space equipped with a topological
$G$-action.
We may apply the above considerations with $\Con(G,V)$ in place
of $\Fun(G,V)$.

\proclaim{Definition 3.2.8} We define the convex $K$-vector space $V_{\con}$
of continuous vectors in $V$ to be the closed subspace
of $\Delta_{1,2}(G)$-invariant elements of $\Con(G,V)$,
equipped with the continuous $G$-action induced by
the right regular action of $G$ on $\Con(G,V)$.
\endproclaim

The formation of $V_{\con}$ is evidently covariantly
functorial in $V$ and contravariantly functorial in $G$.

\proclaim{Proposition 3.2.9} The restriction of
the map $\ev_e$ to $V_{\con}$ is a continuous $G$-equivariant injection
of $V_{\con}$ into $V$.
Its image consists of the subspace of vectors $v$ in $V$ for which
the orbit map $\orbit_v$ is continuous.
\endproclaim
\demo{Proof}
Since the map $\ev_e:\Con(G,V) \rightarrow V$ is continuous, so is
its restriction to the closed subspace $V_{\con}$ of $\Con(G,V)$.
It is $G$-equivariant by lemma~3.2.1.
By corollary~3.2.5, an element $\phi\in \Con(G,V)$ is fixed under the action
of $\Delta_{1,2}(G)$ if and only if it is of the form $\orbit_v$ for
some $v \in V$. Thus
an element $v$ of $V$ lies in $\ev_e(V_{\con})$ if and only if
$\orbit_v$ is continuous.
\qed\enddemo

Now suppose that the $G$-action on $V$ is continuous, so that
$V_{\con}$ maps continuously and bijectively onto $V$.
The map $v\mapsto \orbit_v$ then yields an embedding
$\orbit: V \rightarrow \Con(G,V)$ of $K$-vector
spaces, which is a section to the map $\ev_e:\Con(G,V) \rightarrow V$.

\proclaim{Proposition 3.2.10} If $V$ is a Hausdorff convex $K$-vector
space, equipped with a continuous $G$-action,
then the map $\orbit$ is continuous
and $G$-equivariant, when
$\Con(G,V)$ is equipped with the right regular $G$-action,
and restricts to a topological isomorphism of $V$ onto the
closed subspace $V_{\con}$ of $\Con(G,V)$.
\endproclaim
\demo{Proof}
Lemma~3.2.1 and corollary~3.2.5 taken together imply that
the map $\orbit$ is $G$-equivariant, with image equal
to $V_{\con}$.  To show that
$\orbit$ is continuous, it suffices to show that 
for each compact subset $C$ of $G$ and
each neighbourhood $M$ of zero in $V$ 
there is a neighbourhood
$M'$ of zero in $V$ such that for all $v \in M'$, the function $\orbit_v$ takes
$C$ into $M$.
In other words, regarding $C$ as fixed for a moment,
for every neighbourhood $M$ of zero,
there should be a neighbourhood $M'$ of zero such that $C M'
\subset M$.  This is precisely the condition that the action of
$C$ be equicontinuous. Since $C$ is compact, lemma~3.1.4 shows
that it does act continuously.  Thus the map $\orbit$ is indeed continuous,
and provides a
continuous inverse to the the continuous map $V_{\con} \rightarrow V$
induced by $\ev_e$.  
\qed\enddemo

\proclaim{Lemma 3.2.11}
The isomorphism of lemma~3.2.4 induces a topological isomorphism of
$\Con(G,V)$ onto itself.
\endproclaim
\demo{Proof}
We begin by noting that since the $G$-action on $V$ is assumed
to be continuous, the function $\tilde{\phi}: G \rightarrow V$
is continuous if the function $\phi$ is.  (Here we are using the
notation introduced in the statement of lemma~3.2.4.)

Fix a compact subset $C$ of $G$ and an open neighbourhood $M$
of zero in $V$.  Since $C$ acts equicontinuously on $V$,
by lemma~3.1.4, we may find an open neighbourhood $M'$ of zero such
that $C M' \subset M.$  Then
the map $\phi \mapsto \tilde{\phi}$ takes $\{ \phi \, | \,
\phi(C) \subset M'\}$ into the set $\{\phi \, | \, \phi(C) \subset M\}.$
Thus $\phi \mapsto \tilde{\phi}$ is continuous.  
It is equally straightforward to check that its inverse is continuous,
and thus the lemma is proved.
\qed\enddemo

Similarly, one may show that
the map $\phi \mapsto (\tilde{\phi}: g \mapsto
g\phi(g))$ induces an isomorphism between $\Con(G,V)$ and $\Con(G,V)$.
(As observed following the proof of lemma~3.2.4, this map intertwines
the $\Delta_{1,3}(G)$-action and the right regular action on its
source and target.)

The following results will be required in subsequent subsections.

\proclaim{Lemma 3.2.12} Let $G$ be a group
and $V$ be a Hausdorff convex $K$-vector space equipped with
a topological $G$-action.  If $W$ is an $FH$-subspace of $V$
that is invariant under the $G$-action on $V$ then the topological $G$-action
on $W$ induces a topological $G$-action on $\overline{W}$.
\endproclaim
\demo{Proof} This follows from proposition~1.1.2~(ii).
\qed
\enddemo

\proclaim{Definition 3.2.13} If $G$ is a topological group
and $V$ is a Hausdorff convex $K$-vector space equipped
with a topological $G$-action then we say that a $G$-invariant $FH$-subspace
$W$ of $V$ is continuously invariant if the induced $G$-action
on $\overline{W}$ is continuous.
\endproclaim

\proclaim{Proposition 3.2.14} Let $G$ be a compact
topological group and $V$ be a Hausdorff convex $K$-vector space equipped
with a topological $G$-action.
If $W$ is an $FH$-subspace of $V$
then among all the continuously $G$-invariant $FH$-subspaces of $W$
is one that contains all the others, 
which we denote by $\core{G}{W}$.
If $W$ is a $BH$-subspace of $V$, then so is $\core{G}{W}$.
\endproclaim
\demo{Proof}
The continuous injection $\overline{W} \rightarrow V$ induces a
continuous injection $\Con(G,\overline{W}) \rightarrow \Con(G,V)$.
Define $\overline{U}$ to be the preimage under this map of the
closed subspace $V_{\con}$ of $\Con(G,V)$.  Then $\overline{U}$
is a closed subspace of $\Con(G,\overline{W}),$ and thus a Fr\'echet
space.  By construction it is invariant under the right regular
action of $G$ on $\Con(G,\overline{W}),$ and thus endowed with
a continuous $G$-action.  We have continuous injections
$\overline{U} \rightarrow V_{\con} \rightarrow V$; let $U$ denote
the image of $\overline{U}$ under the composite of these maps.

By construction, $U$ is a continuously $G$-invariant $FH$-subspace
of $W$.  Furthermore, $U$ is equal to the set of vectors
$w\in W$ for which $G w \subset W$ and $\orbit_w: G \rightarrow \overline{W}$
is continuous, and so $U$ contains any continuously $G$-invariant
$FH$-subspace of $W$.  Thus we set $\core{G}{W} := U$.
To complete the proof of the proposition,
note that if $\overline{W}$ is in fact a Banach space, then
so also is $\overline{U}$, and so $U$ is a $BH$-subspace of $V$.
\qed\enddemo 

The following result is also useful.

\proclaim{Proposition 3.2.15}
If $G$ is a compact topological group and $V$ is a semi-complete
Hausdorff convex $K$-vector space equipped with a continuous
$G$-action, then any $BH$-subspace of $V$ is contained in a
$G$-invariant $BH$-subspace.
\endproclaim
\demo{Proof}
The discussion preceding proposition~1.1.12 shows that
it suffices to prove that any bounded subset of $V$
may be embedded in $G$-invariant bounded subset of $V$.
Equivalently, if $A$ is a bounded subset of $V$, we must
show that the set $G A$ of all $G$-translates
of elements of $A$ is bounded.  This follows from lemma~3.1.4,
which shows that $G$ acts equicontinuously.
\qed\enddemo

\section{3.3}
Let $\G$ be an affinoid rigid analytic group, and suppose that
the group of $L$-rational points $G := \G(L)$ is Zariski dense 
in $\G$.  Let $V$ be a Hausdorff convex $K$-vector space equipped with a
topological $G$-action.  Our goal in this subsection
is to define the convex space of
$\G$-analytic vectors in $V$.

We begin by supposing that $V$ is a $K$-Banach space
equipped with a topological action of $G$.  We can repeat the
discussion that begins subsection~3.2 in the context
of $\An(\G,V)$.
In particular, the functoriality of $\An(\G,V)$ in $V$ induces a
topological action of $G$ on $\An(\G,V)$ that
we refer to as the pointwise action of $G$, which commutes
with both the left and right regular actions. 
Since the
left and right regular actions also commute one with the
other, we obtain an action of $G\times G \times G$ on
$\An(\G,V),$ the first factor acting via the pointwise action,
the second via the left regular action, and the third via
the right regular action.   
As in subsection~2.1 we let
$\Delta_{1,2}: G \rightarrow G\times G \times G$ denote
the map $g \mapsto (g,g,1)$.  Then $\Delta_{1,2}$ induces
a topological $G$-action on $\An(\G,V)$ that commutes with the right regular
action of $G$.

\proclaim{Definition 3.3.1} If $V$ is a $K$-Banach space equipped
with a topological action of $G$, then define the Banach space $V_{\G-\an}$
of $\G$-analytic vectors
in $V$ to be the closed subspace of $\Delta_{1,2}(G)$-invariant elements
of $\An(\G,V)$, equipped with the right regular action of $G$. 
\endproclaim

The formation of $V_{\G-\an}$ is clearly covariantly functorial
in $V$ and contravariantly functorial in $\G$.  Note that the 
$G$-action on $V_{\G-\an}$ is continuous, since the right
regular $G$-action on $\An(\G,V)$ is so (by proposition~3.1.6).

Proposition~2.1.20 yields a continuous injection
$$\An(\G,V) \rightarrow \Con(G,V)\tag 3.3.2$$
that is compatible with the
$G\times G\times G$-action on source and target.

\proclaim{Proposition 3.3.3} The subspace $V_{\G-\an}$ of 
$\An(\G,V)$ is equal to the preimage under~(3.3.2) of the closed
subspace $V_{\con}$ of $\Con(G,V)$.  Evaluation at $e$ yields
a continuous injection $V_{\G-\an} \rightarrow V$.
The image of $V_{\G-\an}$ in $V$ consists of the subspace
of those vectors $v$ for which the orbit map $\orbit_v$ is 
(the restriction to $G$ of) an element of $\An(\G,V)$.
\endproclaim
\demo{Proof} This is an immediate consequence of the definitions.
\qed\enddemo

\proclaim{Proposition 3.3.4} If the Banach space
$\An(\G,V)$ is equipped with its right regular $G$-action,
then the natural map $\An(\G,V)_{\G-\an} \rightarrow \An(\G,V)$
of proposition~3.3.3 is an isomorphism.
\endproclaim
\demo{Proof} We consider the rigid analytic analogue
of the diagrams~(3.2.6) and~(3.2.7).
Namely, let $\alpha: \G \times \G \rightarrow \G$
denote the rigid analytic map $\alpha: (g_1,g_2) \mapsto g_2g_1,$
and let $\beta: \G \rightarrow \G\times \G$ denote the
rigid analytic map $\beta: g \mapsto (e,g).$
Then $\alpha$ is faithfully flat, $\beta$ is a closed
immersion, and the composite $\alpha\circ \beta$ is equal to the identity map
from $\G$ to itself.

Passing to spaces of rigid analytic functions, we obtain 
morphisms of Banach spaces
$\alpha^*: \An(\G,K) \rightarrow \An(\G\times \G,K)$ and 
$\beta^*: \An(\G\times \G) \rightarrow \An(\G,K),$
the first of which is a closed embedding, and the second of which is a
strict surjection.  There is a natural isomorphism $$\imath:
\An(\G\times \G,K) \iso \An(\G,K) \cotimes \An(\G,K) = \An(\G, \An(\G,K)).$$
Via this isomorphism, the action of $\Delta_{1,2}(G)$ on
$\An(\G,\An(\G,K))$ induces an action of $\Delta_{1,2}(G)$ on
$\An(\G\times \G,K),$  given explicitly by the formula
$$(g\cdot \phi)(g_1,g_2) = \phi(g^{-1} g_1, g_2 g).$$

As $G$ is dense in $\G$, we deduce that $\alpha^*$ identifies
$\An(\G,K)$ with the $\Delta_{1,2}(G)$-invariant subspace of
$\An(\G\times \G,K)$.  Also, since $\ev_e$ is equal to the composite 
$\beta^*\circ \imath^{-1},$ and since $\beta^*\circ \alpha^*$ is
the identity on $\An(\G,K)$, we conclude that
$\ev_e:\An(\G,K)_{\G-\an} = \An(\G, \An(\G,K))^{\Delta_{1,2}(G)}
\rightarrow \An(\G,K)$ is an isomorphism.
Tensoring this isomorphism through by $V$, and completing, we conclude
similarly that $\ev_e: \An(\G,V)_{\G-\an} \rightarrow \An(\G,V)$ is
an isomorphism, as required.
\qed\enddemo

\proclaim{Proposition 3.3.5} If $W \rightarrow V$ is a $G$-equivariant
closed embedding of $K$-Banach spaces equipped with topological $G$-actions
then the diagram
$$\xymatrix{W_{\G-\an} \ar[d]^-{(3.3.3)}\ar[r] & V_{\G-\an} \ar[d]^-{(3.3.3)} \\
W \ar[r] & V}$$ is Cartesian. In particular, the
induced map $W_{\G-\an} \rightarrow V_{\G-\an}$ is 
a closed embedding of Banach spaces.
\endproclaim
\demo{Proof} The closed embedding $W \rightarrow V$ induces a
closed embedding of Banach spaces
$\An(\G,W) \rightarrow \An(\G,V)$, which by proposition~2.1.23
identifies $\An(\G,W)$
with the subspace of functions in $\An(\G,V)$ that are $W$-valued.
Passing to $\Delta_{1,2}(G)$-invariants yields the proposition.
\qed
\enddemo

\proclaim{Corollary 3.3.6} If $V$ is a $K$-Banach space equipped with
a topological $G$-action then the natural map $(V_{\G-\an})_{\G-\an}
\rightarrow V_{\G-\an}$ is an isomorphism.
\endproclaim
\demo{Proof}
Proposition~3.3.4 shows that the natural map $\An(\G,V)_{\G-\an}
\rightarrow \An(\G,V)$ is an isomorphism.  Taking this
into account, the corollary follows by applying proposition~3.3.5 to
the closed embedding $V_{\G-\an} \rightarrow \An(\G,V)$.
\qed
\enddemo

\proclaim{Proposition 3.3.7} If $V$ is a $K$-Banach space equipped with
a topological $G$-action then the continuous injection
$\An(\G,V) \rightarrow \Con(G,V)$
of proposition~2.1.20 induces an isomorphism
$\An(\G,V) \iso \Con(G,V)_{\G-\an}.$  (The analytic vectors of the
target are computed with respect to the right regular $G$-action.)
\endproclaim
\demo{Proof}
The functoriality of the construction of the space of
analytic vectors yields
a map $\An(\G,V)_{\G-\an} \rightarrow \Con(G,V)_{\G-\an},$ which
when composed with the inverse of the isomorphism of proposition~3.3.4
yields a map $\An(\G,V) \rightarrow \Con(G,V)_{\G-\an}$.  We will
show this map to be an isomorphism.

Consider the commutative diagram
$$\xymatrix{\An(\G,V) \ar[dd] \ar[r]_-{\alpha^*}^-{\sim} &
\An(\G\times\G,V)^{\Delta_{1,2}(G)} \ar[r]^-{\sim} \ar[dd] &
(\An(\G,K)\cotimes \An(\G,V))^{\Delta_{1,2}(G)} \ar[d] \\
&& (\An(\G,K)\cotimes \Con(G,V))^{\Delta_{1,2}(G)} \ar[d] \\
\Con(G,V) \ar[r]_-{\alpha^*}^-{\sim} &
\Con(G\times G,V)^{\Delta_{1,2}(G)} \ar[r]^-{\sim} &
(\Con(G,K) \cotimes \Con(G,V))^{\Delta_{1,2}(G)} ,} $$
in which $\alpha^*$ denotes the map induced on function spaces
by the map $\alpha: (g_1,g_2) \mapsto g_2g_1,$  and all of whose
horizontal arrows are isomorphisms.
By definition, the space $\Con(G,V)_{\G-\an}$ is equal to the space
$(\An(\G,K)\cotimes \Con(G,V))^{\Delta_{1,2}(G,V)}$, and
this diagram shows that this space may be identified with
those $V$-valued functions on $G\times G$ that are 
of the form $(g_1,g_2) \mapsto f(g_2g_1)$ for some $f\in \Con(G,V),$
and that are analytic with respect to the variable $g_1$.  By
restricting such a function to $G\times e \subset G\times G,$
we see that $f$ must be an analytic function, and thus that
the map $\An(\G,V) \rightarrow \Con(G,V)_{\G-\an}$ is a continuous
bijection.  Since its source and target are Banach spaces,
the open mapping theorem shows that it is an isomorphism, as claimed.
\qed
\enddemo

\proclaim{Corollary 3.3.8} If $V$ is a $K$-Banach space equipped with a
topological $G$-action then the natural map
$(V_{\con})_{\G-\an} \rightarrow V_{\G-\an}$ is an isomorphism.
\endproclaim
\demo{Proof} By construction there is a $G$-equivariant closed embedding
$V_{\con} \rightarrow \Con(G,V),$ and so by propositions~3.3.5 and~3.3.7
a Cartesian diagram
$$\xymatrix{(V_{\con})_{\G-\an} \ar[d] \ar[r] & \An(\G,V) \ar[d] \\
V_{\con} \ar[r] & \Con(G,V).}$$
The definition of $V_{\con}$ then implies that
$(V_{\con})_{\G-\an}$ maps isomorphically onto the subspace of
$\Delta_{1,2}(G)$-invariants in $\An(\G,V)$; that is, $V_{\G-\an}$.
\qed\enddemo

The next two lemmas are included for later reference.  For both
lemmas, we assume given an open affinoid subgroup $\H$ of $\G$,
equal to the Zariski closure of its group $H:= \H(L)$ of
$L$-valued points.   If $V$ is a $K$-Banach space
equipped with a representation of $H$, then we may restrict the
left and right regular action of $G$ on $\An(\G,V)$ to actions
of $H$.  Combining these actions with the pointwise action of $H$
arising from the action of $H$ on $V$, we obtain an action
of $H\times H \times H$ on $\An(\G,V)$ (the first factor acting
via the pointwise action, the second factor by the left regular
representation, and the third factor by the right regular
representation). The formation of $\An(\G,V)$ with its
$H\times H \times H$-action is obviously functorial in $V$.

\proclaim{Lemma 3.3.9} If $V$ is a $K$-Banach space equipped
with a topological $H$-action, then the natural map
$\An(\G,V_{\H-\an})^{\Delta_{1,2}(H)} \rightarrow 
\An(\G,V)^{\Delta_{1,2}(H)}$ is an isomorphism.
\endproclaim
\demo{Proof}
Directly from the definitions, this is the map
$$\An(\G, \An(\H,V)^{\Delta_{1,2}(H)})^{\Delta_{1,2}(H)}
\rightarrow \An(\G,V)^{\Delta_{1,2}(H)} \tag 3.3.10 $$
induced by restricting the evaluation map
$\ev_e: \An(\H,V) \rightarrow V$ 
to $\An(\H,V)^{\Delta_{1,2}(H)}$.
There is a canonical isomorphism 
$\An(\G,\An(\H,V)) \iso \An(\G \times \H,V),$
which induces isomorphisms
$$\multline \An(\G,\An(\H,V)^{\Delta_{1,2}(H)})^{\Delta_{1,2}(H)} \\
\iso
\{f \in \An(\G \times \H, V) \, | \,
f(h_1g,h_2) = f(g, h_2 h_1) \text{ and } \\
f(g,h_1 h_2) = h_1 f(g, h_2) \text{ for all } g \in G, h_1,h_2 \in H\}
\endmultline$$
and
$$\An(\G,V)^{\Delta_{1,2}(H)} \iso
\{f \in \An(\G,V) \, | \, f(h g) = h f(g) \text{ for all } g \in G, 
h \in H \}. $$
With respect to these descriptions of its source and target,
we may reinterpret the map~(3.3.10) as the map 
$$\multline \{f \in \An(\G \times \H, V) \, | \,
f(h_1g,h_2) = f(g, h_2 h_1) \text{ and }
f(g,h_1 h_2) = h_1 f(g, h_2) \text{ for all } \\ g \in G, h_1,h_2 \in H\}
\rightarrow
\{f \in \An(\G,V) \, | \, f(h g) = h f(g) \text{ for all } g \in G, 
h \in H \},\endmultline $$
given by restricting to $\G \times e \subset \G \times \H$.  
Described this way, it is immediate that~(3.3.10) 
induces an isomorphism between its source and target.
\qed\enddemo

\proclaim{Lemma 3.3.11} The evaluation map $\ev_e: \An(\H,K) \rightarrow K$
induces an isomorphism $\An(\G, \An(\H,K))^{\Delta_{1,2}(H)}
\iso \An(\G,K)$.  
\endproclaim
\demo{Proof}
This is proved in an analogous manner to the previous lemma.
\qed\enddemo

\proclaim{Proposition~3.3.12} If $U$ and $V$ are two $K$-Banach
spaces, each equipped with a topological $G$-action,
then the map
$U_{\G-\an} \cotimes_{K} V_{\G-\an} \rightarrow U\cotimes_{K} V$
(obtained by taking the completed tensor product of the
natural continuous injections
$U_{\G-\an} \rightarrow U$ and $V_{\G-\an} \rightarrow V$)
factors through the natural continuous injection
$(U\cotimes_{K} V)_{\G-\an} \rightarrow U\cotimes_{K} V$
(where $U\cotimes_{K} V$
is equipped with the diagonal $G$-action).
\endproclaim
\demo{Proof}
If $\phi_1$ is an element of $\An(\G,U)$ and $\phi_2$
is an element of $\An(\G,V)$, then the map
$g \mapsto \phi_1(g) \otimes \phi_2(g)$
yields an element of $\An(\G,U\cotimes_K V)$.
Thus we obtain a continuous bilinear map
$\An(\G,U) \times \An(\G,V) \rightarrow \An(\G, U\cotimes_K V),$
which induces a continuous map
$\An(\G,U)\cotimes \An(\G,V) \rightarrow \An(\G,U\cotimes_K V).$
This restricts to a map
$$
\An(\G,U)^{\Delta_{1,2}(G)}
\cotimes_K \An(\G,V)^{\Delta_{1,2}(G)} 
\rightarrow \An(\G,U\cotimes_K V)^{\Delta_{1,2}(G)}$$
(where the $G$-action on $U\cotimes_K V$ is taken to be
the diagonal action).  Taking into account definition~3.3.1,
this proves the proposition.
\qed\enddemo

We now allow $V$ to be an arbitrary convex $K$-vector space equipped
with a topological action of $G$, and extend definition~3.3.1
to this more general context.

\proclaim{Definition 3.3.13} If $V$ is a locally convex Hausdorff
topological $K$-vector space equipped
with a topological action of $G$, define the convex $K$-vector space
$V_{\G-\an}$ of $\G$-analytic vectors in $V$ to be the locally convex inductive
limit $V_{\G-\an} = \ilim{W}\overline{W}_{\G-\an},$ where $W$ runs
over all the $G$-invariant $BH$-subspaces of $V$, and each
$\overline{W}$ is equipped with
the $G$-action provided by lemma~3.2.12. 
\endproclaim

The formation of $V_{\G-\an}$ (as a convex $K$-vector space with topological
$G$-action) is covariantly functorial in
$V$ (by proposition~1.1.7) and contravariantly functorial in $\G$.
Since $V_{\G-\an}$ is defined as the locally convex inductive
limit of a family of Banach spaces, it is both barrelled and bornological.

Proposition~3.3.3 provides for each $G$-invariant $BH$-subspace
$W$ of $V$ a continuous injection
$\overline{W}_{\G-\an} \rightarrow \overline{W}$, which we may compose
with the continuous injection $\overline{W} \rightarrow V$ to obtain
a continuous injection $$\overline{W}_{\G-\an} \rightarrow V.\tag 3.3.14$$
Taking the limit over all such $W$ then yields a continuous injection
$$V_{\G-\an} \rightarrow V .\tag 3.3.15.$$
As a particular consequence,
we conclude that $V_{\G-\an}$ is Hausdorff.

The following result justifies the designation of $V_{\G-\an}$ as
the space of $\G$-analytic vectors in $V$.

\proclaim{Theorem 3.3.16}
If $V$ is a Hausdorff convex $K$-vector
space equipped with a topological $G$-action then
there is a (uniquely determined) continuous
injection $V_{\G-\an} \rightarrow \An(\G,V)$ such that the diagram
$$\xymatrix{V_{\G-\an} \ar[r]^-{\text{(3.3.15)}}\ar[d] & V \ar[d]^-{\orbit} \\
\An(\G,V) \ar[r] & \Fun(G,V)}\tag 3.3.17$$ is Cartesian on the level of abstract
vector spaces.  Thus the image of the continuous injection~(3.3.15)
contains precisely those vectors $v\in V$ for which the orbit map
$\orbit_v$ is (the restriction to $G$ of) an element of $\An(\G,V)$. 
\endproclaim
\demo{Proof} 
Since the map~(3.3.15) and the natural map
$\An(\G,V) \rightarrow \Fun(G,V)$ are injective (the latter
since we have assumed that $G$ is Zariski dense in $\G$),
there is clearly at most one such map which makes~(3.3.17) 
commute.

If $W$ is a $G$-invariant $BH$-subspace of $V$, then
by construction there is a closed embedding
$\overline{W}_{\G-\an} \rightarrow \An(\G,\overline{W}).$
Composing with the injection $\An(\G,\overline{W}) \rightarrow
\An(\G,V)$, and then
passing to the inductive limit over all $\overline{W}$,
yields a continuous injection
$V_{\G-\an} \rightarrow \An(\G,V)$
that makes~(3.3.17) commute.
It remains to be shown that~(3.3.17) is Cartesian.

Let $V_1$ denote the subspace of $V$ consisting
of those vectors $v\in V$ for which $\orbit_v$ is given by an
element of $\An(\G,V)$.
We must show that the image of~(3.3.15) is equal to $V_1$.

The commutativity of~(3.3.17) shows that image of~(3.3.15) lies
in $V_1$.
Conversely, let $v$ be a vector in $V_1$.   The orbit map $\orbit_v:
G \rightarrow V$ 
is thus given by an element of $\An(\G,V),$ and so by
an element of $\An(\G,\overline{W})$ for some
$BH$-subspace $W$ of $V$.
The continuous map
$\An(\G,\overline{W}) \rightarrow \Con(G,\overline{W}) \rightarrow \Con(G,V)$
is $G$-equivariant, if we endow all the spaces appearing with the
right regular $G$-action.
Let $\overline{W}_1$ denote the
preimage in $\An(\G,\overline{W})$ 
of the $G$-invariant closed subspace
$V_{\con}$ of $\Con(G,V)$;
then $\overline{W}_1$ is a $G$-invariant closed subspace of
$\An(\G,\overline{W}).$
Let $W_1$ denote the image of $\overline{W}_1$
under the continuous map
$\overline{W}_1\rightarrow V_{\con} \buildrel \ev_e \over
\longrightarrow V;$ 
by construction $W_1$ is a $G$-invariant $BH$-subspace of $V$.
Again by construction, $W_1$ is contained in $V_1$ and contains
the element $v$.
We claim that the natural map $(\overline{W}_1)_{\G-\an} \rightarrow
\overline{W}_1$ is an isomorphism.  Granting this, it follows that
$v$ is contained in the image of $(\overline{W}_1)_{\G-\an}$ in $V_1$
and the commutativity of~(3.3.17) is proved.

We now prove the claim. Observe that
by proposition~3.3.4 the natural map
$\An(\G,\overline{W})_{\G-\an} \rightarrow \An(\G,\overline{W})$
is an isomorphism.
Our claim follows upon
applying proposition~3.3.5 to the closed immersion $\overline{W}_1 \rightarrow
\An(\G,\overline{W}).$
\qed\enddemo

\proclaim{Proposition~3.3.18}
If $V$ is a Fr\'echet space, then
the map $V_{\G-\an} \rightarrow \An(\G,V)$ provided by
proposition~3.3.16 is a closed embedding.
In particular, $V_{\G-\an}$ is again a Fr\'echet space.
\endproclaim
\demo{Proof}
Proposition~2.1.13~(ii) shows that $\An(\G,V)$ is a Fr\'echet space,
and the proof of that proposition, together with proposition~3.2.15,
shows that, as $W$ ranges
over all $G$-invariant $BH$-subspaces of $V$, the images of the maps
$\An(\G,\overline{W}) \rightarrow \An(\G,V)$ are cofinal
in the directed set of all $G$-invariant $BH$-subspaces of $\An(\G,V)$.
Since the closed subspace $\An(\G,V)^{\Delta_{1,2}(G)}$
of $\An(\G,V)$ is again a Fr\'echet space, and hence ultrabornological,
we deduce from proposition~1.1.12
that the natural map $\ilim{W} \An(\G,\overline{W})^{\Delta_{1,2}(G)}
\rightarrow \An(\G,V)^{\Delta_{1,2}(G)}$ is an isomorphism.
This proves the proposition.
\qed\enddemo

\proclaim{Proposition 3.3.19}
If $V$ is a Hausdorff convex $K$-vector
space equipped with a topological $G$-action then
the action of $G$ on
$V_{\G-\an}$ is continuous.
\endproclaim
\demo{Proof} If $W$ is a $G$-invariant $BH$-subspace
of $V$ then the action of $G$ on $\overline{W}_{\G-\an}$ is continuous.
This implies that the action of $G$ on $V_{\G-\an}$ is
separately continuous.  Since $V_{\G-\an}$ is barrelled, the
$G$-action on $V_{\G-\an}$ is continuous, as claimed.
\qed
\enddemo

The following result shows that in definition~3.3.13 we may restrict
our attention to continuously $G$-invariant $BH$-subspaces of $V$.

\proclaim{Proposition 3.3.20} If $V$ is a Hausdorff convex $K$-vector
space equipped with a topological $G$-action then 
the natural map $\ilim{U}\overline{U}_{\G-\an} \rightarrow
V_{\G-\an},$ in which the inductive limit is taken over all continuously
$G$-invariant $BH$-subspaces $U$ of $V$, is an isomorphism.
\endproclaim
\demo{Proof}
Let $W$ be a $G$-invariant $BH$-subspace of $V$.  Then the image $U$ of
$\overline{W}_{\con}$ in $W$ is a continuously $G$-invariant
$BH$-subspace of $V$ that is contained in $W$.  Since
$\overline{U}$ is isomorphic to $\overline{W}_{\con}$ by construction,
corollary~3.3.8 shows that the map $\overline{U}_{\G-\an}
\rightarrow \overline{W}_{\G-\an}$ induced by the inclusion $U\subset W$
is an isomorphism.  This proves the proposition.
\qed\enddemo

\proclaim{Corollary 3.3.21} If $V$ is a Hausdorff convex $K$-vector space
of $LF$-type (respectively of $LB$-type), 
and if $V$ is equipped with a topological $G$-action,
then $V_{\G-\an}$ is an $LF$-space (respectively an $LB$-space.).
\endproclaim
\demo{Proof}
Suppose $V=\bigcup_{n=1}^{\infty} V_n$, where $V_1\subset V_2 \subset
\cdots \subset V_n \subset \cdots $ is an increasing sequence
of $FH$-subspaces of $V$.  Proposition~1.1.10 shows that the sequence
$\{V_n\}_{n\geq 1}$ is cofinal in the directed set of $FH$-subspaces
of $V$.  Applying proposition~3.2.14, define $W_n := \core{G}{V_n}$;
then the sequence $\{W_n\}_{n\geq 1}$ is cofinal in the directed set
of continuously $G$-invariant $FH$-subspaces of $V$.
Proposition~3.3.20 shows that $V_{\G-\an} \iso \ilim{W}
(\overline{W})_{\G-\an},$ where $W$ runs through all the
continuous $G$-invariant $FH$-subspaces of $V$.  Any such
$W$ lies in one of the $W_n$, and so we see that
$V_{\G-\an} \iso \ilim{n} (\overline{W}_n)_{\G-\an}.$
Proposition~3.3.18 shows that $(\overline{W}_n)_{\G-\an}$
is a Fr\'echet space for each value of $n$, and thus that
$V_{\G-\an}$ is an $LF$-space.

If in fact each of the $V_n$ is a $BH$-subspace of $V_n$,
then so is each $W_n$.  Thus $(\overline{W}_n)_{\G-\an}$ is
a Banach space for each value of $n$, and in this case we
see that $V_{\G-\an}$ is even an $LB$-space.
\qed\enddemo

\proclaim{Proposition 3.3.22}
If $V$ is a Hausdorff convex $K$-vector
space equipped with a topological $G$-action then
the continuous injection
$(V_{\G-\an})_{\G-\an} \rightarrow V_{\G-\an}$ is in fact
a topological isomorphism.
\endproclaim
\demo{Proof}
Let $W$ range over all $G$-invariant $BH$-subspace of $V$.  The functoriality
of the formation of $\G$-analytic vectors
induces a commutative diagram
$$\xymatrix{\ilim{W}(\overline{W}_{\G-\an})_{\G-\an} \ar[r]\ar[d] &
(V_{\G-\an})_{\G-\an} \ar[d] \\
\ilim{W} \overline{W}_{\G-\an} \ar[r] & V_{\G-\an},}$$
in which the inductive limits are taken over all such $W$.
Corollary~3.3.6 implies that the left-hand vertical map is an
isomorphism, while by definition the lower horizontal arrow
is an isomorphism.  Since the other arrows are continuous injective maps,
we see that the right-hand vertical arrow is also a topological isomorphism,
as claimed.
\qed\enddemo

\proclaim{Proposition 3.3.23}
Let $V$ be a Hausdorff convex $K$-vector space
equipped with a topological $G$-action, and let $W$ be a
$G$-invariant closed subspace of $V$.  If $W_1$ denotes the preimage
of $W$ in $V_{\G-\an}$ under the natural map $V_{\G-\an} \rightarrow V$
of~(3.3.15) (so that $W_1$ is a $G$-invariant closed subspace of $V_{\G-\an}$),
then the
map $W_{\G-\an} \rightarrow V_{\G-\an}$ (induced by functoriality of
the analytic vectors) induces a continuous
bijection of $W_{\G-\an}$ onto $W_1$.  If furthermore $V_{\G-\an}$
is either a Fr\'echet space or of compact type,
then this continuous bijection is a topological isomorphism.
\endproclaim
\demo{Proof}
The first statement of the corollary
follows from propositions~2.1.23 and theorem~3.3.16.

To prove the second, we first observe that proposition~3.3.22
yields an isomorphism $(V_{\G-\an})_{\G-\an} \iso V_{\G-\an}$,
and so the first statement of the proposition, applied to
the closed embedding $W_1 \rightarrow V_{\G-\an},$ induces
a continuous bijection $(W_1)_{\G-\an} \rightarrow W_1.$
As $V_{\G-\an}$ is either a Fr\'echet space or a space of compact type,
the same is true of its closed subspace $W_1$, and
hence by corollary~3.3.21, $(W_1)_{\G-\an}$ is an $LF$-space.
This bijection is thus a topological isomorphism.
This isomorphism fits into the sequence of continuous injections
$(W_1)_{\G-\an} \rightarrow W_{\G-\an} \rightarrow W_1 \rightarrow W,$
(the first of these being induced by the third and functoriality of
the formation of $\G$-analytic vectors).  
We deduce that the map $W_{\G-\an} \rightarrow W_1$
is also a topological isomorphism, as required.
\qed\enddemo

\proclaim{Proposition 3.3.24} If 
$V$ is either a Fr\'echet space or a convex space of
compact type, 
then there is a natural isomorphism
$\Con(G,V)_{\G-\an} \iso \An(\G,V)$.
(Here $\Con(G,V)$ and $\An(\G,V)$ are regarded as $G$-representations
via the right regular $G$-action.)
\endproclaim
\demo{Proof}
Propositions~2.1.6 and~3.2.15 together imply
that as $W$ ranges over all $G$-invariant $BH$-subspaces of $V$,
the images of the maps $\Con(G,\overline{W}) \rightarrow
\Con(G,V)$ are cofinal in the set of all $G$-invariant $BH$-subspaces
of $\Con(G,V)$.  Thus $\Con(G,V)_{\an} = \ilim{W} \Con(G,\overline{W})_{\G-\an}
\iso \ilim{W} \An(\G,\overline{W}) = \An(\G,V)$ (the isomorphism
being provided by proposition~3.3.7).
\qed\enddemo

\proclaim{Corollary 3.3.25} If $V$ is either a Fr\'echet space or a convex
space of compact type, then the natural injection $\An(\G,V)_{\G-\an}
\rightarrow \An(\G,V)$ is a topological isomorphism.
(Here $\An(\G,V)$ is regarded as a $G$-representation
via the right regular $G$-action.)
\endproclaim
\demo{Proof}
Proposition~3.3.24 shows that
the injection $\An(\G,V)_{\G-\an} \rightarrow \An(\G,V)$
is obtained by passing to $\G$-analytic vectors in the injection
$\Con(G,V)_{\G-\an} \rightarrow \Con(G,V).$ 
The corollary follows from proposition~3.3.22.
\qed\enddemo

\proclaim{Corollary 3.3.26} If $V$ is either a Fr\'echet space
or a convex space of compact type, then there is a natural
isomorphism $\La(G,V)_{\G-\an} \iso \An(\G,V)$.
(Here $\La(G,V)$ and $\An(\G,V)$ are regarded as $G$-representations
via the right regular $G$-action.)
\endproclaim
\demo{Proof}
If we pass to $\G$-analytic vectors in the sequence of continuous
injections $\An(\G,V) \rightarrow \La(G,V) \rightarrow \Con(G,V)$
we obtain a sequence of continuous injections
$\An(\G,V)_{\G-\an} \rightarrow \La(G,V)_{\G-\an} \rightarrow
\Con(G,V)_{\G-\an}.$ As we observed in the proof of corollary~3.3.25,
the composite of these bijections is an isomorphism,
and thus each of these bijections is itself an isomorphism.
Applying proposition~3.3.24 we see that $\La(G,V)_{\G-\an} \iso
\An(\G,V),$ as required.
\qed\enddemo

\proclaim{Proposition 3.3.27} If $V$ is a Hausdorff convex $K$-vector
space, and there is an isomorphism $\ilim{n} V_n \iso V$,
where $\{V_n\}_{n\geq 1}$ is an inductive sequence of Hausdorff
$LB$-spaces with injective transition maps,
then there is a natural isomorphism
$\ilim{n} (V_n)_{\G-\an} \iso V_{\G-\an}.$
\endproclaim
\demo{Proof}
It follows from proposition~1.1.10 that any $BH$-invariant subspace
of $V$ lies in the image of the map $V_n\rightarrow V$ for
some sufficiently large value of $n$.  This yields the proposition.
\qed\enddemo

\section{3.4} In this subsection we suppose that $\G$ is a
rigid analytic group defined over $L$, which is $\sigma$-affinoid
as a group, in the sense that $\G = \bigcup_{n=1}^{\infty} \G_n,$ where
each $\G_n$ is an admissible affinoid open subgroup of $\G$.
We write $G := \G(L),$ and for each $n$, write $G_n := \G_n(L),$
so that $G = \bigcup_{n=1}^{\infty} G_n.$  We also assume that
$G_n$ is Zariski dense in $\G_n$ for each $n \geq 1$,
and thus that $G$ is Zariski dense in $\G.$ 

\proclaim{Definition 3.4.1} If $V$ is a Hausdorff convex $K$-vector space
equipped with a topological $G$-action, then we define the convex $K$-vector space
$V_{\G-\an}$ of $\G$-analytic vectors in $V$ to be the projective limit
$\plim{n} V_{\H-\an}$, where the projective limit is taken over all
admissible affinoid open subgroups $\H$ of $\G$.
\endproclaim

If $\H_1 \subset \H_2$ is an inclusion of admissible affinoid open
subgroups of $\G$, then the natural map $V_{\H_2-\an} \rightarrow
V_{\H_1-\an}$ is an injection (since when composed with the natural
injection $V_{\H_1-\an} \rightarrow V$ it yields the natural injection
$V_{\H_2-\an} \rightarrow V$).  Thus the transition maps in
the projective system of definition~3.4.1 are all injections.
The remarks preceding definition~2.1.18 show that we may in fact
restrict this
projective limit to be taken over the subgroups $\G_n$ of $\G$.

The formation of $V_{\G-\an}$ is evidently covariantly functorial
in $V$ and contravariantly functorial in $\G$.

\proclaim{Proposition 3.4.2} There is a natural continuous $G$-action
on $V$, as well as a natural continuous $G$-equivariant
injection $V_{\G-\an} \rightarrow V.$
\endproclaim
\demo{Proof} We regard $V_{\G-\an}$ as the projective limit
$\plim{n} V_{\G_n-\an}.$  If $m \leq n$ then the continuous $G_n$-action
on $V_{\G_n-\an}$ provided by proposition~3.3.19 restricts to a continuous
$G_m$-action on $V_{\G_n-\an}$.  Passing to the projective limit,
we obtain a continuous $G_m$-action on $V_{\G-\an}$.  
Since $G$ is the union of its subgroups $G_m,$ we obtain
a continuous $G$-action on $V.$
The $G$-equivariant injection $V_{\G-\an} \rightarrow V$ is obtained
by composing the natural projection $V_{\G-\an} \rightarrow V_{\G_n-\an}$
(for some choice of $n$)
with the continuous injection $V_{\G_n-\an} \rightarrow V$.  (The resulting
continuous injection is obviously independent of the choice of $n$.)
\qed\enddemo

\proclaim{Theorem 3.4.3}
If $V$ is a Hausdorff convex $K$-vector
space equipped with a topological $G$-action then
there is a (uniquely determined) continuous
injection $V_{\G-\an} \rightarrow \An(\G,V)$ such that the diagram
$$\xymatrix{V_{\G-\an} \ar[r]^-{\text{(3.4.2)}}\ar[d] & V \ar[d]^-{\orbit} \\
\An(\G,V) \ar[r] & \Fun(G,V)}$$ is Cartesian on the level of abstract
vector spaces.  In particular, the image of the continuous injection~(3.4.2)
contains precisely those vectors $v\in V$ for which the orbit map
$\orbit_v$ is (the restriction to $G$ of) an element of $\An(\G,V)$. 
\endproclaim
\demo{Proof} 
For each value of $n$, theorem~3.3.16 yields a diagram of
continuous maps
$$\xymatrix{V_{\G_n-\an} \ar[r]\ar[d] & V \ar[d]^-{\orbit} \\
\An(\G_n,V) \ar[r] & \Fun(G_n,V),}$$ Cartesian on the level of abstract
$K$-vector spaces.  These diagrams are compatible with respect
to the maps $\G_{n} \rightarrow \G_{n+1},$ and
passing to the projective limit in $n$ yields the diagram of
the theorem.
\qed\enddemo

\proclaim{Proposition 3.4.4}
If $V$ is a Fr\'echet space, then
the map $V_{\G-\an} \rightarrow \An(\G,V)$ provided by
theorem~3.4.3 is a closed embedding.
In particular, $V_{\G-\an}$ is again a Fr\'echet space
(since the remark following definition~2.1.18 shows that $\An(\G,V)$
is a Fr\'echet space).
\endproclaim
\demo{Proof}
Proposition~3.3.18 shows that the natural map
$V_{\G_n-\an} \rightarrow \An(\G_n,V)$ is a closed embedding,
for each $n\geq 1$.  Since the projective limit of closed
embeddings is a closed embedding, the result follows.
\qed\enddemo

\proclaim{Corollary~3.4.5} If $V$ is a $K$-Fr\'echet space
equipped with a topological action of $G$, 
then there exists a natural $G$-equivariant topological isomorphism
$V_{\G-\an} \iso (\An(\G,K)\cotimes_K V)^{\Delta_{1,2}(G)}.$
\endproclaim
\demo{Proof}
Proposition~2.1.19 yields a natural topological isomorphism
$\An(\G,V) \iso \An(\G,K)\cotimes_K V.$  The claim now
follows from proposition~3.4.4.
\qed\enddemo

\proclaim{Corollary~3.4.6} If $\G$ is strictly $\sigma$-affinoid,
and if $V$ is a nuclear Fr\'echet space
equipped with a topological $G$-action,
then $V_{\G-\an}$ is again a nuclear Fr\'echet space.
\endproclaim
\demo{Proof}
Since $\G$ is strictly $\sigma$-affinoid, the space
$\An(\G,K)$ is a nuclear Fr\'echet space.
Proposition~2.1.19 and corollary~3.4.5 show that $V_{\G-\an}$ is a
closed subspace of a nuclear Fr\'echet space, and hence is
again a nuclear Fr\'echet space.
\qed\enddemo

\proclaim{Corollary 3.4.7} If $V$ is a $K$-Fr\'echet space equipped
with a topological $G$-action, then the natural map
$(V_{\G-\an})_{\G-\an} \rightarrow V_{\G-\an}$ is a topological
isomorphism.
\endproclaim
\demo{Proof}
Since $V$ is a Fr\'echet space, the same is true of
$V_{\G-\an}$, by proposition~3.4.4.  Corollary~3.4.5
thus yields an isomorphism
$$\multline
(V_{\G-\an})_{\G-\an} \iso (\An(\G,K)\cotimes_K V_{\G-\an})^{\Delta_{1,2}(G)} \\
= ((\plim{n} \An(\G_n,K)) \cotimes_K
(\plim{n} V_{\G_n-\an}))^{\Delta_{1,2}(G)}.\endmultline \tag 3.4.8$$

Proposition~3.3.18 shows that each $V_{\G_n-\an}$ is a
Fr\'echet space, and we obtain isomorphisms
$$\multline ((\plim{n} \An(\G_n,K)) \cotimes_K
(\plim{n} V_{\G_n-\an}))^{\Delta_{1,2}(G)} \\
\iso
(\plim{n} (\An(\G_n,K)\cotimes_K V_{\G_n-\an}))^{\Delta_{1,2}(G)}
\\
\iso \plim{n} (\An(\G_n,K) \cotimes_K V_{\G_n-\an})^{\Delta_{1,2}(G_n)} 
\iso \plim{n} V_{\G_n-\an} = V_{\G-\an}.\endmultline \tag 3.4.9$$
(The first isomorphism is provided by proposition~1.1.29,
the second is evident,
and the third is provided by propositions~2.1.13~(ii), 3.3.18
and~3.3.22.)
Composing the isomorphism provided by~(3.4.9) with that
provided by~(3.4.8) establishes the corollary.
\qed\enddemo

\proclaim{Corollary 3.4.10}
Let $V$ be a Hausdorff convex $K$-vector space
equipped with a topological $G$-action, and let $W$ be a
$G$-invariant closed subspace of $V$.  If $W_1$ denotes the preimage
of $W$ in $V_{\G-\an}$ under the natural map $V_{\G-\an} \rightarrow V$
(a $G$-invariant closed subspace of $V_{\G-\an}$), then the
map $W_{\G-\an} \rightarrow V_{\G-\an}$ (induced by functoriality of
the formation of analytic vectors) induces a continuous
bijection of $W_{\G-\an}$ onto $W_1$.  If furthermore $V_{\G-\an}$ is a
Fr\'echet space, then this map is even a topological isomorphism.
\endproclaim
\demo{Proof}
The first claim follows from theorem~3.4.3 and proposition~2.1.23.
The second claim is proved in an analogous manner to the
second claim of proposition~3.3.23, by appealing to proposition~3.4.4
and corollary~3.4.7.
\qed\enddemo

\proclaim{Proposition 3.4.11} If $V$ is either a $K$-Fr\'echet space or a
convex $K$-vector space of compact type, then there is a natural
isomorphism $\Con(G,V)_{\G-\an} \iso \An(\G,V).$
(Here $\Con(G,V)$ and $\An(\G,V)$ are regarded as $G$-representations
via the right regular $G$-action.)
\endproclaim
\demo{Proof}
Proposition~3.3.24 yields an isomorphism
$\An(\G_n,V) \iso \Con(G_n,V)_{\G_n-\an}$ for each $n\geq 1$.
The naturality of these
isomorphisms implies that they are compatible with respect to
the maps $\G_n \rightarrow \G_{n+1},$ and passing to the projective
limit in $n$ yields the proposition.
\qed\enddemo

\proclaim{Proposition 3.4.12} If $V$ is a Hausdorff convex $K$-vector space
equipped with a topological action of $G$,
and if $\H$ is an affinoid rigid analytic
subgroup of $\G$ that is normalised by $G$ (so that $G$ acts naturally
on $V_{\H-\an}$),
then the natural map $(V_{\H-\an})_{\G-\an} \rightarrow V_{\G-\an}$
is a topological isomorphism.
\endproclaim
\demo{Proof}
Recall that $\G$ may be written as an increasing union of a sequence
of affinoid open subgroups $\{\G_n\}_{n\geq 1}$.  We may choose
this sequence of subgroups so that $\H = \G_1$.  
If $n\geq 1,$ then there are continuous injections
$V_{\G_n-\an} \rightarrow V_{\H-\an} \rightarrow V,$
and hence continuous injections $(V_{\G_n-\an})_{\G_n-\an}
\rightarrow (V_{\H-\an})_{\G_n-\an} \rightarrow V_{\G_n-\an}$.
Proposition~3.3.22 shows that the composite of these maps is
a topological isomorphism, and thus in particular so is the second of these,
the map
$$(V_{\H-\an})_{\G_n-\an} \rightarrow V_{\G_n-\an}.\tag 3.4.13$$
Passing to the projective limit of the maps~(3.4.13) as $n$ tends
to infinity, the proposition follows.
\qed\enddemo

\proclaim{Proposition 3.4.14} Suppose that $\H$ is a rigid analytic
affinoid group defined over $L$, containing $\G$ as
a rigid analytic subgroup, and assume that $H := \H(L)$
normalises $\G$.
If $V$ is
equipped with a topological $H$-action, then the natural map
$(V_{\G-\an})_{\H-\an} \rightarrow V_{\H-\an}$ is a 
topological isomorphism.
\endproclaim
\demo{Proof} Since $\G$ is contained in $\H$, we have natural
continuous injections $V_{\H-\an} \rightarrow V_{\G-\an}
\rightarrow V.$  As $H$ normalises $\G,$ the $H$-action on $V$
lifts to a $H$-action on $V_{\G-\an}$, and both injections are
$H$-equivariant.  Passing to $\H$-analytic vectors yields continuous
injections
$(V_{\H-\an})_{\H-\an} \rightarrow (V_{\G-\an})_{\H-\an} 
\rightarrow V_{\H-\an}$.
Proposition~3.3.22 shows that the composite of these two maps is
a topological isomorphism.  The same is thus true of each of these
maps separately. 
\qed\enddemo

\proclaim{Corollary 3.4.15}
Suppose that $\J \subset \G$ is
an inclusion of $\sigma$-affinoid groups, which factors
as $\J \subset \G_1 \subset \G,$ and such that $G$ normalises $\J$.
If $V$ is a Hausdorff convex $K$-vector space 
equipped with a topological
$G$-action then the natural map $(V_{\J-\an})_{\G-\an} \rightarrow
V_{\G-\an}$ is a topological isomorphism.
\endproclaim
\demo{Proof}
By proposition~3.4.14, the natural map
$(V_{\J-\an})_{\G_n-\an} \rightarrow V_{\G_n-\an}$ is an isomorphism
for each $n\geq 1$.  Passing to the projective limit in $n$
yields the corollary.
\qed\enddemo

\section{3.5}
In this subsection we let $G$ denote a locally $L$-analytic
group.  We begin by introducing the notion of an analytic
open subgroup of $G$.

Suppose that $(\phi,H,\H)$ is a chart of $G$ --
thus $H$ is a compact open subset of $G$, $\H$
is an affinoid rigid analytic space over $L$ isomorphic to
a closed ball, and $\phi$ is an isomorphism $\phi:H \iso \H(L)$ 
-- with the additional property that $H$ is a
subgroup of $G$.
Since $H$ is Zariski dense in $\H$, there is at most one
rigid analytic group structure
on $\H$ giving rise to the group structure on $H$.
If such a rigid analytic group structure exists on $\H$, we
will refer to the chart $(\phi,H,\H)$
as an analytic open subgroup of $G$.
Usually, we will suppress the isomorphism $\phi,$ and simply
refer to an analytic open subgroup $H$ of $G$, and write
$\H$ to denote the corresponding rigid analytic group determined by $H$.

The analytic open subgroups of $G$ form a directed set in an obvious
fashion:  if $H' \subset H$ is an inclusion of open subgroups of $G$
each of which is equipped with
the structure of an analytic open subgroup of $G$,
then we say that it is an inclusion of analytic open subgroups if
it lifts to a rigid analytic map $\H' \rightarrow \H.$  (Since $H'$
and $H$ are Zariski dense in $\H'$ and $\H$ respectively, such a lift
is uniquely determined, and is
automatically a homomorphism of rigid analytic groups.)
Forgetting the chart structure yields an order-preserving
map from the directed set
of analytic open subgroups of $G$ to the set of all open
subgroups of $G$.
Since the group structure on $G$ is locally analytic by
assumption, the image of this map is
cofinal in the directed set of all open subgroups of $G$.
(For example, see \cite{\SERLG, LG~4 \S 8};
standard open subgroups -- in the terminology of that
reference -- provide examples of analytic open subgroups.)

If $H$ is an analytic open subgroup of $G$ and $g$ is an element of $G$,
then multiplication by $g$ induces an isomorphism $H \iso g H.$
Thus $g H$ is naturally a chart of $G$.  Just as we let $\H$
denote the rigid analytic structure induced on $H$ by its structure
as a chart of $G$, we will let $g\H$
denote the rigid analytic space induced on $g H$ by its structure
as a chart of $G$.  (In fact, if $\phi: H \iso \H(L)$ describes the
chart structure on $H$, then the chart structure on $g H$ is
describe by the map $\phi': g H \iso \H(L)$ defined as
$\phi': g h \mapsto \phi(h)$.
Thus the rigid analytic  structure on $g H$ is obtained not by changing
the rigid analytic space $\H$, but rather by changing the map
$\phi$.  However, it is easier and more suggestive in the
exposition to suppress the maps $\phi$ and $\phi'$, and to instead
use the notation $g \H$ to denote the rigid analytic structure on $g H$.)
If $g$ ranges over a set of right coset representatives of $H$ in $G$,
then the collection of charts $\{g \H\}_{g \in G/H}$ forms an analytic
partition of $G$.  If $G$ is compact, then the set of such analytic partitions
obtained by allowing $H$ to run over all analytic open
subgroups of $G$ is cofinal in the set of all analytic partitions of $G$.

If $V$ is a Hausdorff locally convex topological $K$-vector
space equipped with a topological action of $G$ then for each analytic
open subgroup $H$ of $G$ we may form the convex $K$-vector space $V_{\H-\an}$
of $\H$-analytic vectors in $V$, equipped with its $H$-invariant
continuous injection $V_{\H-\an} \rightarrow V$.
Since the formation of this space and injection is functorial
in $H$ these spaces form an inductive
system:   if $H' \subset H$ is an inclusion of analytic open
subgroups of $G$ then there is a continuous $H'$-equivariant
injection $V_{\H-\an} \rightarrow V_{\H'-\an}$, compatible
with the injections of each of the source source and target 
into $V$.  Passing to the locally convex inductive limit
we obtain a continuous injection $\ilim{H} V_{\H-\an}
\rightarrow V$.

\proclaim{Lemma 3.5.1} Let $V$ be a Hausdorff convex $K$-vector space equipped
with a topological action of $G$.
If $g$ is an element of $G$ and $H$ is an analytic open
subgroup of $G$ (so that $g H g^{-1}$ is a second
such subgroup) then the automorphism of $V$
induced by the action of $g$ lifts in a unique fashion to an isomorphism
$V_{\H-\an} \iso V_{g\H g^{-1}}$.
\endproclaim
\demo{Proof}
The meaning of the lemma is that there is a unique way to fill
in the top horizontal arrow of the following diagram so as to
make it commute:
$$\xymatrix{V_{\H-\an}\ar[r]\ar[d] & V_{g\H g^{-1} - \an} \ar[d] \\
V \ar[r]^-{g \cdot } & V .}$$
The uniqueness is clear, since the vertical arrows are injections.
To see the existence, let $\phi: H \rightarrow g H g^{-1}$
denote the conjugation map $h \mapsto  g h g^{-1}$.  Then
the continuous linear map from $V$ to itself given by the
action of $g$ is equivariant with respect to $\phi$.  The
functoriality of the construction of analytic vectors
shows that the action of $g$ on $V$ lifts to a map
$V_{\H-\an} \rightarrow V_{g\H g^{-1}-\an}.$  
\qed\enddemo

\proclaim{Corollary 3.5.2} In the preceding situation there is
an action of $G$ on the inductive system $\{V_{\H-\an}\}$,
and hence on the inductive limit $\ilim{H} V_{\H-\an}$,
uniquely determined by the condition that the continuous
injection $\ilim{H} V_{\H-\an} \rightarrow V$ be $G$-equivariant.
\endproclaim
\demo{Proof}
The uniqueness is clear, since all the transition maps in
the inductive system $\{V_{\H-\an}\}$ are injective.
Lemma~3.5.1 shows that if $g$ is an element of $G$ and
$H$ is any analytic open subgroup of $G$
then the action of $g$ on $V$ lifts to a map
$V_{\H-\an} \rightarrow V_{g\H g^{-1}-\an}.$
By functoriality
of the construction of analytic vectors, these combine to yield
the required action of $G$ on the inductive system $\{V_{\H-\an}\}$.
\qed\enddemo

\proclaim{Definition 3.5.3} Suppose that 
$V$ is a locally convex
Hausdorff $K$-topological vector space equipped with a
topological action of $G$.  We define the locally convex
space $V_{la}$ of locally analytic vectors in $V$ to be the
locally convex inductive limit
$$V_{\la} := \ilim{\H} V_{\H-\an},$$ where $H$ runs
over all the analytic open subgroups of $\G$ (and $\H$
denotes the rigid analytic group corresponding to $H$).
If we wish to emphasise the role of the locally analytic
group $G$ then we will write $V_{G-\la}$ in place of $V_{\la}$.
\endproclaim

The space $V_{\la}$ is equipped with a continuous injection
into $V$, and thus is Hausdorff.
Lemma~3.5.1 shows that $V_{\la}$ 
is equipped with a continuous action of $G$ with respect to which
its injection into $V$ is $G$-equivariant.
The construction of $V_{\la}$ (as a locally convex
topological $K$-vector space with $G$-action)
is obviously contravariantly functorial in $V$ and covariantly
functorial in $G$.

Theorem~3.5.7 below justifies our designation of $V_{\la}$ as
the space of locally analytic vectors of $V$.
(See the remarks following that theorem for a comparison
of our definition of the space of locally analytic vectors
with the definition given in \cite{\SCHTBD, \S 3}
and \cite{\SCHTNEW, p.~38}.)

\proclaim{Lemma 3.5.4}
If $H$ is an open subgroup of $G$ and 
$V$ is a Hausdorff convex $K$-vector space equipped with a topological
$G$-action then the natural $H$-equivariant morphism
$V_{H-\la} \iso V_{G-\la}$ is an isomorphism.
\endproclaim
\demo{Proof}
This follows immediately from the fact that the directed set of analytic
open subgroups of $G$ that are contained in $H$ is cofinal in
the directed set of analytic open subgroups of $G$.
\qed\enddemo

\proclaim{Proposition 3.5.5} If $V$ is a Hausdorff convex $K$-vector space
equipped with a topological action of $G$ then the continuous injection
$(V_{\la})_{\la} \rightarrow V_{\la}$ is a topological isomorphism.
\endproclaim
\demo{Proof}
If we fix an analytic open subgroup $H$ of $G$ and an
$H$-invariant $BH$-subspace $W$ of $V$, then
the image of $\overline{W}_{\H-\an}$ in $V_{\la}$ is an $H$-invariant
$BH$-subspace of $V_{\la}$, and
corollary~3.3.6 shows that the natural map $(\overline{W}_{\H-\an})_{\H-\an}
\rightarrow \overline{W}_{\H-\an}$ is an isomorphism.
Taking the limit over all such $W$ and $H$ we obtain continuous injections 
$$\ilim{W,H} (\overline{W}_{\H-\an})_{\H-\an}
\rightarrow (V_{\la})_{\la} \rightarrow
V_{\la} = \ilim{W,H} \overline{W}_{\H-\an}.$$
Our preceding remarks imply that the composite of these injections
is a topological isomorphism, and consequently the second injection
is also a topological isomorphism.
\qed\enddemo

\proclaim{Proposition 3.5.6}
If $V$ is a Hausdorff convex $K$-vector space of $LF$-type,
(respectively $LB$-type),
then $V_{\la}$ is an $LF$-space
(respectively an $LB$-space).
\endproclaim
\demo{Proof}
Let $\{H_n\}_{n\geq 1}$ denote a cofinal sequence of 
of analytic open subgroups of $\G$. Then
$V_{\la} = \ilim{n}V_{\H_n-\an}.$  
Corollary~3.3.21 shows that each
$V_{\H_n-\an}$ is either an $LF$-space or an $LB$-space
(depending on the hypothesis on $V$), and
thus the same is true of $V_{\la}$.
\qed\enddemo

\proclaim{Theorem 3.5.7}  
If $V$ is a Hausdorff convex $K$-vector space equipped with a
topological action of $G$ then
there is a (uniquely determined) continuous injection $V_{\la} \rightarrow
\La(G,V)$ such that the diagram
$$\xymatrix{V_{\la} \ar[d] \ar[r] & V \ar[d]^-{\orbit} \\
\La(G,V) \ar[r] & \Fun(G,V)}\tag 3.5.8 $$ is Cartesian on the level
of abstract vector spaces.
Thus $V_{\la}$ maps via a continuous bijection
onto $\La(G,V)^{\Delta_{1,2}(G)}$,
and the image of the continuous
injection $V_{\la} \rightarrow V$ is equal to the
subspace of $V$ consisting of those
vectors $v\in V$ for which $\orbit_v$ is an element of
$\La(G,V)$.
\endproclaim
\demo{Proof}
Since the natural maps $V_{\la} \rightarrow V$ and $\La(G,K)
\rightarrow \Fun(G,V)$  are injective, there is at most one
map $V_{\la} \rightarrow V$ making~(3.5.8) commute.
If $H$ is an analytic open subgroup of $G$ then
theorem~3.3.16 shows that the association of
$\orbit_{v | H}$ to any $v\in V$
yields a continuous map
$V_{\H-\an} \rightarrow \An(\H,V).$
From this we deduce that the association of $\orbit_{v | g H}$ to $v \in V$
induces a continuous map $V_{\H-\an} \rightarrow \An(g\H,V)$ for any
$g \in G$.
Indeed, multiplication on the left by $g^{-1}$ induces an analytic isomorphism
$g\H \rightarrow \H,$ and thus an isomorphism $\An(\H,V) \iso \An(g \H,V)$.
The continuous action of $g$ on $V$ induces an isomorphism
$\An(g \H,V) \iso \An(g \H, V).$  The composite
$$\xymatrix{ V_{\H-\an} \ar[r]^-{v \mapsto \orbit_{v | \H} } &
\An(\H,V) \ar[r]^-{\sim} & \An(g \H , V) \ar[r]^-{v \mapsto gv} &
\An(g \H, V)}$$
is thus a continuous map from $V_{\H-\an}$ to $\An(g \H, V)$ which
is immediately seen to be equal to the map $v \mapsto \orbit_{v | g \H}.$
If we take the product of these maps as $g$ ranges over a set
of right coset representatives for $g \in G$, we obtain a continuous
map $V_{\H-\an} \rightarrow \prod_{g \in G/H} \An(g \H, V),$ given
by $g \mapsto g v$.
Composing with the natural map $\prod_{g \in G/H} \An(g \H,V)
\rightarrow \La(G,V)$ (recall that $\{ g \H\}_{g \in G/H}$ forms
an analytic partition of $G$)
and then passing to the inductive limit over all analytic
open subgroups $H$ of $G$, we obtain the continuous map
$V_{\la} \rightarrow \La(G,V)$ which makes~(3.5.8) commutes.

We turn to showing that~(3.5.8) is Cartesian.
If $v$ is an element of $V$ for which $\orbit_v$
lies in $\La(G,V)$  
then in particular there is an analytic open subgroup
$H$ of $G$ such that $\orbit_{v | H}$ lies in $\An(\H,V)$.  Theorem~3.3.16
shows that $v$ then lies in the image of $V_{\H-\an}$, and so in
particular in the image of $V_{\la}$.  This proves that~(3.5.8)
is Cartesian.
\qed\enddemo

In the papers \cite{\SCHTBD} and \cite{\SCHTNEW}, the
space of locally analytic vectors in a continuous representation of $G$
on a Banach space $V$ is {\it defined} to be the closed
subspace $\La(G,V)^{\Delta_{1,2}(G)}$ of $\La(G,V)$.
(In these references, the authors use the terminology
``analytic vectors'' rather than ``locally analytic vectors'',
and write $V_{\an}$ rather than $V_{\la}$.)
Thus the topology imposed on the space $V_{\la}$ in these
references is in general coarser than the topology that we
impose on $V_{\la}$.

\proclaim{Proposition 3.5.9} If $V$
is a convex space of compact type equipped with a topological
$G$-action, then the map $V_{\la} \rightarrow \La(G,V)$
of theorem~3.5.7 is a closed embedding.  In particular, $V_{\la}$
is again of compact type.
\endproclaim
\demo{Proof}
Let $H$ be a compact open subgroup of $G$.  It is observed
in \cite{\SCHTBD, \S 3} that the natural map
$\La(H,V)^{\Delta_{1,2}(H)} \rightarrow \La(G,V)^{\Delta_{1,2}(G)}$
is a topological isomorphism.  Thus, replacing $G$ by $H$ if necessary,
we may suppose without loss of generality that $G$ is compact.
Proposition~2.1.28 then shows that $\La(G,V)$ is of compact type,
while proposition~3.5.6 shows that $V_{\la}$ is an $LB$-space.
The image of $V_{\la}$ under the map of theorem~3.5.7
is a closed subspace of $\La(G,V)$, and thus is of compact type.
Theorem~1.1.17 now shows that $V_{\la}$ maps isomorphically onto
this subspace.
\qed\enddemo

\proclaim{Proposition 3.5.10}
Let $V$ be a Hausdorff convex $K$-vector space
equipped with a topological $G$-action, and let $W$ be a
$G$-invariant closed subspace of $V$.  Let $W_1$ denote the preimage
of $W$ in $V_{\la}$ under the natural map $V_{\la} \rightarrow V,$
which is a $G$-invariant closed subspace of $V_{\la}$.  Then the
map $W_{\la} \rightarrow V_{\la}$ (induced by functoriality of
the formation of locally analytic vectors) induces a continuous
bijection of $W_{\la}$ onto $W_1$. 
If $V_{\la}$ is of compact type,
then this bijection is even a
topological isomorphism.
\endproclaim
\demo{Proof}
The first statement follows from proposition~2.1.27 and theorem~3.5.7.

We now prove the second statement,
assuming that $V_{\la}$ is of compact type. 
Proposition~3.5.5 implies that the natural map
$(V_{\la})_{\la} \rightarrow V_{\la}$ is a topological isomorphism.
The first statement of the corollary, applied to the closed embedding
$W_1 \rightarrow V_{\la}$, yields a continuous bijection
$(W_1)_{\la} \rightarrow W_1$. 
As $V_{\la}$ is assumed to be of compact type, its closed subspace $W_1$
is also of compact type.  Proposition~3.5.9
implies that $(W_1)_{\la}$ is again of compact type, and thus that this
continuous bijection is in fact a topological isomorphism.
This isomorphism fits into the sequence of continuous injections
$(W_1)_{\la} \rightarrow W_{\la} \rightarrow W_1 \rightarrow W,$
(the first of these being induced by the third and functoriality of
the formation of locally analytic vectors).  We conclude that
the second of these maps is also a topological isomorphism,
as required.
\qed\enddemo

\proclaim{Proposition 3.5.11}
If $G$ is compact,
and if $V$ is either a $K$-Fr\'echet space or a convex $K$-vector space of
compact type, then there is a natural isomorphism
$\Con(G,V)_{\la} \iso \La(G,V)$.
(Here the spaces
$\Con(G,V)$ and $\La(G,V)$ are regarded as $G$-representations
via the right regular $G$-action.)
\endproclaim
\demo{Proof} 
If $H$ is an
analytic open subgroup of $G$ then it is of finite index
in $G$, and we may write $G = \coprod_{i\in I} g_i H,$ where
$I$ is some finite index set.  Then
$\Con(G,V) \iso \prod_{i\in I} \Con(g_i H,V),$ and so proposition~3.3.24
shows that
$$\Con(G,V)_{\H-\an} \iso \prod_{i\in I} \Con(g_i H,V)_{\H-\an}
\iso \prod_{i\in I} \An(g_i \H, V).$$  As $H$ ranges over all
analytic open subgroups of $G$, the partitions $\{g_i H\}_{i\in I}$
are cofinal in the set of all analytic partitions of $G$.
Thus we see that
$$\Con(G,V)_{\la} = \ilim{H}\Con(G,V)_{\H-\an} \iso
\ilim{H} \prod_{i \in I} \An(g_i \H,V) \iso \La(G,V) .\qed$$
\enddemo

The next result is originally due to Feaux de Lacroix
\cite{\FETH, bei.~3.1.6}.

\proclaim{Proposition 3.5.12}
If $G$ is compact, then for any
Hausdorff convex $K$-vector space, the natural map
$\La(G,V)_{\la} \rightarrow \La(G,V)$ is a topological isomorphism.
(Here $\La(G,V)$ is regarded as a $G$-representation via the
right regular $G$-action.)
\endproclaim
\demo{Proof} 
Since $G$ is compact, there is a natural isomorphism
$$\ilim{W} \La(G,\overline{W}) \iso \La(G,V),$$ where $W$ runs over all
$BH$-subspaces of $V$.  The functoriality of the formation
of locally analytic vectors gives rise to the commutative diagram
$$\xymatrix{\ilim{W} \La(G,\overline{W})_{\la} \ar[d]\ar[r] &
\La(G,V)_{\la} \ar[d] \\
\ilim{W}  \La(G,\overline{W})_{\la} \ar[r]^{\sim} & \La(G,V) .}\tag 3.5.13$$
Propositions~3.5.5 and~3.5.11 together show that for each $W$,
the natural map $\La(G,\overline{W})_{\la} \rightarrow \La(G,\overline{W})$
is a topological isomorphism.   
The left-hand vertical arrow of~3.5.13 is thus an isomorphism.
Since the lower horizontal arrow is also an isomorphism,
and the other arrows are continuous injections,
we see that the right-hand vertical arrow is a topological isomorphism,
as required.
\qed\enddemo

\proclaim{Proposition 3.5.14} If $V$ is a Hausdorff convex $K$-vector
space, and there is an isomorphism $\ilim{n} V_n \iso V$,
where $\{V_n\}_{n\geq 1}$ is an inductive sequence of Hausdorff
$LB$-spaces with injective transition maps,
then there is a natural isomorphism
$\ilim{n} (V_n)_{\la} \iso V_{\la}.$
\endproclaim
\demo{Proof}
It follows from proposition~1.1.10 that any $BH$-invariant subspace
of $V$ lies in the image of the map $V_n\rightarrow V$ for
some sufficiently large value of $n$.  This yields the proposition.
\qed\enddemo

\proclaim{Proposition 3.5.15} If $U$ and $V$ are Hausdorff $LB$-spaces,
each equipped with a topological action of $G$,
for which $U_{\la}$ and $V_{\la}$ are spaces of compact type,
then the natural map
$U_{\la} \cotimes_K V_{\la} \rightarrow U\cotimes_K V$
factors as a composite of continuous $G$-equivariant maps
$$U_{\la} \cotimes_K V_{\la} \rightarrow
(U \cotimes_K V)_{\la} \rightarrow U\cotimes_K V,$$
where the second arrow is the natural injection.
(Here $U\cotimes_K V$ is equipped with the diagonal action of $G$.)
\endproclaim
\demo{Proof}
Each of $U$ and $V$ has a cofinal sequence of $BH$-subspaces,
say $\{W_{1,n}\}_{n \geq 1}$ and $\{W_{2,n}\}_{n \geq 2}$
respectively.  If $\{H_n\}_{n \geq 1}$ is a cofinal sequence
of analytic open subgroups of $G$,
then $U_{\la} \iso \ilim{n} (\overline{W}_{1,n})_{\H_n-\an}$
and $V_{\la} \iso \ilim{n} (\overline{W}_{2,n})_{\H_n-\an}.$
Proposition~1.1.32 yields an isomorphism
$$U_{\la} \cotimes_K V_{\la} \iso
\ilim{n} (\overline{W}_{1,n})_{\H_n-\an} \cotimes_K 
(\overline{W}_{2,n})_{\H_n-\an}. \tag 3.5.16$$
Proposition~3.3.12 yields a map
$$(\overline{W}_{1,n})_{\H_n-\an} \cotimes_K (\overline{W}_{2,n})_{\H_n-\an}
\rightarrow
(\overline{W}_{1,n}\cotimes_K \overline{W}_{2,n})_{\H_n-\an}$$
for each $n$, while functoriality of the formation of analytic
vectors yields a map
$$(\overline{W}_{1,n}\cotimes_K \overline{W}_{2,n})_{\H_n-\an}
\rightarrow (U\cotimes_K V)_{\H_n-\an}$$ for each $n$.
Altogether we obtain a map
$$\ilim{n}
(\overline{W}_{1,n})_{\H_n-\an} \cotimes_K (\overline{W}_{2,n})_{\H_n-\an}
\rightarrow (U\cotimes_K V)_{\la}.\tag 3.5.17$$
Composing~(3.5.17) with~(3.5.16) yields a map
$U_{\la}\cotimes_K V_{\la} \rightarrow (U\cotimes_K V)_{\la}.$
Composing this map with the injection
$(U \cotimes_K V)_{\la} \rightarrow U\cotimes_K V,$
yields the map 
$U_{\la} \cotimes_K V_{\la} \rightarrow U\cotimes_K V.$
Thus we have obtained the required factorisation.
\qed\enddemo

\section{3.6}  In this subsection we recall the definition of a locally
analytic representation of a locally $L$-analytic group on a convex
$K$-vector space, and establish some basic properties of this notion.
We also introduce the related notion of an analytic representation
of the group of $L$-valued points of a $\sigma$-affinoid rigid analytic group.

\proclaim{Definition 3.6.1} Let $\G$ be a $\sigma$-affinoid rigid analytic
group defined over $L$, such that $G := \G(L)$ is Zariski dense
in $\G$.  A $\G$-analytic representation of $G$ consists
of a barrelled Hausdorff convex $K$-vector space $V$ equipped with a
topological $G$-action,
having the property that the natural injection $V_{\G-\an} \rightarrow V$
is in fact a bijection.
\endproclaim

The justification for this terminology is provided by theorem~3.3.16.
The requirement that $V$ be barrelled is made so as to ensure that
the $G$-action on $V$ is continuous.  (Indeed, since the $G$-action
on $V_{\G-\an}$ is continuous, we see that the $G$-action on $V$
is certainly separately continuous.
Since $V$ is barrelled, it follows that the $G$-action
on $V$ is continuous.)

If $V$ is an arbitrary Hausdorff convex $K$-vector space equipped with
a topological $G$-action, and if $\G$ is affinoid,
then $V_{\G-\an}$ is barrelled, and
so proposition~3.3.22 shows that $V_{\G-\an}$ affords a 
$\G$-analytic representation of $G$.  If $V$ is a $K$-Fr\'echet
space, then for any $\sigma$-affinoid $\G$, corollary~3.4.7
shows that $V_{\G-\an}$ affords a $\G$-analytic representation of $G$.

Let us remark that the notion of $\G$-analytic representation
is closely related to the notion of a uniformly analytic
(``gleichm\"assig analytisch'') representation of $G$ defined in
\cite{\FETH, def.~3.1.5}.  Indeed, any such representation
of $G$ is $\H$-analytic for some sufficiently small analytic
open subgroup $H$ of $G$.

\proclaim{Proposition 3.6.2}
Let $\G$ be a $\sigma$-affinoid rigid analytic group such that $G := \G(L)$ 
is Zariski dense in $\G$, let
$V$ be a Hausdorff convex $K$-vector space
equipped with a $\G$-analytic representation,
and let $W$ be a closed $G$-invariant subspace of $V$.
If $W$ (respectively $V/W$) is barrelled,
then $W$ (respectively $V/W$) is again an analytic representation of  $\G$.
\endproclaim
\demo{Proof} Proposition~3.3.23 shows that the natural map
$W_{\G-\an} \rightarrow W$ is a bijection. 
On the other hand, functoriality of the formation
of $\G$-analytic vectors induces a commutative diagram
$$\xymatrix{ V_{\G-\an} \ar[r]\ar[d] & (V/W)_{\G-\an} \ar[d]\\
V \ar[r] & V/W.}$$
Since the bottom horizontal arrow is surjective, while both vertical arrows
are injective, we see that if
the left-hand vertical arrow is a bijection, the same is true
of the right-hand vertical arrow.
\qed\enddemo

\proclaim{Theorem 3.6.3}
If $\G$ is an affinoid rigid analytic group such that $G := \G(L)$ 
is Zariski dense in $\G$ and $V$ is a Hausdorff $LF$-space
equipped with a $\G$-analytic representation
of $G$, then the natural map $V_{\G-\an} \rightarrow V$ is a $G$-equivariant
topological isomorphism.
\endproclaim
\demo{Proof} Corollary~3.3.21 shows that $V_{\G-\an}$ is also
an $LF$-space.
Thus the map $V_{\G-\an} \rightarrow V$
is a continuous $G$-equivariant bijection between $LF$-spaces,
and so an isomorphism, by theorem~1.1.17.
\qed\enddemo

The preceding theorem is the analogue
for analytic representations of proposition~3.2.11.
Taking into account theorem~3.3.16, it shows that if $V$ is a $\G$-analytic
representation on a Hausdorff $LF$-space, then there is a
{\it continuous} orbit map $\orbit: V \rightarrow \An(\G,V)$.

The next result provides a variant in the rigid analytic setting of
the untwisting map of lemma~3.2.4.

\proclaim{Lemma 3.6.4} If $V$ is a $\G$-analytic representation of $G$
on a Hausdorff $LF$-space, then there is a
natural isomorphism $\An(\G,V) \iso \An(\G,V)$  which intertwines 
the $\Delta_{1,2}(G)$-action on its source and the left regular $G$-action
on its target.
\endproclaim
\demo{Proof}
The natural map $V_{\G-\an} \rightarrow V$ is a topological isomorphism,
by theorem~3.6.3.  Proposition~3.3.21 and its proof
then show that $V$ may be written as a $G$-equivariant inductive limit (with
injective transition maps)
of Fr\'echet spaces equipped with a $\G$-analytic action of $G$,
say $V \iso \ilim{n} V_n$.
There is then a corresponding isomorphism
$\An(\G,V) \iso \ilim{n} \An(\G,V_n)$
(since any map from a Banach space to $V$ factors through some $V_n$).
Thus it suffices to prove the lemma with $V$ replaced by each $V_n$
in turn, and thus we assume for the remainder of the proof that $V$
is a Fr\'echet space.

Proposition~2.1.13~(ii) gives an isomorphism
$\An(\G,V) \iso \An(\G,K) \cotimes_K V,$ while proposition~3.3.18
gives a closed embedding
$V \rightarrow \An(\G,V).$  Together, these give a closed embedding
$$\multline \An(\G,V) \buildrel (2.1.13) \over \longrightarrow
\An(\G,K) \cotimes_K V
\buildrel \id \cotimes (3.3.18) \over \longrightarrow
\An(\G,K) \cotimes_K \An(\G,V) \\
\buildrel \id \cotimes~(2.1.13) \over \longrightarrow
\An(\G,K) \cotimes_K \An(\G,K) \cotimes_K V.\endmultline\tag 3.6.5$$
Now the antidiagonal embedding $a: \G \rightarrow \G \times \G$
(defined by $a: g \mapsto (g,g^{-1})$) induces a continuous map
$$\An(\G,K) \cotimes_K \An(\G,K) \iso \An(\G\times \G,K)
\buildrel a^* \over \longrightarrow \An(\G,K).$$
Composing $a^*\cotimes \id$ with~(3.6.5), and with the inverse of the
isomorphism of proposition~2.1.13~(ii), yields a continuous map
$\An(\G,V) \rightarrow \An(\G,V).$ 

From its construction, one easily checks that the automorphism
of $\An(\G,V)$ constructed in the preceding paragraph
takes a function $\phi$ to the function $\tilde{\phi}: g \mapsto
g^{-1} \phi(g).$  That is, we have shown that the automorphism
of $\Fun(G,V)$ provided by lemma~3.2.4 leaves 
$\An(\G,V)$ invariant.  Similar considerations show that the inverse
of this automorphism likewise leaves $\An(\G,V)$ invariant.
Thus the automorphism of $\Fun(G,V)$ provided by lemma~3.2.4
restricts to a topological automorphism of $\An(\G,V)$.
Lemma~3.2.4 shows that this automorphism intertwines the action
of $\Delta_{1,2}(G)$ and the left regular $G$-action on
$\An(\G,V)$.
\qed\enddemo

Similarly, we may find an isomorphism from $\An(\G,V)$ to itself
that intertwines
the $\Delta_{1,3}(G)$-action and the right regular $G$-action
on $\An(\G,V)$.

\proclaim{Proposition 3.6.6}
Let $\G$ be an affinoid rigid analytic group such that $G := \G(L)$ 
is Zariski dense in $\G$.
If $V$ is a Hausdorff convex $K$-vector space
equipped with a topological $G$-action, and if $W$ is a finite dimensional
$\G$-analytic representation of $G$, then there is a natural isomorphism
$(V \otimes_K W)_{\G-\an} \iso V_{\G-\an} \otimes_K W.$  (Here the tensor
product is equipped with the diagonal $G$-action.)
\endproclaim
\demo{Proof} 
Suppose first that $V$ is a $K$-Banach space.  If
$\phi \in \An(\G, V \otimes_K W)$
then we denote by $\tilde{\phi}$ the function
$g \mapsto (\id \otimes g^{-1}) \phi(g).$  (Here $\id \otimes g^{-1}$
is denoting the indicated
automorphism of $V \otimes_K W$.)  Since the $G$-action on
$W$ is assumed to be $\G$-analytic, the function $\tilde{\phi}$ again
lies in $\An(\G,V\otimes_K W),$ and the association of $\tilde{\phi}$
to $\phi$ induces an isomorphism of $\An(\G, V\otimes_K W)$ with
itself.  (The inverse isomorphism takes a function $\tilde{\phi}$ to
the function $\phi: g \mapsto (\id \otimes g) \phi(g)$.)
In fact, it simplifies our discussion if we use the canonical isomorphism
$\An(\G, V\otimes_K W ) \iso \An(\G, V) \otimes_K W$ to regard $\phi
\mapsto \tilde{\phi}$ as an
isomorphism $\An(\G, V\otimes_K W) \iso \An(\G,V) \otimes_K W$.
One then checks that this isomorphism intertwines the $\Delta_{1,2}(G)$-action
on the source with tensor product of the $\Delta_{1,2}(G)$-action
and the trivial action on the target.
Passing to $\Delta_{1,2}(G)$-invariants on the source,
we obtain the desired isomorphism
$(V \otimes_K W)_{\G-\an} \iso V_{\G-\an} \otimes_K W.$

Now suppose that $V$ is arbitrary.  Proposition 1.1.8 shows that
as $U$ ranges over all $BH$-subspaces of $V$, the tensor products
$U\otimes_K W$ are cofinal among the $BH$-subspaces of $V\otimes_K W$.
Combining this observation with the result already proved for Banach spaces,
and the fact that tensor product with $W$ commutes with inductive limits,
we obtain a natural isomorphism
$$\multline
(V\otimes_K W)_{\G-\an}
\iso \ilim{U} (\overline{U} \otimes_K W)_{\G-\an}
\iso \ilim{U} (\overline{U}_{\G-\an} \otimes_K W) \\
\iso (\ilim{U} \overline{U}_{\G-\an}) \otimes_K W
= V_{\G-\an} \otimes_K W,\endmultline$$
as required. 
\qed\enddemo

\proclaim{Corollary 3.6.7}
Let $\G$ be an affinoid rigid analytic group such that $G := \G(L)$ 
is Zariski dense in $\G$.
If $V$ is a Hausdorff convex $K$-vector space
equipped with a $\G$-analytic representation of $G$,
and if $W$ is a finite dimensional
$\G$-analytic representation of $G$, then the diagonal $G$-action
on $V\otimes_K W$ makes this tensor product a $\G$-analytic representation
of $G$.
\endproclaim
\demo{Proof}
Since $V$ is barrelled and $W$ is finite dimensional, the tensor product
$V\otimes_K W$ is certainly barrelled.  The corollary is now an immediate
consequence of proposition~3.6.6.
\qed\enddemo

\proclaim{Proposition 3.6.8}
If $V$ is a Hausdorff convex $K$-vector space
equipped with a topological $G$-action, and if
there is an isomorphism $\ilim{i\in I} V_i \iso V,$
where $\{V_i\}_{i \in I}$ is a $G$-equivariant inductive system
of Hausdorff $K$-vector spaces, each equipped with a $\G$-analytic
action of $G$,
then the $G$-action on $V$ is again $\G$-analytic.
\endproclaim
\demo{Proof}
Functoriality of the formation of $\G$-analytic vectors
yields the commutative diagram
$$\xymatrix{ \ilim{i\in I} (V_i)_{\G-\an} \ar[d]\ar[r] & V_{\G-\an} \ar[d] \\
\ilim{i\in I} V_i \ar[r] & V.}$$
The left-hand vertical arrow and lower horizontal arrow
are both continuous bijections, by assumption, and thus
so is the right-hand vertical arrow.
\qed\enddemo

Having introduced analytic representations, we now consider
locally analytic representations.

\proclaim{Definition 3.6.9} Let $G$ be a locally $L$-analytic group.
A locally analytic representation of $G$ is a topological action of $G$
on a barrelled Hausdorff locally convex topological $K$-vector space
for which the natural map $V_{\la} \rightarrow V$ is a bijection.
\endproclaim

This definition is taken from \cite{\SCHTAN, p.~12}.
(That it is equivalent with the definition given in this reference
follows from theorem~3.5.7.)
The requirement that $V$ be barrelled is made so as to ensure that the
$G$-action on $V$ is continuous.  (See the discussion following
definition~3.6.1.)

Note that if $\G$ is an affinoid rigid analytic group such that $G := \G(L)$ and
$G$ is Zariski dense in $\G$, then a $\G$-analytic representation of
$G$ is in particular a locally analytic representation of $G$.

If $V$ is an arbitrary convex $K$-vector space equipped with a topological
action of $G$ then $V_{\la}$ is barrelled, and proposition~3.5.5 then
shows that $V_{\la}$ is equipped with a locally analytic representation
of $G$.  Thus if $G$ is compact and
if $V$ is either a Fr\'echet space or a convex space of compact type, then,
by proposition~3.5.11,
$\La(G,V)$ is a locally analytic representation of $G$, which
by propositions~2.1.6 and~3.5.6 is an $LB$-space.

If $V$ is a locally analytic representation of $G$
then there is induced on $V$ an
action of the universal enveloping algebra $U(\lie{g})$ of the Lie
algebra $\lie{g}$ of $G$ by continuous linear operators \cite{\SCHTAN, 
p.~13}.   (The continuity of this action again depends on the
fact that $V$ is barrelled.)  This action is functorial in $V$.

\proclaim{Proposition 3.6.10} Suppose that $L = \Q_p$
and that $K$ is local.
Then any continuous action of a locally $L$-analytic
group $G$ on a finite dimensional $K$-vector space is locally analytic.
\endproclaim
\demo{Proof} Replacing $G$ by a compact open subgroup if necessary,
we may assume that $G$ is compact.
Let $n$ denote the dimension of $V$ over $K$, and
choose a basis of $V$ over $K$.  The $G$-action on $V$ then determines
a continuous homomorphism of locally $\Q_p$-analytic groups
$G \rightarrow \GL_n(K)$,  
which is necessarily locally $\Q_p$-analytic
\cite{\SERLG, LG~5.42, thm.~2}.
This implies that the action of $G$ on $V$ is analytic.
\qed
\enddemo

\proclaim{Proposition 3.6.11} If $H$ is an open subgroup of the locally
$L$-analytic group $G$, and $V$ is a barrelled Hausdorff convex $K$-vector
space equipped with a topological $G$-action,
then $V$ is a locally analytic representation of $G$ if and only
if it is a locally analytic representation of $H$.
\endproclaim
\demo{Proof} 
This follows from lemma~3.5.4.
\qed\enddemo

\proclaim{Theorem 3.6.12} If $G$ is a locally $L$-analytic group
and $V$ is a locally analytic representation of $G$
on a Hausdorff $LF$-space,
then the natural map $V_{\la} \rightarrow V$
is a $G$-equivariant topological isomorphism.  
\endproclaim
\demo{Proof} Proposition~3.5.6 shows that $V_{\la}$ is again
an $LF$-space.
Thus the map $V_{\la} \rightarrow V$
is a continuous $G$-equivariant bijection between $LF$-spaces.
Theorem~1.1.17 implies that it is a topological isomorphism.  
\qed\enddemo

The preceding theorem is the analogue
for locally analytic representations of proposition~3.2.11.
Taking into account theorem~3.5.7, it shows that if $V$ is a locally analytic
representation on an $LF$-space then there is a
{\it continuous} orbit map $\orbit: V \rightarrow \La(G,V).$

\proclaim{Corollary 3.6.13} A $K$-Fr\'echet space $V$ equipped with
a topological action of $G$ is a locally analytic representation of $G$
if and only if there is an
analytic open subgroup $H$ of $G$ such that the 
natural map $V_{\H-\an} \rightarrow V$ is an isomorphism.
\endproclaim
\demo{Proof} If such an $H$ exists then proposition~3.6.11 shows
that the $G$-action on $V$ gives rise to a locally analytic representation
of $G$, since by assumption $V$ affords an $\H$-analytic, and so
locally analytic, representation of $H$.
Conversely, suppose that $V$ is endowed with a locally analytic
representation of $G$.  Let $H_n$ be a cofinal sequence of analytic
open subgroups of $G$.  Theorem~3.6.12 shows that there is
an isomorphism $\ilim{n}V_{\H_n-\an} \iso V,$ and propositions~1.1.7
and~3.3.18
then imply that the natural map $V_{\H_n-\an} \rightarrow V$
is an isomorphism for some $n,$ as required.
\qed\enddemo

(In the case of $K$-Banach spaces, this result
is equivalent to \cite{\FETH, kor.~3.1.9}.)

\proclaim{Lemma 3.6.14} Let $V$ be a Hausdorff convex $K$-vector space
equipped with a locally analytic representation of $G$,
and let $W$ be a closed $G$-invariant subspace of $V$.
If $W$ (respectively $V/W$) is barrelled, 
then $W$ (respectively $V/W$) is also a locally analytic $G$-representation.
\endproclaim
\demo{Proof} Proposition~3.5.10 implies that the natural
map $W_{\la} \rightarrow W$ is a bijection.
On the other hand, functoriality of the formation
of locally analytic vectors induces a commutative diagram
$$\xymatrix{ V_{\la} \ar[r]\ar[d] & (V/W)_{\la} \ar[d]\\
V \ar[r] & V/W.}$$
Since the bottom horizontal arrow is surjective, while both vertical arrows
are injective, we see that if
the left-hand vertical arrow is a bijection, the same is true
of the right-hand vertical arrow.
\qed\enddemo

\proclaim{Proposition 3.6.15} If $V$ is a Hausdorff convex $K$-vector space
equipped with a topological $G$-action, and if $W$ is a finite dimensional
locally analytic representation of $G$, then there is a natural isomorphism
$(V \otimes_K W)_{\la} \iso V_{\la} \otimes_K W.$
\endproclaim
\demo{Proof} Since $W$ is finite dimensional, we may find
an analytic open subgroup $H$ of $G$ for which $W$ is $\Bbb H$-analytic.
Let $\{H'\}$ denote the directed set of analytic open subgroups of $H$;
this set is then cofinal in the directed set of analytic open subgroups
of $G$.  Thus we obtain isomorphisms
$$\multline
(V\otimes_K W)_{\la} \iso \ilim{H'} (V\otimes_K W)_{\H'-\an}
\iso \ilim{H'} (V_{\H'-\an} \otimes_K W ) \\
\iso (\ilim{H'} V_{\H'-\an})\otimes_K W \iso V_{\la} \otimes_K W,\endmultline$$
in which the second isomorphism is provided by proposition~3.6.6, and
the third isomorphism follows from the fact that tensor product with
the finite dimensional space $W$ commutes with inductive limits.
\qed\enddemo

\proclaim{Corollary 3.6.16}
If $V$ is a Hausdorff convex $K$-vector space
equipped with a locally analytic representation of $G$,
and if $W$ is a finite dimensional
locally analytic representation of $G$, then the diagonal $G$-action
on $V\otimes_K W$ makes this tensor product a locally-analytic representation
of $G$.
\endproclaim
\demo{Proof}
Since $V$ is barrelled and $W$ is finite dimensional, the tensor product
$V\otimes_K W$ is certainly barrelled.  The corollary is now an immediate
consequence of proposition~3.6.15.
\qed\enddemo

\proclaim{Proposition 3.6.17}
If $V$ is a Hausdorff convex $K$-vector space
equipped with a topological $G$-action, and if
there is an isomorphism $\ilim{i\in I} V_i \iso V,$
where $\{V_i\}_{i \in I}$ is a $G$-equivariant inductive system
of Hausdorff $K$-vector spaces, each equipped with a locally analytic
action of $G$,
then the $G$-action on $V$ is again locally analytic.
\endproclaim
\demo{Proof}
This is proved in the same manner as proposition~3.6.8.
\qed\enddemo

\proclaim{Proposition 3.6.18} If $U$ and $V$ are compact type
spaces, each equipped with a locally analytic representation
of $G$, then the diagonal $G$-action on the completed tensor product
$U\cotimes_K V$ is again a locally analytic
representation of $G$.
\endproclaim
\demo{Proof}
This follows from proposition~3.5.15, and the fact that
the natural maps $U_{\la} \rightarrow U$ and $V_{\la} \rightarrow V$
are both isomorphisms, by theorem~3.6.12.
\qed\enddemo

Suppose now that $E$ is a local subfield of $L$. 
As discussed in subsection~2.3,
we may restrict scalars from $L$ to $E$,
and so regard the locally $L$-analytic group $G$
as a locally $E$-analytic group.  
Thus if $G$ acts on a convex $K$-vector space $V,$ we can speak of the space of
locally analytic vectors in $V$ with respect to the action of $G$ regarded
as a locally $E$-analytic group.  (That is,
with respect to $\Res^L_E G$.)  We will refer to this space as the space
of locally $E$-analytic vectors of $V$, and denote it by $V_{E-\la}$. 
We conclude this subsection with some discussion about the relations between
locally $L$-analytic $G$-representations and locally $E$-analytic
representations.

If $\lie{g}$ denotes the Lie algebra of $G$ (a Lie algebra
over $L$), then the Lie algebra of $\Res^L_E G$ is also equal to
$\lie{g},$ but now regarded as a Lie algebra over $E$.
The Lie algebra of $(\Res^L_E G)_{/L}$ is thus equal to
$L\otimes_E \lie{g}$.  
The natural map $(\Res^L_E G)_{/L} \rightarrow G$ given
by~(2.3.2) induces a map on Lie algebras, which is just the
natural map $L\otimes_E \lie{g} \rightarrow \lie{g}$ given
by the $L$-linear structure on $\lie{g}$.  Let $(L\otimes_E\lie{g})^0$
denote the kernel of this map.

Observe that if $V$ is a Hausdorff convex $K$-vector space equipped
with a locally $E$-analytic action of $G$, then the $\lie{g}$-action
on $V$ extends to a $K\otimes_E \lie{g}$-action, and so in particular
to an $L\otimes_E \lie{g}$-action.

\proclaim{Proposition 3.6.19}
If $V$ is a barrelled Hausdorff convex $K$-vector space equipped with a locally $E$-analytic $G$-action,
then this action is locally $L$-analytic if and only if $V$ is annihilated by
the Lie subalgebra
$(L\otimes_E \lie{g})^0$ of $L\otimes_E \lie{g}$.
\endproclaim
\demo{Proof}
If $\H$ is an analytic open subgroup of $G$, let $\H_0$ denote
the restriction of scalars of $\H$ from $L$ to $E$.
Then $V_{\la} \iso \ilim{W,H} \overline{W}_{\H-\an},$ 
where $W$ ranges over the $G$-invariant $BH$-subspaces of $V$, 
and $H$ ranges over the analytic open subgroups of $G$,
while $V_{E-\la} \iso \ilim{W,H} \overline{W}_{\H_0-\an}.$ 
By assumption, the map $V_{E-\la} \rightarrow V$ is bijective, and
we must show that the map $V_{\la} \rightarrow V$ is also bijective
if and only if $V$ is annihilated by $(L\otimes_E \lie{g})^0$.
For this, it suffices to show, for each fixed $W$ and
$H$, that the natural map
$$\overline{W}_{\H-\an} \rightarrow \overline{W}_{\H_0-\an}$$
is an isomorphism if and only if $\overline{W}_{\H_0-\an}$ is
annihilated by $(L\otimes_E \lie{g})^0$.
This last claim follows from the fact that
$\An(\H,\overline{W})$ is equal to the subspace
$\An(\H_0,\overline{W})^{(L\otimes_E \lie{g})^0}$ of
$\An(\H_0,\overline{W})$ consisting of functions annihilated
by $(L\otimes_E \lie{g})^0$ (under the right regular action),
together with the defining isomorphisms
$\overline{W}_{\H-\an} \iso \An(\H,\overline{W})^{\Delta_{1,2}(H)}$
and
$\overline{W}_{\H_0-\an} \iso \An(\H_0,\overline{W})^{\Delta_{1,2}(H)}.$
\qed\enddemo

\proclaim{Corollary~3.6.20}
If $V \rightarrow W$ is a continuous injection of $G$-representations
on convex $K$-vector spaces, such that $V$ is locally $E$-analytic
and $W$ is locally $L$-analytic, then $V$ is necessarily locally $L$-analytic.
\endproclaim
\demo{Proof}
Proposition~3.6.19 shows that $(L\otimes_E \lie{g})^0$ acts trivially on $V$.  Thus this
Lie algebra also acts trivially on $W$, and by that same proposition, we are done.
\qed\enddemo

\proclaim{Proposition~3.6.21}
Suppose that the Lie algebra of $G$ is semi-simple.
If $0 \rightarrow U \rightarrow V \rightarrow W \rightarrow 0$ 
is a short exact sequence of locally $E$-analytic representations of $G$
on convex $K$-vector spaces,
such both $U$ and $W$ are locally $L$-analytic,
then $V$ is locally $L$-analytic.
\endproclaim
\demo{Proof}
Proposition~3.6.19 shows that $(L\otimes_E \lie{g})^0$ annihilates both $U$ and $W$.
Since $\lie{g}$ is semi-simple by assumption, the same is true of $L\otimes_E \lie{g},$
and hence of its subalgebra $(L\otimes_E \lie{g})^0$.  Thus the extension $V$ of $W$ by
$U$ is also annihilated by $(L\otimes_E \lie{g})^0,$ and appealing again to proposition~3.6.19,
we are done.
\qed\enddemo

The preceding result does not hold in general if $G$ has non-semi-simple
Lie algebra.  For example, suppose that $G = \Cal O_L$,
and suppose that $\psi$ is a $K$-valued character of $\Cal O_L$
that is locally $E$-analytic, but note locally $L$-analytic.
The two dimensional representation of $\Cal O_L$ defined
by the matrix $\pmatrix 1 & \psi \\ 0 & 1 \endpmatrix$
is then an extension of the trivial representation of $\Cal O_L$
by itself, but is not locally $L$-analytic.

\head 4. Smooth, locally finite, and locally algebraic vectors \endhead

\section{4.1} In this subsection we develop some connections
between the notions of smooth and locally finite representations,
and the notion of locally analytic representations.
These connections have already been treated in \cite{\SCHTUF};
our goal is to recall some points from this reference, and to present
some minor additions to its discussion.

\proclaim{Proposition-Definition 4.1.1}
Let $G$ be a topological group and $V$
a $K$-vector space equipped with a $G$-action.

(i) We say that a vector $v \in V$ is smooth if there is an
an open subgroup of $G$ which fixes $v$. 

(ii) The set of all smooth vectors in $V$ forms a
$G$-invariant $K$-linear subspace of $V$, which we denote by
$V_{\sm}$.

(iii) We say that $V$ is a smooth representation of $G$
if every vector in $V$ is smooth; that is, if $V_{\sm} = V$.

(iv) We say that $V$ is an admissible smooth representation of $G$
if it is a smooth representation of $V$ having the additional
property that for any open subgroup $H$ of $G$, the subspace of $H$-fixed
vectors in $V$ is finite dimensional over $K$.
\endproclaim
\demo{Proof}
The only statement to be proved is that of~(ii), stating that
$V_{\sm}$ is a $G$-invariant $K$-subspace of $V$, and this is clear.  
\qed\enddemo

The formation of $V_{\sm}$ is covariantly functorial
in the $G$-representation $V$.

\proclaim{Lemma 4.1.2} If $G$ is a compact topological group
and $V$ is a smooth representation of $G$ then $V$ is semi-simple.
If furthermore $V$ is irreducible, then $V$ is finite dimensional.
\endproclaim
\demo{Proof} This is well-known, and is in any case easily proved.
(It easily reduces to the fact that
representations of finite groups on vector spaces over a field
of characteristic zero are semi-simple.)
\qed\enddemo

\proclaim{Proposition-Definition 4.1.3}
Let $G$ be a locally $L$-analytic group, and let
$V$ be a Hausdorff convex $K$-vector space equipped
with a locally analytic representation of $G$.

(i)
A vector $v\in V$ is said to be $U(\lie{g})$-trivial if it is
annihilated by the action of $\lie{g}$.

(ii) The subset of $V$ consisting of all $U(\lie{g})$-invariant
vectors is a closed $G$-invariant $K$-linear subspace of $V$, which we denote
by $V^{\lie{g}}$.

(iii) The representation $V$ is
said to be $U(\lie{g})$-trivial if every vector in it is $U(\lie{g})$-trivial;
that is, if $V^{\lie{g}} = V$.
\endproclaim
\demo{Proof}
All that is to be proved is that for any locally analytic
representation $V$ of $G$, the set of $U(\lie{g})$-trivial vectors
$V^{\lie{g}}$ is a closed $G$-invariant subspace of $V$.
That $V^{\lie{g}}$ is closed
follows from the fact that the $U(\lie{g})$-action on $V$
is via continuous operators.  That it is $G$-invariant
follows from the formula $X g v = g \Ad_{g^{-1}}(X) v$.
\qed
\enddemo

The following lemma and its corollaries relate definitions~4.1.1
and~4.1.3.

\proclaim{Lemma 4.1.4}
Suppose that $G$ is equal to the
group of $L$-valued points of a connected affinoid rigid analytic group
$\G$ defined over $L$, and that furthermore $G$ is Zariski dense in $\G$.
Let $V$ be a $\G$-analytic representation of $G$.
Then a vector $v$ in $V$ is $U(\lie{g})$-trivial if
and only if $v$ is stabilised by $G$.
\endproclaim
\demo{Proof}
We may find a $G$-invariant $BH$-subspace $W$ of $V$ such that
$v$ lies in the image of the natural map $\overline{W}_{\G-\an}\rightarrow V$.
Consider the tautological $G$-equivariant map
$\overline{W}_{\G-\an} \rightarrow \An(\G,\overline{W})$
(the target being equipped with the right regular $G$-action).
If $v$ is $U(\lie{g})$-trivial, then its image lies
in the $U(\lie{g})$-trivial subspace of $\An(\G,\overline{W}),$
which consists precisely of the constant functions.
Thus it is $G$-invariant. Conversely, if $v$ is $G$-invariant,
then it is certainly $U(\lie{g})$-trivial.
\qed\enddemo

\proclaim{Corollary 4.1.5} 
Suppose that $G$ is equal to the
group of $L$-valued points of a connected affinoid rigid analytic group
$\G$ defined over $L$, and that furthermore $G$ is Zariski dense in $\G$.
If $V$ is a Hausdorff convex $K$-vector
space equipped with a topological $G$-action then the image of
the composite
$(V_{\G-\an})^{\lie{g}} \subset V_{\G-\an} \rightarrow V$
is equal to the closed subspace of $V$ consisting of $G$-fixed vectors.
\endproclaim
\demo{Proof}  Since $V_{\G-\an}$ is a $\G$-analytic representation
of $G$, lemma~4.1.4 shows that $(V_{\G-\an})^{\lie{g}}$
is the subspace of $V_{\G-\an}$ consisting of $G$-fixed vectors.
It remains to be shown that the continuous injection $V_{\G-\an}
\rightarrow V$ induces a bijection on the corresponding subspaces
of $G$-fixed vectors.
That this is so follows from proposition~3.3.3, since a vector $v \in V$
is $G$-fixed if and only if the orbit map $\orbit_v$ is constant, and constant
maps are certainly rigid analytic.
\qed\enddemo

\proclaim{Corollary 4.1.6}
If $G$ is a locally $L$-analytic group and $V$ is a Hausdorff convex $K$-vector
space equipped with a topological action of $G$ then the image of the composite
$(V_{\la})^{\lie{g}} \subset V_{\la} \rightarrow V$ is equal to the subspace
$V_{\sm}$ of $V$.
\endproclaim
\demo{Proof} If we let $H$ run over all the analytic open subgroups
of $G$ then $V_{\la} = \ilim{\H} V_{\H-\an},$ and so there
is a continuous bijection $\ilim{\H} (V_{\H-\an})^{\lie{g}}
\rightarrow (V_{\la})^{\lie{g}}.$  Since the directed set of analytic
open subgroups of $G$ maps cofinally to the directed set of all open
subgroups of $G$, the corollary follows from corollary~4.1.5.
\qed
\enddemo

It is important to note that if we endow $V_{\sm}$ with the topology
it inherits as a subset of $V$ then the natural map
$V_{\la}^{\goth g} \rightarrow V_{\sm},$ which is a $G$-equivariant isomorphism
of abstract $K$-vector spaces, need not be a topological isomorphism.
For example, if $V = \Con(G,K),$ then by proposition~3.5.11
$\Con(G,K)_{\la} \iso \La(G,K)$.  The subspace $\Con(G,K)_{\sm}$ consists of the
locally constant $G$-valued functions.  This is necessarily a closed
subspace of $\La(G,K)$ (being equal to $\La(G,K)^{\lie{g}}$), but
is dense as a subspace of $\Con(G,K)$.

\proclaim{Corollary 4.1.7}
If $G$ is a locally $L$-analytic group and $V$ is a Hausdorff convex
$K$-vector space equipped with a
locally analytic representation of $G$ then a vector $v$ in
$V$ is $U(\lie{g})$-trivial if and only if it is smooth.
In particular, $V_{\sm}$ is the underlying $K$-vector space of the
closed subspace $V^{\lie{g}}$ of $V$.
\endproclaim
\demo{Proof} 
If $V$ is locally analytic then by definition the natural map
$V_{\la} \rightarrow V$ is a continuous bijection, and so the
induced map $(V_{\la})^{\lie{g}} \rightarrow V^{\lie{g}}$ is
also a continuous bijection.  The result then follows from
corollary~4.1.6.
\qed\enddemo

In the remainder of this subsection we discuss the
class of locally finite $G$-representations, which contains
the smooth representations as a subclass.

\proclaim{Proposition-Definition 4.1.8}
Let $G$ be a topological group and $V$
a $K$-vector space equipped with a $G$-action.

(i) We say that a vector $v \in V$ is locally finite if there is an
open subgroup $H$ of $G$, and a finite dimensional $H$-invariant
subspace $W$ of $V$ that contains $v$, such that the $H$-action
on $W$ is continuous (with respect to the natural Hausdorff
topology on the finite dimensional $K$-vector space $W$).

(ii) The set of all locally finite vectors in $V$ forms a
$G$-invariant $K$-linear subspace of $V$, which we denote by
$V_{\lf}$ (or by $V_{G-\lf}$, if we wish to emphasise the group $G$).

(iii) We say that $V$ is a locally finite representation of $G$
if every vector in $V$ is locally finite; that is, if $V_{\lf} = V$.
\endproclaim
\demo{Proof}
The only statement to be proved is that of~(ii), stating that
$V_{\lf}$ is a $G$-invariant $K$-subspace of $V$, and this is clear.  
\qed\enddemo

\proclaim{Proposition 4.1.9}
Let $G$ be a compact topological group, and let $V$ be a locally finite
$G$-representation.

(i) If $W$ is any finite-dimensional $G$-invariant subspace of $V$,
then the $G$-action on $W$ is continuous.

(ii) The natural map $\ilim{W} W \rightarrow V$, where $W$ ranges
over the directed set of
all finite dimensional $G$-invariant subspaces of $V$,
is an isomorphism.
\endproclaim
\demo{Proof}
Suppose that $W$ is a finite dimensional $G$-invariant subspace of $V$.
If $\{w_1,\ldots, w_n\}$ denotes a basis of $W$, then for
each of the vectors $w_i$ we may find an open subgroup $H_i$
of $G$, and a finite dimensional $H_i$-invariant
subspace $W_i$ of $V$ that contains $w_i$, and on which the
$H_i$-action is continuous.  If we let $H = \bigcap_i H_i,$
then the natural map $\bigoplus_i W_i \rightarrow V$ is
an $H$-equivariant map whose image contains $W$, and on whose
source the $H$-action is continuous.  It follows that the
$H$-action on $W$ is continuous.  Since $H$ is open in $G$,
and since any group action on a finite dimensional vector space
is topological, lemma~3.1.3 shows that
the $G$-action on $W$ is also continuous.  This proves~(i).

To prove~(ii), we must show that each vector $v \in V$ is
contained in a finite dimensional $G$-invariant subspace
of $V$.  By assumption, we may find an open subgroup $H$ of
$G$ and a finite dimensional $H$-invariant subspace $W$ of $V$
containing $v$.  Since $G$ is compact, the index of $H$ in $G$
is finite, and so the $G$-invariant subspace of $V$ spanned
by $W$ is again finite dimensional.  This proves~(ii).
\qed\enddemo

\proclaim{Proposition-Definition 4.1.10}
Let $G$ be a locally $L$-analytic group, and let
$V$ be a Hausdorff convex $K$-vector space equipped
with a locally analytic representation of $G$.

(i)
A vector $v\in V$ is said to be $U(\lie{g})$-finite if it is
contained in a $U(\lie{g})$-invariant finite-dimensional subspace
of $V$.

(ii) The set of all $U(\lie{g})$-finite vectors in $V$ forms a
$G$-invariant $K$-linear subspace of $V$, which we denote by
$V_{\lie{g}-\lf}$.

(iii) The representation $V$ is
said to be $U(\lie{g})$-locally finite
if every vector in it is $U(\lie{g})$-finite;
that is, if $V_{\lie{g}-\lf} = V$.

(iv) The representation $V$ is said to $U(\lie{g})$-finite
if there is an ideal of cofinite dimension in $U(\lie{g})$
which annihilates $V$.
\endproclaim
\demo{Proof}
The only statement to be proved is that of~(ii),
which follows from the formula $X g v = g \Ad_{g^{-1}}(X) v$.
\qed\enddemo

The notion of $U(\lie{g})$-finite representation was
introduced and studied in \cite{\SCHTUF, \S 3}.

The following results relate definitions~4.1.8 and~4.1.10.

\proclaim{Lemma 4.1.11} If $G$ is a locally $L$-analytic group,
and if $V$ is a locally analytic representation of $V$,
then a vector $v \in V$ is $U(\lie{g})$-finite if and only
if it locally finite under the action of $G$.
\endproclaim
\demo{Proof}
If $v$ lies in a finite dimensional $G$-invariant subspace of $V$,
then $v$ is certainly $U(\lie{g})$-finite. Conversely, if
$v$ is $U(\lie{g})$-finite, let $U$ denote a finite dimensional
$U(\lie{g})$-invariant subspace of $V$ containing $v$.
If we choose the open subgroup $H$ of $G$ to be sufficiently small,
then we may lift the $U(\lie{g})$-action on $U$ to an action
of $H$.  The space $\Hom(U,V) = \check{U} \otimes_K V$ is
then a locally analytic $H$-representation (by corollary~3.6.16),
and the given inclusion of $U$ into $V$ gives a $U(\lie{g})$-fixed
point of this space. By corollary~4.1.7, replacing $H$ by a 
smaller open subgroup if necessary, we may assume that this point is
in fact $H$-invariant, and thus that $U$ is an $H$-invariant subspace
of $V$.  Thus the $U(\lie{g})$-finite vector $v$ is $G$-locally finite,
as claimed.
\qed\enddemo

\proclaim{Proposition 4.1.12}
If $G$ is a locally $\Q_p$-analytic group, if $K$ is local,
and if $V$ is a Hausdorff convex $K$-vector
space equipped with a topological action of $G$,
then the image of the composite
$(V_{\la})_{\lie{g}-\lf} \subset V_{\la} \rightarrow V$ is equal to the subspace
$V_{\lf}$ of $V$.
\endproclaim
\demo{Proof}  Since $V_{\la}$ is a locally analytic representation
of $G$, lemma~4.1.11 shows that $(V_{\la})_{\lie{g}-\lf}$
is equal to the space of $(V_{\la})_{\lf}$.
It remains to be shown that the continuous injection $V_{\la}
\rightarrow V$ induces a bijection on the corresponding subspaces
of locally finite vectors.

Choose a compact open subgroup $H$ of $G$. 
If $v \in V$ is a locally finite vector, then by proposition~4.1.9,
we may find 
a finite dimensional $H$-invariant subspace $W$ of $V$ containing
$v$ on which $H$ acts continuously.  By proposition~3.6.10,
$W$ is in particular a $BH$-subspace of $V$ equipped with
a locally analytic action of $H$, and so the inclusion
of $W$ into $V$ factors through the natural map
$V_{\la}\rightarrow V$.
\qed
\enddemo

Just in the case of smooth vectors,
it is important to note that if we endow $V_{\lf}$ (respectively
$(V_{\la})_{\lie{g}-\lf}$) with the topology
it inherits as a subspace of $V$ (respectively $V_{\la}$),
then the natural map
$(V_{\la})_{\goth g-\lf} \rightarrow V_{\lf},$
which is a $G$-equivariant isomorphism
of abstract $K$-vector spaces, need not be a topological isomorphism.

The preceding result implies that if $V$ is a $K$-vector
space equipped with a locally finite action of a locally
$\Q_p$-analytic group $G$, and if
we endow $V$ with its finest convex topology, then
the natural map $V_{\la} \rightarrow V$ is a topological
isomorphism, and in particular $V$ becomes a locally analytic
$U(\lie{g})$-locally finite locally analytic representation of $G$.

\proclaim{Proposition 4.1.13}
Suppose that $G$ is a locally $\Q_p$-analytic group.
Let $0 \rightarrow U \rightarrow V \rightarrow W \rightarrow 0$
be a short exact sequence of continuous $G$-representations
on topological $K$-vector spaces.
If $U$ and $W$ are locally finite, then the same is true of
$V$.
\endproclaim
\demo{Proof}
Let $v$ be a vector in $V$, with image $w$ in $W$.
Replacing $G$ by a compact open subgroup if necessary,
we may assume that $G$ is compact,
and so find a finite-dimensional $G$-subrepresentation of $W$
that contains $w$.  Pulling back our extension along the
embedding of this subspace into $W$, we may thus assume
that $W$ is finite-dimensional.  Choose a finite-dimensional
subspace $V'$ of $V$ which contains $v$ and projects
isomorphically onto $W$.

If we tensor our short exact sequence through with
the contragredient $\check{W}$ of $W$ we obtain
a short exact sequence
$$0 \rightarrow U\otimes_K \check{W} \rightarrow V\otimes_K \check{W}
\rightarrow W\otimes_K \check{W} \rightarrow 0.$$
Our choice of $V'$ gives rise to an element $\phi \in V\otimes_K \check{W}$,
which projects onto 
the identity map from $W$ to itself, regarded as an
element of $W\otimes_K \check{W}$.  Since this identity map
is $G$-equivariant, the element $\phi$ gives rise
to a continuous 1-cocycle $\sigma: G \rightarrow U\otimes_K \check{W}$.
Since $G$ is compact and $\Q_p$-locally analytic,
it is in particular topologically
finitely generated.  Since $U$ is a locally finite representation of
$G$, we may find a finite-dimensional $G$-invariant
subrepresentation $U' \subset U$ such that $\sigma$ takes values
in $U'\otimes_K \check{W}$.  

Now consider the short exact sequence
$0 \rightarrow U/U' \rightarrow V/U' \rightarrow W \rightarrow 0$.
By construction, the subspace $V'$ is $G$-invariant
when regarded as a subspace of $V/U'$.
Thus the subspace $U'\bigoplus V'$ of $V$ is  a finite dimensional
$G$-equivariant subrepresentation of $V$ that contains $v$.
\qed\enddemo

\proclaim{Proposition 4.1.14} If $G$ is a compact locally
$\Q_p$-analytic group with semi-simple Lie algebra, and if $K$ is local,
then the category of locally finite $G$-representations is semi-simple,
and any irreducible locally finite $G$-representation is finite dimensional.
\endproclaim
\demo{Proof}
We begin by considering a surjection $V \rightarrow W$ of
finite dimensional continuous representations of $G$.
We must show that this surjection can be split in a $G$-equivariant fashion.
Applying $\Hom(W, \text{--}\,)$ to the given surjection yields
the surjection of finite dimensional continuous
representations $ \Hom(W,V) \rightarrow \Hom(W,W) $.
Proposition~3.6.10 shows that both of these finite-dimensional representations
are in fact locally analytic, and so are
equipped with a natural $\goth g$-action.  Since
$\goth g$ is semi-simple, passing to $\goth g$-invariants yields
a surjection $\Hom_{\goth g}(W,V) \rightarrow \Hom_{\goth g}(W,W).$
This is a surjection of smooth representations of the compact group
$G$.  Since by lemma~4.1.2 smooth representations of a compact group are
semi-simple,
we may find a $G$-invariant map from $W$ to $V$ lifting the identity
map from $W$ to itself.

Now suppose that $V \rightarrow W$ is any surjection of locally
finite $G$-representations.  Taking into account
proposition~4.1.9,
the result of the preceding paragraph
shows that we may find a $G$-equivariant section to this map
over any finite dimensional subrepresentation of $W$.
A simple application of Zorn's lemma now allows us to split this
surjection over all of $W$.

Finally, since $G$ is compact, proposition~4.1.9 implies that
any irreducible locally finite representation of $G$ is finite.
\qed\enddemo

\proclaim{Corollary 4.1.15} If $G$ is a compact locally $\Q_p$-analytic
group with semi-simple Lie algebra, and if $K$ is local,
then any continuous $G$-representation
on a topological $G$-vector space which may be written as
an extension of smooth $G$-representations is again smooth.
\endproclaim
\demo{Proof} Replacing $G$ by a compact open subgroup, if necessary,
we may suppose that $G$ is compact.  Since any smooth representation
is locally finite, the corollary is now seen to be an immediate
consequence of propositions~4.1.13 and~4.1.14.
\qed
\enddemo

\proclaim{Corollary 4.1.16} If $G$ is a locally $\Q_p$-analytic
group with semi-simple Lie algebra, and if $K$ is local, then
any continuous $G$-representation on a topological vector space
which may be written as
an extension of two $G$-representations equipped with trivial
$G$-action is itself equipped with trivial $G$-action.
\endproclaim
\demo{Proof}
Corollary~4.1.15 shows that any such extension is smooth.
By lemma~4.1.2 it is semi-simple.  The corollary follows.
\qed\enddemo

The preceding three results are false for more general locally analytic
groups, as is illustrated by the two-dimensional
representation of the group $\Z_p$ which sends $t \in \Z_p$
to the matrix $\pmatrix 1 & t \\ 0 & 1\endpmatrix.$

\section{4.2} In this subsection we let $\G$ denote a connected reductive
linear algebraic group defined over $L$, and let $G$ denote an open
subgroup of $\G(L)$.
In this situation, $G$ is Zariski dense in $\G$
(since $\G$ is unirational, by \cite{\BOR, thm.~18.2}).  We let $\Aff(\G,K)$
denote the affine ring of $\G$ over $K$.  Since $G$ is Zariski dense
in $\G,$ the restriction of elements of $\Aff(\G,K)$ to $G$
induces an injection $\Aff(\G,K) \rightarrow \La(G,K)$.

Let $\Cal R$ denote the category of finite dimensional algebraic 
representations of $\G$ defined over $K$.  Since $\G$ is reductive
the category $\Cal R$ is a semi-simple abelian category.  Passing
to the underlying $G$-representation is a fully faithful functor,
since $G$ is Zariski dense in $\G$.

\proclaim{Definition 4.2.1} Let $V$ be a $K$-vector space
equipped with an action of $G$ and $W$ be an object of $\Cal R$. 
We say that a vector $v$ in $V$ is locally $W$-algebraic
if there exists an open subgroup $H$ of $G$, a natural number $n$,
and an $H$-equivariant
homomorphism $W^n \rightarrow V$ whose image contains the vector $v$.
We say that the $G$-representation $V$ is locally
$W$-algebraic if each vector of $V$ is locally $W$-algebraic.
\endproclaim

Note that $W$-locally algebraic vectors are in particular $G$-locally
finite, and thus a $W$-locally algebraic representation is a $G$-locally
finite representation.   Thus, if we equip a $W$-locally algebraic
$V$ with its finest convex topology,
then $V$ becomes a locally analytic $G$-representation, which
is $U(\lie{g})$-finite.  (It is annihilated by the annihilator
in $U(\lie{g})$ of the finite dimensional representation $W$.)

For the next several results, we fix an object $W$ of $\Cal R.$

\proclaim{Proposition-Definition 4.2.2} If $V$ is a $K$-vector space equipped
with an action of $G$, then we let $V_{W-\lalg}$ denote the $G$-invariant
subspace of $W$-locally
analytic vectors of $V$.
\endproclaim
\demo{Proof}
We must show that $V_{W-\lalg}$ is a $G$-invariant subspace of $V$.
It is clear that the $W$-locally algebraic vectors do form
a subspace of $V$.
If $v$ lies in $V_{W-\lalg}$, then by assumption there is an open subgroup
$H$ of $G$, a natural number $n$, and an $H$-equivariant map $\phi:W^n
\rightarrow V$ with $v$ in its image. Let $g$ be an element of $G$.
Then $g\phi(g^{-1} \text{--}\,): W^n \rightarrow V$
is $g H g^{-1}$-equivariant and has $g v$ in its image.  Thus
$V_{W-\lalg}$ is $G$-invariant.
\qed\enddemo

The formation of $V_{W-\lalg}$ is covariantly functorial in $V$.
Let $A_1$ denote the image of $U(\lie{g})$ in $\End(W)$, and let $A_2$
denote the image of the group ring $K[G]$ in $\End(W)$.  Since each
of $A_1$ and $A_2$ act faithfully on the semi-simple module $W$, they
are both semi-simple $K$-algebras.
Write $$B: = \End_{A_1}(W) = \End_{\G}(W) = \End_G(W) = \End_{A_2}(W).$$
(The second equality holds because $\G$ is connected, the third equality
holds because $G$ is Zariski dense in $\G$, and the fourth equality
holds by definition of $A_2$.)  The ring $B$ is a semi-simple
$K$-algebra.   The double commutant theorem implies that
$A_1 = A_2 = \End_B(W)$, and from now on we write $A$ in place of
$A_1$ or $A_2$.   
The left $B$-action on $W$ induces a right $B$-action on
the dual space $\check{W}$ (the transpose action).

\proclaim{Proposition 4.2.3}
(i) There is a natural isomorphism $\check{W}\otimes_B W \iso \check{A}$.

(ii) The is a natural
injection $\check{A} \rightarrow \Fun(G)$,
such that the composite
$\check{W}\otimes_B W \iso \check{A} \longrightarrow \Fun(G)
\buildrel \ev_e \over \longrightarrow K$ is the natural evaluation
map $\check{W}\otimes_B W \rightarrow K$.
\endproclaim
\demo{Proof} As already observed, the double commutant theorem implies
that the faithful action of $A$ on $W$ induces an isomorphism
$A \rightarrow \End_B(W);$ dualising this yields the isomorphism
of part~(i).   Dualising the surjection $K[G] \rightarrow A$ induces an
injection of $\check{A}$ into the space of functions $\Fun(G)$ on $G$;
this is the map of part~(ii).  The claim of part~(ii) is easily checked.
\qed\enddemo

Since $G$ acts on $W$ through an algebraic representation
of $\G$, it is immediate that the elements of $\Fun(G)$
in the image of $\check{A}$ via the map of part~(ii) of the preceding
proposition are the restriction
to $G$ of elements of $\Aff(\G)$. 

If $U$ is a smooth $G$-representation on a right $B$-module,
then $U\otimes_B W$ is certainly a locally $W$-algebraic representation of $G$.
Note that we may recover $U$ via the natural isomorphism
$$\Hom(W, U\otimes_B W)_{\sm} \iso U \otimes_B \End_{\lie{g}}(W) =
U \otimes_B B \iso U.$$
There is an important converse to this construction.
Note that if $V$ is a $K$-vector space equipped with a $G$-action,
then $\Hom(W,V)$ 
is equipped with its natural $G$-action, and a commuting right $B$-action
induced by the $G$-equivariant left $B$-action on $W$.
The subspace of smooth vectors $\Hom(W,V)_{\sm}$ is a $G$-invariant
$B$-submodule, on which the $G$-action is smooth.

\proclaim{Proposition 4.2.4} If $V$ is a $K$-vector space equipped 
with a $G$-representation,
then the evaluation map $\Hom(W,V)_{\sm}\otimes_B W
\rightarrow V$ is an injection, with image equal to $V_{W-\lalg}$.
\endproclaim
\demo{Proof}
By definition, the space $\Hom(W,V)_{\sm}$ is the inductive limit of
the spaces $\Hom_H(W,V)$ as $H$ runs over all open subgroups of $G$. 
Thus $\Hom(W,V)_{\sm} \otimes_B W$ is the inductive limit of the spaces
$\Hom_H(W,V) \otimes_B W.$ 

If $\sum_{i=1}^n \phi_i \otimes w_i$ is an element of $\Hom_H(W,V)
\otimes_B W$, then its image in $V$ under the evaluation map
is equal to $\sum_{i=1}^n \phi_i(w_i).$  Let $\phi = 
(\phi_1,\ldots, \phi_n) \in \Hom_H(W,V)^n = \Hom_H(W^n,V),$
and let $w = (w_1,\ldots,w_n) \in W^n$.  Then $\phi(w) = v,$ and
thus $v$ lies in $V_{W-\lalg}$.
Conversely, if $v\in V_{W-\lalg}$ then there exists some such $H$ and an element
$\phi = (\phi_1,\ldots,\phi_n) \in \Hom_H(W^n,V)$ such that 
$v$ lies in the image of $\phi$.  If $w = (w_1,\ldots, w_n)$ is such that
$\phi(w) = v,$ then $\sum_{i=1}^n \phi_i \otimes w_i \in
\Hom_H(W,V)\otimes_B W$ maps to $v$ under the evaluation map
$\Hom_H(W,V) \otimes_B W \rightarrow V$.  
This shows that $V_{W-\lalg}$ is equal to the image of the evaluation map.

It remains to be shown that for each open subgroup $H$ of $G$
the evaluation map $\Hom_H(W,V) \otimes_B W \rightarrow V$ is
injective.  Let $\Fun(H)$ denote the space of $K$-valued functions
on $H$.  Applying proposition~4.2.3 with $H$ in place of $G$,
we obtain an $H\times H$-equivariant
injection $\check{W}\otimes_B W \iso \check{A} \longrightarrow \Fun(H)$.
Tensoring through by $V$ yields an $H\times H \times H$-equivariant injection
$$ \Hom(W,V) \otimes_B W =
V\otimes \check{W} \otimes_B W \rightarrow V\otimes \Fun(H).$$
Now taking $\Delta_{1,2}(H)$-equivariant invariants yields an injection
$\Hom_H(W,V)\otimes_B W \longrightarrow (V\otimes \Fun(H))^{\Delta_{1,2}(H)}
\iso V,$ as required.  (The final isomorphism is provided by the
map $\ev_e$.  That this composite is the natural evaluation map
follows from part~(ii) of proposition~4.2.3.)
\qed\enddemo

\proclaim{Proposition 4.2.5} If $W_1$ and $W_2$ are two objects
of $\Cal R$ then the map $V_{W_1-\lalg}\oplus V_{W_2-\lalg} \rightarrow V$
induced by the inclusions $V_{W_1-\lalg} \rightarrow V$ and
$V_{W_2-\lalg} \rightarrow V$ has image equal to $V_{W_1\oplus W_2 -\lalg}$.
If furthermore $\Hom_{\G}(W_1,W_2)=0,$ then the resulting surjection
$V_{W_1-\lalg} \oplus V_{W_2-\lalg} \rightarrow V_{W_1\oplus W_2-\lalg}$
is in fact an isomorphism.
\endproclaim
\demo{Proof} Obviously both $V_{W_1-\lalg}$ and $V_{W_2-\lalg}$ are contained
in $V_{W_1\oplus W_2-\lalg}.$
Furthermore, the
definition immediately implies
that any element of $V_{W_1\oplus W_2-\lalg}$ may be written
as the sum of an element of $V_{W_1-\lalg}$ and an element of $V_{W_2-\lalg}$.
Thus the image of $V_{W_1-\lalg} \oplus V_{W_2-\lalg}$
is equal to $V_{W_1\oplus W_2-\lalg}$.

Now suppose that $\Hom_{\G}(W_1,W_2) = 0$. 
We must show that the intersection $V_{W_1-\lalg} \bigcap V_{W_2-\lalg}$
is equal to zero. 
If $v$ is 
a vector lying in this intersection then we may find open subgroups $H_i$
of $G$ and $H_i$-equivariant maps
$\phi_i:W_i^{n_i} \rightarrow V$ (for $i=1,2$) that contain $v$ in
their image.  Let $H = H_1\bigcap H_2$, and let $U$ denote the $H$-invariant
subspace of $V$ spanned by $v$.  By our choice of $H$, $U$ is an $H$-equivariant
subquotient of both $W_1^{n_1}$ and $W_2^{n_2}$.
Since the category $\Cal R$ is semi-simple, we may find an $H$-equivariant
surjection $\Cal W_1^{n_1} \rightarrow U$ and an $H$-equivariant embedding
$U \rightarrow \Cal W_2^{n_2}$.  Since $\Hom_{H}(W_1,W_2) = \Hom_{\G}(W_1,
W_2) = 0$ (by assumption, together with the fact that $H$ is Zariski dense),
we conclude that $U=0$, hence that $v=0$, and hence that $V_{W_1-\lalg}
\bigcap V_{W_2-\lalg} = 0.$ 
\qed\enddemo

\proclaim{Proposition-Definition 4.2.6}
Let $V$ be a $K$-vector space equipped with
an action of $G$.

(i)
We say that a vector $v \in V$ is locally algebraic
if it is $W$-locally algebraic for some object $W$ of $\Cal R$.

(ii)
The set of all locally algebraic vectors of $V$ is
a $G$-invariant subspace of $V$, which we by denote
$V_{\lalg}$.

(iii)  We say that 
$V$ that is a locally algebraic representation of $G$
if every vector of $V$ is locally algebraic; that is,
if $V_{\lalg} = V$.
\endproclaim
\demo{Proof}
That $V_{\lalg}$ is a subspace of $V$ follows from proposition~4.2.5.
\qed
\enddemo

The formation of $V_{\lalg}$ is covariantly functorial in $V$.

The class of locally algebraic representations
of $G$ is introduced in the appendix \cite{\PRA} to \cite{\SCHTUF}.
We will see below that our definition of this notion coincides with
that of \cite{\PRA} (when $\Bbb G$ is reductive, as we are assuming
here).

Let $\hat{\G}$
denote a set of isomorphism class representatives for
the irreducible objects of $ \Cal R$.

\proclaim{Corollary 4.2.7} If $V$ is equipped with an action
of $G$, then the natural map $\bigoplus_{W \in \hat{\G}}
V_{W-\lalg} \rightarrow V_{\lalg}$ is an isomorphism.
\endproclaim
\demo{Proof} Since $\Cal R$ is a semi-simple abelian category,
any object $W$ of $\Cal R$ is isomorphic to a direct sum of powers
of elements of $\hat{\G}$.  Given this, the result follows
immediately from proposition~4.2.5.
\qed\enddemo

\proclaim{Proposition 4.2.8}
Let $V$ be an irreducible
locally algebraic representation of $G$.  Then there is an element
$W$ of $\hat{G}$, and, writing $B = \End_{\G}(W),$ an irreducible
smooth representation of $V$ on a right $B$-module $U$ (here irreducibility
is as a $B$-module), such that $V$ is isomorphic to the tensor product
$U\otimes_B W.$
Conversely, given such a $W$ and $U$, the tensor product $U\otimes_B W$
is an irreducible locally algebraic representation of $G$.
\endproclaim
\demo{Proof}
Suppose first that $V$ is an irreducible locally algebraic representation
of $G$.  Then corollary~4.2.7 shows that $V = V_{W-\lalg}$ for some element
$W$ of $\hat{\G}$.  Proposition~4.2.4 then yields an isomorphism
$\Hom(W,V)_{\sm} \otimes_B W \iso V.$  This shows that $V$ is isomorphic
to the tensor product over $B$ of a smooth representation and $W$.  Clearly
$\Hom(W,V)_{\sm}$ must be irreducible, since $V$ is.

Conversely, given such an element $W$ of $\hat{\G}$, and
an irreducible smooth representation on a $B$-module $U$,
consider the tensor product $V = U\otimes_B W$.  We obtain
a natural isomorphism of $B$-modules $U \iso \Hom(W,V)_{\sm}$,
and so the irreducibility of $U$ implies the irreducibility of
$V$.
\qed\enddemo

The following result is an analogue, for locally algebraic
representations, of proposition~4.1.14.

\proclaim{Corollary 4.2.9} Suppose that $G$ is compact.
If $V$ is a locally algebraic representation
of $G$,
then $V$ is semi-simple as a representation of $G$,
and each irreducible summand of $V$ is finite dimensional.
\endproclaim
\demo{Proof} Corollary~4.2.7 shows that $V$ is isomorphic
to the direct sum $\bigoplus_{W\in \hat{\G}} V_{W-\lalg},$ and so we may
restrict our attention to $V$ being $W$-algebraic, for some
irreducible algebraic representation $W$ of $\G$.  
Let $B = \End_{\G}(W)$.  Proposition~4.2.4 shows that the natural map
$\Hom(W,V)_{\sm} \otimes_B  W \rightarrow V$ is an isomorphism.
Since $\Hom(W,V)_{\sm}$ is a smooth representation of the compact
group $G$, lemma~4.1.2 shows that it is a
direct sum of finite dimensional irreducible representations of
$G$ on a $B$-module.
Thus it suffices to show that for any finite
dimensional irreducible smooth representation $S$ of $G$ on a $B$-module,
the tensor product
representation $S\otimes_B W$ is semi-simple.  In fact, proposition~4.2.8
shows that such a representation is irreducible.
\qed\enddemo

The preceding corollary shows that our definition of
locally algebraic representation agrees with that of \cite{\PRA}.
Indeed, the condition of definition~4.2.6~(iii) amounts to
condition~(2) of the definition of \cite{\PRA}, and corollary~4.2.9
shows that condition~(2) implies condition~(1) of that definition.
Our proposition~4.2.8
is thus a restatement of \cite{\PRA, thm.~1} (and our method of proof
coincides with that of \cite{\PRA}).  Note however that \cite{\PRA}
does not consider the possibility that $\End_{\G}(W)$ might be larger than
the ground field $K$.

Suppose now that $V$ is a locally convex $K$-vector
space of compact type equipped with a topological $G$-action,
let $W$ be an object of $\Cal R$, and write $B = \End_G(W)$.
The evaluation map
$\Hom(W,V)_{\sm} \otimes_B W \rightarrow V$
of proposition~4.2.4 is then a continuous $G$-equivariant map
(if we equip $\Hom(W,V)_{\sm}$ with its topology as a subspace
of $\Hom(W,V)$), which that proposition
shows to be a continuous bijection
between its source and the subspace $V_{W-\lalg}$ of its target.

\proclaim{Proposition 4.2.10} If $V$ is a locally convex $K$-vector
space of compact type equipped with a locally analytic representation of $G$,
then for any $W \in \Cal R,$ the evaluation map of proposition~4.2.4
is a closed embedding,
and induces a topological isomorphism between its source
and $V_{W-\lalg}$ (which is thus a closed subspace of~$V$).
\endproclaim
\demo{Proof}
Since $V$ is a locally analytic representation of $G$,
the space $\Hom(W,V) \iso V \otimes_K \check{W}$ is again a
locally analytic representation of $G$, by corollary~3.6.16.
Thus $\Hom(W,V)_{\sm} = \Hom(W,V)^{\lie{g}}$ is a closed
subrepresentation of $\Hom(W,V),$
and so
$\Hom(W,V)_{\sm} \otimes_B W$ is a closed subrepresentation
of $V\otimes_K \check{W} \otimes_B W.$
As in the proof of proposition~4.2.4, we see that the closed embedding
$\Hom(W,V)_{\sm} \otimes_B W \rightarrow V \otimes_K \check{W}
\otimes_B W$ 
has image lying in
$(V \otimes_K \check{W} \otimes_B W)^{\Delta_{1,2}(\lie{g})}
= (V\otimes_K \check{A})^{\Delta_{1,2}(\lie{g})}.$
We claim that the evaluation map 
$$(V\otimes_K \check{A})^{\Delta_{1,2}(\lie{g})} \rightarrow V
\tag 4.2.11$$
is a closed embedding.  Once this is known,
the proposition follows from proposition~4.2.4.

Since $V$ is a locally analytic representation on a space of compact
type, we may write $V \iso \ilim{n} V_n,$ where each 
$V_n$ is a Banach space equipped with a $\H_n$-analytic representation,
for some open analytic subgroup $H_n$ of $G$, and the transition
maps are compact and injective.
Thus (by corollary~3.6.7 and lemma~4.1.4) we obtain an isomorphism
$$ (V\otimes_K \check{A})^{\Delta_{1,2}(\lie{g})}
\iso \ilim{n} (V_n \otimes_K \check{A})^{\Delta_{1,2}(H_n)},$$
again with compact and injective transition maps.
The elements of $\check{A}$, when regarded as functions on $\G$,
are algebraic, and so in particular, they restrict to analytic
functions on each $\H_n$.  Thus for each $n$ we
have a closed embedding $\check{A} \rightarrow \La(\H_n,K)$, 
and hence a closed embedding
$$(V_n \otimes_K \check{A})^{\Delta_{1,2}(H_n)}
\rightarrow \La(\H_n,V_n)^{\Delta_{1,2}(H_n)} \iso V_n,$$
the isomorphism being provided by the evaluation map (and the
fact that $V_n$ is an analytic $\H_n$-representation).
Since the transition maps for varying $n$ on either side of
this closed embedding are compact and injective, and
since 
$$(V_{n+1} \otimes_K \check{A})^{\Delta_{1,2}(H_{n+1})} \bigcap V_n
= (V_n \otimes_K \check{A})^{\Delta_{1,2}(H_n)}$$
(as is easily seen),
the map obtained
after passing to the locally convex inductive limit in $n$
is again a closed embedding.
Thus the evaluation map~(4.2.11) 
is a closed embedding, as claimed.
\qed\enddemo

In the situation of the preceding proposition,
one finds in particular that
for each $W \in \hat{\G}$, the space $V_{W-\lalg}$ is again
a locally analytic representation of $G$ on a space of compact type.
Hence the direct sum $\bigoplus_{W \in \hat{\G}} V_{W-\lalg}$ (equipped
with its locally convex direct sum topology)
is again a locally analytic $G$-representation
on a convex $K$-vector space of compact type.
The isomorphism of corollary~4.2.7 thus allows us to
regard $V_{\lalg}$ as a locally analytic $G$-representation
on a convex $K$-vector space of compact type.  By construction,
the natural map $V_{\lalg} \rightarrow V$ is continuous.

Note that $V_{\lalg}$ is contained in $V_{\lie{g}-\lf}$.
If $\G$ is semi-simple and simply connected, then
the proof of \cite{\SCHTUF, prop.~3.2} shows that in fact
there is equality $V_{\lalg} = V_{\lie{g}-\lf}$.  However, the example
of \cite{\SCHTUF, p.~120} shows that this is not true for arbitrary 
reductive groups~$\G$. 

\head 5. Rings of distributions \endhead

\section{5.1}
In this subsection we describe the convolution structure
on the space of distributions on a group (whether topological,
rigid analytic, or locally analytic),
and those modules over the ring of distributions
that can be constructed out of an appropriate representation
of the group. These module structures can be related to some simple forms
of Frobenius reciprocity, which we also recall.

Let $G$ be a group acting on a $K$-vector space $V$.  If $W$ is
a second vector space then the
function space $\Fun(G,W)$ is endowed with its commuting left and right
regular $G$-actions,  and the space $\Hom_G(V,\Fun(G,W))$ (defined
by regarding $\Fun(G,W)$ as a $G$-representation via the right 
regular $G$-action) is equipped with a $G$-action (via the left 
regular $G$-action on $\Fun(G,W)$).  We will repeatedly apply
this observation below, and the obvious variants in which
$G$ is a topological, rigid analytic, or locally analytic group,
$V$ and $W$ are convex spaces, $\Fun(G,W)$ is replaced by $\Con(G,W)$,
$\An(\G,W)$, or $\La(G,W)$ as appropriate, and $\Hom$ is replaced
by $\Lin$.

The following result gives the archetypal form of Frobenius
reciprocity for continuous, analytic, and locally analytic
representation.

\proclaim{Proposition 5.1.1} 
(i) Let $G$ be a locally compact topological
group and let $V$ be a Hausdorff convex $K$-vector space equipped
with a continuous action of $G$.
If $W$ is any Hausdorff convex
$K$-vector space then the map $\ev_e: \Con(G,W)
\rightarrow W$ induces a $G$-equivariant topological
isomorphism of convex spaces
$\Lin_{G,b}(V,\Con(G,W)) \iso \Lin_b(V,W).$  

(ii) Let $\G$ be an affinoid rigid analytic group over $L$ such that
$G := \G(L)$ is dense in $\G$ and let $V$ be a Hausdorff convex $LF$-space
equipped with an analytic $G$-representation. 
If $W$ is any Hausdorff convex space
then the map $\ev_e:
\Con(G,W) \rightarrow W$ induces a continuous $G$-equivariant bijection
$\Lin_{G,b}(V,\An(\G,W)) \rightarrow \Lin_b(V,W)$.
If $V$ and $W$ are $K$-Fr\'echet
spaces then this map is even a topological isomorphism.

(iii) Let $G$ be a compact locally $L$-analytic group and let $V$ be
a Hausdorff convex $LF$-space
equipped with a locally analytic $G$-representation.
If $W$ is any Hausdorff convex $K$-vector space
then the map $\ev_e: \Con(G,W) \rightarrow W$ induces a $G$-equivariant
continuous bijection 
$\Lin_{G,b}(V,\La(G,W)) \rightarrow \Lin_b(V,W)$.
If $V$ and $W$ are of compact type then this map is even a topological
isomorphism.
\endproclaim
\demo{Proof} We prove each part in turn.  Thus we first suppose that
$G$ is a locally compact topological group, that $V$ is a Hausdorff convex
space equipped with a continuous $G$-action,
and that $W$ is an arbitrary Hausdorff convex space.

Proposition~2.1.7 shows that the natural map
$\Lin_b(V,W) \rightarrow \Lin_b(\Con(G,V),\Con(G,W))$ is continuous.
Note that its image lies in $\Lin_{G,b}(\Con(G,V),\Con(G,V)).$
Proposition~3.2.10 shows that $\orbit: V \rightarrow \Con(G,V)$
is continuous and $G$-equivariant, and so induces a continuous
map $\Lin_{G,b}(\Con(G,V),\Con(G,W)) \rightarrow \Lin_{G,b}(V,\Con(G,W))$.
Composing this with the preceding map  yields a continuous map
$\Lin_b(V,W) \rightarrow \Lin_{G,b}(V,\Con(G,W)).$  It is easily checked
that this is inverse to the continuous map
$\Lin_{G,b}(V,\Con(G,W)) \rightarrow \Lin_b(V,W)$ induced by $\ev_e$,
and so part~(i) of the proposition is proved.

The proof of part~(ii) is similar.  The functoriality
of the construction of $\An(\G,\text{--}\,)$ yields a natural
map of abstract $K$-vector spaces
$$\Lin_b(V,W) \rightarrow \Lin_{G,b}(\An(\G,V),\An(\G,W)).\tag 5.1.2$$
Since $V$ is an $LF$-space, theorem~3.6.3 shows
that the orbit map $\orbit:V \rightarrow \An(\G,V)$ is a continuous
embedding, and combining this with~(5.1.2) we obtain a
$K$-linear map
$\Lin_b(V,W) \rightarrow \Lin_{G,b}(V,\An(\G,W)),$ which is easily
checked to provide a $K$-linear inverse to the continuous
map $\Lin_{G,b}(V,\Con(G,W)) \rightarrow \Lin_b(V,W)$ induced
by $\ev_e$.  Thus this latter map is a continuous bijection,
as claimed.  If $V$ and $W$ are Fr\'echet spaces then the map~(5.1.2)
is continuous, by proposition~2.1.24, and so in this case
we even obtain a topological isomorphism.  This completes
the proof of part~(ii).

The proof of part~(iii) is similar again.  One uses theorem~3.6.12
and proposition~2.1.31.
\qed\enddemo

Let us put ourselves in the situation of part~(i) of the
preceding proposition,
with $W=K$.
Proposition 1.1.36 shows that passing to the transpose yields a
topological embedding
$$\Lin_{G,b}(V,\Con(G,K)) \rightarrow \Lin_{G,b}(\DCon(G,K)_b,V'_b).\tag 5.1.3$$
Thus part~(i) yields a topological embedding
$V'_b \rightarrow \Lin_{G,b}(\DCon(G,K)_b,V'_b),$
and hence a map
$$\DCon(G,K)_b \times V'_b \rightarrow V'_b,\tag 5.1.4$$
that is $G$-equivariant in the first variable, and 
separately continuous.
If $\mu \in \DCon(G,K)$ and $v'\in V'$ we let $\mu*v'$ denote the image of
the pair $(\mu,v')$.
Tracing through the definitions shows that for any element $v\in V$
we have
$$ \langle \mu * v', v \rangle = \int_G \langle v' , g v \rangle \d \mu(g).
\tag 5.1.5$$
In particular if $g\in G$
then $\delta_g * v' = g^{-1}v'$.

If we take $V$ to be $\Con(G,K)$,
we obtain as a special case of~(5.1.4) a map
$$\DCon(G,K)_b \times \DCon(G,K)_b \rightarrow \DCon(G,K)_b,$$
which we denote by $(\mu,\nu) \mapsto \mu*\nu$,
which has the property that $\delta_g \times \delta_h \mapsto \delta_{h g}.$
As a special case of~(5.1.5) we find that
$$\int_G f(g) \d (\mu * \nu)(g) = \int_G \left(
\int_G f(h g) \d \mu(h) \right) \d \nu(g).  \tag~5.1.6$$

\proclaim{Corollary~5.1.7} Suppose that $G$ is a locally compact
topological group.

(i) The formula~(5.1.6) defines an associative
product on $\DCon(G,K)_b$, which 
is separately continuous in each of its variables.
If $G$ is compact then it is even jointly continuous.

(ii) If $V$ is a Hausdorff convex space equipped with a continuous
$G$-action then the map~(5.1.4) makes $V'_b$ a left $\DCon(G,K)_b$-module,
and this map is separately continuous in each of its variables.
If $G$ is compact then it is jointly continuous.

(iii) If $V$ is as in~(ii), then there is a natural isomorphism
of $K$-vector spaces $\Hom_{\DCon(G,K)}(\DCon(G,K),V') \iso
\Lin_G(V,\Con(G,K)).$
\endproclaim
\demo{Proof}
We have already observed that the statements concerning
separate continuity  are consequences
of part~(i) of proposition~5.1.1.  The associativity statement
of part~(i) can be checked directly from the formula~(5.1.6),
and a similarly direct calculation with the formula~(5.1.5)
shows that~(5.1.4) makes $V'_b$ a left $\DCon(G,K)$-module.

If $G$ is compact then $\DCon(G,K)_b$ is a Banach space.
Since the product~(5.1.4) is defined by the continuous map~(5.1.3),
we see that given any bounded subset $A$ of $\DCon(G,K)$,
and any open subset $U_1$ of $V'_b$,
there is an open subset $U_2$ of $V'_b$ such that
$A\times U_2$ maps into $U_1$ under~(5.1.4).  Since $\DCon(G,K)$
is a Banach space, it has a neighbourhood basis of the origin
consisting of bounded subsets, and thus~(5.1.4) is actually
jointly continuous.  This completes the proof of parts~(i) and~(ii).

To prove part~(iii), note that tautologically
$\Hom_{\DCon(G,K)}(\DCon(G,K),V') \iso V',$ and that
by part~(i) of proposition~5.1.1 there is a $K$-linear isomorphism
$V' = \Lin(V,K) \iso \Lin_G(V,\Con(G,K)).$
\qed\enddemo

In the case where $K$ is local and $V$ is a Banach space,
this result almost coincides with the construction of
\cite{\SCHTIW,\S 2}. The one distinction is that we work with
the strong topology on our dual spaces, while in \cite{\SCHTIW}
the authors work with the bounded weak topology on these dual
spaces.

Analogous results hold in the situations of parts~(ii) and~(iii)
of proposition~5.1.1.  

\proclaim{Corollary~5.1.8} Suppose that $\G$ is an affinoid rigid
analytic group over $L$ such that $G := \G(L)$ is dense in $G$.

(i)  Formula~(5.1.6) defines
an associative product on $\DAn(\G,K)_b$ which
is jointly continuous as a map $\DAn(\G,K)_b \times
\DAn(\G,K)_b \rightarrow \DAn(\G,K)_b.$

(ii) If we are given an analytic $G$-representation on a Hausdorff
convex $LF$-space~$V$,
then formula~(5.1.5) makes $V'_b$ a left $\DAn(\G,K)$-module.
The resulting map
$$\DAn(\G,K)_b \times V'_b \rightarrow V'_b$$
is continuous in its first variable.
If $V$ is a Fr\'echet space then it is jointly continuous.

(iii) If $V$ is as in~(ii), then there is a natural isomorphism
of $K$-vector spaces $\Hom_{\DAn(G,K)}(\DAn(G,K),V') \iso
\Lin_G(V,\An(G,K)).$
\endproclaim
\demo{Proof} The proof of this result follows the same lines
as that of corollary~5.1.7, using part~(ii) of proposition~5.1.1.
In order to get the joint continuity statements of~(i) and~(ii),
one should take into account that $\DAn(\G,K)_b$ is a Banach space,
being the strong dual of the Banach space $\An(\G,K)$.
\qed\enddemo

\proclaim{Corollary~5.1.9} Suppose that $G$ is a 
locally $L$-analytic group.

(i)  Formula~(5.1.6) defines
an associative product on $\DLa(G,K)_b$ which
is separately continuous as a map $\DLa(G,K)_b \times
\DLa(G,K)_b \rightarrow \DLa(G,K)_b.$
If $G$ is compact then it is even jointly continuous.

(ii) If we are given a locally analytic representation of $G$ on a
Hausdorff convex $LF$-space $V$,
then formula~(5.1.5) makes $V'_b$ a left $\DLa(G,K)$-module.
The resulting map
$$\DLa(G,K)_b \times V'_b \rightarrow V'_b$$
is continuous in its first variable.
If $V$ is of compact type then it is separately continuous in each
variable, and if furthermore $G$ is compact then it is jointly continuous.

(iii) If $V$ is as in~(ii), then there is a natural isomorphism
of $K$-vector spaces $\Hom_{\DLa(G,K)}(\DLa(G,K),V') \iso
\Lin_G(V,\La(G,K)).$
\endproclaim
\demo{Proof} Suppose first that $G$ is compact.
The proof of this result follows the same lines
as that of corollary~5.1.7, using part~(iii) of proposition~5.1.1.
In order to get the joint continuity statements of~(i) and~(ii), one takes into
account the facts that because $\La(G,K)$ is of compact type its strong dual
$\DLa(G,K)_b$ is a Fr\'echet space, that if $V$ is of
compact type then $V'_b$ is a Fr\'echet space, and that a separately
continuous map on a product of Fr\'echet spaces is necessarily
jointly continuous.

If $G$ is not compact, one can still check that~(5.1.5) and~(5.1.6)
define an associative product and a left module structure respectively.
In order to see the required continuity properties,
choose a compact open subgroup $H$ of $G$,
and let $G = \coprod g_i H$ be a decomposition of $G$
into right $H$-cosets.  Then $\DLa(G)_b \iso \bigoplus \DLa(g_i H)_b
= \bigoplus \DLa(H)_b *\delta_{g_i}.$ 
Now right multiplication by $\delta_{g_i}$ acts as an automorphism
of $\DLa(G)$ (induced by the automorphism of $G$ given by left multiplication
by $g_i$), and so we see that separate continuity of the multiplication
on $\DLa(G)_b$ follows from the joint continuity of the multiplication
on $\DLa(H)_b$.  Similarly, if $V$ is an $LF$-space
equipped with a locally analytic $G$-action, then the multiplication
by $\delta_{g_i}$ is an automorphism of $V'_b$ (it is just given by the 
contragredient action of $g_i^{-1}$ on $V'_b$), and so the continuity
properties of the $\DLa(G)_b$ action on $V'_b$ can be deduced
from the continuity properties of the $\DLa(H)_b$-action.
Thus parts~(i) and~(ii) are proved for arbitrary $G$.

The proof of part~(iii) proceeds by a similar reduction.
One notes that $$\Hom_{\DLa(G,K)}(\DLa(G,K),V')\iso 
V' \iso \Lin_H(V,\La(H,K)) \iso \Lin_G(V,\La(G,K)),$$
where the last isomorphism can be deduced immediately using
the decomposition $G = \coprod g_i H$ and the corresponding
isomorphism $\La(G,K) \iso \prod \La(g_i H,K).$
\qed\enddemo

Part~(i) of the preceding result is originally due to F\'eaux de Lacroix
\cite{\FEDI, 4.2.1}, and is recalled in \cite{\SCHTAN, \S 2}.
We have included a proof here, since the work of de Lacroix is unpublished.
The case of part~(ii) when $V$ is of compact type is
originally due to Schneider and Teitelbaum \cite{\SCHTAN, cor.~3.3}.

There are alternative approaches to obtaining a multiplicative
structure on the distribution spaces $\DCon(G,K),$ $\DAn(\G,K)$
and $\DLa(G,K)$.  For example, if $G$ is a compact topological group
one can use the isomorphism between $\Con(G,K)\hat{\otimes} \Con(G,K)$
and $\Con(G\times G,K)$ to obtain a coproduct on $\Con(G,K),$
which dualises to yield the product on $\DCon(G,K)$, and similar
approaches are possible in the analytic or locally analytic
cases.  (This is the approach to part~(i) of corollary~5.1.9 explained
in \cite{\SCHTAN, \S 2}.)

Similarly, one can obtain the module structure on $V'$ in other ways.
For example, if $V$ is a Hausdorff convex $K$-vector space equipped with a
continuous action of the topological group $G$, then
thinking of $\DCon(G,K)$ as a group ring, 
it is natural to attempt to extend the contragredient
$G$-action on $V'$ to a left $\DCon(G,K)$-module structure.
This is the approach adopted in \cite{\SCHTIW}, and in \cite{\SCHTAN, \S 3}  
in the locally analytic situation.

We also remark that the product structures we have obtained
on our distribution spaces are opposite to the ones usually considered.
(For example, $\delta_g * \delta_h  = \delta_{h g},$ rather than
$\delta_{g h},$ and $\delta_g * v' = g^{-1} v'$ rather than $g v'$.)
If we compose these structures
with the anti-involution obtained by sending $g \mapsto g^{-1}$,
then we obtain the usual product structures. 

Let us recall from \cite{\SCHTAN}
the following converse to corollary~5.1.9~(ii).

\proclaim{Proposition 5.1.10} If $V$ is a space of compact type,
and $V'_b$ is equipped with a $\DLa(G)$-module structure for which
the corresponding product map $\DLa(G) \times V'_b \rightarrow V'_b$
is separately continuous, then the $G$-action on $V$ obtained
by passing to the transpose makes $V$ a locally analytic representation
of $G$.
\endproclaim
\demo{Proof} This is \cite{\SCHTAN, cor.~3.3}.  We recall a proof.
We may replace $G$ by a compact open subgroup, and thus assume that
$G$ is compact.  Then $\DLa(G)$ and $V'_b$ are both nuclear Fr\'echet
spaces, and the product map $\DLa(G) \times V'_b \rightarrow V'_b$
induces a surjection $\DLa(G) \cotimes_K V'_b \rightarrow V'_b$.
Proposition~1.1.32 shows that the source is a nuclear Fr\'echet space,
and so this surjection is strict, by the open mapping theorem.
Passing to duals, and taking into account propositions~1.1.32 and~2.1.28,
we obtain a closed embedding
$V \rightarrow \La(G,V),$
which is $G$-equivariant, with respect to the transposed $G$-action
on $V$ and the right regular $G$-action on $\La(G,V)$.
Proposition~3.5.11 and lemma~3.6.14 shows that $V$ is a locally analytic
representation of $G$.
\qed\enddemo

We have the following application of Frobenius
reciprocity for locally analytic representations.

\proclaim{Proposition 5.1.11} If $V$ is a $K$-Banach space equipped
with a locally analytic action of the locally $L$-analytic group
$G$ then the contragredient action of $G$ on $V'_b$ is again
locally analytic.
\endproclaim
\demo{Proof} We may and do assume that $G$ is compact.
Frobenius reciprocity (more precisely, proposition~5.1.1~(iii), with $W = K$)
yields a $G$-equivariant isomorphism
of abstract $K$-vector spaces
$$V' \iso \Lin_G(V,\La(G,K)).\tag 5.1.12$$
Since $\La(G,V)$ is reflexive, passing
to the transpose yields a $G$-equivariant isomorphism
$$\Lin_G(V,\La(G,K)) \iso \Lin_G(\DLa(G,K)_b,V'_b).\tag 5.1.13 $$
Proposition~2.2.10 yields a natural isomorphism
$\Lin(\DLa(G,K)_b,V'_b) \iso \La(G,V'_b),$ and hence
a $G$-equivariant isomorphism (of abstract $K$-vector spaces)
$$\Lin_G(\DLa(G,K)_b,V'_b) \iso \La(G,V'_b)^{\Delta_{1,2}(G)}
\iso (V'_b)_{\la}.\tag 5.1.14$$
Working through the definitions,
one checks that the composite of~(5.1.12), (5.1.13), and~(5.1.14)
provides a $K$-linear inverse to the natural map $(V'_b)_{\la}
\rightarrow V'_b,$ and thus that this natural map is an isomorphism
of abstract $K$-vector spaces.  Thus $V'_b$ is equipped with a
locally analytic representation of $G$, as claimed.
\qed\enddemo

A version of this result is proved in the course of proving
\cite{\SCHTBD, prop.~3.8}.

We remark that in the context of the preceding proposition, corollary~3.6.13
guarantees that there is an analytic open subgroup $H$ of $G$ such
that $V$ is $\H$-analytic.  However, it is not the case in general
that $V'_b$ is also $\H$-analytic (although again by corollary~3.6.13
it will be $\H'$-analytic for some analytic open subgroup $H'$ of $H$).

We can use proposition~5.1.11 to give an  example of a non-trivial locally
analytic representation on a complete $LB$-space
that is neither of compact type nor a Banach space.  For this,
let $G$ be a locally $L$-analytic group which is $\sigma$-compact,
let $H$ be an analytic open subgroup of $G$, and let $G = \coprod_i
H g_i$ be the decomposition of $G$ into (countably many) left $H$ cosets.
The right multiplication of $g_i^{-1} H g_i$ on $H g_i$ makes
$\An(\H g_i)$ an analytic $g_i^{-1} \H g_i$-representation on a Banach space,
and so by proposition~5.1.11 the contragredient representation on
$\DAn(\H g_i,K)_b$
is a locally analytic $g_i^{-1} H g_i$-representation.  
The right multiplication of $G$ on itself induces a natural right
$G$-action on the countable direct sum
$\bigoplus_i \DAn(\H g_i,K)_b,$  extending the $g_i^{-1} H g_i$-action
on the $i^{th}$ summand, and thus this direct sum is a locally analytic
$G$-representation. 
Since it is a direct sum of a countable collection
of Banach spaces it is a complete $LB$-space \cite{\TVS, prop.~9,
p.~II.32}.
It is not of compact type.

\proclaim{Theorem 5.1.15}
(i) Let $G$ be a compact topological
group and let $V$ be a Hausdorff $LB$-space
equipped with a continuous $G$-action.  If we equip
$V'$ with the $\DCon(G,K)$-action of corollary~5.1.7, then
any surjection of left $\DCon(G,K)$-modules $\DCon(G,K)^n \rightarrow V'$
(for any natural number $n$) is obtained by dualising a closed $G$-equivariant
embedding $V \rightarrow \Con(G,K)^n$.  In particular,
$V'$ is finitely generated as a left $\DCon(G,K)$-module if and only
if $V$ admits a closed $G$-equivariant embedding into $\Con(G,K)^n$
for some natural number~$n$.

(ii) Let $\G$ be an affinoid rigid analytic group defined
over $L$, and assume that $G := \G(L)$ is Zariski dense
in $\G$.  Let $V$ be a Hausdorff $LB$-space
equipped with a $\G$-analytic $G$-action.  If we equip
$V'$ with the $\DAn(\G,K)$-action of corollary~5.1.8 then
any surjection of left $\DAn(\G,K)$-modules $\DAn(\G,K)^n \rightarrow V'$
for some natural number $n$ arises by dualising a closed $G$-equivariant
embedding $V \rightarrow \An(\G,K)^n$.  In particular,
$V'$ is finitely generated as a left $\DAn(\G,K)$-module if and only
if $V$ admits a closed $G$-equivariant embedding into $\An(\G,K)^n$
for some natural number~$n$.

(iii) Let $G$ be a compact locally $L$-analytic 
group, and let $V$ be a Hausdorff $LB$-space
equipped with a locally analytic $G$-action.
If we equip $V'$ with the $\DLa(G,K)$-action of corollary~5.1.9, then
any surjection $\DLa(G,K)^n \rightarrow V'$
of left $\DLa(G,K)$-modules,
for some natural number $n$, arises by dualising a closed $G$-equivariant
embedding $V \rightarrow \La(G,K)^n$.  In particular,
$V'$ is finitely generated as a left $\DLa(G,K)$-module if and only
if $V$ admits a closed $G$-equivariant embedding into $\La(G,K)$
for some natural number $n$.
\endproclaim
\demo{Proof}
Suppose first that $G$ and $V$ are as in~(i), and that
we are given a surjection of left $\DCon(G,K)$-modules
$$\DCon(G,K)^n_b \rightarrow V'_b.\tag 5.1.16 $$
Since the map $\DCon(G,K)^n_b \times V'_b \rightarrow V'_b$ describing
$V'_b$ as a left $\DCon(G,K)^n_b$-module is continuous in its first variable,
the surjection~(5.1.16) is necessarily continuous.

Write $V = \ilim{n} V_n$ as the locally convex inductive limit
of a sequence of Banach spaces.  Dualising we obtain a continuous
bijection $$V'_b \rightarrow \plim{n} (V_n)'_b.\tag 5.1.17$$
Since each $V_n$ is a Banach space, so is each $(V_n)'_b,$ and
so $\plim{n} (V_n)'_b$ is the projective limit of a sequence of
Banach spaces, and so is a Fr\'echet space. 
The composite of~(5.1.16) and~(5.1.17) is a surjection
whose source and target are Fr\'echet spaces,
and so by the open mapping theorem it is strict.  Thus each map
in this composite must itself be strict.

Dualising the strict surjection~(5.1.16) 
yields a closed $G$-equivariant embedding
$$(V'_b)' \rightarrow (\DCon(G,K)_b^n)' \iso ((\Con(G,K)'_b)'_b)^n.
\tag 5.1.18$$

Let $\imath_i: \DCon(G,K)_b \rightarrow V'_b$ denote the $i^{th}$
component of~(5.1.16), where $1 \leq i \leq n$.
By part~(iii) of corollary~5.1.7 each $\imath_i$ is obtained by dualising
a continuous $G$-equivariant map $V \rightarrow \Con(G,K).$  Taking the
direct sum of these we obtain a continuous $G$-equivariant map
$$V \rightarrow \Con(G,K)^n,\tag 5.1.19$$ and our original
surjection (5.1.16) is obtained by dualising~(5.1.19).
We can embed the maps~(5.1.18) and~(5.1.19) into the following diagram
$$\xymatrix{V \ar[r]\ar[d] & \Con(G,K)^n \ar[d] \\
(V'_b)'_b \ar[r] & ((\Con(G,K)_b')'_b)^n,}$$ in which
the vertical arrows are double duality maps, which commutes
by construction. The bottom arrow is a topological embedding,
and since $V$ and $\Con(G,K)$ are both barrelled, so are
the vertical arrows.  Thus the top arrow is an embedding.
It remains to prove that it has closed image.  But since
$V$ is an $LB$-space by assumption, and normable, since
it can be embedded in the Banach space $\Con(G,K)^n$, proposition~1.1.18
implies that $V$ is a Banach space, and so must have closed image
in $\Con(G,K)^n$. 

Conversely, if we dualise a closed $G$-equivariant embedding
$V \rightarrow \Con(G,K)^n,$ then we certainly obtain a surjection
of left $\DCon(G,K)$-modules $\DCon(G,K)^n \rightarrow V'.$ 
This proves part~(i).

Parts~(ii) and~(iii) are proved in an identical fashion.
(In fact, the proof of part~(iii) is even a little simpler,
since $\DLa(G,K)$ is reflexive.)
\qed\enddemo

\section{5.2}  In this subsection we explain the relationship
between the ring of analytic distributions on an analytic open
subgroup of a locally $L$-analytic group, and certain completions
of universal enveloping algebras.

We begin with some simple rigid analysis.
For any integer $d$, rigid analytic affine
$d$-space $\A^d$ represents the functor
that attaches to any rigid analytic space
$X$ over $L$ the set $\Gamma(X,\Cal O_X)^d$
of $d$-tuples of globally defined rigid analytic
functions on $X$.  If $r$ is any positive real number,
then the closed ball of radius $r$ centred at the origin
(denoted by
$\B_d(r)$) is an open subdomain of $\A^d$.
It represents the functor that attaches to $X$
the set $(\Gamma(X,\Cal O_X)^{|\,|_{\sup} \leq r})^d$
consisting of $d$-tuples of globally defined
rigid analytic functions on $X$, each of whose
sup-norm is bounded above by $r$.

There is a coordinate-free version
of these constructions.
Fix a free finite rank $\Cal O_L$-module $M$,
of rank $d$ say, and write $V = L \otimes_{\Cal O_L} M$.  
The functor
$$X \mapsto \Gamma(X,\Cal O_X) \otimes_L V$$
is represented by a rigid analytic space that we denote
$\A(V)$.  The $L$-valued points of $\A(V)$ are canonically
identified with $V$, as is the tangent space to $\A(V)$ at any $L$-valued
point of $\A(V)$.

If $r$ is a positive real number, then the functor
$$X \mapsto \Gamma(X,\Cal O_X)^{|\,|_{\sup} \leq r}
\otimes_{\Cal O_L} M $$
is represented by an open subdomain of $\A(V),$
that we denote by $\B(M,r).$  The $L$-valued points of 
$\B(M,1)$ are canonically identified with $M$.

Let us return now to the consideration of a locally $L$-analytic group $G$.
Suppose that $\lie{h}$ is an $\Cal O_L$-Lie subalgebra
of $\lie{g}$, free of finite rank as an $\Cal O_L$-module,
such that the natural map $L \otimes_{\Cal O_L} \lie{h} \rightarrow \lie{g}$
is an isomorphism, and such that $[\lie{h},\lie{h}]
\subset a \lie{h}$ for some element $a \in \Cal O_L$.
If $\ord_L(a)$ is sufficiently large, 
then it is proved in \cite{\SERLG, LG Ch.~V, \S~4}
that we can use the Baker-Campbell-Hausdorff formula 
to define a rigid analytic group structure on $\B(\lie{h},1)$,
for which the induced Lie algebra structure on $\lie{g}$
(thought of as the tangent space to the identity
in $\B(\lie{h},1)$) agrees with the given Lie algebra structure on
$\lie{g}.$
We let $\H$ denote $\B(\lie{h},1)$ equipped with
this rigid analytic group structure, and as usual
we let $H$ denote the underlying locally analytic group.
(Thus $H$ and $\lie{h}$ are equal as sets;
$\lie{h}$ denotes this set regarded as an $\Cal O_L$-Lie algebra,
while $H$ denotes this set regarded as a locally $L$-analytic group.)
For any real number $0 < r < 1,$ 
the subdomain $\B(\lie{h},r)$ is in fact an open subgroup of $\B(\lie{h},1)$;
we denote this open subgroup by $\H_r$, and its underlying locally
$L$-analytic group of $L$-rational points by $H_r$.
(Here and below, let us make the convention that we only consider
values of $r$ that arise as the absolute values of elements
in the algebraic closure of $L$.  These elements are dense
in the interval $(0,1)$, 
and have the merit that the corresponding groups $\H_r$ are affinoid
\cite{\BGR, thm.~6.1.5/4}.)

Again under the hypothesis that $\ord_L(a)$ is sufficiently large,
it is proved in \cite{\SERLG, LG~5.35, cor.~2} that
we may construct an embedding of locally $L$-analytic groups
$\exp: H \rightarrow G$,
and thus realise
$H$ as an analytic open subgroup of $G$.
We fix such an identification from now on,
and refer to analytic open subgroups of $G$ constructed in this
manner as good analytic open subgroups.  The discussion of \cite{\SERLG,
LG Ch.~IV, \S 8} shows that the set of good analytic
open subgroups is cofinal in the directed set of all
analytic open subgroups of $G$.

If we fix a basis $X_1,\ldots,X_n$ for $\lie{h}$,
then then we may define canonical coordinates of the
second kind on $H$.  More precisely,
we define a rigid analytic isomorphism 
between $\B_d(1)$ and $\H$ via the map
$$(t_1,\ldots,t_d) \mapsto \exp(t_1 X_1) \cdots \exp(t_d X_d).$$
For any $0 < r < 1$ this map restricts to an isomorphism
between $\B_d(r)$ and $\H_r$.
Via this isomorphism, we obtain the following description
of the Banach space of $K$-valued rigid analytic functions on the
affinoid group $\H_r$:
$$\An(\H_r,K) \iso
\{ \sum_{I=(i_1,\ldots,i_d)} a_I t_1^{i_1}\cdots t_d^{i_d} \quad | \quad
a_I \in K, \, \lim_{|I| \rightarrow \infty} |a_I| r^{|I|} = 0\}.\tag 5.2.1$$
(Here $I$ runs over all multi-indices of length $d$,
and $|I| = i_1  + \ldots + i_d$ for any multi-index $I = (i_1,\ldots,i_d).$)

The right regular representation of $H_r$ on $\An(\H_r,K)$ makes
this latter space a locally analytic $H$-representation,
and thus induces a continuous action of $U(\lie{g})$ on $\An(\H_r,K)$.
Concretely, an element $X \in \lie{g}$ acts on $\An(\H_r,K)$ by the 
endomorphism
$$(X f)(h) = \dfrac{d}{d t} f(h \exp(t X) )_{|t = 0}.$$
(See the discussion of \cite{\SCHTAN, p.~449} -- but note that
we consider rigid analytic functions, rather than locally
analytic functions, and consider $\An(\H_r,K)$ as an $H_r$-representation
via the right regular action rather than the left regular action.)
This in turn yields a morphism
$U(\lie{g})\rightarrow \DAn(\H_r,K),$
defined via
$$\langle X, f \rangle := (X f) (e).$$

In particular, any monomial $X_1^{i_1}\cdots X_d^{i_d}$
yields an element of $\DAn(\H_r,K)$.  
The entire space $\DAn(\H_r,K)$ admits the following
description in terms of such monomials:
$$\DAn(\H_r,K) \iso 
\{\sum_{I = (i_1,\ldots,i_d)} b_I
\dfrac{X_1^{i_1}\cdots X_d^{i_d}}
{i_1! \cdots i_d!}
\quad | \quad 
|b_I| \leq C r^{|I|} \text{ for some } C > 0\}.$$
(Compare the discussion of \cite{\SCHTAN, p.~440}.)
In fact this is an isomorphism not just of $K$-vector spaces,
but of $K$-algebras.  (Note that the $K$-algebra structure on $U(\lie{g})$
yields a $K$-algebra structure on the space of series in
the monomials $X_1^{i_1}\cdots X_d^{i_d}$ appearing on
the right hand side of this isomorphism.)

Actually, it will be important for us to consider not only 
the affinoid groups $\H_r$, 
but also certain corresponding strictly $\sigma$-affinoid
subgroups.
For any value of $r$, we define
$\open{\H}_{r} = \bigcup_{r' < r} \H_{r'}$.
The open subspace $\open{\H}_r$ is an
open subgroup of $\H$, isomorphic as a
rigid analytic space to the open ball of radius $r$.
We write $\open{H}_r := \open{\H}_r(L)$ to denote
the corresponding locally $L$-analytic group of $L$-valued
points.
In the particular case when $r = 1$, we write simply
$\open{\H}$ and $\open{H}$ rather than $\open{\H}_1$
and $\open{H}_1$.

The space $\An(\open{\H},K)$ is defined as the projective limit
$\plim{r < 1} \An(\H_r,K)$; it is a locally analytic representation
of $\open{H}$, and is a nuclear Fr\'echet space (since
$\open{\H}$ is strictly $\sigma$-affinoid).
Its strong dual space $\DAn(\open{\H},K)_b$
is then isomorphic to the inductive limit $\ilim{r < 1}\DAn(\H_r,K)_b.$
It is a space of compact type.
Corollary~5.1.8 shows that each of the spaces
$\DAn(\H_r, K)_b$ is naturally a topological ring.
Passing to the inductive limit, we obtain a ring
structure on $\DAn(\open{\H},K)_b$, and we may conclude
{\it a priori} that the multiplication map
$$\DAn(\open{\H},K)_b \otimes_K \DAn(\open{\H},K)_b \rightarrow
\DAn(\open{\H},K)_b$$
is separately continuous.  However, since $\DAn(\open{\H},K)_b$
is of compact type, it follows from proposition~1.1.31 that
this map is in fact jointly continuous.

The above explicit description of each of the Banach algebras
$\DAn(\open{\H}_r,K)$ yields the following explicit
description of the inductive limit $\plim{r < 1} \DAn(\open{\H}_r,K)$:
$$ \multline \DAn(\open{\H},K) \iso  \\
\{\sum_{I} b_I
\dfrac{X_1^{i_1}\cdots X_d^{i_d}}
{i_1! \cdots i_d!}
\quad | \quad 
|b_I| \leq C r^{|I|} \text{ for some } C > 0 \text{ and some } r < 1\}.
\endmultline$$

For the remainder of this subsection, we assume that $K$ is discretely
valued (so that $\Cal O_K$ is $p$-adically separated).
Under this assumption, the compact type algebra
$\DAn(\open{\H},K)$ admits another description as a direct limit
of Banach algebras, via the technique of partial divided powers
introduced in \cite{\BER}.

We first recall the following elementary lemma.

\proclaim{Lemma~5.2.2} If $n$ is a natural number
and $m$ a non-negative integer, and if we use the
Euclidean algorithm to write
$n = p^m q + r$, with $0 \leq r < p^m$,
then
$$\dfrac{n!}{q!} = (p^m!)^q r! u,$$
where $u$ is a $p$-adic unit.
\endproclaim
\demo{Proof}
This follows directly from the well-known formula,
stating that the $p$-adic ordinal of $n!$ is equal
to $\dfrac{n-s(n)}{p-1},$ where $s(n)$ denotes
the sum of the $p$-adic digits of $n$.
\qed\enddemo

Fix an integer $m \geq 0,$ and consider the
$\Cal O_K$-subring of the universal enveloping algebra $U(\lie{g})$
generated by the monomials 
$\dfrac{X_j^i}{i!}$
for $0 \leq i \leq p^m$ and $1 \leq j \leq d$;
denote this $\Cal O_K$-subalgebra
by $A^{(m)}$.
Lemma~5.2.2 shows that
the ring $A^{(m)}$ admits the following alternative description:
$$A^{(m)} = 
\{\sum_{I} b_I \dfrac{q_1! \cdots q_d! }{i_1! \cdots i_d!}
X_1^{i_1}\cdots X_d^{i_d} \quad | \quad b_I \in \Cal O_K \text{ and }
b_I = 0 \text{ for almost all } I\},$$
where $q_j$ denotes the integral part of the fraction $i_j/p^m$.

We let $\hat{A}^{(m)}$ denote the $p$-adic completion of $A^{(m)},$
and write $\DAn(\open{\H},K)^{(m)}:= K \otimes_{\Cal O_K} \hat{A}^{(m)}.$
Thus $\DAn(\open{\H},K)^{(m)}$ is naturally a $K$-Banach algebra.

\proclaim{Proposition 5.2.3} There is a natural isomorphism
of topological $K$-algebras of compact type
$\ilim{m} \DAn(\open{\H},K)^{(m)} \iso \DAn(\open{\H},K)$
(assuming that $K$ is discretely valued, so that the source
of this isomorphism is defined).
\endproclaim
\demo{Proof}
This follows from the fact that the $p$-adic ordinal of
$q_1!\cdots q_d!$ is asymptotic to $\dfrac{i_1 + \cdots + i_d}{(p-1)p^m}$
as $i_1 + \cdots + i_d \rightarrow \infty$.
\qed\enddemo

\section{5.3} Let $G$ be locally $L$-analytic group,
and let $H$ be a compact open subgroup of $G$.
The main goal
of this subsection is to present a proof of the fact that $\DLa(H,K)$
is a Fr\'echet-Stein algebra.  This result was originally proved
by Schneider and Teitelbaum \cite{\SCHTNEW, thm.~5.1}, by methods related
to those of \cite{\LAZ}.  The proof we give is quite different; it relies
on an extension of the methods used in \cite{\BER} to prove
the coherence of the sheaf of rings $\Cal D^{\dagger}.$  Our approach
also shows that for any good analytic open subgroup $H$ of $G$,
the algebra $\DAn(\open{\H},K)$ is coherent.

\proclaim{Proposition 5.3.1} If $H$ is a good analytic open 
subgroup of $G$ (in the sense of subsection~5.2), then the
$K$-Fr\'echet algebra $\DLa(\open{H},K)_b$
is a weak Fr\'echet-Stein algebra.
\endproclaim
\demo{Proof}
Let $\{r_n\}_{n\geq 1}$ be a strictly decreasing sequence of numbers
lying in the open interval between $0$ and $1$, that converges to 0,
and that lie in
$|\overline{L}^{\times}|$,
and form the corresponding decreasing sequence
$\{\open{\H}_{r_n}\}_{n\geq 1}$
of $\sigma$-affinoid rigid analytic open subgroups of $\H$.
(As in the preceding subsection,
we have written $\open{\H}_{r_n} = \bigcup_{r < r_n} \H_r$.)
Each of the subgroups $\open{\H}_{r_n}$ is normalised by $\open{H}$.

For each $n\geq 1$ there are continuous injections
$$\Con(\open{H},K)_{\open{\H}_{r_n}-\an} \rightarrow
\Con(\open{H},K)_{\open{\H}_{r_{n+1}}-\an}\tag 5.3.2$$
and
$$\Con(\open{H},K)_{\open{\H}_{r_n}-\an} \rightarrow \La(\open{H},K)\tag 5.3.3$$
(the latter being compatible with the former)
which, upon passing to the inductive limit in $n$, yield an isomorphism
$$\ilim{n} \Con(\open{H},K)_{\open{\H}_{r_n}-\an} \iso \La(\open{H},K).\tag 5.3.4$$

For each value of $n\geq 1$,
we write $D(\open{\H}_{r_n},\open{H})$ to denote the strong dual to the space
$\Con(\open{H},K)_{\open{\H}_{r_n}-\an}$.    The restriction map
$$\Con(\open{H},K)_{\open{\H}_{r_n}-\an} \rightarrow
\Con(\open{H}_{r_n},K)_{\open{\H}_{r_n}-\an}
\iso \An(\open{\H}_{r_n},K)$$
(the isomorphism being provided by proposition~3.4.11)
yields a closed embedding
$\DAn(\open{\H}_{r_n},K) \rightarrow D(\open{\H}_{r_n},\open{H}).$  
The topological ring structure
on $\DAn(\open{\H}_{r_n},K)$ extends naturally to a ring structure on
$D(\open{\H}_{r_n},\open{H}),$ and one immediately sees that
$D(\open{\H}_{r_n},\open{H}) = \bigoplus_{h \in \open{H}/\open{H}_{r_n}}
\DAn(\open{\H}_{r_n},K) * \delta_h$. 
(As indicated, the
direct sum ranges over a set of coset representatives of $\open{H}_{r_n}$ in
$\open{H}$.)
Dualising the isomorphism~(5.3.4) yields an
isomorphism of topological $K$-algebras $\DLa(\open{H},K)_b \iso
\plim{n} D(\open{\H}_{r_n},\open{H}).$
We will show that this isomorphism induces a weak Fr\'echet-Stein structure
on $\DLa(\open{H},K)$.

Each of the algebras $D(\open{\H}_{r_n},\open{H})$ is of compact type
as a convex $K$-vector space, and so satisfies condition~(i) of definition~1.2.6.
Since for each $n\geq 1,$ the inclusion $\open{\H}_{r_{n+1}} \subset \open{\H}_{r_n}$
factors through the inclusion $\H_{r_{n+1}} \subset \open{\H}_{r_n},$ the map~(5.3.2)
is compact (as follows from proposition~2.1.17), and thus so is the dual map
$D(\open{\H}_{r_{n+1}},\open{H}) \rightarrow D(\open{\H}_{r_n},\open{H}).$
Thus by lemma~1.1.14,
condition~(ii) of definition~1.2.6 is also satisfied.
Finally, since for each $n\geq 1$ the map~(5.3.3) 
is a continuous injection of reflexive spaces, the dual map
$\DLa(\open{H},K)_b \rightarrow D(\open{\H}_{r_n},\open{H})_b$
has dense image.  Thus condition~(iii) of definition~1.2.6 is satisfied,
and the proposition is proved.
\qed\enddemo

We intend to show that $\DLa(\open{H},K)_b$ is in fact a
Fr\'echet-Stein algebra,
and we begin by abstracting the argument
used to prove \cite{\BER, thm.~3.5.3}.
We suppose that $A$ is a (not necessarily commutative)
$\Z_p$-algebra, which is $p$-torsion free and $p$-adically
separated.  If $\hat{A}$ denotes the $p$-adic completion of $A$,
then $\hat{A}$ is also $p$-torsion free, and the natural map
$A \rightarrow \hat{A}$ is an injection.  Tensoring with $\Q_p$
over $\Z_p,$ we obtain an injection $\Q_p\otimes_{\Z_p} A
\rightarrow \Q_p \otimes_{\Z_p} \hat{A}$.  We regard $A$, $\Q_p\otimes_{\Z_p}
A,$ and $\hat{A}$ as subalgebras of $\Q_p \otimes_{\Z_p} \hat{A}$.
One easily checks that $A$ is saturated in $\hat{A}$, in the sense
that the inclusion $A \rightarrow \left (\Q_p \otimes_{\Z_p} A\right ) \bigcap
\hat{A}$ is an equality.

Multiplication in $\Q_p\otimes_{\Z_p} \hat{A}$ induces
a natural map $\hat{A} \otimes_A \left ( \Q_p \otimes_{\Z_p} A \right )
\rightarrow \Q_p \otimes_{\Z_p} \hat{A},$ which is evidently an isomorphism.
More generally, if $M$ is a left $A$-submodule of $\Q_p \otimes_{\Z_p} A,$ then
multiplication in $\Q_p\otimes_{\Z_p}\hat{A}$ induces
a natural map of left $\hat{A}$-modules
$\imath_M : \hat{A} \otimes_A M \rightarrow \Q_p \otimes_{\Z_p}
\hat{A}.$  The following lemma concerns the image of the map $\imath_M$.

\proclaim{Lemma 5.3.5} Let $M$ be a left $A$-submodule of
$\Q_p \otimes_{\Z_p} A$.  The image of the map $\imath_M$ is contained
in $M + \hat{A}$ (the sum taking place in $\Q_p \otimes_{\Z_p} \hat{A}$).
If $A \subset M$, then image of $\imath_M$ is in fact equal to $M + \hat{A}$.
\endproclaim
\demo{Proof} The $A$-module $M$
is the direct limit of finitely generated $A$-submodules.
The first statement of the lemma is clearly compatible with passing to
direct limits, and so it suffices to prove it for $M$ that are finitely
generated.  We may thus assume that $p^n M \subset A$ for some natural
number $n$.  
The image of $\imath_M$
is spanned as a $\Z_p$-module by the product $\hat{A} M$.
Since $\hat{A} = A + p^n \hat{A}$, we see that this product is
contained in $M + \hat{A},$ as claimed.  This proves the first claim
of the lemma.

To prove the second claim, note that if $1 \in M,$
then $\hat{A} M$ contains both $M$ and $\hat{A}$.  Thus in this
case $M + \hat{A}$ is contained in the image of $\imath_M,$ and
the claimed equality follows.
\qed\enddemo
 
\proclaim{Corollary 5.3.6} If $B$ is a subring of $\Q_p\otimes_{\Z_p} A$
that contains $A$, then the image of $\imath_B: \hat{A}\otimes_A B
\rightarrow \Q_p \otimes_{\Z_p} \hat{A}$ is a subring of
$\Q_p\otimes_{\Z_p} \hat{A}$.
\endproclaim
\demo{Proof} Since $B$ is in particular a left $A$-submodule of
$\Q_p\otimes_{\Z_p} A$ that contains $A$, lemma~5.3.5 implies
that the image of $\imath_B$ is equal to $B + \hat{A}$.
Applying the obvious analogue of lemma~5.3.5  for right $A$-modules,
we find that the same is true of the image of the natural map
$B \otimes_A \hat{A} \rightarrow \Q_p \otimes_{\Z_p} \hat{A}.$ 
It is now clear that $B + \hat{A}$ is a subring of $\Q_p \otimes_{\Z_p}
\hat{A}$.
\qed\enddemo

\proclaim{Lemma 5.3.7} In the situation of lemma~5.3.5, if we write
$N$ to denote the image of $\imath_M$, and 
if we assume that $A \subset M$,
then $M = N \bigcap \left ( \Q_p \otimes_{\Z_p} A \right )$.
\endproclaim
\demo{Proof}
Lemma~5.3.5 shows that $N = M + \hat{A}.$  Taking this
into account,
the lemma follows from the fact that
$A = \left( \Q_p \otimes_{\Z_p} A \right) \bigcap \hat{A},$
together with the assumption that $A \subset M$.
\qed\enddemo  

\proclaim{Lemma 5.3.8} In the situation of corollary~5.3.6,
the natural map $B \rightarrow B + \hat{A}$ of $\Z_p$-algebras
induces an isomorphism on $p$-adic completions.
\endproclaim
\demo{Proof}
Since $\hat{A} = A + p^n \hat{A} \subset B + p^n(B + \hat{A})$ for
any natural number $n$, we see that
the natural map $B/p^n \rightarrow \left (B + \hat{A} \right )/p^n$
is surjective for any natural number $n$.  Lemma~5.3.7 implies
that this map is also an
injection.  Passing to the projective limit in $n$ proves the lemma.
\qed\enddemo

\proclaim{Lemma 5.3.9} In the situation of corollary~5.3.6,
suppose that
$B$ is equipped with an exhaustive increasing filtration
by $\Z_p$-submodules $F_0 \subset F_1 \subset \cdots$ satisfying
the following assumptions:

(i) For each pair $i,j \geq 0,$ we have 
$F_i F_j \subset F_{i+j}$.  (That is, the $\Z_p$-submodules $F_i$
filter $B$ in the ring-theoretic sense.)

(ii) $F_0 = A$.  (Note that when combined with~(i), this implies
that $F_i$ is a two-sided $A$-submodule of $B,$ for each $i \geq 0$.)

(iii)
The associated graded algebra $\Gr^F_{\bullet} B$ is
finitely generated over $A$ ($= \Gr^F_0 B$) by central elements.

If we let $C$ denote the image of $\imath_B$ (a subalgebra of
$\Q_p \otimes_{\Z_p} \hat{A},$ by corollary~5.3.6), and if for
each $i\geq 0$, we let $G_i$ denote the image of $\imath_{F_i}$,
then $G_{\bullet}$ is an exhaustive filtration of $C$
satisfying properties analogous to~(i), (ii) and~(iii) above. Namely:

(i') For each pair $i,j \geq 0,$ we have 
$G_i G_j \subset G_{i+j}$. 

(ii') $G_0 = \hat{A}$. 

(iii')
The associated graded algebra $\Gr^{G}_{\bullet} C$ is
finitely generated over $\hat{A}$ ($= \Gr^{G}_0 C$) by central elements.
\endproclaim
\demo{Proof}
Since $A = F_0 \subset F_i,$ lemma~5.3.5  shows that
the image $G_i$ of $\imath_{F_i}$ is equal to $F_i + \hat{A}.$
It is now immediately checked that $G_{\bullet}$ is
an exhaustive filtration of $C$, satisfying conclusions~(i') and~(ii')
of the lemma.
Lemma~5.3.7 implies that for each $i\geq 1,$ the inclusion
$F_{i-1} \subset \left( F_{i} \bigcap (F_{i-1} + \hat{A})\right )$ is an equality,
and thus that the natural map
$\Gr^{F}_i B \rightarrow \Gr^{G}_i C$ is an isomorphism if $i \geq 1$.
Part~(iii') thus follows from the corresponding assumption~(iii),
and the lemma is proved.
\qed\enddemo

We are now ready to generalise \cite{\BER, thm.~3.5.3}.

\proclaim{Proposition 5.3.10}
Let $A$ be a (not necessarily commutative)
$p$-torsion free and $p$-adically separable left Noetherian 
$\Z_p$-algebra.
If $B$ is a $\Z_p$-subalgebra of $\Q_p\otimes_{\Z_p} A$
that contains $A$, and satisfies conditions~(i), (ii) and~(iii)
of lemma~5.3.9,
then the rings
$\Q_p\otimes_{\Z_p} \hat{A}$ and $\Q_p \otimes_{\Z_p} \hat{B}$
are left Noetherian $\Q_p$-algebras, and
$\Q_p\otimes_{\Z_p} \hat{B}$ is flat 
as a right $\Q_p\otimes_{\Z_p} \hat{A}$-module.
\endproclaim
\demo{Proof}
We begin by noting that assumption~(iii) of lemma~5.3.9
and the Hilbert basis theorem
imply that $\Gr^F_{\bullet} B$ is left Noetherian,
and thus that the same is true of $B$.  Since the
element $p$ is central in $A$ and $B$, it follows that the $p$-adic completions
$\hat{A}$ and $\hat{B}$ are both left Noetherian,
and thus the same is true of
$\Q_p\otimes_{\Z_p} \hat{A}$ and $\Q_p \otimes_{\Z_p} \hat{B}$.
The main point of the proposition, then, is to prove that 
$\Q_p\otimes_{\Z_p} \hat{B}$ is flat as a right
$\Q_p\otimes_{\Z_p} \hat{A}$-module.

As in the statement of lemma~5.3.9, let $C$ denote the image of
$\imath_B$, equipped with its filtration $\hat{F}_{\bullet}$.
Since $\hat{A}$ is left Noetherian, conclusion~(iii') of lemma~5.3.9 
and the Hilbert basis theorem imply that $\Gr^{\hat{F}}_{\bullet} C$
is left Noetherian, and hence that $C$ is left Noetherian.
As $p$ is central in $C$, the Artin-Rees
theorem applies to show that $\hat{C}$ is flat 
as a right $C$-module.  Tensoring with $\Q_p$ over $\Z_p,$ we conclude
that $\Q_p\otimes_{\Z_p} \hat{C}$ is right flat over $\Q_p \otimes_{\Z_p} C$.

Since $A \subset B \subset \Q_p \otimes_{\Z_p} A,$ and since
$C = B + \hat{A},$ it is immediate that the inclusion $\hat{A} \subset
C$ becomes an equality after tensoring through with $\Q_p$ over $\Z_p.$
Lemma~5.3.8 implies that the inclusion $B \subset C$ induces an isomorphism
$\hat{B} \iso \hat{C}$, which of course remains an isomorphism after
tensoring through with $\Q_p$ over $\Z_p$.  Combining these remarks with
the conclusion of the preceding paragraph, we deduce that $\Q_p\otimes_{\Z_p}
\hat{B}$ is right flat over $\Q_p \otimes_{\Z_p} \hat{A}$, as required.
\qed\enddemo

We now present some applications of the preceding proposition.

\proclaim{Proposition 5.3.11} Let $H$ be a good analytic open subgroup
of $G$ (in the sense of subsection~5.2), and suppose that $K$ is
discretely valued. 

(i) For each $m \geq 0,$ the ring $\DAn(\open{\H},K)^{(m)}$ (as
defined in the discussion preceding proposition~5.2.3) is Noetherian.

(ii)
If $0 \leq m_1 \leq m_2,$ then the natural map
$\DAn(\open{\H},K)^{(m_1)} \rightarrow \DAn(\open{\H},K)^{(m_2)}$ is flat.
\endproclaim
\demo{Proof}
Let $\lie{h}$ be the $\Cal O_K$-Lie subalgebra of $\lie{g}$ that gives
rise to the good analytic open subgroup $H$, and let $X_1,\ldots,X_d$
be an $\Cal O_K$-basis for $\lie{h}$.
For any natural number $m$,
we let $A^{(m)}$ denote the $\Cal O_K$-subalgebra of the enveloping algebra
$U(\lie{g})$ of $\lie{g}$ over $K$ generated by the monomials $X_j^i/i!$ for $0 \leq i \leq p^m$
and $1\leq j\leq d$.  (Compare the discussion preceding proposition~5.2.3.)
The algebra $U(\lie{g})$ has its usual graded structure
$U(\lie{g}) = \bigoplus_i U(\lie{g})_i$, where $U(\lie{g})$
denotes the subspace of $U(\lie{g})$ spanned by monomials in the basis elements $X_j$
of total degree equal to $i$.
Since $A^{(m)}$ is generated as an $\Cal O_K$-algebra by monomials,
it is a graded $\Cal O_K$-subalgebra of $U(\lie{g})$, and so
$A^{(m)} = \bigoplus_i A^{(m)}_i,$ where of course
$A^{(m)}_i = A^{(m)} \bigcap U(\lie{g})_i.$
This grading on $A^{(m)}$ induces a corresponding filtration on $A^{(m)}$,
whose $i$th filtered piece is equal to
$\bigoplus_{i'\leq i} A^{(m)}_{i'}.$
The associated graded algebra is a commutative $\Cal O_K$-algebra of finite type.
In particular, it is Noetherian, and thus the same is true of $A^{(m)}$.

Now consider the inclusions
$A^{(m)} \subset A^{(m+1)} \subset \Q_p\otimes_{\Z_p} A^{(m)} = U(\lie{g}).$
We equip $A^{(m+1)}$ with the following filtration:
$$F_i := A^{(m)} (\bigoplus_{i'\leq i} A^{(m+1)}_{i'}).$$
For each $i \geq 1,$ it is clear that
$A^{(m)} A^{(m+1)}_i = A^{(m+1)}_i A^{(m)}$,
and thus that $F_i F_j \subset F_{i+j}$ for any $i,j \geq 0$.
Since $A^{(m+1)}$ is generated over $A^{(m)}$ by finitely many monomials,
the associated graded algebra
$\Gr^{F}_{\bullet} A^{(m+1)}$ is finitely generated over $A^{(m)}$.
Furthermore, these monomials reduce
to central elements in the associated graded algebra.
Indeed, this is a consequence of the standard fact that
$$[A^{(m+1)}_i,A^{(m+1)}_j] \subset A^{(m+1)}_{i+j-1},$$
for any $i,j \geq 0,$ together with the identity
$$[a,\dfrac{X_j^{p^{m+1}}}{p^{m+1}!}] = u \left ( \dfrac{X_j^{p^m}}{p^m!}
\right )^{p-1}  [a, \dfrac{X_j^{p^m}}{p^m!}] \in A^{(m)}$$
(where $u$ is a rational number that is a $p$-adic unit), for
$a \in A^{(m)}$ and $j = 1,\ldots,d$.
Altogether, we see that we are in the situation 
of proposition~5.3.10 (if we take the ring $A$ of that proposition to be
$A^{(m)}$ and the ring $B$ to be $A^{(m+1)}$).  That proposition then
implies that the each of the rings 
$\DAn(\open{\H},K)^{(m)}$ and
$\DAn(\open{\H},K)^{(m+1)}$ is Noetherian, and that the natural map
$\DAn(\open{\H},K)^{(m)} \rightarrow \DAn(\open{\H},K)^{(m+1)}$ is flat.
Both parts of the proposition now follow.
\qed\enddemo

\proclaim{Corollary 5.3.12}  If $H$ is a good analytic open
subgroup of the locally $L$-analytic group $G$,
and if $K$ is discretely valued,
then $\DAn(\open{\H},K)$ is coherent.
\endproclaim
\demo{Proof}
Propositions~5.2.3 and~5.3.11 together imply that $\DAn(\open{\H})$
is equal to the inductive limit of a sequence of Noetherian 
rings with flat transition maps.  This immediately implies that
$\DAn(\open{\H})$ is coherent.
\qed\enddemo

\proclaim{Proposition 5.3.13} Let $H$ be a good analytic open subgroup
of $G$ (in the sense of subsection~5.2), and suppose that $K$ is
discretely valued. 
If $0 < r < 1$ is a real number lying in $|\bar{L}^{\times}|$,
and if $m$ is a sufficiently large
natural number, then the natural map
$\DAn(\open{\H}_r,K) \rightarrow \DAn(\open{\H},K)$ factors
through the natural map $\DAn(\open{\H},K)^{(m)} \rightarrow
\DAn(\open{\H},K)$, and the resulting map
$\DAn(\open{\H}_r,K) \rightarrow \DAn(\open{\H},K)^{(m)}$ is flat.
\endproclaim
\demo{Proof}
After making a change of coefficients from $K$ to some finite
extension of $K$, we may assume that $r$ is equal to the
absolute value of some element $\alpha$ of $K^{\times}$.

As in the proof of proposition~5.3.11, we
let $\lie{h}$ denote the $\Cal O_K$-Lie subalgebra of
$\lie{g}$ that exponentiates to give $H$, and choose
a basis $X_1,\ldots,X_j$ of $\lie{h}$.

For any natural number $m',$
let $A^{(m')}(r)$ denote the $\Cal O_K$-subalgebra of the enveloping algebra
$U(\lie{g})$ generated by the monomials $(\alpha X_j)^i/i!$ for
$1 \leq i \leq p^{m'}$.  If $\hat{A}^{(m')}(r)$ denotes the
$p$-adic completion of this ring,
then
$$\DAn(\open{\H}_r,K) \cong \ilim{m} K\otimes_{\Cal O_K}
\hat{A}^{(m')}(r).$$
(This follows from proposition~5.2.3,
applied to $\H_r$ rather than $\H$.)
Thus, to prove the proposition,
it suffices to exhibit a natural number $m$ such that,
for every $m'\geq 0$, the map
$$K\otimes_{\Cal O_K} \hat{A}^{(m')}(r) \rightarrow \DAn(\open{\H},K)\tag
5.3.14$$
factors as a composite
$$K\otimes_{\Cal O_K} \hat{A}^{(m')}(r) \rightarrow
\DAn(\open{\H},K)^{(m)} \rightarrow
\DAn(\open{\H},K),\tag 5.3.15$$
with the first factor being flat. 

The $K$-algebra $\DAn(\open{\H},K)^{(m)}$ is equal to
the tensor product $K \otimes_{\Cal O_K} \hat{A}^{(m)},$
where $\hat{A}^{(m)}$ denotes the $p$-adic completion of
the $\Cal O_K$-algebra $A^{(m)}$ considered in the discussion
preceding proposition~5.2.3, and again in the proof of
proposition~5.3.11.  Just as with $A^{(m')}(r)$, we regard
$A^{(m)}$ as an $\Cal O_K$-subalgebra of $U(\lie{g})$.

As in proof of lemma~5.2.2, let $s(i)$ denote the sum of
the $p$-adic digits of $i$, for any natural number $i$.
If we fix an integer $i$, and write $i = p^m q + r$
(respectively $i = p^{m'} q + r'$)
with $0 \leq r < p^m$
(respectively $0\leq r' < p^{m'}$),
then
$$\ord_K(q!) = \dfrac{i - s(i)}{(p-1)p^m} - \dfrac{r - s(r)}{(p-1)p^m},$$
while
$$\ord_K(\alpha^i q'!)
=  \ord_K(\alpha) i +
\dfrac{i - s(i)}{(p-1)p^{m'}} - \dfrac{r' - s(r')}{(p-1)p^{m'}}.$$
The values of $r$ and $r'$ are bounded independently of $i.$
Thus, if we choose $m$ so that
$\ord_K(\alpha) \geq 1/(p-1)p^m$
(a condition independent of $m'$),
then we find that 
$$\ord_K(\alpha^i q'!) \geq 	\ord_K(q!) - C,$$
for some constant $C \geq 0$ independent of $i$.
This inequality in turn shows that the product
$A^{(m')}(r) A^{(m)}$ is contained in $\dfrac{1}{p^k}A^{(m)}$
for some sufficiently large natural number $k$.
(Note also that the product $A^{(m')}(r) A^{(m)}$ is obviously equal
to $A^{(m)} A^{(m')}(r)$, and so $A^{(m')}(r) A^{(m)}$ is
an $\Cal O_K$-subalgebra of $U(\lie{g})$.)
If we let $(A^{(m')}(r) A^{(m)})\hat{}$ denote the $p$-adic
completion of $A^{(m')}(r) A^{(m)}$,
then it follows that the natural map
$$K\otimes_{\Cal O_K} (A^{(m')}(r) A^{(m)})\hat{}
\rightarrow K\otimes_{\Cal O_K} \hat{A}^{(m)} = \DAn(\open{\H},K)^{(m)}$$
is a topological isomorphism.
The inclusion
$A^{(m')}(r) \rightarrow 
A^{(m')}(r) A^{(m)}$
also induces a morphism
$\hat{A}^{(m')}(r) \rightarrow 
(A^{(m')}(r) A^{(m)})\hat{}.$
Thus we see that the map~(5.3.14) factors as in~(5.3.15).
It remains to be shown that the first factor of~(5.3.15)
is flat.

Proposition~5.3.11, applied to $\H_r$, shows that
if $m' \leq m''$, then the natural map
$$K\otimes_{\Cal O_K} \hat{A}^{(m')}(r) \rightarrow
K\otimes_{\Cal O_K} \hat{A}^{(m'')}(r)$$ 
is flat.  Thus, in showing that the first factor
of~(5.3.15) is flat, we may,
without loss of generality,
replace $m'$ by any larger natural number.
In particular, we may, and do, assume that $m' \geq m$.

If we equip $U(\lie{g})$ with its usual grading,
then (as was already observed in the proof of proposition~5.3.11)
$A^{(m)}$ is a graded $\Cal O_K$-subalgebra of $U(\lie{g})$; 
we write $A^{(m)} = \bigoplus_{i\geq 0} A^{(m)}_i.$
We equip the product $A^{(m')}(r) A^{(m)}$
with the following filtration:
$$F_i := A^{(m')}(r)(\bigoplus_{i'\leq i} A^{(m)}_{i'}).$$
It is immediate that this filtration satisfies
conditions~(i) and~(ii) of lemma~5.3.9
(taking the algebra $A$ of that lemma to be $A^{(m')}(r),$
and the algebra $B$ to be $A^{(m')}(r) A^{(m)}.$
It also satisfies condition~(iii) of that lemma.
Indeed, the evident inclusion
$$[\alpha \lie{h}, \lie{h}]
= \alpha [\lie{h},\lie{h}] \subset \alpha \lie{h}$$
yields the inclusion
$$[A^{(m')}(r),A^{(m')}_i] \subset F_{i-1},$$
from which condition~(iii) of lemma~5.3.9 follows
directly.  (It is here that we take into account
that $m \leq m'$.)
Proposition~5.3.10 now implies 
that the natural map
$$K\otimes_{\Cal O_K} \hat{A}^{(m')}(r)
\rightarrow K\otimes_{\Cal O_K} (A^{(m')}(r) A^{(m)})\hat{}$$
is flat.  This is precisely the first arrow of~(5.3.15),
and thus the proof of the proposition is completed.
\qed\enddemo

In the context of the preceding proposition,
the homomorphism 
$\DAn(\open{\H}_r,K) \rightarrow \DAn(\open{\H},K)$
extends to a homomorphism
$$D(\open{\H}_r,\open{H}) \rightarrow \DAn(\open{\H},K)\tag 5.3.16$$
(in the notation of the proof of proposition~5.3.1;
this homomorphism is dual to the natural map
$\An(\open{\H},K) \rightarrow \Con(\open{H},K)_{\open{\H}_r-\an}.$)
We wish to strengthen proposition~5.3.13
so that is applies to the map~(5.3.16).

The following lemma provides the necessary bootstrap. 
It provides an analogue in a simple noncommutative situation
of the following commutative algebra fact:  if the composite
$A \rightarrow B \rightarrow C$ is a flat morphism of
commutative rings, and if $B$ is finite \'etale over $A$,
then $B\rightarrow C$ is also flat.  

\proclaim{Lemma 5.3.17}
Let $A$ be an associative $K$-algebra, and let $B$ be an extension
of $A$, satisfying the following hypotheses:

(i)  There exists an element
$x \in B$ such that $B$
is generated as a ring by $A$ together with $x$.

(ii) There is a unit $a \in A^{\times}$ such that $x^n = a$.

(iii) The element $x$ of $B$ (which is a unit of $B$, by~(ii))
normalises $A$; that is, $x A x^{-1} = A$.

If $C$ is an extension of $B$ that is flat as a left (respectively right)
$A$-module,
then $C$ is flat as a left (respectively right) $B$-module.
\endproclaim
\demo{Proof}
If $M$ is a right $B$-module, then we define a $C$-linear automorphism
$\sigma$ of $M\otimes_A C$ as follows:
$$\sigma(m \otimes c) = m x^{-1} \otimes x c.$$
The formation of $\sigma$ is evidently functorial in $M$.
It follows from hypothesis~(iii) in the statement of the lemma
that $\sigma$ is well-defined, and hypothesis~(ii) implies
that $\sigma^n = 1$.
Hypothesis~(i) implies that $M\otimes_B C$ is isomorphic to
the coinvariants of $M\otimes_A C$ with respect to the cyclic
group $\langle \sigma \rangle$ generated by $\sigma$.

Suppose that $C$ is flat as a left $A$-module.  Since passing to
$\langle \sigma \rangle$-coinvariants is exact (since $K$
is of characteristic zero), we see that the formation
of $M\otimes_B C$ is an exact functor of $M$.  Thus $C$
is also flat as a left $B$-module.

The case when $C$ is flat as a right $A$-module is proved by
an analogous argument.
\qed\enddemo

\proclaim{Proposition 5.3.18}  
Let $H$ be a good analytic open subgroup of $G$,
and suppose that $K$ is discretely valued.
If $0 < r < 1$ lies in $|\bar{L}^{\times}|$,
and if $m$ is a sufficiently large
natural number, then map~(5.3.16)
factors
through the natural map $\DAn(\open{\H},K)^{(m)} \rightarrow
\DAn(\open{\H},K)$, and the resulting map
$D(\open{\H}_r,\open{H}) \rightarrow \DAn(\open{\H},K)^{(m)}$ is flat.
\endproclaim
\demo{Proof}
The ring $D(\open{\H}_r,\open{H})$ is free of finite rank as
a $\DAn(\open{\H}_r,K)$-module (on either side), and so in
particular is finitely generated as a ring over $\DAn(\open{\H}_r,K)$.
Since $\DAn(\open{\H},K)$ is the direct limit of its subrings
$\DAn(\open{\H},K)^{(m)}$, and since, by proposition~5.3.13,
the image of the restriction of~(5.3.16) to  $\DAn(\open{\H}_r,K)$
lies in $\DAn(\open{\H},K)^{(m)}$ for some sufficiently
large value of $m$, we see that $D(\open{\H}_r,\open{H})$
also maps into $\DAn(\open{\H},K)^{(m)}$ for some (possibly
larger) value of $m$.

Again by proposition~5.3.13, the composite of the sequence of maps
$$\DAn(\open{\H}_r,K) \rightarrow D(\open{\H}_r,\open{H})
\rightarrow \DAn(\open{\H}_r,K)^{(m)}$$
is flat.  Our goal is to show that the second of these maps is
flat.

Since $\open{H}$ is a pro-$p$-group containing $\open{H}_r$
as a normal open subgroup, we may find a filtration by normal subgroups
$$\open{H}_r = G_0 \subset G_1
\subset \cdots \subset G_n = \open{H},$$
such that each 
quotient $G_{i+1}/G_i$ is cyclic.
For each $i \geq 0$, let $D(\open{\H}_r,G_i)$ denote the
dual to the space $\Con(G_i,K)_{\open{\H}_r-\an},$
equipped with its natural ring structure.
For each value of $i,$ the inclusion $D(\open{\H}_r,G_i) \subset
D(\open{\H}_r,G_{i+1})$ then satisfies the hypothesis of
the inclusion $A \subset B$ of lemma~5.3.17.  Taking
$\DAn(\open{\H}_r,K)^{(m)}$
to be the ring $C$ of that lemma, and arguing by induction,
starting from case $i = 0$ (and noting that
$D(\open{\H}_r,G_0) = D(\open{\H}_r,\open{H}_r) = \DAn(\open{\H}_r,K)$)
we find that 
$\DAn(\open{\H}_r,K)^{(m)}$
is flat over
$D(\open{\H}_r,G_n) = D(\open{\H}_r,\open{H}),$
as required.
\qed\enddemo

We now use the preceding results to give a new proof
of \cite{\SCHTNEW, thm.~5.1}.

\proclaim{Corollary 5.3.19} 
If $K$ is discrete,
and if $H$ is a compact open subgroup
of $G$, then the nuclear Fr\'echet algebra
$\DLa({H},K)$ is a Fr\'echet-Stein algebra.
\endproclaim
\demo{Proof}
Any compact open subgroup of $G$ contains as a normal
subgroup a good analytic open
subgroup of $G$, and so it suffices to establish that
$\DLa(\open{H},K)$ is a Fr\'echet-Stein algebra for
any good analytic open subgroup $\H$ of $G$.
(Compare the second half of the argument in
step 1 of the proof of \cite{\SCHTNEW, thm.~5.1}.)

Let $\pi$ denote the uniformiser of $\Cal O_K$,
and write $r = |\pi|$.  Consider the decreasing sequence
$\{\open{\H}_{r^n}\}_{n \geq 0}$ of affinoid subgroups of $\H$.
Fix a value of $n$, and 
apply proposition~5.3.18 to the inclusion
$\open{\H}_{r^{n+1}} \rightarrow \open{\H}_{r^n}$.
We see that we may find an integer $m_n$
such that the map
$D(\open{\H}_{r^{n+1}},\open{H}_{r^n}) \rightarrow \DAn(\open{\H}_{r^n},K)$
factors as
$$D(\open{\H}_{r^{n+1}},\open{H}_{r^n}) \rightarrow
\DAn(\open{\H}_{r^n},K)^{(m_n)}
\rightarrow
\DAn(\open{\H}_{r^n},K),$$
with both arrows being flat.  (The flatness of the first arrow 
is the conclusion of proposition~5.3.18; that of the second
arrow follows from propositions~5.2.3 and~5.3.11.)
Tensoring with $D(\open{\H}_{r^{n+1}},\open{H})$
over $D(\open{\H}_{r^{n+1}},\open{H}_{r^n}),$
we obtain the sequence of maps
$$\multline
D(\open{\H}_{r^{n+1}},\open{H}) \rightarrow
D(\open{\H}_{r^{n+1}},\open{H}) \otimes_{D(\open{\H}_{r^{n+1}},\open{H}_{r^n})}
\DAn(\open{\H}_{r^n},K)^{(m_n)} \\
\rightarrow D(\open{\H}_{r^n},\open{H}).\endmultline \tag 5.3.20$$
The tensor product
$D(\open{\H}_{r^{n+1}},\open{H}) \otimes_{D(\open{\H}_{r^{n+1}},\open{H}_{r^n})}
\DAn(\open{\H}_{r^n},K)^{(m_n)}$
is naturally a Banach algebra.
More precisely,
we have the isomorphism
$$D(\open{\H}_{r^{n+1}},\open{H})
\iso
\bigoplus_{h \in \open{H}/\open{H}_{r^n}}
\delta_h * D(\open{\H}_{r^{n+1}},\open{H}_{r^n}),$$ 
and hence a corresponding isomorphism
$$D(\open{\H}_{r^{n+1}},\open{H}) \otimes_{D(\open{\H}_{r^{n+1}},\open{H}_{r^n})}
\DAn(\open{\H}_{r^n},K)^{(m_n)}
\iso
\bigoplus_{h \in \open{H}/\open{H}_{r^n}}
\delta_h * 
\DAn(\open{\H}_{r^n},K)^{(m_n)}.$$
The right hand side has a natural topological algebra structure
(since $\open{H}$ normalises $\open{\H}_{r^n}$),
and is furthermore a Noetherian Banach algebra,
since the same is true of
$\DAn(\open{\H}_{r^n},K)^{(m_n)}$.
With this Banach algebra structure on the middle term,
we find that~(5.3.20) becomes a sequence of continuous
flat homomorphisms of topological $K$-algebras.

If we now allow $n$ to vary, and take into account the isomorphism
$\DLa(\open{H},K) \iso \plim{n} D(\open{\H}_{r^n},\open{H})$
established in the course of proving proposition~5.3.1,
as well as the sequences of flat maps provided by~(5.3.20),
we obtain an isomorphism
$$\DLa(\open{H},K) \iso \plim{n} 
D(\open{\H}_{r^{n+1}},\open{H}) \otimes_{D(\open{\H}_{r^{n+1}},\open{H}_{r^n})}
\DAn(\open{\H}_{r^n},K)^{(m_n)}.$$
This isomorphism describes $\DLa(\open{H},K)$ as the projective limit of
Noetherian Banach algebras related by flat transition maps, 
and so establishes that it is a Fr\'echet-Stein algebra.
\qed\enddemo

In fact, we will require a more general result than the
preceding corollary.

\proclaim{Definition 5.3.21}  Assume that $K$ is discretely
valued.  We say that a 
Fr\'echet algebra over $K$ admits an
integral Fr\'echet-Stein structure 
if we may find a projective sequence $\{\Cal A_n\}_{n\geq 1}$ of $p$-adically
separated and complete and $p$-torsion free Noetherian $\Cal O_K$-modules
with flat transition maps such that there is an isomorphism
$A \iso \plim{n} K\otimes_{\Cal O_K} \Cal A_n$
(where each of the tensor products $K \otimes_{\Cal O_K} \Cal A_n$
is equipped with its natural $K$-Banach algebra structure,
defined by taking $\Cal A_n$ to be the unit ball).
\endproclaim

Note that any Fr\'echet algebra satisfying the condition of the preceding
definition is certainly a Fr\'echet-Stein algebra, since each Banach algebra
$K\otimes_{\Cal O_K} \Cal A_n$ is Noetherian, and the transition maps
$K\otimes_{\Cal O_K} \Cal A_{n+1} \rightarrow
K\otimes_{\Cal O_K} \Cal A_n$ are flat.

For example, if $\X$ is a $\sigma$-affinoid rigid analytic space
over $K$, then $\An(\X,K)$ admits an integral Fr\'echet-Stein 
structure.  (This follows from the main result of \cite{\MEH}; see also
\cite{\BOL}.  The point is if we write $\X = \bigcup_{n\geq 1} \X_n$
as the union of an increasing sequence of affinoid open subdomains,
then the flattening results of these
papers allow us to find integral models for the flat maps
$\An(\X_{n+1},K) \rightarrow \An(\X_n,K).$)

\proclaim{Proposition 5.3.22}
Suppose that $K$ is discretely valued, and that $H$ is a compact open
subgroup of $G$.
If $A$ is a $K$-Fr\'echet algebra that admits an integral Fr\'echet-Stein
structure, then the completed tensor product
$A\cotimes_K \DLa(H,K)$ is a Fr\'echet-Stein algebra.
\endproclaim
\demo{Proof}
Let $H$ be an good analytic open subgroup of $G$,
and let $A \rightarrow B$ be 
a flat map of Noetherian $p$-adically
complete $\Cal O_K$-algebras.
If $0 < r < 1$ is an element of $|L^{\times}|$,
then an argument analogous to that used to prove
proposition~5.3.13, with $U(\lie{g})$ replaced by $A\otimes_{\Cal O_K}
U(\lie{g})$, and all of the various $\Cal O_K$-subalgebras of
$U(\lie{g})$
replaced by the corresponding $A$-subalgebras,
shows that the natural map
$$A \cotimes_{\Cal O_K} \DAn(\open{\H}_r,K)
\rightarrow A\cotimes_{\Cal O_K} \DAn(\open{\H},K)$$
factors through the natural map
$$A\cotimes_{\Cal O_K} \DAn(\open{\H},K)^{(m)} \rightarrow 
A\cotimes_{\Cal O_K} \DAn(\open{\H},K),$$
for some sufficiently large
natural number $m$,
and that the resulting map
$$A \cotimes_{\Cal O_K} \DAn(\open{\H}_r,K)
\rightarrow A\cotimes_{\Cal O_K} \DAn(\open{\H},K)^{(m)}$$
is flat. 
An argument analogous to that used to prove
proposition~5.3.18 allows us to strengthen this result
to show that the natural map
$$A \cotimes_{\Cal O_K} D(\open{\H}_r,\open{H})
\rightarrow A\cotimes_{\Cal O_K} \DAn(\open{\H},K)$$
factors through the natural map
$$A\cotimes_{\Cal O_K} \DAn(\open{\H},K)^{(m)} \rightarrow 
A\cotimes_{\Cal O_K} \DAn(\open{\H},K),$$
for some sufficiently large
natural number $m$,
and that the resulting map
$$A \cotimes_{\Cal O_K} D(\open{\H}_r,\open{H})
\rightarrow A\cotimes_{\Cal O_K} \DAn(\open{\H},K)^{(m)}$$
is flat. 

Each of the rings 
$A\otimes_{\Cal O_K} A^{(m)}$
and 
$B\otimes_{\Cal O_K} A^{(m)}$
is Noetherian (as follows by an argument
analogous to that used in the proof of proposition~5.3.11
to show that $A^{(m)}$ is Noetherian),
and that the map
$A\otimes_{\Cal O_K} A^{(m)}
\rightarrow 
B\otimes_{\Cal O_K} A^{(m)}$
is flat
(since it is a base-change of the flat map $A \rightarrow B$).
Passing to $p$-adic completions thus yields a flat map
$A\cotimes_{\Cal O_K} A^{(m)}
\rightarrow 
B\cotimes_{\Cal O_K} A^{(m)}$.
Finally, tensoring with $K$ over $\Cal O_K$,
we find that the map
$A\cotimes_{\Cal O_K} \DAn(\open{\H},K)^{(m)}
\rightarrow
B\cotimes_{\Cal O_K} \DAn(\open{\H},K)^{(m)}$
is flat.  Combining this result with that of
the preceding paragraph, we find that
the natural map
$$A \cotimes_{\Cal O_K} D(\open{\H}_r,\open{H})
\rightarrow B\cotimes_{\Cal O_K} \DAn(\open{\H},K)$$
factors through the natural map
$$B\cotimes_{\Cal O_K} \DAn(\open{\H},K)^{(m)} \rightarrow 
B\cotimes_{\Cal O_K} \DAn(\open{\H},K),$$
for some sufficiently large natural number $m$,
and that the resulting map
$$A \cotimes_{\Cal O_K} D(\open{\H}_r,\open{H})
\rightarrow B\cotimes_{\Cal O_K} \DAn(\open{\H},K)^{(m)}$$
is flat. 
The proposition now follows from this result,
in the same way that corollary~5.3.19 follows
from proposition~5.3.18.
\qed\enddemo

\head 6. Admissible locally analytic representations \endhead

\section{6.1} Fix a locally $L$-analytic group $G$.
In this subsection we introduce the notion of an admissible locally
analytic $G$-representation.

\proclaim{Definition 6.1.1} Let $V$ be a locally analytic
representation of $G$.  We say that $V$ is admissible  if it
an $LB$-space, and if for
a cofinal sequence of
analytic open subgroups $H$ of $G$,
the space $(V_{\H-\an})'$ is finitely
generated as a left $\DAn(\H,K)$-module. 
\endproclaim

Note that if $V$ is an $LB$-space, then corollary~3.3.21
implies that $V_{\H-\an}$ is 
also an $LB$-space, and so theorem~5.1.15 implies that $(V_{\H-\an})'$ is
finitely generated as a $\DAn(\H,K)$-module if and only if
the space $V_{\H-\an}$ of $\H$-analytic vectors
admits a closed embedding into $\An(\H,K)^n$ for some natural number $n$.

We will show below that a locally analytic $G$-representation $V$
is admissible if and only if $V'_b$ is a coadmissible topological
$\DLa(H,K)$ for one (or equivalently any) analytic open subgroup
$H$ of $G$.  This will show that our definition of an admissible
locally analytic representation of $G$ is equivalent to that of
\cite{\SCHTNEW, def., p.~33}.

\proclaim{Proposition 6.1.2}  If $V$ is an admissible  locally
analytic representation of $G$ then $(V_{\H-\an})'$ is finitely generated as 
a left $\DAn(\H,K)$-module for every analytic open subgroup $H$ of $G$.
In particular, for every such $H$, the space $V_{\H-\an}$ is a Banach space.
\endproclaim
\demo{Proof} Let $H$ be an analytic open subgroup of $G$.
By assumption we may find an analytic open subgroup $H' \subset H$
such that $V_{\H'-\an}'$ is finitely generated over $\DAn(\H',K)$,
or equivalently, by the remarks following definition~6.1.1,
such that there is an $H'$-equivariant
closed embedding $V_{\H'-\an} \rightarrow \An(\H',K)^n$ for some natural
number $n$.  

We now construct a closed embedding
$$\multline
V_{\H-\an} = \An(\H,V)^{\Delta_{1,2}(H)} \rightarrow
\An(\H,V)^{\Delta_{1,2}(H')} \iso \An(\H,V_{\H'-\an})^{\Delta_{1,2}(H')} \\
\rightarrow \An(\H,\An(\H',K)^n)^{\Delta_{1,2}(H')} \iso \An(\H)^n.
\endmultline $$
(Here the isomorphisms are provided by lemmas~3.3.9 and~3.3.11 
respectively.)  The existence of this closed embedding in turn
implies that $(V_{\H-\an})'$ is finitely generated over $\DAn(\H,K)$,
as required.
\qed\enddemo

\proclaim{Proposition 6.1.3}  If $V$ is an admissible  locally
analytic representation of $G$ then $V$ is of compact type.
\endproclaim
\demo{Proof} Replacing $G$ by a compact open subgroup of itself if
necessary, we assume that $G$ is compact.
Let $\{H_m\}_{m \geq 1}$ be a cofinal sequence
of normal analytic open subgroups of $G$, chosen so that
the rigid analytic map $\H_{m+1} \rightarrow \H_m$ is
relatively compact (in the sense of definition~2.1.15)
for each $m\geq 1$. Fix one such value of $m$, and
write $H_m = \prod_{i\in I} h_i H_{m+1}.$  
Proposition~6.1.2 yields an $H_m$-equivariant closed embedding
of $K$-Banach spaces
$V_{\H_{m+1}-\an}
 \rightarrow \prod_{i\in I} \An(h_i \H_{m+1},K)^n$ for some $n$.
Passing to $\H_m$-analytic vectors and applying proposition~3.3.23
yields a Cartesian diagram
$$\xymatrix{V_{\H_m-\an} \ar[r]\ar[d] & \An(\H_m,K)^n \ar[d] \\
V_{\H_{m+1}-\an} \ar[r] & \prod_{i\in I} \An( h_i \H_{m+1} , K)^n}$$
in which the horizontal arrows are closed embeddings.
Since the right-hand vertical arrow is compact, by proposition~2.1.16,
and since the
horizontal arrows are closed embeddings, we see that the left-hand
vertical arrow is also compact.  Thus $V_{\la} \iso \ilim{n}V_{\H_m-\an}$
is of compact type.
From theorem~3.6.12 we conclude that the natural map  
$V_{\la} \rightarrow V$ is a topological isomorphism, and thus
that $V$ is of compact type.
\qed\enddemo

\proclaim{Proposition 6.1.4} If $V$ is an admissible  locally analytic
representation of $G$ and $W$ is a $G$-invariant closed subspace of $G$
then $W$ is an admissible  locally analytic representation of $G$.
\endproclaim
\demo{Proof} Proposition~6.1.3 shows that $V$ is of compact type, and
thus that $W$ is also of compact type.  In particular, $W$ is an $LB$-space.

Let $H$ be an analytic open subgroup of $G$.  Proposition~6.1.2 shows
that $V_{\H-\an}$ is a Banach space, and so corollary~3.3.18
shows that the natural map $W_{\H-\an} \rightarrow V_{\H-\an}$ 
is a closed embedding. Hence
$W_{\H-\an}'$ is a quotient of $V_{\H-\an}',$ and so is finitely
generated over $\DAn(\H,K)$ (since this is true of $V_{\H-\an}'$ by
assumption).
\qed\enddemo

\proclaim{Proposition 6.1.5}
If $V$ is an admissible  locally analytic representation of $G$
and $W$ is a finite dimensional locally analytic representation
of $G$ then the tensor product $V\otimes_K W$ is again an admissible 
locally analytic representation of $G$.
\endproclaim
\demo{Proof}
From corollary~3.6.16 we conclude that $V\otimes_K W$ is a locally analytic
representation of $G$.  By assumption $V$ is an $LB$-space and $W$
is finite dimensional, and so their tensor product is clearly also an
$LB$-space.

Since $W$ is finite dimensional we may find an analytic open subgroup
$H$ of $G$ such that $W$ is an $\H$-analytic representation of $H$.  
If $H'$ is an analytic open subgroup of $H$ then proposition~3.6.6
yields an isomorphism
$(V\otimes_K W)_{\H'-\an} \iso V_{\H'-\an} \otimes_K W.$
Proposition~6.1.2 yields an $H'$-equivariant
closed embedding $V_{\H'-\an} \rightarrow
\An(\H',K)^n$ for some $n$, and thus we obtain a closed embedding
$V_{\H'-\an} \otimes_K W \rightarrow \An(\H',K)^n \otimes_K W$
(both source and target being equipped with the diagonal action of $G$).
The remark following the proof of lemma~3.6.4
yields an isomorphism
$\An(\H',K)^n \otimes_K W \iso \An(\H',K)^n\otimes_K W,$
in which the source is equipped with the diagonal action of $G$,
and the target is equipped with the tensor product of the right regular
$G$-action on $\An(\H',K)^n$ and the trivial $G$-action on $W$.  Thus,
if $d$ denotes the dimension of $W$, we obtain altogether a $G$-equivariant
closed embedding $(V\otimes_K W)_{\H'-\an} \rightarrow \An(\H',K)^{nd}.$
This shows that $V\otimes_K W$ is admissible, as required.
\qed\enddemo

\proclaim{Lemma 6.1.6} If $V$ is an admissible  locally analytic 
representation of $G$, then for any good analytic open subgroup
$H$ of $\G$ (in the sense of subsection~5.2),
the space $V_{\open{\H}-\an}$ of $\open{\H}$-analytic
vectors in $V$ is a nuclear Fr\'echet space.
\endproclaim
\demo{Proof}
Write $\open{\H} = \bigcup_{n=1}^{\infty} \H_{r_n},$
where $\{r_n\}_{n\geq 1}$ is a
strictly increasing sequence of real numbers belonging to
$|\overline{L}^{\times}|$
and converging to 1.
Then
$\{\H_{r_n}\}_{n\geq 1}$ is an increasing sequence of analytic open
subgroups, with the property that each of the maps
$\H_{r_n} \rightarrow \H_{r_{n+1}}$ is
relatively compact,
and that each $\H_{r_n}$ is normalised by $H$.
This latter hypothesis implies that if $V$ is an
$H$-representation, then for each value of $n,$
the space $V_{\H_{r_n}-\an}$ is naturally an $H$-representation.
As in subsection~5.2, write $\open{H} := \open{\H}(L)$, and similarly, write
$H_{r_n} := \H_{r_n}(L)$ for each $n\geq 1$.

In the course of proving proposition~6.1.3, it was observed that each
of the maps $V_{\H_{r_{n+1}}-\an} \rightarrow V_{\H_{r_n}-\an}$ is 
compact.  Passing to the projective limit over $n$, we find that
$V_{\open{\H}-\an} = \plim{n} V_{\H_{r_n}-\an}$ is a projective limit
of a sequence of Banach spaces with compact transition maps,
and so is a nuclear Fr\'echet space.
\qed\enddemo

\proclaim{Proposition 6.1.7} If $V$ is a locally analytic representation
of $G$ on an $LB$-space, then the following are equivalent.

(i) $V$ is an admissible locally analytic representation.

(ii) For every good analytic open subgroup $\H$ of $\G,$
there exists a natural number $n$ and an
$\open{H}$-equivariant closed embedding
$V_{\open{\H}-\an} \rightarrow \An(\open{\H},K)^n.$

(iii) For a cofinal sequence of good analytic open subgroups $\H$ of $\G,$
there exists a natural number $n$ and an
$\open{H}$-equivariant closed embedding
$V_{\open{\H}-\an} \rightarrow \An(\open{\H},K)^n.$
\endproclaim
\demo{Proof}
Suppose first that $V$ is admissible, and that $\H$ is a good analytic
open subgroup of $G$.
As in the proof of lemma~6.1.6, we write $\open{\H} = \bigcup_{n=1}^{\infty}
\H_{r_n}.$ 
Proposition~6.1.2 yields an $H_{r_1}$-equivariant closed embedding
$V_{\H_{r_1}-\an} \rightarrow \An(\H_{r_1},K)^m$ for some natural number $m$.
Since $\open{H}$ acts on $V_{\H_{r_1}-\an},$ this induces in a natural way
an $\open{H}$-equivariant closed embedding
$V_{\H_{r_1}-\an} \rightarrow \Con(\open{H},K)_{\H_{r_1}-\an}^m.$

Passing to $\open{\H}$-analytic vectors,
and appealing to proposition~3.4.4 and corollary~3.4.10,
we obtain an $\open{H}$-equivariant closed embedding
$$(V_{\H_{r_1}-\an})_{\open{\H}-\an} \rightarrow
(\Con(\open{H},K)_{\H_{r_1}-\an}^m)_{\open{\H}-\an}.
\tag 6.1.8$$
Proposition~3.4.12 (and the fact that $V$ is assumed to be an $LB$-space)
implies that
the natural map $(V_{\H_{r_1}-\an})_{\open{\H}-\an} \rightarrow
V_{\open{\H}-\an}$ is an isomorphism.
Similarly, the natural map
$(\Con(\open{H},K)_{\H_{r_1}-\an}^m)_{\open{\H}-\an} \rightarrow
\Con(\open{H},K)_{\open{\H}-\an}^m$ is an isomorphism.
Thus we may rewrite~(6.1.8) as an $\open{H}$-equivariant
closed embedding
$$V_{\open{\H}-\an} \rightarrow \Con(\open{H},K)_{\open{\H}-\an}^m.
\tag 6.1.9$$
Since $\H$ was an arbitrary analytic open subgroup of $G$,
we find that~(i) implies~(ii).

It is obvious that~(ii) implies~(iii), and so we turn to showing
that~(iii) implies~(i).
Suppose that $\H$ is a good analytic open subgroup of $G$ for which
we have a $\open{H}$-equivariant closed embedding
$V_{\open{\H}-\an} \rightarrow \An(\open{\H},K)^m$ for some natural
number $m$.  Since $H$ acts on the source (a consequence of the
fact that $H$ normalises $\open{\H}$) this lifts to
an $H$-equivariant closed embedding
$V_{\open{\H}-\an} \rightarrow \Con(H,K)^m_{\open{\H}-\an}.$
Passing to $\H$-analytic vectors, and appealing to
Propositions~3.3.18 and~3.3.23,
we obtain a closed $H$-equivariant embedding
$$(V_{\open{\H}-\an})_{\H-\an} \rightarrow
(\Con(H,K)_{\open{\H}-\an})_{\H-\an}^m.\tag 6.1.10$$

A consideration of propositions~3.3.7 and~3.4.14
shows that the natural map
$$\An(\H,K)^m \rightarrow \Con(H,K)^m_{\open{\H}-\an}$$
induces an isomorphism
$$\An(\H,K)^m \iso (\Con(H,K)^m_{\open{\H}-\an})_{\H-\an},$$
while proposition~3.4.14 also implies that the natural map
$(V_{\open{\H}-\an})_{\H-\an} \rightarrow V_{\H-\an}$ is an
isomorphism. 
Thus~(6.1.10) can be rewritten as a closed embedding
$V_{\H-\an} \rightarrow \An(\H,K)^m.$
Since $H$ was one of a cofinal sequence of analytic
open subgroups of $G$, we see that $V$ satisfies the conditions
of definition~6.1.1, and so is an admissible representation of $G$.
\qed\enddemo

Let $H$ be a good analytic open subgroup of $G$.  
If as in the proof of lemma~6.1.6 we write
$\open{\H} = \bigcup_{n=1}^{\infty} \H_{r_n}$,
then 
$\DAn(\open{\H},K)_b \iso \ilim{n} \DAn(\H_{r_n},K)_b.$
As observed in subsection~5.2, $\DAn(\open{\H},K)_b$
is a naturally a topological $K$-algebra of compact type.

Suppose now that $V$ is an admissible locally analytic representation
of $G$.
As observed in the proof of lemma~6.1.6, we see that
$V_{\open{\H}-\an} = \plim{n} V_{\H_n-\an}$ is a projective
limit of Banach spaces with compact transition maps.  Passing to duals,
we find that
$(V_{\open{\H}-\an})'_b = \ilim{n} (V_{\H_n-\an})'_b.$
Corollary~5.1.8 shows that
for each $n$, the Banach space $(V_{\H_n-\an})'_b$ is a left
$\DAn(\H_n,K)_b$-module, and that the multiplication map
$\DAn(\H_n,K)_b \times (V_{\H_n-\an})'_b \rightarrow (V_{\H_n-\an})'_b$
is jointly continuous.
Passing to the locally convex inductive limit in $n$, we obtain a map
$$\DAn(\open{\H},K)_b \times (V_{\open{\H}-\an})'_b \rightarrow
(V_{\open{\H}-\an})'_b,\tag 6.1.12$$
which makes $(V_{\open{\H}-\an})'_b$
a left $\DAn(\open{\H},K)$-module,
and which is {\it a priori} separately continuous.
Since each of $\DAn(\open{\H},K)_b$ and $(V_{\open{\H}-\an})'_b$
is of compact type, it follows from proposition~1.1.31 that in
fact this map is jointly continuous, that is, that
$(V_{\open{\H}-\an})'_b$ is a topological
$\DAn(\open{\H},K)_b$-module.

\proclaim{Lemma 6.1.13} If $V$ is an admissible locally
analytic $G$-representation, and if $\H$ is a good analytic open
subgroup of $G$, 
then the topological $\DAn(\open{\H},K)_b$-module $(V_{\open{\H}-\an})'_b$ is
finitely generated (in the sense of definition~1.2.1~(iii)).
\endproclaim
\demo{Proof}
Dualising the $\open{H}$-equivariant
closed embedding $V_{\open{\H}-\an} \rightarrow \An(\open{\H},K)^m$
of part~(ii) of proposition~3.4.2 yields a surjection
$\DAn(\open{\H},K)_b^m \rightarrow (V_{\open{\H}-\an})'_b$.
Since this is a surjection of spaces of compact type, 
it is an open map, and so $(V_{\open{\H}-\an})'_b$ is finitely generated
in the strict sense of definition~1.2.1~(iii).
\qed\enddemo

Suppose now that
$J \subset H$ is an inclusion of good analytic open subgroups of $G$,
with the property that $\open{H}$ normalises $\open{\J}$. 
Lemma~6.1.13 shows
that the spaces $(V_{\open{\H}-\an})'_b$ and $(V_{\open{\J}-\an})'_b$ are
modules over the convex $K$-algebras 
$\DAn(\open{\H},K)_b$ and $\DAn(\open{\J},K)_b$, respectively.
Let us simplify our notation a little,
and denote these topological rings
by $D(\open{\H})$ and $D(\open{\J})$, respectively.
The inclusion $\open{\J} \rightarrow \open{\H}$ induces
a continuous ring homomorphism $D(\open{\J}) \rightarrow D(\open{\H})$.  

There is another ring that we will need to consider.
We let $D(\open{\J},\open{H})$ denote the strong dual to the space
$\Con(\open{H},K)_{\open{\J}-\an}$.    The restriction map
$\Con(\open{H},K)_{\open{\J}-\an} \rightarrow \Con(\open{J},K)_{\open{\J}-\an}
\iso \An(\open{\J},K)$ (the isomorphism is provided by proposition~3.4.5)
yields a closed embedding
$D(\open{\J}) \rightarrow D(\open{\J},\open{H}).$  
The topological ring structure
on $D(\open{\J})$ extends naturally to a ring structure on
$D(\open{\J},\open{H}),$ and one immediately sees that
$D(\open{\J},\open{H}) = \bigoplus_{h \in \open{H}/\open{J}}
D(\open{\J}) * \delta_h$. 
(As indicated, the
direct sum ranges over a set of coset representatives of $\open{J}$ in
$\open{H}$.)
Dualising the continuous injection
$\An(\open{\H},K) \rightarrow \Con(\open{H},K)_{\open{\J}-\an}$ yields
a continuous ring homomorphism $D(\open{\J},\open{H}) \rightarrow D(\open{\H}),$
which extends the homomorphism
$D(\open{\J}) \rightarrow D(\open{\H})$.
(The notions introduced in this paragraph are similar to
those introduced in the proof of proposition~5.3.1.)

Since $\open{H}$ normalises $\open{\J}$, there is a natural action of
$\open{H}$ on $V_{\open{\J}-\an},$
and hence the $D(\open{\J})$-module structure on
$(V_{\open{\J}-\an})'_b$ extends in a natural fashion to a
$D(\open{\J},\open{H})$-module structure. 
The continuous $H$-equivariant injection
$V_{\open{\H}-\an} \rightarrow V_{\open{\J}-\an}$ induces
a continuous map
$$(V_{\open{\H}-\an})'_b \rightarrow (V_{\open{\J}-\an})'_b,\tag 6.1.14$$
which is compatible with the natural map
$D(\open{\J},\open{H}) \rightarrow D(\open{\H})$
and the topological module structures on its source and target.
The map~(6.1.14) thus induces a continuous map
of $D(\open{\H})$-modules
$$ D(\open{\H}) \cotimes_{D(\open{\J},\open{H})} (V_{\open{\J}-\an})'_b
\rightarrow (V_{\open{\H}-\an})'_b .  \tag 6.1.15$$

The following result deals with a more general situation.

\proclaim{Proposition 6.1.16}
If $M$ is a Hausdorff convex $K$-vector space of compact type equipped with a
finitely generated topological
$D(\open{\J},\open{H})$-module structure,
then there is a natural $\open{H}$-equivariant isomorphism
$$(D(\open{\H})\cotimes_{D(\open{\J},\open{H})} M)'_b
\iso (M'_b)_{\open{\H}-\an}.$$
\endproclaim
\demo{Proof}
Since $M$ is of compact type, its
strong dual $M'_b$ is a nuclear Fr\'echet space.
Corollary~3.4.5 thus
implies that there is a natural isomorphism
$$(M'_b)_{\open{\H}-\an} \iso
(\An(\open{\H},K)\cotimes_K M'_b)^{\Delta_{1,2}(\open{H})}.\tag 6.1.17$$

Proposition~1.2.5, and the fact that compact type spaces
are hereditarily complete, implies that the tensor product
$$D(\open{\H}) \cotimes_{D(\open{\J},\open{H})} M \tag 6.1.18$$
is the quotient of the tensor product
$D(\open{\H}) \otimes_{K} M$
by the closure of its subspace spanned by expressions of
the form $\mu*\nu \otimes m - \mu \otimes \nu * m,$
with $\mu \in D(\open{\H}),$ $\nu \in D(\open{\J},\open{\H})$,
and $m \in M$.
Since the delta functions $\delta_h$, for $h \in \open{H},$
span a strongly dense subspace of $D(\open{\J},\open{H})$ 
(see the discussion following definition~2.2.3),
we see that~(6.1.18) can also be described as the quotient of
$D(\open{\H}) \otimes_{K} M$
by the closure of its subspace spanned by expressions of the
form $\mu*\delta_h \otimes m - \mu \otimes \delta_h * m,$
with $\mu \in D(\open{\H}),$ $h \in \open{H},$
and $m \in M$.  

This description of~(6.1.18), together with proposition~1.1.32~(ii),
implies that there is a natural isomorphism
$$(D(\open{\H})\cotimes_{D(\open{\J},\open{H})} M)'_b \iso
(\An(\open{\H},K)\cotimes_K M'_b)^{\Delta_{1,2}(\open{H})}.$$
Thus~(6.1.17) yields the required isomorphism.
\qed\enddemo

\proclaim{Corollary 6.1.19}
The morphism~(6.1.15) is a topological isomorphism.
\endproclaim
\demo{Proof}
Proposition~6.1.16 yields an isomorphism
$(D(\open{\H})\cotimes_{D(\open{\J},\open{H})} (V_{\open{\J}-\an})'_b)'_b
\iso (V_{\open{\J}-\an})_{\open{\H}-\an}.$
By corollary~3.4.15, the natural map
$(V_{\open{\J}-\an})_{\open{\H}-\an} \rightarrow V_{\open{\H}-\an}$
is an isomorphism.
It is easily checked that dualising the composite of these two
isomorphisms yields the map~(6.1.15).
Thus this map is an isomorphism, as claimed.
\qed\enddemo

\proclaim{Theorem 6.1.20} 
Let $H$ be a good analytic open subgroup of $G$.
Passage to the dual induces an anti-equivalence of categories
between the category of admissible locally analytic $\open{H}$-representations
(with morphisms being continuous $\open{H}$-equivariant maps) and
the category of coadmissible locally convex topological
modules over the Fr\'echet-Stein algebra $\DLa(\open{H},K)_b$.
\endproclaim
\demo{Proof}
Suppose first that $V$
is equipped with an admissible locally analytic representation
of $\open{H}$.
Let $\{r_n\}$ be a strictly decreasing
sequence of real numbers lying in the open
interval between 0 and 1, that converges to 0.
Then $V \iso \ilim{n} V_{\open{\H}_{r_n}-\an},$ and so
$V'_b  \iso \plim{n} (V_{\open{\H}_{r_n}-\an})'_b.$
In the proof of proposition~5.3.1, it is proved that
the isomorphism
$\DLa(\open{H},K) \iso \plim{n} D(\open{\H}_{r_n},\open{H})$
exhibits $\DLa(\open{H},K)$ as a weak Fr\'echet-Stein structure.
We will show that $V'_b$ is coadmissible with
respect to this weak Fr\'echet-Stein structure.

For each $n\geq 1,$ lemma~6.1.13 shows that
$V_{\open{\H}_{r_n}-\an}$ is finitely generated
over $D(\open{\H}_{r_n}),$ and so in particular over
$D(\open{\H}_{r_n},\open{H})$,  
while corollary~6.1.19 yields an isomorphism
$D(\open{\H}_{r_n},\open{H})\cotimes_{D(\open{\H}_{r_{n+1}},\open{H})}
(V_{\open{\H}_{r_{n+1}}-\an})'_b \iso (V_{\open{\H}_{r_n}-\an})'_b.$
Thus the projective system $\{V_{\open{\H}_{r_n}-\an}\}_{n\geq 1}$ satisfies
the conditions of definition~1.2.8, showing that $V'_b$ is a coadmissible
$\DLa(\open{H},K)$-module.

Conversely, suppose that $M = \plim{n} M_n$ is a coadmissible
$\DLa(\open{H},K)$-module, so that $M_n$ is a Hausdorff 
finitely generated locally convex
$D(\open{\H}_{r_n},\open{H})$-module, for each $n \geq 1.$
For each such value of $n$, we obtain a surjection
of spaces of compact type $D(\open{\H}_{r_n},\open{H})^m
\rightarrow M_n$, for some natural number $m$. Dualising
this yields a closed embedding of nuclear Fr\'echet spaces
$(M_n)'_b \rightarrow \Con(\open{H},K)^m_{\open{\H}_{r_n}-\an}.$
Passing to $\H_{r_n}$-analytic vectors,
and applying propositions~3.3.23 and~3.4.14
we obtain a closed embedding
$$((M_n)'_b)_{\H_{r_n}-\an}
\rightarrow \Con(\open{H},K)^m_{\H_{r_n}-\an}.\tag 6.1.21$$

Condition~(ii) of definition~1.2.8, together with
proposition~6.1.16, yields, for each $n\geq 1$, a natural isomorphism
$(M_n)'_b \iso ((M_{n+1})'_b)_{\open{\H}_{r_n}-\an}.$
Passing to $\H_{r_n}$-analytic vectors, and applying proposition~3.4.14,
we obtain an isomorphism
$$((M_n)'_b)_{\H_{r_n}-\an} \iso ((M_{n+1})'_b)_{\H_{r_n}-\an}.$$
Since $M'_b \iso \ilim{n} (M_n)'_b,$
we deduce that the natural map
$$((M_n)'_b)_{\H_{r_n}-\an} \rightarrow (M'_b)_{\H_{r_n}-\an}$$
is an isomorphism.  In light of the closed embedding~(6.1.21),
we conclude that $M'_b$ is an admissible locally analytic
$\open{H}$-representation, as claimed.
\qed\enddemo

The next result shows that our definition of an admissible locally
analytic representation of $G$ agrees with that of \cite{\SCHTNEW,
p.~33}.

\proclaim{Corollary 6.1.22} If $G$ is a locally $L$-analytic group, and
if $V$ is a convex $K$-vector space of compact type equipped with a locally
analytic representation of $G$, then $V$ is admissible if and only
if $V'_b$
is a coadmissible module with respect to the natural $\DLa(H,K)$-module
structure on $V'$, for some (or equivalently, every) compact open subgroup
$H$ of $G$.
\endproclaim
\demo{Proof}
Let $J$ be a good analytic open subgroup of $G$ contained in $H$.
The algebra $\DLa(H,K)$ is free of finite rank (on both sides)
over $\DLa(\open{J},K)$, and so $V'_b$ is coadmissible
as a $\DLa(H,K)$-module if and only if it is so as a $\DLa(\open{J},K)$-module.
Similarly, $V$ is an admissible locally analytic $G$-representation if and only
if it is an admissible locally analytic $\open{J}$-representation.
The corollary thus follows from theorem~6.1.20.
\qed\enddemo

The following result is a restatement of \cite{\SCHTNEW, prop.~6.4}.

\proclaim{Corollary 6.1.23} If $G$ is a locally $L$-analytic group,
then the category of admissible locally analytic $G$-representations
and continuous $G$-equivariant morphisms is closed under passing to 
closed subobjects and Hausdorff quotients.  Furthermore, any morphism
in this category is necessarily strict.   Consequently, this
category is abelian.
\endproclaim
\demo{Proof}
This is a consequence of the general properties of coadmissible modules
over Fr\'echet-Stein algebras, summarised in theorem~1.2.11 and the
remarks that follow it.
\qed\enddemo

\section{6.2}
The following class of locally analytic representations was
introduced by Schneider and Teitelbaum in \cite{\SCHTAN}.

\proclaim{Definition 6.2.1}
A locally analytic representation $V$ of $G$ is called strongly
admissible  if it is an $LB$-space, and if its dual space $V'$ is 
finitely generated as a left $\DLa(H,K)$-module for one (and hence
every) compact open subgroup $H$ of $G$.
\endproclaim

Note that theorem~5.1.15~(iii) implies that a locally analytic
representation of $G$ on an $LB$-space $V$ satisfies the condition of
definition~6.2.1 if and only if 
$V$ admits a closed embedding into $\La(H,K)^n$ for some natural
number $n$, and one (or equivalently, every) compact open subgroup
of $G$.  In particular, any such $V$ is necessarily of compact type.

We also remark that the strongly admissible  locally analytic representations of
$G$ are precisely the continuous
duals of those $\DLa(H,K)$-modules which are analytic in the sense
of \cite{\SCHTUF, p.~112}.

\proclaim{Proposition 6.2.2} If $V$ is a strongly
admissible locally analytic representation of $G$, then $V$
is admissible.
\endproclaim
\demo{Proof}
Since $V$ is strongly admissible, it is an $LB$-space by assumption.
If $H$ is an analytic open subgroup of $G$ then
by assumption we may find a surjection of left $\DLa(G,K)$-modules
$\DLa(G,K)^n \rightarrow V'$ for some natural number $n$, which
by theorem~3.3.1 arises by dualising a $G$-equivariant closed embedding
$V \rightarrow \La(G,K)^n.$  
Passing to $\H$-analytic vectors and applying proposition~3.3.23
and corollary~3.3.26
we obtain an $H$-equivariant closed embedding 
$V_{\H-\an} \rightarrow \An(\H,K)^n$.
It follows that $(V_{\H-\an})'$ is finitely generated as a left $\DAn(\H,K)$-module.
\qed\enddemo

A basic method for producing strongly admissible locally analytic representations of $G$ is by
passing to the locally analytic vectors of continuous representations
in $G$ that satisfy the admissibility condition introduced in \cite{\SCHTIW}.
We begin by reminding the reader of that definition.

\proclaim{Proposition-Definition 6.2.3}  
A continuous $G$-action on a Banach space $V$ is said to be an
admissible continuous representation of $G$ (or an admissible Banach space
representation of $G$) if $V'$ is finitely
generated as a left $\DCon(H,K)$-module for one (and hence every)
compact open subgroup $H$ of $G$.  (This is equivalent, by theorem~5.1.15~(i),
to the existence of
a closed $H$-equivariant embedding $V \rightarrow \Con(H,K)^n$
for some $n\geq 0$.)
\endproclaim
\demo{Proof} We must check that the definition is independent
of the choice of compact open subgroup $H$ of $G$.  This follows
immediately from the fact that if $H_1 \subset H_2$ is an inclusion
of compact open subgroups of $H$, and if we write $H_2 = \coprod h_i H_1$
as a finite disjoint union of right $H_1$ cosets, then
$\DCon(H_2,K) = \bigoplus \DCon(H_1,K)\delta_{h_i}$, and so
$\DCon(H_2,K)$ is finitely generated as a left $\DCon(H_1,K)$-module.
\qed\enddemo

Note that in \cite{\SCHTIW}, the authors restrict their attention
to the case of local $K$.  In this case, it follows
from \cite{\SCHTIW, lem.~3.4} that definition~6.2.3 
agrees with the notion of admissibility introduced in that
reference.  (See also proposition~6.5.6 below.)

\proclaim{Proposition 6.2.4}
If $V$ is a Banach space equipped
with an admissible continuous representation of $G$ 
then $V_{\la}$ is a strongly
admissible locally analytic representation of $G$.
\endproclaim
\demo{Proof} Since $V$ is a Banach space, proposition~3.5.6 implies
that $V_{\la}$ is an $LB$-space.  By assumption, we may find
a surjection $\DCon(H,K)^n \rightarrow V'$ for some natural number $n$
and some compact open subgroup $H$ of $G$,
and so a closed $H$-equivariant embedding $V\rightarrow \Con(H,K)^n.$
Now by proposition~3.5.11 there is an isomorphism
$\Con(H,K)_{\la} \iso \La(H,K).$  Since $\La(H,K)$ is
of compact type, we obtain by proposition~3.5.10 a closed embedding
$V_{\la} \rightarrow \La(H,K)^n.$  Dualising this yields a
surjection $\DLa(H,K)^n \rightarrow (V_{\la})',$ proving that
$V_{\la}$ is strongly admissible.  
\qed\enddemo

When $K$ is local,
the preceding result has been established independently by
Schneider and Teitelbaum \cite{\SCHTNEW, thm.~7.1}.
(As was noted in the remark following the
proof of theorem~3.5.7, the topology with which
these authors equip $V_{\la}$ is {\it a priori} coarser
than the topology with which we equip that space.
In the case of $V$ being a continuous admissible representation,
however, both topologies are of compact type, since they
each underly a strongly admissible locally analytic
representation of $G$.  Thus the topologies coincide.)

The remainder of this subsection is devoted to establishing
some fundamental properties of admissible continuous $G$-representations.
In the local case, all these results (other than propositions~6.2.6
and~6.2.7)
are contained in \cite{\SCHTIW}.
Using a little more functional analysis, we are able to extend
these results to the case where $K$ is discretely valued.

In fact, the first three results require no assumption on $K$ at all.

\proclaim{Proposition 6.2.5} If $V$ is a $K$-Banach space
equipped with an admissible continuous representation of $G$
and $W$ is a closed $G$-invariant $K$-subspace of $V$,
then the continuous $G$-representation on $W$ is also
admissible.
\endproclaim
\demo{Proof} This follows immediately from the fact that
$W'$ is a quotient of $V'$.
\qed\enddemo

\proclaim{Proposition 6.2.6} If $V$ is a $K$-Banach space equipped with
an admissible continuous representation of $G$ and $W$ is a finite-dimensional
continuous representation of $G$ then $V\otimes_K W,$ equipped with
the diagonal action of $G$, is an admissible continuous representation
of $G$.
\endproclaim
\demo{Proof}
Replacing $G$ by a compact open subgroup if necessary, we may assume
that $G$ is compact and that there is a closed $G$-equivariant embedding
$V \rightarrow \Con(G,K)^n$ for some $n\geq 0$.  Tensoring with $W$ over
$K$ yields a closed embedding
$V\otimes_K W \rightarrow \Con(G,K)^n \otimes_K W.$  The remark following
the proof of lemma~3.2.11 shows that there is a natural isomorphism
$\Con(G,K)^n \otimes_K W \iso \Con(G,K)^n \otimes_K W$ which intertwines
the diagonal $G$-action on the source with the tensor product of the
right regular $G$-action and the trivial $G$-action on the target.
If $d$ denotes the dimension of $W$, then we may thus find a
$G$-equivariant closed embedding
$V\otimes_K W \rightarrow \Con(G,K)^{nd}.$  This proves that $V\otimes_K W$
is equipped with an admissible continuous representation of $G$.
\qed\enddemo

\proclaim{Proposition 6.2.7} If $0 \rightarrow U \rightarrow V
\rightarrow W \rightarrow 0$ is a $G$-equivariant exact sequence of $K$-Banach
spaces equipped with continuous $G$-actions, and
each of $U$ and $W$ is an admissible continuous $G$-representation,
then $V$ is also an admissible continuous $G$-representation.
\endproclaim
\demo{Proof}
Without loss of generality, we may assume that $G$ is compact.
Passing to strong duals, and taking into account corollary~5.1.7,
we obtain the short exact sequence of $\DCon(G,K)$-modules
$0 \rightarrow W'_b \rightarrow V'_b \rightarrow U'_b \rightarrow 0.$
If each of the end terms is finitely generated over $\DCon(G,K)$,
then the same is true of the middle term.
\qed\enddemo

As in \cite{\SCHTNEW}, 
the remainder of our results depend on the fact that
$\DCon(G,\Cal O_K)$ is Noetherian.   When $K$ is local,
this is provided by \cite{\LAZ, V.2.2.5}.  In fact, it
is easily extended to the case where $K$ is discretely valued.

\proclaim{Theorem 6.2.8}
If $G$ is a compact locally $L$-analytic
group and $K$ is discretely valued, then the $K$-Banach algebra $\DCon(G,K)$,
and its unit ball $\DCon(G,\Cal O_K)$, are both Noetherian rings.
\endproclaim
\demo{Proof}
It follows from \cite{\LAZ, V.2.2.4} that
$\DCon(G,\Z_p)$ is Noetherian.
Indeed, that reference shows that $G$ contains an open subgroup
$H$ such that $\DCon(H,\Z_p)$ is equipped with
an exhaustive filtration $F^{\bullet} \DCon(H,\Z_p)$,
with respect to which it is complete, whose associated graded ring 
$\roman{Gr}_F^{\bullet} \DCon(H,\Z_p)$ is isomorphic
to the enveloping algebra over the ring $\Gamma:= \F_p[\epsilon]$
of a finite rank Lie algebra over $\Gamma$.

If we endow $\Cal O_K$ with its $p$-adic filtration, then we may endow
the completed tensor product
$\Cal O_K \cotimes_{\Z_p} \DCon(H,\Z_p)$
with the (completion of the) tensor product filtration.
The associated graded ring is then isomorphic to the tensor product
$(\Cal O_K/p)\otimes_{\F_p} \roman{Gr}_F^{\bullet} \DCon(H,\Z_p),$
and thus is isomorphic to the enveloping algebra over the ring
$(\Cal O_K /p)[\epsilon]$ of a finite rank Lie algebra over this ring.
In particular, it is Noetherian, and hence the same is true of the
completed tensor product
$\Cal O_K \cotimes_{\Z_p} \DCon(H,\Z_p)$
itself.  This ring is naturally isomorphic to $\DCon(H,\Cal O_K)$,
and $\DCon(G,\Cal O_K)$ is free of finite rank as a module over
this latter ring.  Thus the theorem is proved.
\qed\enddemo

\proclaim{Proposition 6.2.9} If $K$ is discretely valued, then any
continuous $G$-equivariant map between admissible continuous
$G$-representations is strict.
\endproclaim
\demo{Proof}
We may replace $G$ by a compact open subgroup, and thus suppose
that $G$ is compact.  If $\phi: V \rightarrow W$ is a continuous
$G$-equivariant map between two admissible continuous $G$-representations,
then to show that $\phi$ is strict, it suffices to show that
$\phi': V'_b \rightarrow W'_b$ is strict \cite{\TVS, cor.~3, p.~IV.30}.
The proposition now follows from
theorem~6.2.8 and proposition~1.2.4.
\qed\enddemo

\proclaim{Proposition 6.2.10} Suppose that $G$ is compact,
and that $K$ is discretely valued.
The association of $V'$ to $V$ induces
an anti-equivalence between the category of admissible continuous
$G$-representations (with morphisms being continuous $G$-equivariant maps)
and the category of finitely generated $\DCon(G,K)$-modules.
\endproclaim
\demo{Proof}
Since the transpose of a continuous linear map between Hausdorff convex
$K$-vector spaces vanishes if and only if the map itself does, we see that
the passage from $V$ to $V'$ is faithful.
We turn to proving that this functor is full.

Suppose that $V$ and $W$ are admissible continuous $G$-representations,
and that we are given a $\DCon(G,K)$-linear map
$$W' \rightarrow V'.\tag 6.2.11 $$
We must show that this map arises as the transpose of a continuous
$G$-equivariant map $V \rightarrow W$.

Let us choose surjections of $\DCon(G,K)$-modules
$\DCon(G,K)^m \rightarrow V'$ and
$\DCon(G,K)^n \rightarrow W',$ for some $m,n\geq 0.$
Then~(6.2.11) may be lifted to yield a commutative square
of $\DCon(G,K)$-linear maps
$$\xymatrix{\DCon(G,K)^n_b \ar[d]\ar[r] & \DCon(G,K)^m_b \ar[d] \\
W'_b \ar[r] & V'_b .} \tag 6.2.12 $$
Part~(iii) of corollary~5.1.7 shows that, if we remove the lower horizontal
arrow from~(6.2.12), then the resulting diagram arises by dualising a
diagram of continuous $G$-equivariant maps of the form 
$$\xymatrix{V \ar[d] & W \ar[d] \\ \Con(G,K)^m \ar[r] & \Con(G,K)^n.}
\tag 6.2.13$$
Furthermore, part~(i) of theorem~5.1.15
implies that the vertical arrows in~(6.2.13)
are closed embeddings.

The fact that the lower horizontal arrow of~(6.2.12)
can be filled in implies that the lower horizontal arrow of~(6.2.13)
restricts to a $G$-equivariant map $\phi: V \rightarrow W.$
By construction, the transpose of $\phi$ is equal to the morphism~(6.2.11)
with which we began, and so the functor under consideration is full,
as claimed.

In order to prove that our functor is essentially surjective,
it suffices to show
that any surjection of $\DCon(G,K)$-modules 
$$\DCon(G,K)^n \rightarrow M \tag 6.2.14$$ arises by dualising a closed 
$G$-equivariant embedding $V \rightarrow \Con(G,K)^n$.
For this, it suffices in turn to show that the kernel of~(6.2.14) is closed
in $\DCon(G,K)_s^n$.  For if we then let $M_s$ denote $M$
equipped with the topology obtained by regarding it as a quotient
of $\DCon(G,K)_s^n,$ we obtain $V$ as the closed subspace
$(M_s)'$ of $\Con(G,K)^n \iso (\DCon(G,K)_s^n)'$.

We must show that any $\DCon(G,K)$-submodule $N$ of $\DCon(G,K)^n_s$
is closed.  
Since $\DCon(G,K)$ is the dual of the Banach space $\Con(G,K)$,
it suffices to show that the intersection of $N$ with the unit
ball of $\DCon(G,K)_b^n$ is closed in $\DCon(G,K)_s$
\cite{\TVS, cor.~3, p.~IV.25}.

As in the statement of theorem~6.2.8,
let $\DCon(G,\Cal O_K)$ denote the unit ball of the Banach
algebra $\DCon(G,K)_b.$
If $\DCon(G,\Cal O_K)_s$ denotes $\DCon(G,\Cal O_K)$
equipped with the topology induced from $\DCon(G,K)_s$,
then $\DCon(G,\Cal O_K)_s$ is c-compact 
\cite{\SCHNA, lem.~12.10}. 

Passing from $\DCon(G,K)$ to $\DCon(G,K)^n,$
we find that $\DCon(G,\Cal O_K)^n_s$ is the unit ball of $\DCon(G,K)_b^n$,
equipped with the topology induced from $\DCon(G,K)_s^n$.
The intersection
$N \cap \DCon(G,\Cal O_K)^n_s$
is a $\DCon(G,K)$-submodule
of $\DCon(G,\Cal O_K)^n_s$, and we must show that it is closed.

Since the ring $\DCon(G,\Cal O_K)$ is Noetherian, by theorem~6.2.8, 
the $\DCon(G,\Cal O_K)$-module
$N \cap \DCon(G,\Cal O_K)^n_s$
is finitely generated.
If $n_1',\ldots,
n_r'$ are generators of this module, then it is the image
of $\DCon(G,\Cal O_K)_s^r$ under the map
$$(\mu_1,\ldots,\mu_r) \mapsto \mu_1 * v_1' + \cdots + \mu_r * v_r'.
\tag 6.2.15 $$
Now~(6.2.15) is the restriction to $\DCon(G,\Cal O_K)_s^r$
of a corresponding map $\DCon(G,K)_s^r \rightarrow \DCon(G,K)^n,$
and part~(iii) of corollary~5.1.7 shows that this map arises as the
dual of a continuous map $\Con(G,K)^n \rightarrow \Con(G,K)^r,$
and so is continuous.  Thus~(6.2.15) is continuous.  Its source
being c-compact , it has c-compact, and hence closed, image.
Thus we conclude that $N \cap \DCon(G,\Cal O_K)^n_s$
is closed in $\DCon(G,\Cal O_K)^n_s,$ as required.
\qed\enddemo

\proclaim{Corollary 6.2.16}  If $K$ is discretely valued,
then the category of admissible continuous 
$G$-representations (with morphisms being continuous
$G$-equivariant maps) is abelian.
\endproclaim
\demo{Proof}
If $G$ is compact, then this follows from proposition~6.2.10.
Indeed, the ring $\DCon(G,K)$ is Noetherian by theorem~6.2.8,
and so the category of finitely generated $\DCon(G,K)$-modules
is abelian.  The general case follows easily from this one.
\qed\enddemo

We remark 
that proposition~6.2.9 is also a formal consequence of
this result.  As already noted,
in the case where $K$ is local,
proposition~6.2.10, and hence corollary~6.2.16 and
proposition~6.2.9, are established in \cite{\SCHTIW}.

\section{6.3} In this subsection we develop some connections
between the notions of smooth and locally algebraic representations,
and of admissible locally analytic representations.

\proclaim{Lemma 6.3.1} 
Suppose that $G$ is equal to the
group of $L$-valued points of an affinoid rigid analytic group
$\G$ defined over $L$, and that furthermore $G$ is Zariski dense in $\G$.
If $\An(\G,K)$ is regarded as a locally analytic representation of $G$
via the right regular $G$-action
then $\An(\G,K)^{\lie{g}}$ is one-dimensional, consisting precisely
of the constant functions.
\endproclaim
\demo{Proof} Lemma~4.1.4 shows that $\An(\G,K)^{\lie{g}}$ consists
of those rigid analytic function on $\G$ that are stabilised by
the right regular action of $G$. Since $G$ is Zariski dense in $\G,$
these are precisely the constant functions.
\qed\enddemo

\proclaim{Proposition 6.3.2} If $V$ is an admissible smooth representation
of a locally $L$-analytic group $G$, then,
when equipped with its finest convex topology,
$V$ becomes an admissible locally analytic representation of $G$
that is $U(\lie{g})$-trivial.
Conversely, if $V$ is an admissible locally analytic representation
of $G$ that is $U(\lie{g})$-trivial then the topology on $V$ is its
finest locally convex topology, and the $G$-representation on the
$K$-vector space underlying $V$ is an admissible smooth representation.
\endproclaim
\demo{Proof} First suppose that $V$ is an admissible smooth
representation of $G$.  Let $H_n$ be a cofinal sequence of
normal analytic open subgroups of $G$. Then $V = \bigcup_n V^{H_n}.$  Since
each $V^{H_n}$ is finite dimensional by assumption, we see
that $V$, when equipped with its finest convex topology,
is isomorphic to the locally convex inductive
limit $\ilim{n} V^{H_n},$ and so is of compact type.
Corollary~4.1.7 implies that $(V_{\la})^{\lie{g}}$ maps bijectively
onto $V,$ and thus in particular that $V$ is locally analytic
and that $V^{\lie{g}} = V$.  Corollary~4.1.5 then
shows that $V_{\H_n-\an} = V^{H_n}$.  Since by assumption this
is a finite dimensional space, its dual is certainly finitely generated
over $\DAn(\H_n,K)$.

Conversely, let $V$ be an admissible locally analytic representation
that is $U(\lie{g})$-trivial. It follows from corollary~4.1.5 that
the $G$-representation on $V$ is smooth.  We have to show that it is
furthermore an admissible smooth representation, and also that the 
topology on $V$ is its finest convex topology.

As above, let $H_n$ denote a cofinal
sequence of normal analytic open subgroups of $V$.
Theorem~3.3.1 shows that for each $n$ we may find an
embedding of $V_{\H_n-\an}$ into $\An(\H,K)^m$, for some natural number $m$,
and hence an embedding
$(V_{\H_n-\an})^{\lie{g}} \rightarrow (\An(\H,K)^m)^{\lie{g}}.$
Lemma~6.3.1 shows that the target of this injection is
finite dimensional, and thus that $(V_{\H_n-\an})^{\lie{g}}$
is also finite dimensional.

By assumption
$(V_{\H_n-\an})^{\lie{g}}$ is equal to $V_{\H_n-\an}$ for each $n$,
while corollary~4.1.5 shows that
$(V_{\H_n-\an})^{\lie{g}}$ maps bijectively onto the subspace
of $H_n$-fixed vectors in $V$.  Thus on the one hand we find that
the subspace of $H_n$-fixed vectors in $V$ is finite dimensional,
proving the the smooth $G$-representation on $V$ is admissible,
while on the other hand, theorem~3.6.12 provides an
isomorphism $\ilim{n} (V_{\H_n-\an})^{\lie{g}} =
\ilim{n} V_{\H_n-\an} = V_{\la} \iso V$, so that
$V$ is the inductive limit of a family of finite dimensional spaces,
and so is endowed with its finest convex topology.
\qed\enddemo

Thus there is an equivalence of categories between the category of admissible
smooth representations of $G$ and the category of
$U(\lie{g})$-trivial admissible locally
analytic representations of $G$.
(This result, and the next,
have been proved independently by Schneider and Teitelbaum
\cite{\SCHTNEW, thm.~6.5}.)

\proclaim{Corollary 6.3.3} If $V$ is an admissible locally analytic
representation of a locally $L$-analytic group $G$,
then the closed subspace $V^{\goth g}$ of $V$ is
endowed with its finest convex topology, and the $G$-representation on
$V^{\lie{g}}$ is an admissible smooth representation of $G$.
\endproclaim
\demo{Proof} Since $V^{\lie{g}}$ is a closed $G$-invariant subspace
of $V$,
proposition~6.1.4 implies
that $V^{\lie{g}}$ is again an admissible representation of $G$, which is
$U(\lie{g})$-trivial by construction.  The corollary follows
from proposition~6.3.2.
\qed\enddemo

For the sake of completeness, we also
recall the following definition and result.

\proclaim{Definition 6.3.4}
If $G$ is a locally $L$-analytic group then
we say that a smooth representation $V$ of $G$ is a
strongly admissible smooth representation if for one
(or equivalently, every) compact open subgroup $H$ of $G$,
there is an $H$-equivariant embedding of $V$ into the
a finite direct sum of copies of the space of locally constant
functions on $H$.
\endproclaim

This is immediately checked to be equivalent to the
definition of \cite{\SCHTUF, p.~115}.

\proclaim{Proposition 6.3.5} If $V$ is a strongly admissible locally
analytic representation of a locally $L$-analytic group $G$, then
the closed subspace $V^{\goth g}$ of $V$ is endowed with its finest convex
topology, and the $G$-representation on $V^{\lie{g}}$ is a strongly admissible
smooth representation of $G$.
Conversely, if we endow any strongly admissible smooth representation of $G$
with its finest convex topology, we obtain a strongly admissible locally
analytic representation of $G$.
\endproclaim
\demo{Proof}
This is \cite{\SCHTUF, prop.~2.1}.
\qed\enddemo

We now suppose that $\G$ is a reductive linear algebraic
group defined over $L$, and that $G$ is an open subgroup of $\G(L)$.
As in subsection~4.2, we let $\Cal R$ denote the category
of finite dimensional algebraic representations of $\G$
defined over $K$.  We will also use the other definitions
and notations introduced in that subsection.

\proclaim{Proposition 6.3.6} 
If $V$ is an admissible locally analytic 
representation of $G$ and $W$ is an object of $\Cal R$,
then $V_{W-\lalg}$ is a closed subspace of $V$, and the topology
induced on $V_{W-\lalg}$ by the topology on $V$ is its
finest convex topology.  In particular, $V_{W-\lalg}$ is again
an admissible locally analytic representation of $G$.
\endproclaim
\demo{Proof}
Proposition~4.2.10 shows that $V_{W-\lalg}$ is
a closed subspace of $V$, topologically isomorphic to 
$\Hom(W,V)_{\sm} \otimes_B W$ (where $B = \End_G(W)$).
Since $V$ is admissible, proposition~6.1.5 (applied to the tensor product
$V\otimes_K \check{W}$)
implies that $\Hom(W,V)$ is also admissible.
Proposition~6.3.3 shows that the space
$\Hom(W,V)_{\sm}$ is equipped with its finest convex topology,
and thus the same is true of the tensor product
$\Hom(W,V)_{\sm} \otimes_B W.$
The proposition follows.  
\qed\enddemo

\proclaim{Corollary 6.3.7} 
If $V$ is an admissible locally analytic
representation of $G$, and if $V$ is also locally algebraic as a
representation of $G$, then 
$V$ is equipped with is finest convex topology.
\endproclaim
\demo{Proof} 
Consider the natural map $$\bigoplus_{W\in \hat{\G}} V_{W-\lalg} \rightarrow V.
\tag 6.3.8$$
Proposition~6.3.6 shows that each summand in the source of~(6.3.8)
is of compact type when endowed with its finest convex topology.
Since $\hat{\G}$ is countable, we see that
if we endow the source with its finest convex topology then it is again
a space of compact type.  Corollary~4.2.7 then
shows that~(6.3.8) is a continuous bijection between spaces of compact type,
and thus that~(6.3.8) is a topological isomorphism.  
\qed\enddemo

In light of the preceding corollary, we make the following definition.

\proclaim{Definition 6.3.9}
If $V$ is a locally algebraic representation of $G$
that becomes admissible as a locally analytic representation 
when equipped with its finest convex topology,
then we say that $V$ is an admissible locally algebraic representation.
If $V$ is furthermore locally $W$-algebraic, for some object $W$ of $\Cal R$,
then we say that $V$ is an admissible locally $W$-algebraic representation.
\endproclaim

\proclaim{Proposition 6.3.10} If $V$ is an admissible locally
$W$-algebraic representation of~$G$,
then $V$ is isomorphic to a representation
of the form $U \otimes_B W,$ where $B$ denotes the semi-simple $K$-algebra
$\End_{\G}(W),$ and $U$ is an admissible smooth representation
of $G$ over $B$, equipped with its finest convex topology.
Conversely, any such tensor product is an admissible locally 
$W$-algebraic representation of $G$.
\endproclaim
\demo{Proof}  The first statement of the proposition
was observed in the course of proving proposition~6.3.6.
The converse statement follows from the remark preceding
the statement of proposition~4.2.4, together with
propositions~6.3.2 and 6.1.5.
\qed\enddemo

\proclaim{Proposition 6.3.11} If $V$ is a locally algebraic representation
of $G$, and if, following
proposition~4.2.4 and corollary~4.2.7, we write 
$V \iso \bigoplus_n U_n \otimes_{B_n} W_n,$
where $W_n$ runs over the elements of $\hat{\G}$, $B_n = \End_{\G}(W_n),$
and $U_n$ is a smooth representation of $G$ over $B_n$,
then $V$ is an admissible locally algebraic representation of $G$
if and only if the direct sum $\bigoplus_n U_n$ is an admissible
smooth representation of $G$.
\endproclaim
\demo{Proof} Let $H$ be an analytic open subgroup of $G$.
Propositions~3.3.27 and~3.6.6 and lemma~4.1.4 yield an isomorphism
$$V_{\H-\an} \iso \bigoplus_n (U_n)^H \otimes_{B_n} W_n$$
(where each side is equipped with its finest convex topology).
If $\bigoplus_n U_n$ is admissible, then the space
$\bigoplus_n (U_n)^H $ is finite dimensional,
and so the same is true of $V_{\H-\an}$.  Since $H$ was arbitrary,
we see that $V$ is an admissible locally analytic representation
of $G$.

On the other hand, if $\bigoplus_n U_n$ is not an admissible
smooth representation of $G$, then
we may find $H$ such that $\bigoplus_n (U_n)^H$ is not
finite dimensional.  In other words, we may find an infinite sequence
$\{n_i\}_{i \geq 1}$ such that each space $(U_{n_i})^H$ is non-zero,
and hence such that $\bigoplus_i W_{n_i}$ is a subspace
of $V_{\H-\an}$, necessarily closed, since $V_{\H-\an}$
is equipped with its finest convex topology.  
However, the direct sum $\bigoplus_i W_{n_i}$ (an infinite
direct sum of algebraic representations) cannot
embed as a closed subspace of a finite direct sum $\An(\H,K)^m$.
Thus neither can $V_{\H-\an}$, and so $V$ is not an admissible
locally analytic representation of $G$.
\qed\enddemo

\section{6.4} 
We begin this subsection by proving some results concerning abelian
locally $L$-analytic groups.

\proclaim{Proposition 6.4.1} If $Z$ is an abelian locally
$L$-analytic group, then the following conditions on $Z$
are equivalent:

(i) $Z$ is topologically finitely generated.

(ii) For every compact open subgroup $Z_0$ of $Z$,
the quotient $Z/Z_0$ is a finitely generated abelian
group.

(iii) $Z$ contains a compact open subgroup $Z_0$ such
that the quotient $Z/Z_0$ is a finitely generated abelian
group.

(iv) $Z$ contains a (necessarily unique) maximal compact
open subgroup $Z_0$, and the quotient $Z/Z_0$ is a free
abelian group of finite rank.

(v) $Z$ 
is isomorphic to a product $Z_0 \times \Lambda$, where $Z_0$
is a compact locally $L$-analytic group
and $\Lambda$ is a free abelian group of finite rank.
\endproclaim
\demo{Proof}
If $Z$ is topologically finitely generated, and if $Z_0$
is a compact open subgroup of $Z$, then the quotient
$Z/Z_0$ is a topologically finitely generated discrete
group, and hence is a finitely generated abelian group.
Thus~(i) implies~(ii), while~(ii) obviously implies~(iii).

Suppose now that $Z$ satisfies~(iii),
and hence contains a compact open subgroup $Z_0'$ such that
the quotient $Z/Z_0'$ is a finitely generated abelian group.
Let $Z_0$ be the preimage in $Z$ of the torsion subgroup of $Z/Z_0'$.
Since this subgroup is finite, we see that $Z_0$ is a compact
open subgroup of $Z$.  Since any compact open subgroup of $Z$
has finite, and hence torsion, image in the discrete subgroup
$Z/Z_0,$ we see that $Z_0$ is furthermore 
the maximal compact open subgroup of $Z$.  Since $Z/Z_0$ is
isomorphic to the quotient of $Z/Z_0'$ by its torsion subgroup,
it is a free abelian group of finite rank.  Thus~(iii) implies~(iv).
(Note that if $Z_0'$ and $Z_0'$ are two compact open subgroups
of $Z$, then the product $Z_0'Z_0''$ is again a compact open
subgroup of $Z$, containing both $Z_0'$ and $Z_0''$.
Thus if $Z$ contains a maximal compact open subgroup, it is
necessarily unique.)

That~(iv) implies~(v) is obvious.
Finally, suppose that $Z$ satisfies~(v).  Since $Z_0$ contains
a finite index subgroup isomorphic to $\Cal O_L^d$ for some
natural number $d$, we see that $Z_0$ is
topologically finitely generated.  The same is then
true of $Z \iso Z_0 \times \Lambda$, and so $Z$ satisfies~(i). 
This completes the proof of the proposition.
\qed\enddemo

\proclaim{Definition 6.4.2}
If $Z$ is an abelian locally $L$-analytic group,
and $X$ is a rigid analytic space over $L$,
then we let $\hat{Z}(X)$ denote the group
of characters on $Z$ with values
in $\An(X,L)^{\times}$, with the property that
for any admissible open affinoid subspace $X_0$ of $X$,
the induced character $\hat{Z}(X) \rightarrow \An(X,L)^{\times}
\rightarrow \An(X_0,L)^{\times}$ (the second arrow
being induced by restricting functions from $X_0$ to $X$)
is locally analytic when regarded as a $\An(X_0,L)$-valued function
on $Z$.
\endproclaim

\proclaim{Lemma 6.4.3} Let $X$ be a rigid analytic space over $L$,
and let $\{X_i\}_{i \in I}$ be an admissible affinoid open cover of $X$.
If $\phi: Z \rightarrow \An(X,L)^{\times}$ is a character
such that each of the composites
$Z \buildrel \phi \over \longrightarrow \An(X,L)^{\times} \longrightarrow
\An(X_i, L)^{\times}$ (the second arrow being induced by restriction to
$X_i$) is a locally analytic $\An(X_i,L)$-valued function on $Z$,
then $\phi$ lies in $\hat{Z}(X)$.
\endproclaim
\demo{Proof}
Let $X_0$ be an admissible affinoid open in $X$, and
write $X_{0,i} := X_0 \bigcap X_i,$ for each $i\in I$.
These sets form an admissible open cover of $X_0$. Since $X_0$
is affinoid, it is quasi-compact in its Grothendieck topology,
and so we find a finite subset $I'$ of $I$ such that
$\{X_{0,i}\}_{i\in I'}$ covers $X_0$.  
The open immersions $X_{0,i} \rightarrow X_0$ for $i\in I'$
induce a surjection $\coprod_{i\in I} X_{0,i} \rightarrow X_0.$
This in turn induces a closed embedding
$\An(X_0,L) \rightarrow \bigoplus_{i\in I} \An(X_{0,i}).$
Thus to prove that the composite
$Z \rightarrow \An(X,L)^{\times} \rightarrow \An(X_0,L)^{\times}$
is locally analytic, it suffices (by proposition~2.1.23)
to show that each of the
composites 
$Z \rightarrow \An(X,L)^{\times} \rightarrow \An(X_{0,i},L)^{\times}$
is locally analytic.   But each such composite factors through the
composite
$Z \rightarrow \An(X,L)^{\times} \rightarrow
\An(X_i,L)^{\times},$ which is locally analytic by assumption.
This proves the lemma.
\qed\enddemo

The preceding lemma implies that the formation of $\hat{Z}(X)$ is local in
the Grothendieck topology on the rigid analytic space $X$.
We observe in the following corollary that it also implies
that $\hat{Z}$ is a functor on the category of rigid analytic spaces
(and so consequently is a sheaf, when this category is made into a site
by equipping each rigid analytic space with its Grothendieck topology).

\proclaim{Corollary 6.4.4} The formation of $\hat{Z}(X)$ is contravariantly
functorial in the rigid $L$-analytic space $X$.
\endproclaim
\demo{Proof}  Let $f: X \rightarrow Y$ be a rigid analytic morphism
between rigid $L$-analytic spaces, and suppose that $\phi:Z \rightarrow 
\An(Y,L)^{\times}$ is an element of $\hat{Z}(Y)$.  Composing $\phi$
with pullback morphism induced by $f$ we obtain a character
$\phi': Z \rightarrow \An(X,L)^{\times}$.  We must show that
this character lies in $\hat{Z}(X)$.  Fix an admissible open cover
$\{Y_i\}_{i\in I}$ of $Y$, and choose an admissible open cover
$\{X_j\}_{j \in J}$ of $X$ that refines the open cover
$\{f^{-1}(Y_i)\}_{i \in I}$.
Then if $X_j$ lies in $f^{-1}(Y_i)$, the map
$f$ induces a continuous map of Banach spaces
$\An(Y_i,L) \rightarrow \An(X_j,L)$, and consequently
the composite
$Z \buildrel \phi' \over \longrightarrow \An(X,L)^{\times}
\longrightarrow \An(X_j,L)^{\times}$
(which is equal to the composite
$Z \buildrel \phi \over \longrightarrow \An(Y,L)^{\times}
\longrightarrow \An(Y_i,L)^{\times} \longrightarrow \An(X_j,L)^{\times}$)
is a locally analytic $\An(X_j,L)$-valued function on $Z$.
Lemma~6.4.3 now implies that $\phi'$ lies in $\hat{Z}(X)$.
\qed\enddemo

In the cases of interest to us, the functor $\hat{Z}$ is representable,
as we now show.

\proclaim{Proposition 6.4.5} If $Z$ is an abelian
locally $L$-analytic group
satisfying the equivalent conditions of proposition~6.4.1,
then the functor $\hat{Z}$ of definition~6.4.2
is representable
by a strictly $\sigma$-affinoid rigid analytic space over $L$.
\endproclaim
\demo{Proof}
Proposition~6.4.1~(v) allows us to write $Z \iso Z_0 \times
\Lambda,$  with $Z_0$ compact and $\Lambda$ a free abelian
group of finite rank.
Thus there is an isomorphism of functors
$\hat{Z} \iso \hat{Z}_0 \times \hat{\Lambda}.$
The functor $\hat{\Lambda}$ is clearly representable
by the $\sigma$-affinoid rigid analytic space $\check{\Lambda}
\otimes_{\Z} \G_m$ (where $\check{\Lambda}$ denotes
the dual of the finitely generated free abelian group $\Lambda$,
and $\G_m$ denotes the multiplicative group
over $L$, regarded as a rigid analytic space); the point is
that since $\Lambda$ is discrete, the condition of locally analyticity
in the definition of $\hat{\Lambda}$ is superfluous.
Thus it remains
to show that the $\hat{Z}_0$ is representable by a strictly $\sigma$-affinoid
rigid analytic space.

Since $Z_0$ is compact and locally $L$-analytic, it contains
a finite index subgroup isomorphic to $\Cal O_L^d$.  
In \cite{\SCHTPF} it is shown that $\hat{\Cal O}_L$ is
representable on affinoid spaces
by a twisted form of the open unit disk, which
is in particular strictly $\sigma$-affinoid. 
Lemma~6.4.3 shows that this space in fact represents the functor
$\hat{\Cal O}_L$ on the category of
all rigid analytic spaces.
Restricting characters from $Z_0$ to $\Cal O_L^d$ then
realises $\hat{Z}_0$ as a finite faithfully flat cover
of $\hat{\Cal O}_L^d$, which must also be strictly $\sigma$-affinoid.
\qed\enddemo

If $Z$ satisfies the equivalent criterion of proposition~6.4.1,
then we let $\hat{Z}$ denote the rigid
analytic space over $L$ constructed in proposition~6.4.5, that
parameterises the rigid analytic characters of $L$.
(Yoneda's lemma assures us that $\hat{Z}$ is determined up to natural 
isomorphism.)
Since $\hat{Z}$ is strictly $\sigma$-affinoid, the algebra
of rigid analytic functions $\An(\hat{Z},K)$ is a nuclear Fr\'echet
algebra.

\proclaim{Proposition~6.4.6} If $Z$ is an abelian
locally $L$-analytic group satisfying the equivalent conditions
of proposition~6.4.1, then there is a natural continuous injection
of topological $K$-algebras
$\DLa(Z,K)_b \rightarrow \An(\hat{Z},K).$
This map has dense image, and if $Z$ is compact, it is even
an isomorphism.
\endproclaim
\demo{Proof}
If $\mu \in \DLa(Z,K)_b$, then integrating against  
$\mu$ induces a map
$\hat{Z}(K) \rightarrow K$.  We will show that this
is obtained by evaluating an element of $\An(\hat{Z},K)$,
and that the resulting map $\DLa(Z,K)_b \rightarrow \An(\hat{Z},K)$
is continuous, with dense image.
If we write $Z = Z_0 \times \Lambda$,
as in proposition~6.4.1~(v), then
$\DLa(Z,K)_b \iso \DLa(Z_0,K)_b\otimes_K K[\Lambda]$,
while $\An(\hat{Z},K) \iso \An(\hat{Z}_0,K)_b\cotimes_K
\An(\hat{\Lambda},K)_b,$
and it suffices to prove the proposition with $Z$ replaced by $Z_0$
and $\Lambda$ in turn.

Since $Z_0$ is compact, it
contains a finite index subgroup isomorphic to $\Cal O_L^d$,
and so our claim is easily reduced to the case where
$Z_0 = \Cal O_L^d$.  In this case,
it follows from \cite{\SCHTPF, thm.~2.3}
that integrating characters against distributions yields a
topological isomorphism
$\DLa(Z_0,K)_b \iso \An(\hat{Z}_0,K).$

Choosing a basis 
$z_1,\ldots, z_r$ for $\Lambda$ yields an isomorphism $\hat{Z} \iso 
\G_m^r$.
Thus $\DLa(\Lambda,K) \iso K[\Lambda] = K[z_1,z_1^{-1}, \ldots, z_r,z_r^{-1}]$,
the target being equipped with its
finest complex topology (the variable $z_i$ in the target
of this isomorphism corresponds to the delta function
$\delta_{z_i}$ in its source),
while
$\An(\hat{Z},K) \iso \An(\G_m^r,K) \iso K\{\{ z_1, z_1^{-1},
\ldots, z_r,z_r^{-1}\}\}$
(the variable $z_i$ in the target of this isomorphism
corresponds to the rigid analytic function that takes
a character on $\hat{Z}$ to its value at $z_i$; that is,
to the delta function $\delta_{z_i}$).
Thus we see that integrating characters against distributions
yields a map $\DLa(\Lambda,K) \rightarrow \An(\hat{\Lambda},K)$,
which corresponds, with respect to the preceding isomorphisms,
to the usual inclusion
$K[z_1,z_1^{-1}, \ldots, z_r,z_r^{-1}] \subset
K\{\{ z_1, z_1^{-1},\ldots,z_r,z_r^{-1}\}.$
In particular, this map has dense image. 
\qed\enddemo

The proof of the preceding proposition shows that if $Z$ is not
compact, then the map $\DLa(Z,K)_b \rightarrow \An(\hat{Z},K)$
is not an isomorphism.  Rather,
$\An(\hat{Z},K)$ is a completion of $\DLa(Z,K)_b$
with respect to a certain continuous (and locally convex)
metric structure on $\DLa(Z,K)_b$.

If $V$ is a convex $K$-vector space of compact type, equipped
with a locally analytic action of $Z$, then proposition~5.1.9~(ii)
shows that $V'$ is a module over the ring $\DLa(Z,K)_b$, and that
the multiplication map $\DLa(Z,K)_b \times V'_b \rightarrow V'_b$
is separately continuous.

\proclaim{Proposition 6.4.7}
Let $V$ be a convex $K$-vector space of compact type,
equipped with a topological action of the topologically finitely
generated abelian locally $L$-analytic group $Z$.
The following conditions are equivalent:

(i) The $Z$-action on $V$ extends to a $\An(\hat{Z},K)$-module
structure for which the multiplication map $\An(\hat{Z},K)\times V
\rightarrow V$ is separately continuous.  (Such an extension is
unique, if it exists.)

(ii) The contragredient $Z$-action on $V'_b$ extends to
a topological $\An(\hat{Z},K)$-module structure on $V'_b$.
(Such an extension is unique, if it exists.)

(iii) The $Z$-action on $V$ is locally analytic,
and may write $V$
as a union of an increasing sequence of $BH$-subspaces, each
invariant under $Z$.
\endproclaim
\demo{Proof}
If we write $Z = Z_0 \times \Lambda,$ as in
proposition~6.4.1~(v), then $\An(\hat{Z},K)
\iso \An(\hat{Z}_0,K) \cotimes_K \An(\hat{\Lambda},K)$,
and it suffices to prove the above proposition separately
for the two cases $Z = Z_0$ and $Z = \Lambda$.

If $Z = Z_0,$ then proposition~6.4.6 yields
an isomorphism of nuclear Fr\'echet algebras
$\DLa(Z_0,K)_b \iso \An(\hat{Z}_0,K)$.
The equivalence of conditions~(i) and~(ii) follows from
the equivalence of~(i) and~(ii) of proposition~1.2.14.
These conditions are also equivalent to the $Z_0$-action
on $V$ being locally analytic, as follows from corollary~5.1.9
and proposition~5.1.10.  Finally, since $Z_0$ is compact,
any compact type space $V$ equipped with a locally analytic
$Z_0$-representation is the inductive limit of a sequence
of $Z_0$-invariant $BH$-subspaces (by proposition~3.2.15).

Now consider the case $Z = \Lambda$.  Let $z_1,\ldots,z_r$
be a basis for $\Lambda$.  Then
$$\multline \An(\hat{\Lambda},K) \iso K\{\{z_1,z_1^{-1},\ldots, z_r,z_r^{-1}\}\}
\\ \iso \plim{n}
K\langle\langle
p^n z_1, p^n z^{-1}, \ldots, p^n z_r, p^n z_r^{-1}\rangle \rangle.\endmultline$$
(The second isomorphism arises from writing the rigid analytic
space $\check{\Lambda} \otimes_{\Z} \G_m \iso \G_m^r$ as the union
of the affinoid subdomains $|p|^n \leq |z_i| \leq |p|^{-n}$
($i = 1,\ldots, r$) as $n$ tends to infinity.  For each $n \geq ,$
we have written
$K\langle\langle
p^n z_1, p^n z^{-1}, \ldots, p^n z_r, p^n z_r^{-1}\rangle \rangle$
to denote the Tate algebra of rigid analytic functions on the corresponding
subdomain.)
It follows from proposition~1.2.14 that conditions~(i) and~(ii) of
the proposition are equivalent, and that each in turn is equivalent
to the following condition:

(iii') There is an isomorphism $V \iso \ilim{n} V_n$, where
each $V_n$ is a Banach module over the Tate algebra
$K\langle\langle
p^n z_1, p^n z^{-1}, \ldots, p^n z_r, p^n z_r^{-1}\rangle \rangle,$
and the maps $V_n \rightarrow V$ are $Z$-equivariant.
(Note that $Z$ embeds as a group of units in the Tate algebra
$K\langle\langle
p^n z_1, p^n z^{-1}, \ldots, p^n z_r, p^n z_r^{-1}\rangle \rangle.$)
 
Clearly~(iii') implies~(iii).  We must show that the converse holds.
For this, it suffices to show that if $W$ is a Banach space
equipped with a topological action of $\Lambda$, then for some sufficiently
large value of $n$, the Banach space $W$ admits a (necessarily unique)
structure of topological 
$K\langle\langle
p^n z_1, p^n z^{-1}, \ldots, p^n z_r, p^n z_r^{-1}\rangle \rangle$-module
extending its $\Lambda$-module structure.    If we fix a norm on $W$,
then we may certainly find some $n$ such that
$|| p^n z_i^{\pm 1} || \leq 1$ for $i = 1,\ldots, r$.  (Here we are 
considering the norm of $p^n z_i^{\pm 1}$ as an operator on the normed
space $W$.) It is then clear
that $W$ does indeed admit the required structure of topological
$K\langle\langle
p^n z_1, p^n z^{-1}, \ldots, p^n z_r, p^n z_r^{-1}\rangle \rangle$-module.
\qed\enddemo

\proclaim{Lemma~6.4.8}
If $V$ is a convex $K$-vector space of compact type,
equipped with a locally analytic action of $Z$
that satisfies the equivalent conditions of proposition~6.4.7,
then the same is true of any $Z$-invariant closed subspace of $V$.
\endproclaim
\demo{Proof}
It is clear that if $V$ satisfies condition~(i) of proposition~6.4.7,
then the separately continuous $\An(\hat{Z},K)$-module structure
on $V$ restricts to a separately continuous $\An(\hat{Z},K)$-module
structure on any $Z$-invariant closed subspace of $W$
(since $K[Z]$ is dense in $\An(\hat{Z},K)$).
\qed\enddemo

We now suppose that $G$ is a locally $L$-analytic group,
whose centre $Z$ satisfies the equivalent conditions of proposition~6.4.1.

If $V$ is a
convex $K$-vector space of compact type,
equipped with a locally analytic $G$-representation, 
then by proposition~5.1.9~(ii), $V'_b$ is equipped with a
topological $\DLa(G,K)_b$-module structure.  If $V$ furthermore satisfies the
hypothesis of proposition~6.4.7, then $V'_b$ is also equipped
with a topological $\An(\hat{Z},K)$-module structure. 

Since $Z$ is central in $G$, the $\DLa(G,K)_b$-action and
$\An(\hat{Z},K)$-action
on $V'_b$ commute with one another, and so together they
induce the structure of a topological
$\An(\hat{Z},K) \cotimes \DLa(G,K)_b$-module
on $V'_b$.
If $H$ is any compact open subgroup
of $G$, then we may restrict this to a
$\An(\hat{Z},K)\cotimes_K \DLa(H,K)_b$-module structure on $V'_b$.
By proposition~5.3.22 (and the remark following definition~5.3.21,
which shows that $\An(\hat{Z},K)$ admits an integral Fr\'echet-Stein structure),
this completed tensor product is a Fr\'echet-Stein
algebra.

\proclaim{Definition 6.4.9}  Let $V$ be a convex $K$-vector space
of compact type, equipped with a locally analytic action of $G$.
We say that $V$ is an essentially admissible locally analytic representation
of $G$ if $V$ satisfies the conditions of proposition~6.4.7,
and if for one (equivalently, every) compact open subgroup $H$ of
$G$, the dual $V'_b$ is a coadmissible module when endowed with its natural
module structure over the Fr\'echet-Stein
algebra $\An(\hat{Z},K) \cotimes_K \DLa(H,K)_b$.
\endproclaim

Note that in the situation of definition~6.4.9,
we have embeddings of $\DLa(Z,K)$ into each factor
of the tensor product $\An(\hat{Z},K)\cotimes_K \DLa(H,G),$
and the action of $\DLa(Z,K)$ on $V'$ obtained by regarding 
$\DLa(Z,K)$ as a subalgebra of either of these factors coincides.
(By construction, either of these actions is the natural $\DLa(Z,K)$-action
on $V'$ induced by the locally analytic $Z$-representation on $V$.)
Thus the $\An(\hat{Z},K) \cotimes_K \DLa(G,K)$-action on
$G$ factors through the quotient algebra $\An(\hat{Z},K)\cotimes_{\DLa(Z,K)}
\DLa(G,K)$.  

Similarly, 
if $H$ is a compact open subgroup of $G$, and if we write $Z_0 = H \cap Z$,
then the $\An(\hat{Z},K) \cotimes_K \DLa(G,K)$-action on $V'$
factors through the quotient algebra
$\An(\hat{Z},K)\cotimes_{\DLa(Z_0,K)} \DLa(H,K)$.  Thus $V$ is essentially
admissible if and only if it is coadmissible for the Fr\'echet-Stein algebra
$\An(\hat{Z},K)\cotimes_{\DLa(Z_0,K)} \DLa(H,K)$.  In particular,
if the centre $Z$ of $G$ is compact, then (taking into account
proposition~6.4.6) $V$ is essentially admissible
if and only if it is admissible.

In general, we have the following proposition.

\proclaim{Proposition 6.4.10} Any admissible locally analytic
representation of $G$
is an essentially admissible locally analytic representation of $G$.
\endproclaim
\demo{Proof}
Let $\{H_n\}_{n\geq 1}$ be a cofinal sequence of analytic open
subgroups of $G$.   If $V$ is an admissible locally analytic 
$G$-representation, then let $V_n$ denote the $BH$-subspace of $V$
obtained as the image of the natural map $V_{\H_n-\an} \rightarrow V$.
Clearly $V = \bigcup_{n=1}^{\infty} V_n$, and each $V_n$ is invariant
under $Z$, since multiplication by elements of $Z$ commutes with
the $H_n$ action on $V$.  Thus $V$ satisfies the condition~(i) of
proposition~6.4.7, and so
is naturally a $\An(\hat{Z},K) \cotimes_K \DLa(H,K)$-module.
It is certainly coadmissible with respect to this Fr\'echet-Stein
algebra, since it is coadmissible even with respect to $\DLa(H,K)$.
\qed\enddemo

If $Z$ is not compact, then the converse to proposition~6.4.10
is false.  For example,
if $G = Z$ is a free abelian group of rank $r$, then a locally
analytic $G$-representation
is admissible if and only if it is finite dimensional.  On the other
hand, the essentially admissible locally analytic $G$-representations
correspond (by passing to the dual) to coherent sheaves on the
rigid analytic space $\G_m^r$.  (Among all such sheaves, the admissible
representations correspond to those skyscraper sheaves whose support
is finite.)

More generally, if $G = Z$ is abelian, then the essentially
admissible locally analytic $Z$-representations correspond (by
passing to the dual) to coherent sheaves on the rigid analytic
space $\hat{Z}$.  Such a sheaf corresponds to an admissible locally
analytic $Z$-representation if and only if its pushforward to
$\hat{Z}_0$ (where $Z_0$ denotes the maximal compact subgroup of $Z$)
is again coherent.

\proclaim{Proposition 6.4.11} If $G$ is a locally $L$-analytic group,
then the category of essentially admissible locally analytic $G$-representations
and continuous $G$-equivariant morphisms is closed under passing to 
closed subobjects and Hausdorff quotients.  Furthermore, any morphism
in this category is necessarily strict.   Consequently, this
category is abelian.
\endproclaim
\demo{Proof}
This is a consequence of the general properties of coadmissible modules
over Fr\'echet-Stein algebras, summarised in theorem~1.2.11 and the
remarks that follow it.
\qed\enddemo

\proclaim{Definition 6.4.12} If $V$ is a Hausdorff convex $K$-vector space
satisfying the equivalent conditions of proposition~6.4.7,
and if $J$ denotes an ideal in the ring $\An(\hat{Z},K),$
then we let $V^J$ denote the closed $G$-invariant subspace of
$V$ that is annihilated by the ideal $J$.
\endproclaim

\proclaim{Proposition 6.4.13} If $V$ is an essentially admissible
locally analytic representation of $G$, and if $J$ is a finite-codimension
ideal in $\An(\hat{Z},K)$,
then $V^J$ is an admissible locally analytic
representation of $G$.
\endproclaim
\demo{Proof}
If $\overline{J}$ denotes the closure of $J$ in 
$\An(\hat{Z},K)$, then $V^J = V^{\overline{J}},$
and so, replacing $J$ by $\overline{J},$
we may assume that $J$ is closed.

There is a natural isomorphism
$(V^J)'_b \iso (\An(\hat{Z},K)/J) \cotimes_{\An(\hat{Z},K)} V'_b.$
If we fix a compact open subgroup $H$ of $G$ then
we see that 
since $V'_b$ is coadmissible with respect to
$\An(\hat{Z}, K)\cotimes_K \DLa(H,K),$
its quotient $(V^J)'_b$ is coadmissible with respect to
$(\An(\hat{Z},K)/J) \cotimes_K \DLa(H,K)$.
Since $\An(\hat{Z},K)/J$ is finite dimensional, it is
in fact coadmissible with respect to $\DLa(H,K)$.
Thus $V^{J}$ is an admissible
locally analytic representation of $G$.
\qed\enddemo

If $\chi$ is an element of $\hat{Z}(K)$, then we let $J_{\chi}$ denote
the maximal ideal of $\An(\hat{Z},K)$ consisting of functions that vanish at
$\chi$,
and we write $V^{\chi}$ in place of $V^{J_{\chi}}$.
Note that $V^{\chi}$ is the closed $G$-invariant
subspace of $V$
consisting of vectors that transform via $\chi$ under the action of $Z$.

\proclaim{Corollary 6.4.14} If $V$ is an essentially admissible
locally analytic representation of $G$, and if $\chi$ is an element
of $\hat{Z}(K)$,
then $V^{\chi}$ is an admissible locally analytic
representation of $G$.
\endproclaim
\demo{Proof}
Apply the previous proposition to the ideal $J_{\chi}$.
\qed\enddemo

We close this subsection with the following result.
In conjunction with proposition~6.4.13, it allows
one in certain situations to reduce questions regarding morphisms from
admissible locally algebraic representations into essentially admissible
locally analytic representations to the case in which
the target is assumed admissible.

\proclaim{Proposition 6.4.15}  If $V$ is an admissible locally
algebraic representation of $G$, and if we let $J$ run through the
directed system of ideals of finite codimension in $\An(\hat{Z},K)$,
then the natural map
$\ilim{J} V^J \rightarrow V$ is an isomorphism.
\endproclaim
\demo{Proof}
If we use proposition~4.2.4 and corollary~4.2.7 to write
$V \iso \bigoplus_n U_n \otimes_{B_n} W_n,$
where $W_n$ runs over the elements of $\hat{\G}$, $B_n = \End_{\G}(W_n),$
and $U_n$ is a smooth representation of $G$ over $B_n$, then we
see that it suffices to prove the proposition for each
of the direct summands $U_n \otimes_{B_n} W_n$.  Since $\An(\hat{Z},K)$
acts on $W_n$ through a finite dimensional quotient, we further reduce
ourselves to proving the proposition for an admissible smooth
representation $U_n$.   We may write $U_n \iso \ilim{H} U_n^H$,
where $H$ runs over the directed set of compact open subgroups of $G$.
Each space $U_n^H$ is finite dimensional and $Z$-invariant.   Thus
$\An(\hat{Z},K)$ acts on it through a finite dimensional quotient,
and the proposition is proved.
\qed\enddemo

\section{6.5}
We conclude this section with a discussion of invariant lattices
in convex $K$-vector spaces equipped with a $G$-action.  We suppose
throughout this subsection that $K$ is discretely valued.

Recall that if $V$ is a $K$-vector space,
then a lattice in $V$ is an $\Cal O_K$-submodule of $V$ that spans
$V$ over $K$ \cite{\SCHNA, def., p.~7}.
We say that two lattices $M_1$ and $M_2$ in $V$
are commensurable if there exists an element $\alpha \in K^{\times}$
such that $\alpha M_2 \subset M_1 \subset \alpha^{-1} M_2.$
Commensurability is obviously an equivalence relation.

We will also consider seminorms on $V$.  (All seminorms and norms on $V$
are understood to be non-archimedean.)  We say that two seminorms
$q_1$ and $q_2$ are commensurable if there exists a positive
real number $c$ such that $c q_1(v) \leq q_2(v) \leq c^{-1} q_1(v)$
for all $v \in V$.

There is a natural bijection between the set
of commensurability classes of lattices in $V$ and commensurability
classes of seminorms on $V$, obtained by passing from a lattice to its
corresponding gauge, and from a seminorm to its corresponding unit ball.
(See \cite{\SCHNA, lem.~2.2}.)  

Let $\pi$ be a uniformiser of $\Cal O_K$.  (Recall that we are assuming
that $K$ is discretely valued.)   We say that an $\Cal O_K$-module $M$
(and so in particular, a lattice $M$ in $V$) is separated if
$\bigcap_{n \geq 1} \pi^n M = 0$.  The property of being separated is
an invariant of a commensurability class of lattices.
Also, the property of being a norm is an invariant of a commensurability
class of seminorms.  Under the above bijection, the commensurability
classes of separated lattices correspond to the commensurability classes of
norms.

If $M$ is a lattice in $V$, then we let $V_M$ denote $V$, endowed with
the locally convex topology induced by the gauge of $M$,
and let $\hat{V}_M$ denote the Hausdorff
completion of $V_M$.  Note that $\hat{V}_M$ is naturally a
$K$-Banach space, and depends up to canonical isomorphism only on
the commensurability class of $M$.  There is a natural $K$-linear
map $V \rightarrow \hat{V}_M$, which is an injection if and only if
$M$ is separated (so that its gauge is a norm).

If $G$ acts on $V$, and if a lattice in $V$ 
is $G$-invariant, then the same
is true of its gauge, while if a seminorm is $G$-invariant,
then the same is true of its unit ball.  Thus a commensurability class
of lattices in $V$ contains a $G$-invariant lattice if and only
if the corresponding commensurability class of seminorms contains
a $G$-invariant seminorm.

\proclaim{Lemma 6.5.1} If $V$ is a $K$-vector space equipped
with a $G$-action, and $M$ is a lattice in $V$, then
$G$ preserves the commensurability class of $M$ (that is,
the lattice $g M$ is commensurable with $M$, for each $g \in G$)
if and only if the $G$-action on $V_M$ is topological.
\endproclaim
\demo{Proof}
Given that the topology on $V_M$ is defined via the gauge of $M$,
it is clear that if $G$ preservers the commensurability class
of $M$, then each element of $G$ acts as a continuous automorphism
of $V_M$.
Conversely, if each element of $G$ acts as a continuous automorphism
of $V_M$, then for any $g \in G$, there is an inclusion
$g M \subset \alpha M$ for some $\alpha \in K^{\times}$.
Combining this inclusion with the analogous inclusion for $g^{-1}$
shows that $g M$ is commensurable with $M$, for each $g \in G$.
\qed\enddemo

\proclaim{Lemma 6.5.2} If $V$ is a $K$-vector space equipped
with a $G$-action, and $M$ is a lattice in $V$, then the following
are equivalent:

(i) There exists a $G$-invariant lattice that is commensurable
with $M$.

(ii) There exists a $G$-invariant seminorm that defines the
topology of $V_M$.

(iii) $G$ acts as a group of equicontinuous operators
on $V_M$.
\endproclaim
\demo{Proof}  The gauge of a $G$-invariant lattice is
a $G$-invariant seminorm, and so~(i) implies~(ii), while~(ii)
clearly implies~(iii).
If~(iii) holds, then there is $\alpha \in K^{\times}$ such
that $g M \subset \alpha M$ for all $g \in G$.  If we let
$M'$ be the lattice generated by the lattices $g M$ for all
$g \in G$, then $M \subset M' \subset \alpha M$, and so $M'$
is commensurable with $M$.  Since $M'$ is $G$-invariant
by construction, we see that~(iii) implies~(i).
\qed\enddemo

If $V$ is furthermore a topological $K$-vector space,
then the property of being an open subset of $V$ is
an invariant of a commensurability class of lattices, while
the property of being continuous is an invariant of a commensurability
class of seminorms.  Under the above bijection, the commensurability
classes of open lattices correspond to the commensurability classes
of continuous seminorms.  The open lattices are thus precisely
those for which the natural map $V \rightarrow V_M$ is continuous.

\proclaim{Lemma 6.5.3} Suppose that $G$ is compact.
If $V$ is a topological $K$-vector space
equipped with a continuous action of $G$, and if $M$ is an open
lattice in $V$, then $G$ acts continuously on $V_M$ if and only if
there is a $G$-invariant lattice that is commensurable with $M$.
\endproclaim
\demo{Proof}
If $G$ acts continuously on $V_M$, then since $G$ is compact,
lemma~3.1.4 shows that
it acts as an equicontinuous set of operators.  From lemma~6.5.2,
we conclude that there exists a $G$-invariant lattice
commensurable with $M$.

Conversely, if there exists a $G$-invariant lattice that
is commensurable with $M$, then $G$ acts as a group of equicontinuous
operators on $V_M$.  Thus conditions~(ii) and~(iii) of lemma~3.1.1
hold with respect to the $G$-action on $V_M$.  Since $M$ is open in $V$,
the natural map $V \rightarrow V_M$ is continuous.  Since the
$G$-action on $V$ is continuous (and since $V$ and $V_M$ are equal
as abstract $K$-vector spaces), we see that condition~(i) of lemma~3.1.1
is also satisfied, and so the $G$-action on $V_M$ is continuous.
\qed\enddemo

\proclaim{Lemma 6.5.4} If $V$ is a topological $K$-vector space on which
$G$ acts as an equicontinuous group of operators (for example,
if $G$ is compact and acts continuously on $V$),
then any open lattice in $V$ contains a $G$-invariant open
sublattice, and there is a unique maximal such $G$-invariant open sublattice.
\endproclaim
\demo{Proof}
If $M$ is a given open lattice in $V$, then by assumption
there exists an open sublattice $M'$ of $M$ such $G M' \subset M.$
The sublattice of $M$ spanned by $G M'$ is then a $G$-invariant
open sublattice of $M$.  If we consider the sublattice of $M$
spanned by all invariant open sublattices, then it is obviously
the unique maximal invariant open sublattice.
\qed\enddemo

\proclaim{Lemma 6.5.5} If $V$ is a $K$-Banach space on which
$G$ acts as an equicontinuous group of operators (for example,
if $G$ is compact and acts continuously on $V$),
then any separated open lattice in $V$ is contained in a $G$-invariant
separated open lattice of $V$,
and there is a unique minimal $G$-invariant separated open lattice of $V$
containing it.
\endproclaim
\demo{Proof}
Let $M$ denote the given separated open sublattice of $V$.
We may choose the norm on the Banach space $V$
so that $M$ is the unit ball of $V$.
Since $G$ acts equicontinuously on $V$, there is an $\alpha \in K^{\times}$
such that $G \alpha M \subset M.$  Then $G M \subset \alpha^{-1} M,$
and so we see that the lattice spanned by $G M$ is a $G$-invariant
separated open lattice in $V$ that contains $M$.  Taking the intersection
of all such lattices yields the minimal $G$-invariant separated open
lattice containing $M$.
\qed\enddemo

Let $k = \Cal O_K/\pi$ denote the residue field of $\Cal O_K$.

\proclaim{Lemma 6.5.6} If $V$ is a convex $K$-vector space equipped with
a continuous $G$-action, and $M$ is an open $G$-invariant lattice in $V$,
then the induced $G$-action on the $k$-vector space $M /\pi $
is smooth.
\endproclaim
\demo{Proof} Lemma~6.5.3 implies that $G$ acts continuously on
$V_M$, and thus that we may replace $V$ by $V_M$ in our considerations.
If we equip $M/\pi$ with the topology induced by regarding it
as a subquotient of $V_M,$ then the $G$-action on $M/\pi$ is
continuous.  On the hand, the topology obtained on $M/\pi$ in
this fashion is the discrete topology.  Thus the $G$-action on
$M/\pi$ is smooth, as claimed.
\qed\enddemo

The following result slightly extends \cite{\SCHTIW, lem.~3} (which treats
the case when $K$ is local).

\proclaim{Proposition 6.5.7} If $V$ is a convex $K$-vector space
equipped with a continuous $G$-action, and $M$ is an open
$G$-invariant lattice in $V$, then the following are equivalent:

(i) The Banach space $\hat{V}_M$, when equipped with its natural
$G$-action, becomes an admissible continuous representation of $G$.

(ii) The smooth representation of $G$ on $M/ \pi$ is admissible.
\endproclaim
\demo{Proof}
Replacing $V$ with $\hat{V}_M$, and $M$ with the closure of
its image in $\hat{V}_M$ (which does not disturb the isomorphism
class of the $G$-representation $M /\pi$), we may assume that
$V$ is a Banach space.  Also, replacing $G$ by an appropriate
pro-$p$ open subgroup if necessary, we may assume that $G$
is a pro-$p$ group.

As in theorem~6.2.8, let $\DCon(G,\Cal O_K)$ denote the unit
ball of $\DCon(G,K)$ (that is, the polar of the unit ball
of $\Con(G,\Cal O_K)$ with respect to the sup norm on this
latter space).  Corollary~5.1.7 shows that $V'_b$ is naturally
a topological $\DCon(G,K)$-module, and so also a $\DCon(G,\Cal O_K)$-module.
If we equip $V$ with the norm obtained as the gauge of $M$,
then the polar $M'$ of $M$ is a bounded open
$\Cal O_K$-sublattice of $(V')_b$, which is also a
$\DCon(G,\Cal O_K)$-submodule of $V'$.

Let $\goth{a}$ denote the ideal of $\DCon(G,\Cal O_K)$ generated
by the uniformiser $\pi$ of $\Cal O_K$, together with the
augmentation ideal $I_G$ of $\DCon(G,\Cal O_K)$.  In the case where
$\Cal O_K = \Z_p$, the $\goth{a}$-adic filtration of 
$\DCon(G,\Z_p)$ is cofinal with the filtration of
$\DCon(G,\Z_p)$ considered in \cite{\LAZ},
and hence $\DCon(G,\Z_p)$ is $\goth{a}$-adically complete.
Writing $\DCon(G,\Cal O_K) = \Cal O_K \cotimes_{\Z_p} \DCon(G,\Z_p)$,
we find that $\DCon(G,\Cal O_K)$ is $\goth{a}$-adically complete.

The quotient $M'/\goth{a}$ is isomorphic to the $\Cal O_K/\pi$-linear
dual of the space $(M/\pi)^G$.  Thus $M'/\goth{a}$ is finite dimensional
over $\Cal O_K/\pi$ if and only if $(M/\pi)^G$ is finite dimensional
over $\Cal O_K/\pi$.  On the other hand,
since $\DCon(G,\Cal O_K)$ is $\goth{a}$-adically complete, 
the quotient $M'/\goth{a}$ is finite dimensional over $\Cal O_K/\pi$
if and only if $M'$ is finitely generated over $\DCon(G,\Cal O_K)$.

Thus if $M/\pi$ is an admissible smooth representation of $G$,
then we see that $M'$ is finitely generated over $\DCon(G,\Cal O_K)$,
and hence that $V'$ is finitely generated over $\DCon(G,K)$; that is,
$V$ is an admissible continuous representation of $G$.
Conversely, if $V$ is an admissible continuous representation
of $G$, then replacing $G$ by an open subgroup $H$, it follows
that $V'$ is finitely generated over $\DCon(H,K)$.
Hence $M'$ is finitely generated over $\DCon(H,\Cal O_K)$,
and so $(M/\pi)^H$ is finite dimensional over $\Cal O_K/\pi$.
It follows that $M/\pi$ is an admissible smooth representation
of $G$.
\qed\enddemo

\proclaim{Definition 6.5.8} If $V$ is a convex $K$-vector space equipped with
a continuous $G$-action, we say that a $G$-invariant open lattice
in $V$ is admissible if it satisfies the equivalent conditions
of proposition~6.5.7.
\endproclaim

The following result provides something of a converse to
proposition~6.2.4.

\proclaim{Proposition 6.5.9} If $V$ is an essentially
admissible locally analytic representation of $G$, and
if $G$ is compact, then the following are equivalent:

(i) There is an admissible $G$-invariant separated open 
lattice in $V$.

(ii) There is a continuous $G$-equivariant injection $V \rightarrow W,$
with $W$ a $K$-Banach space equipped with an admissible continuous 
representation of $G$.

(iii) $V$ is a strongly admissible locally analytic $G$-representation.
\endproclaim
\demo{Proof}
Suppose that $V$ contains an admissible $G$-invariant separated open 
lattice $M$.  Proposition~6.5.7 shows that the $G$-action
on $\hat{V}_M$ makes this space an admissible continuous $G$-representation.
Since $M$ is separated and open, the natural map $V \rightarrow \hat{V}_M$
is continuous and injective, and it is certainly $G$-equivariant.
Thus~(i) implies~(ii).

Suppose that~(ii) holds, so that we are given a continuous
$G$-equivariant injection $V \rightarrow W$ with $W$ a $K$-Banach
space on which $G$ acts via an admissible continuous representation.
Lemma~6.5.3 shows that we may choose a $G$-invariant norm
on $W$.   Proposition~6.5.7 shows that
the unit ball of this norm will be an admissible $G$-invariant separated open 
lattice in $W$, and hence its
preimage in $V$ will be an admissible $G$-invariant separated open lattice
in $V$.  Thus~(ii) implies~(i).

If we continue to assume given a $G$-equivariant injection $V \rightarrow W$
as in~(ii), then passing to locally analytic vectors yields
a continuous injection $V \rightarrow W_{\la}$.
Proposition~6.2.4 shows that $W_{\la}$ is a strongly admissible
locally analytic representation of $G$.  Proposition~6.4.11 
shows that $V$ embeds as a closed subspace of $W_{\la}$, and
so $V_{\la}$ is itself strongly admissible.  Thus~(ii) also
implies~(iii).

Finally suppose that~(iii) holds, so that $V$ is strongly admissible.
The remark following definition~6.2.1 shows that we may find a
closed embedding $V \rightarrow \La(G,K)^n$ for some integer $n$.
Composing this map with the continuous injection $\La(G,K)^n \rightarrow
\Con(G,K)^n$ yields a continuous injection $V \rightarrow \Con(G,K)^n$.
Thus~(iii) implies~(ii).
\qed\enddemo

\head 7. Representations of certain product groups \endhead

\section{7.1}
Throughout this subsection, we assume that
$\Gamma$ is a locally compact topological group,
and that $\Gamma$ admits a countable neighbourhood
basis of the identity consisting of open subgroups.

If $V$ is a $K$-vector space equipped with a representation
of $\Gamma$, and if $H$ is a subgroup of $\Gamma$, then we
let $V^H$ denote the subspace of $V$ consisting of $H$-fixed
vectors.  As usual, we say that $V$ is a smooth representation
of $\Gamma$ if and only if each vector of $V$ is fixed by
some compact open subgroup of $V$; equivalently, $V$ is a smooth
representation if and only if the natural map
$\ilim{H}V^H \rightarrow V$ is an isomorphism (the inductive limit
being taken over all the compact open subgroups of $\Gamma$).

If $V$ is a smooth representation of $\Gamma$,
then we let $\pi_H: V \rightarrow V^H$ denote the operator
given by ``averaging over $H$'': if $v \in V$ is fixed by
the subgroup $H'$ of $H$,
then
$\pi_H(v) := [H: H']^{-1} \sum_{h \in H/H'} h v.$

We will consider a strengthening of the notion of a smooth representation,
when $V$ is assumed to be a Hausdorff locally convex topological $K$-vector
space, and the $\Gamma$-action on $V$ is topological.  In this
context, we let $V^H$ denote the space of $H$-fixed vectors,
equipped with the subspace topology.  Thus each $V^H$ is a closed
subspace of $V$.

\proclaim{Definition 7.1.1}
If $V$ is an arbitrary Hausdorff convex $K$-vector space equipped with
a topological $\Gamma$-action, then we define
$V_{\strict}$ to be the locally convex inductive limit
$$V_{\strict} := \ilim{H} V^H,$$
where $H$ runs over the compact open subgroups of $\Gamma$.
\endproclaim

The map $V \mapsto V_{\strict}$ is a covariant functor from 
the category of Hausdorff convex $K$-vector spaces equipped
with topological $\Gamma$-actions (with morphisms being 
continuous and $\Gamma$-equivariant) to itself.

Our assumption on $\Gamma$ shows that the inductive limit in
Definition~7.1.1 can be replaced by a countable inductive limit.
Also, the transition maps are closed embeddings.  Thus the vector
space $V_{\strict}$ is a strict inductive limit, in the sense of
\cite{\TVS, p.~II.33}.  In particular, it is Hausdorff, and is complete
if $V$ is (since each $V^H$ is then complete).
There is a natural continuous injection
$V_{\strict} \rightarrow V$, 
and the $\Gamma$-action on $V$ is smooth if and only if
this map is a bijection.

\proclaim{Definition 7.1.2} If $V$ is a Hausdorff locally convex $K$-vector
space equipped with a topological action of $\Gamma$, then we say
that $V$ is a strictly smooth representation of $\Gamma$ if
the natural map
$V_{\strict} \rightarrow V$ is a topological isomorphism.
\endproclaim

In particular,
a strictly smooth $\Gamma$-representation on a convex $K$-vector
space $V$ is certainly a topological smooth $\Gamma$-action.
In general, the converse need not be true.

\proclaim{Proposition 7.1.3}  Let $V$ be a Hausdorff locally convex $K$-vector
space equipped with a topological action of $\Gamma$.
If either $V$ is of compact type, or if for a cofinal
sequence of compact open subgroups $H$ of $\Gamma$, the closed subspaces
$V^H$ of $\Gamma$ are of compact type,
then $V_{\strict}$ is of compact type.
\endproclaim
\demo{Proof}
If $V$ is of compact type, then so is its closed subspace $V^H$,
for each compact open subgroup $H$ of $\Gamma$.  Thus the first hypothesis
implies the second.  Let $\{H_n\}_{n\geq 1}$ be a cofinal sequence
of compact open subgroups of $\Gamma$ for which $V^{H_n}$ is of compact type.
The isomorphism $\ilim{n} V^{H_n} \iso V_{\strict}$ shows that $V$
is the locally convex inductive limit of a sequence of spaces of
compact type, with injective transition maps,
and so is itself of compact type.
\qed\enddemo

\proclaim{Corollary 7.1.4}
Let $V$ be a Hausdorff locally convex $K$-vector space
equipped with a smooth topological action of $\Gamma$.  
If $V$ is of compact type,
then $V$ is a strictly smooth $\Gamma$-representation.
Conversely, if the $\Gamma$-action on $V$ is strictly smooth,
and if for a cofinal collection of compact open subgroups $H$ of $\Gamma$
the subspaces $V^H$ of $V$ are of compact type,
then $V$ is of compact type.
\endproclaim
\demo{Proof}
By assumption, the $\Gamma$-action on $V$ is smooth,
and hence the natural map
$$V_{\strict} \rightarrow V\tag 7.1.5$$
is a continuous bijection.   Also, either hypothesis on $V$ implies,
by lemma~7.1.3, that $V_{\strict}$ is of compact type.

If $V$ is of compact type, then~(7.1.5) is a continuous bijection
between spaces of compact type, and hence is a topological isomorphism.
Thus $V$ is a strictly smooth $\Gamma$-representation.  Conversely,
if the $\Gamma$-action on $V$ is strictly smooth,
then~(7.1.5) is a topological isomorphism, and so we conclude that
$V$ is of compact type.
\qed\enddemo

\proclaim{Proposition 7.1.6} If $V$ is a Hausdorff convex $K$-vector space
equipped with a strictly smooth action of $\Gamma$,
then for each compact open subgroup $H$ of $\Gamma$,
the operator $\pi_H: V \rightarrow V^H $ is continuous.
\endproclaim
\demo{Proof} If $H'$ is an open subgroup of $H$,
then the formula given above for the restriction of $\pi_H$ to
$V^{H'}$,  and the assumption that the $\Gamma$-action
on $V$ is topological, shows that $\pi_H$ induces
a continuous map $V^{H'} \rightarrow V^H$.
Passing to the inductive limit over $H'$, and using the
fact that the $\Gamma$-action on $V$ is strictly smooth,
we see that the proposition follows.
\qed\enddemo

\proclaim{Corollary 7.1.7}
If $V$ is a Hausdorff convex $K$-vector 
space equipped with a strictly smooth action of $\Gamma$,
then for each compact open subgroup $H$ of $\Gamma$,
we obtain a topological direct sum decomposition
$V^H \bigoplus \ker \pi_H \iso V.$
\endproclaim
\demo{Proof}
This follows from the fact that the operator
$\pi_H$ on $V$ is idempotent, with image equal to $V^H$,
and is continuous, by proposition~7.1.6.
\qed\enddemo

\proclaim{Corollary 7.1.8}
If $V$ is a Hausdorff convex $K$-vector 
space equipped with a strictly smooth action of $\Gamma$,
then $V$ is barrelled if and only if $V^H$ is barrelled for
each compact open subgroup $H$ of $\Gamma$.
\endproclaim
\demo{Proof}
If each $V^H$ is barrelled, then the isomorphism
$\ilim{H} V^H \iso V$ shows that $V$ is barrelled.
Conversely, if $V$ is barrelled, and if $H$ is a compact open
subgroup of $\Gamma$, then corollary~7.1.7
shows that $V^H$ is a topological direct summand of $V$.
It follows that $V^H$ is barrelled \cite{\TVS, cor.~2, p.~III.25}.
\qed\enddemo

We now explain how the countable inductive limit that implicitly appears in
Definition~7.1.1 can be replaced by a countable direct sum.
Let us fix a cofinal decreasing sequence $\{H_n\}_{n\geq 1}$ of 
compact open subgroups of the identity of $\Gamma$.  
If $V$ is a Hausdorff convex $K$-vector space,
equipped with a topological action of $\Gamma$,
then for each value of $n > 1$,
the averaging map $\pi_{H_{n-1}}$, when restricted to
$V^{H^n},$ induces a continuous map $\pi_n: V^{H_n} \rightarrow
V^{H_{n-1}}$.

Set $V_1 = V^{H_1},$ and $V_n = \ker \pi_n$ for $n > 1$.
Each of these spaces is a closed subspace of $V$.
For each $n \geq 1,$ the natural map $\bigoplus_{m=1}^n V_m \rightarrow
V^{H_n}$ is a topological isomorphism, and passing to the
inductive limit in $n$, we obtain a topological isomorphism
$$\bigoplus_{n=1}^{\infty} V_n \iso V_{\strict}.  \tag 7.1.9$$ 
Thus $V_{\strict}$ can be written not only as a strict inductive limit
of closed subspaces of $V$,
but in fact as a countable direct sum of closed subspaces of $V$.

If $W$ is a second Hausdorff convex $K$-vector space equipped with
a topological $\Gamma$-action,
and if $f: V \rightarrow W$ is a continuous $\Gamma$-equivariant map,
then $f$ induces continuous maps $f_n: V_n \rightarrow W_n$
for each $n\geq 1,$ and the map $V_{\strict} \rightarrow W_{\strict}$
induced by $f$ can be recovered as the direct sum
of the $f_n$.  

We now state some consequences of this analysis.

\proclaim{Lemma 7.1.10} If $V$ is a Hausdorff convex $K$-vector space
equipped with a strictly smooth action of $\Gamma,$
then the $\Gamma$-action on $V$ is continuous.
\endproclaim
\demo{Proof}
We apply lemma~3.1.1.  Condition~(i) of that result
holds for any smooth action of $\Gamma$, and by definition a strictly
smooth $\Gamma$-action on $V$ is a topological action, so that condition~(ii)
holds.  We must verify that condition~(iii) holds.  Let $H$ be a compact
open subgroup of $\Gamma$, and choose the cofinal sequence $\{H_n\}_{n \geq 1}$
of compact open subgroups appearing in the preceding discussion to be
normal open subgroups of $H$.  Each $V_n$ (as defined in that discussion)
is then invariant under the action of $H$, and the $H$-action on
$V_n$ factors through the finite quotient $H/H_n$ of $H$.   In particular,
the action of $H$ on $V_n$ is continuous, for each $n \geq 1$.  It follows
that the action of $H$ on $V \leftiso \bigoplus_{n \geq 1} V_n$ is
continuous, and so in particular the required condition~(iii) holds.
\qed\enddemo

\proclaim{Lemma 7.1.11} If $V$ is a Hausdorff
convex $K$-vector space equipped with a topological action of
$\Gamma,$ then the natural map $(V_{\strict})_{\strict} \rightarrow 
V_{\strict}$ is an isomorphism.  Equivalently, the $\Gamma$-representation
$V_{\strict}$ is strictly smooth.
\endproclaim
\demo{Proof}
This follows from the isomorphism~(7.1.9), together
with the evident isomorphism
$ \bigoplus_{m=1}^n V_m \iso (\bigoplus_{m=1}^{\infty} V_m)^{H_n}.$
\qed
\enddemo

\proclaim{Lemma 7.1.12} If $V$ is a Hausdorff convex $K$-vector space
equipped with a strictly smooth action of $\Gamma$,
and if $W$ is a $\Gamma$-invariant subspace of $V$,
then $W$ is closed in $V$ if and only if $W^H$ is closed in
$V^H$ for each compact open subgroup $H$ of $\Gamma$.
Furthermore, if these equivalent conditions hold,
then the $\Gamma$-actions on each of $W$ and $V/W$ are again strictly smooth,
and for each compact open subgroup $H$ of $\Gamma$,
the natural map $V^H/W^H \rightarrow (V/W)^H$ is a topological isomorphism.
\endproclaim
\demo{Proof}
Certainly, if $W$ is a closed subspace of $V$,
then $W^H = V^H \bigcap W$ is a closed subspace of $V^H$,
for each compact open subgroup $H$ of $\Gamma$.
Conversely, suppose that the latter is true.
Then for each $n \geq 1$
we can define $V_n$ and $W_n$ as in the discussion preceding
lemma~7.1.10,
and we see that $W_n$ is a closed subspace of $V_n$ for each value of $n$.
Thus $W \leftiso \bigoplus_{n = 1}^{\infty} W_n
\subset \bigoplus_{n=1}^{\infty} V_n \iso V$
is a closed subspace of $V$.  There is also an isomorphism
$V/W \leftiso \bigoplus_{n=1}^{\infty} V_n/W_n.$
(See proposition~8 and corollary~1 of
\cite{\TVS, p.~II.31}.)
These direct sum decompositions show that each of
$W$ and $V/W$ is again a strictly smooth representation of $\Gamma$.
Finally, if we fix a compact open subgroup $H$ of $\Gamma$,
and choose the sequence $\{H_n\}_{n\geq 1}$ so that $H_1 = H$,
then the same direct sum decompositions show that
$V^H/W^H \iso (V/W)^H$, as required.
\qed\enddemo

\proclaim{Lemma 7.1.13}
Suppose that $V$ and $W$ are a pair of Hausdorff convex
$K$-vector spaces, each
equipped with a strictly smooth action of $\Gamma$,
and that we are given
a continuous $\Gamma$-equivariant map $V \rightarrow W$.
This map 
is strict if and only if
the same is true of the induced maps $V^H \rightarrow W^H$,
for each compact open subgroup $H$ of $\Gamma$.
\endproclaim
\demo{Proof}
This result follows by an argument analogous to the one used
to prove lemma~7.1.12.
\qed\enddemo

The following lemma provides an analysis of
the ``strictly smooth dual'' of a strictly
smooth representation of $\Gamma$.

\proclaim{Lemma 7.1.14}  If $V$ is a Hausdorff convex $K$-vector
space equipped with a strictly smooth action of $\Gamma$,
let $V'_{b}$ 
be equipped with the
contragredient representation of $\Gamma$.

(i) For each compact open subgroup $H$ of $\Gamma$,
there is a natural isomorphism $(V'_b)^H \iso (V^H)'_b$.

(ii)
There is a natural $\Gamma$-equivariant isomorphism
$(V'_b)_{\strict} \iso \ilim{H} (V^H)'_b.$

(iii)
The natural map
$V \rightarrow (V'_b)'_b$ induces
a natural $\Gamma$-equivariant injective map
of strictly smooth $\Gamma$-representations
$V \rightarrow (((V'_b)_{\strict})'_b)_{\strict},$
which is a topological embedding if the same is true
of the first map.

(iv) If $V$ is of compact type, or is a nuclear Fr\'echet space,
then the map
$V \rightarrow (((V'_b)_{\strict})'_b)_{\strict}$
of~(iv) is an isomorphism.
\endproclaim
\demo{Proof}
Let $H$ be a compact open subgroup of $\Gamma$.  Lemma~7.1.6
provides a direct sum decomposition $V^H \bigoplus \ker \pi_H
\iso V.$ Taking strong duals yields an isomorphism
$V'_b \iso (V^H)'_b \bigoplus (\ker \pi_H)'_b.$
Thus to prove~(i), it suffices to show that $((\ker \pi_H)')^H = 0.$
If $v'$ is an $H$-invariant linear functional defined on $\ker \pi_H,$
then $v' \circ \pi_H = v',$ and so $v' = 0$, as required.  Thus~(i)
is proved.  Note that~(ii) is an immediate consequence of~(i).

To prove~(iii), we first construct the required map.
Note that the natural map
$(V'_b)_{\strict} \rightarrow V'_b$ induces a continuous map
$(V'_b)'_b \rightarrow ((V'_b)_{\strict})'_b$.  Composing
with the double duality map $V \rightarrow (V'_b)'_b$
yields a map $V \rightarrow ((V'_b)_{\strict})'_b,$
and hence (since $V_{\strict} \iso V$) a map
$$V \rightarrow (((V'_b)_{\strict})'_b)_{\strict},\tag 7.1.15$$
as required.

Parts~(i) and~(ii),
applied to first to $V$ and then to $(V'_b)_{\strict},$
yield an isomorphism
$\ilim{H} ((V_H)'_b)'_b \iso (((V'_b)_{\strict})'_b)_{\strict}.$
The map~(7.1.15) admits a more concrete description,
in terms of this isomorphism and the isomorphism $\ilim{H} V^H \iso V$:
it is the locally convex
inductive limit over the compact open subgroups $H$ of $\Gamma$
of the double duality maps
$$V^H \rightarrow ((V^H)'_b)'_b.\tag 7.1.16$$
If the double duality map $V \rightarrow (V'_b)'_b$ is an embedding,
then corollary~7.1.7 implies that the same is true for each of these
maps.  Since the strict inductive limit of embeddings is an embeddings,
this completes the proof of part~(iii).

To prove part~(iv), note that if $V$ is of compact type (respectively,
is a nuclear $K$-Fr\'echet space), then the same is true of
each of its closed subspaces $V^H$.  In particular, each of these
spaces is reflexive, and so each of the double duality maps~(7.1.16)
is a topological isomorphism.  Passing to the inductive limit in $H$,
we see that the same is true of the map~(7.1.15).
\qed\enddemo

Note that if $V$ is an admissible smooth representation of $\Gamma$,
equipped with its finest convex topology, then 
$V$ becomes a compact type convex $K$-vector space, equipped with
a smooth topological action of $\Gamma$, necessarily strictly
smooth, by corollary~7.1.4.  
The preceding lemma thus implies that
the strictly smooth dual $(V'_b)_{\strict}$ of $V$
is then isomorphic to the smooth dual to $V$,
equipped with its finest convex topology. 

We close this subsection by introducing a notion that we
will require below.

\proclaim{Definition 7.1.17}
Let $A$ be a Fr\'echet-Stein $K$-algebra, in the sense of definition~1.2.10.
A coadmissible $(A,\Gamma)$-module
consists of a topological left $A$-module $M$, equipped with a topological
action of $\Gamma$ that commutes with the $A$-action on $M$,
such that:

(i) The $\Gamma$-action on $A$ is strictly smooth.

(ii) For each compact open subgroup $H$ of $\Gamma$,
the closed $A$-submodule $M^H$ of $M$ is a coadmissible $A$-module
in the sense of definition~1.2.8.
\endproclaim

In particular, since any coadmissible $A$-module is a
$K$-Fr\'echet space,
we see that $M \iso \ilim{H} M^H$ is a strict inductive limit
of Fr\'echet spaces.  (In fact, the discussion preceding
lemma~7.1.10 shows that it is even a countable direct sum
of Fr\'echet spaces.)

\proclaim{Lemma~7.1.18} If $A$ is a locally convex topological $K$-algebra,
and if $M$ is a topological left $A$-module equipped with a topological
action of $\Gamma$ that commutes with the $A$-action on $M$,
then $M_{\strict}$ is naturally a topological left $A$-module.
\endproclaim
\demo{Proof}
It is clear by functoriality of the formation of $M_{\strict}$
that the $A$-module structure on $M$ induces an $A$-module structure
on $M_{\strict}$.  We must show that the multiplication map
$A \times M_{\strict} \rightarrow M_{\strict}$ is continuous.

If we adopt the notation introduced preceding lemma~7.1.10,
then the isomorphism~(7.1.9) shows that it suffices to prove
that the multiplication map $A \times M_n \rightarrow M_n$ is
jointly continuous for each $n \geq 1$.  Since $M_n$ is an
$A$-invariant subspace of $M$, this follows from the
assumption that the multiplication map $A \times M \rightarrow M$
is jointly continuous (that is, that $M$ is a topological left
$A$-module).
\qed\enddemo

\section{7.2}
As in the preceding subsection, in this subsection we suppose that
$\Gamma$ is a locally compact topological group,
admitting a countable neighbourhood
basis of the identity consisting of open subgroups.
We also suppose that $G$ is a locally $L$-analytic group.
In this subsection we introduce the notions
of admissible continuous, admissible locally analytic,
and essentially admissible locally analytic
representations of the topological group $G \times \Gamma$. 

\proclaim{Definition~7.2.1} Let $V$ be a locally convex Hausdorff
$K$-vector space equipped with a topological action of $G\times \Gamma$.
We say that $V$ is an admissible continuous representation of 
$G\times \Gamma$ if:

(i) For each compact open subgroup $H$ of $\Gamma$ the closed subspace
$V^H$ of $H$-invariant vectors in $V$ is a Banach space, and
the $G$-action on $V^H$ is an admissible continuous representation of $G$.

(ii) The $\Gamma$-action on $V$ is strictly smooth, in the sense
of Definition~7.1.2.
(That is, the natural map $\ilim{H} V^H \rightarrow V$ is a topological
isomorphism, the locally convex inductive limit being taken over
all compact open subgroups of $\Gamma$.)
\endproclaim

It follows from proposition~6.2.5 that it suffices to
check condition~(i) of this definition for a cofinal collection
of compact open subgroups $H$ of $G$.  
The discussion preceding lemma~7.1.10 shows that any 
$V$ satisfying the conditions of Definition~7.2.1 is in fact isomorphic
to a countable direct sum of Banach spaces.

\proclaim{Proposition~7.2.2} 
(i) If $V$ is an admissible continuous representation of $G\times \Gamma$,
then the $G \times \Gamma$-action on $V$ is continuous.

(ii) If $V$ is an admissible continuous representation of $G\times \Gamma$,
and if $W$ is a closed $G\times \Gamma$-invariant subspace of $V$,
then the $G\times \Gamma$-action on $W$ is also an admissible continuous
representation.

(iii) Any
$G\times \Gamma$-equivariant morphism of admissible continuous representations
of $G\times \Gamma$ is strict, with closed image.

(iv)
The category of admissible
continuous representations of $G\times \Gamma$ is abelian.

(v)
If $H$ is a compact open subgroup of $\Gamma$, then the
functor that takes $V$ to the closed subspace of $H$-invariants $V^H$
is an exact functor from the category of admissible continuous representations
of $G\times \Gamma$ to the category of admissible continuous
$G$-representations.
\endproclaim
\demo{Proof}
If $V$ is an admissible continuous representation of $G\times \Gamma$,
then $V$ is in particular a direct sum of Banach spaces, and so is
barrelled.  Since $V$ is isomorphic to a locally convex
inductive limit of admissible
continuous representations of $G$, the $G$-action on $V$ is separately
continuous, and hence continuous.  On the other hand,
lemma~7.1.10 shows that the $\Gamma$-action on $V$ is continuous. 
Thus~(i) is proved.

Now suppose that we are given a
$G \times \Gamma$-invariant subspace $W$ of $V$.
If $H$ is a compact open subgroup of $\Gamma,$ then by assumption
$V^H$ is an admissible continuous representation of $G$.  Since
$W^H = V^H \bigcap W,$ we see that $W^H$ is a closed $G$-invariant
subspace of the admissible continuous $G$-representation $V^H,$
and so by proposition~6.2.5 is again an admissible continuous $G$-representation.
We are also assuming that $V$ is
a strictly smooth $\Gamma$-representation, and so it follows 
from lemma~7.1.12 that the same is true of $W$.  This proves~(ii).

Part~(iii) follows from lemmas~7.1.12 and~7.1.13 and
the fact that $G$-equivariant morphisms of continuous
admissible $G$-representations are strict with closed image,
by proposition~6.2.9.

In light of what we have already proved, and the fact
that the direct sum of two admissible continuous representations
of $G\times \Gamma$ is again such a representation,
to prove part~(iv)
we need only show that if $W$ is a closed subrepresentation
of an admissible continuous representation $V$ of $G \times \Gamma$,
then the quotient $V/W$ is again an admissible continuous
representation of $G \times \Gamma$.  Lemma~7.1.12 shows
that the representation of $\Gamma$ on this quotient is strictly
smooth.  Furthermore, the isomorphism $V^H/W^H \iso (V/W)^H,$
and the fact that by corollary~6.2.16 the category of admissible continuous representations
of $G$ is abelian, shows that $(V/W)^H$ is
an admissible representation of $G$, for each compact open subgroup
$H$ of $\Gamma$.  Thus $V/W$ is an admissible continuous
representation of $G \times \Gamma$, as required.

Part~(v) follows from
the fact that passing to $H$-invariants yields an exact functor
on the category of  smooth $\Gamma$-representations. 
\qed\enddemo

If $V$ is a Hausdorff convex $K$-vector space equipped with a topological action
of $G\times \Gamma$, then we will let $V_{\la}$ denote the
space of $G$-locally analytic vectors attached to $V$, as in definition~3.5.3.
The functoriality of the construction
of locally analytic vectors shows
that the $\Gamma$-action on $V$ lifts to a $\Gamma$-action on $V_{\la},$
and thus that $V_{\la}$ is equipped 
with a topological $G\times \Gamma$-action,
uniquely determined by the requirement that the natural continuous injection
$V_{\la} \rightarrow V$ be $G\times \Gamma$-equivariant.

\proclaim{Definition 7.2.3}
We say that a topological action of $G\times \Gamma$
on a Hausdorff convex $K$-vector space $V$ is
a locally analytic $G \times \Gamma$-representation if
$V$ is barrelled,
if the natural map $V_{\la} \rightarrow V$ is a bijection,
and if the $\Gamma$-action on $V$ is strictly smooth.
\endproclaim

If $V$ is equipped with a locally analytic $G \times \Gamma$-representation,
then the fact that the $G$-action on $V$ is locally
analytic, together with lemma~7.1.10, implies that the $G \times \Gamma$-action
on $V$ is continuous.
Note that if $V$ is of compact type, then corollary~7.1.4 implies
that in the preceding definition,
we may replace the condition that $\Gamma$ act via a strictly
smooth representation on $V$ with the apparently weaker condition
that the action of $G$ on $V$ be smooth.

\proclaim{Lemma 7.2.4}
If $V$ is Hausdorff convex $K$-vector space equipped with
a topological action of $G \times \Gamma$,
such that the $\Gamma$-action is strictly smooth,
then the following are equivalent:

(i) The $G \times \Gamma$-action on $V$ is locally analytic.

(ii) The $G$-action on $V$ is locally analytic in the sense
of definition~3.6.9.

(iii) For each compact open subgroup $H$ of $G$,
the $G$-action on $V^H$ is locally analytic.
\endproclaim
\demo{Proof}
Note that since the $\Gamma$-action on $V$ is assumed to
be strictly smooth, conditions~(i) and~(ii) are 
equivalent by definition.

Corollary~7.1.8 shows that if $V$ is equipped with
a strictly smooth $\Gamma$-action, then $V$ is barrelled
if and only if $V^H$ is barrelled for each compact open
subgroup $H$ of $G$.  It follows from lemma~3.6.14
that~(ii) implies~(iii).

Conversely, if~(iii) holds,
then $(V^H)_{\la} \rightarrow V^H$ is a continuous bijection,
for each compact open subgroup $H$ of $\Gamma$.   The sequence
of isomorphisms
$$V_{\la} \leftiso (\ilim{H} V^H)_{\la} \iso \ilim{H} (V^H)_{\la}$$
(the first following from the assumption that the $\Gamma$-action
on $V$ is strictly smooth, and the second from proposition~3.5.14)
then shows that the natural map $V_{\la} \rightarrow V$ is a continuous
bijection.  Again taking into account corollary~7.1.8, we find that~(iii)
implies~(ii).
\qed\enddemo

\proclaim{Proposition 7.2.5}
If $V$ is a Hausdorff convex $K$-vector space, equipped with
a topological $G\times \Gamma$-action, such that the $\Gamma$-action
is strictly smooth, then for any compact open subgroup $H$ of $\Gamma$
there is a natural isomorphism
$(V^H)_{\la} \iso (V_{\la})^H$.
\endproclaim
\demo{Proof}
Corollary~7.1.7 shows that $V^H$ is a topological direct summand of $V$.
Thus $(V^H)_{\la}$ is naturally identified
with a topological direct summand of $V_{\la}$,
which clearly must equal $(V_{\la})^H$.
\qed\enddemo

\proclaim{Corollary 7.2.6}
If $V$ is a Hausdorff convex $K$-vector space, equipped with
a topological $G \times \Gamma$-action, such that the
$\Gamma$-action is strictly smooth,
then $V_{\la}$ is a locally analytic $G \times \Gamma$-representation.
\endproclaim
\demo{Proof}
Propositions~3.5.5 and~7.2.5 show that for each compact open subgroup $H$ of $\Gamma$,
the space of $H$-invariants $(V_{\la})^H$
is a locally analytic $G$-representation.
By lemma~7.2.4, the corollary will
be proved if we show that the $\Gamma$-action on $V_{\la}$
is strictly smooth.

Let $\{H_n\}_{n \geq 1}$ be a cofinal sequence of compact open
subgroups of $G$.
By assumption there is a natural isomorphism $\ilim{n \geq 1} V^{H_n} \iso V$.
Passing to locally analytic vectors, and again appealing to
proposition~7.2.5,
we obtain isomorphisms
$$(V_{\la})_{\strict} \iso
\ilim{n\geq 1}(V_{\la})^{H_n} \leftiso \ilim{n \geq 1} (V^{H_n})_{\la}
\iso (\ilim{n \geq 1} V^{H_n})_{\la}
\iso V_{\la}$$
(the second-to-last isomorphism following from proposition~3.5.14,
and the last from the fact that $V$ is assumed to be a strictly
smooth $\Gamma$-representation),
and so the $\Gamma$-action on $V_{\la}$ is strictly smooth,
as claimed.
\qed\enddemo

We now extend the notions of admissible and essentially admissible representations
of $G$ 
to the context of representations of $G \times \Gamma$.  Whenever we speak of essentially
admissible locally analytic representations of $G$ or of $G \times \Gamma$,
we assume that the centre of $G$ is topologically finitely generated, so
that the notion of an essentially admissible locally analytic representation
of $G$ is in fact defined.

\proclaim{Definition~7.2.7}
Let $V$ be a locally convex Hausdorff
$K$-vector space equipped with a locally analytic representation
of $G\times \Gamma$.
We say that $V$ is an (essentially)
admissible locally analytic representation of 
$G\times \Gamma$ if
for each compact open subgroup $H$ of $\Gamma$ the closed subspace
$V^H$ of $H$-invariant vectors in $V$ is an (essentially)
admissible locally analytic
representation of $G$.
\endproclaim

Lemma~7.1.4 shows that if $V$ is equipped with an (essentially)
admissible locally analytic representation of $G \times \Gamma$,
then $V$ is of compact type.

\proclaim{Proposition 7.2.8}
(i) If $V$ is an (essentially) admissible locally analytic representation
of $G\times \Gamma$, then the $G \times \Gamma$-action on $V$
is continuous.

(ii) If $V$ is an (essentially)
admissible locally analytic representation of $G\times \Gamma$, and if
$W$ is a closed $G\times \Gamma$-invariant subspace of $V$, then
the $G\times \Gamma$-action on $W$ is also an 
(essentially) admissible locally analytic
representation.

(iii) Any
$G\times \Gamma$-equivariant morphism of (essentially)
admissible locally analytic representations
of $G\times \Gamma$ is strict, with closed image.

(iv)
The category of (essentially) admissible locally
analytic representations of $G\times \Gamma$ is abelian.

(v)
If $H$ is a compact open subgroup of $\Gamma$, then the
functor that takes $V$ to the closed subspace of $H$-invariants $V^H$
is an exact functor from the category of (essentially)
admissible locally analytic representations
of $G\times \Gamma$ to the category of (essentially)
admissible locally analytic representations of $G$.
\endproclaim
\demo{Proof}
This is proved in an analogous manner to 
proposition~7.2.2, once one
recalls corollary~6.1.23 and proposition~6.4.11.
\qed\enddemo

We can give an alternative characterisation of (essentially)
admissible locally analytic representations of $G \times \Gamma$
in terms of their strictly smooth duals.

Suppose that $V$ is a convex $K$-vector space of compact type, equipped
with a locally analytic representation of $G \times \Gamma$.
The strictly smooth
dual $(V'_b)_{\strict}$ of $V$ is then a topological left $\DLa(H,K)$-module
(by corollary~6.1.22 together with lemma~7.1.18),
for any compact open subgroup $H$ of $G$, equipped with a commuting
strictly smooth action of $\Gamma$ (the contragredient action).

\proclaim{Lemma 7.2.9} Suppose that the centre $Z$ of $G$ is
topologically finitely generated.  If $V$ is as in the
preceding discussion, then the $Z$-representation on $V$ satisfies
the equivalent conditions of proposition~6.4.7 if and only if
the same is true for the closed $Z$-subrepresentation $V^H$
of $V$, for each compact open subgroup $H$ of $G$.
\endproclaim
\demo{Proof}
If $V$ is the union of an increasing sequence of $Z$-invariant
$BH$-subspaces, then so is any closed subspace of $V$,
and so in particular, so is each $V^H$.
Conversely, if each $V^H$ is the union of an increasing sequence
of $Z$-invariant $BH$-subspaces, then letting $H$ run through
a cofinal sequence $\{H_n\}_{n \geq 1}$ of compact open subgroups
of $G$, and writing $V = \bigcup V^{H_n}$, we find that the
same is true of $V$.
\qed\enddemo

If the centre $Z$ of $G$ is topologically finitely generated,
and if $V$ satisfies the conditions of proposition~6.4.7,
then lemmas~7.1.18 and~7.2.9, together with 
that proposition, 
show that the strictly smooth dual $(V'_b)_{\strict}$ of $V$
is a topological left
$\An(\hat{Z},K)\cotimes_K \DLa(H,K)$-module,
for any compact open subgroup $H$ of $G$,
equipped with a commuting strictly smooth action of $\Gamma$.

We now state the criterion on $(V'_b)_{\strict}$ for $V$
to be an (essentially) admissible representation of $G \times \Gamma$.

\proclaim{Proposition 7.2.10} 
If $V$ is a convex $K$-vector space of compact type, equipped with
a locally analytic $G \times \Gamma$-action (respectively,
a locally analytic $G \times \Gamma$-action, whose restriction
to the centre $Z$ of $G$ satisfies the conditions of 
proposition~6.4.7),
then $V$ is an admissible (respectively, essentially admissible)
locally analytic $G \times \Gamma$-representation
if and only if $(V'_b)_{\strict}$ is a coadmissible
$(\DLa(H,K), \Gamma)$-module (respectively,
a coadmissible $(\An(\hat{Z},K)\cotimes_K \DLa(H,K),\Gamma)$-module),
for one, or equivalently every,
compact open subgroup $H$ of $G$, in the sense of definition~7.1.17.
\endproclaim
\demo{Proof}
Let $H$ be a compact open subgroup of $\Gamma$.
Lemma~7.1.14 shows that the space of invariants
$((V'_b)_{\strict})^H$ ($= (V'_b)^H$)
is naturally isomorphic to $(V^H)'_b$.  Thus, by the definition
of admissibility (respectively, essential admissibility)
of a locally analytic $G$-representation, the
space $((V'_b)_{\strict})^H$
is a coadmissible $\DLa(H,K)$-module (respectively,
a coadmissible $\An(\hat{Z},K)\cotimes_K \DLa(H,K)$-module),
equipped with its canonical topology,
if and only if $V^H$ is an admissible (respectively, essentially
admissible) locally analytic representation of $G$.
Taking into account the definition of coadmissibility of
a topological $(\DLa(H,K),\Gamma)$-module (respectively,
$(\An(\hat{Z},K)\cotimes_K \DLa(H,K),\Gamma)$-module),
the proposition follows.
\qed\enddemo

\proclaim{Proposition~7.2.11} If $V$ is an admissible continuous
representation of $G\times \Gamma$, then $V_{\la}$ is an admissible
locally analytic representation of $G$.
\endproclaim
\demo{Proof}
It follows from corollary~7.2.6 that $V_{\la}$ is a locally analytic
representation of $G \times \Gamma$.  It remains to show
that for each compact open subgroup $H$ of $G$,
the space $(V_{\la})^H$ is an admissible locally analytic 
representation of $G$.  Proposition~7.2.5 shows that this
space is naturally isomorphic to $(V^H)_{\la}$,
which proposition~6.2.4 shows
is an admissible locally analytic representation of $G$.
\qed\enddemo

\proclaim{Proposition~7.2.12} If $V$ is an admissible 
locally analytic representation of $G$,
and if $W$ is a finite dimensional locally analytic representation
of $G$, then the tensor product $V \otimes_K W$, regarded as a
$G\times \Gamma$-representation via the diagonal action of $G$
and the action of $\Gamma$ on the left factor, is an admissible
locally analytic representation
of $G \times \Gamma$.
\endproclaim
\demo{Proof}
As a $\Gamma$-representation, the tensor product $V \otimes_K W$
is isomorphic to a finite direct sum of copies of $V$,
and so is strictly smooth, since $V$ is assumed to be so.
If $H$ is a compact open subgroup of $\Gamma,$ then
the natural map
$V^H \otimes_K W \rightarrow (V \otimes_K W)^H$ is evidently
a $G$-equivariant topological isomorphism.  Since $V^H$ is
assumed to be an admissible (respectively, essentially admissible)
locally analytic representation of $G$, we see by proposition~6.1.5
that the same is true
of $V^H \otimes_K W$.  Thus $V\otimes_K W$ satisfies the conditions
of Definition~7.2.7.
\qed\enddemo

Now suppose that $\G$ is a linear reductive algebraic group defined
over $L$, and that $G$ is an open subgroup of
$\G(L)$.  We can extend the notion of a locally
algebraic representation of $G$, discussed in subsection~4.2,
to representations of $G\times \Gamma$.  As in that subsection we let
$\Cal R$ denote the semi-simple abelian category of finite dimensional
representations of $\G$.  The functor from $\Cal R$ to the
category of $G$-representations is fully faithful.

Fix an object $W$ of $\Cal R$.
If $V$ is a representation of $G \times \Gamma$ on an abstract $K$-vector
space, then we may in particular
regard $V$ is a representation of $G$, and so (following 
proposition-definition~4.2.2) construct the space $V_{W-\lalg}$
of locally $W$-algebraic vectors of $V$.
Since the formation of locally $W$-algebraic vectors is functorial in $V$,
we see that $V_{W-\lalg}$ is a $G \times \Gamma$-invariant subspace
of $V$.

Similarly, we may define the space $V_{\lalg}$ of all local algebraic
vectors of $V$, as in definition~4.2.6.  This is again a
$G \times \Gamma$-invariant subspace of $V$.  If $\hat{G}$ denotes
a set of isomorphism class representatives of irreducible objects
of $\Cal R$, then  corollary~4.2.7 implies that the natural map
$\bigoplus_{W \in \hat{\G}} V_{W-\lalg} \rightarrow V_{\lalg}$
is an isomorphism.

\proclaim{Proposition~7.2.13} If $V$ is an admissible locally
analytic representation of $G\times \Gamma$, and if $W$ is an
object of $\Cal R$, then $V_{W-\lalg}$ is a closed $G\times \Gamma$-invariant
subspace of~$V$.
\endproclaim
\demo{Proof}
It is evident that the natural map $(V^H)_{W-\lalg} \rightarrow
(V_{W-\lalg})^H$ is an isomorphism.  
Thus proposition~4.2.10
shows that $(V_{W-\lalg})^H$
is a closed subspace of $V^H$ for each compact open subgroup
$H$ of $\Gamma$.  The proposition follows from lemma~7.1.12.
\qed\enddemo

\proclaim{Proposition 7.2.14}
If $V$ is an admissible locally analytic
representation of $G\times \Gamma$ that is also locally algebraic, then
$V$ is equipped with its finest convex topology.
\endproclaim
\demo{Proof}
If $H$ is a compact open subgroup of $\Gamma$ then $V^H$ is an
admissible locally analytic representation of $G$ that is also
locally algebraic as a representation of $G$. It follows from
corollary~6.3.7 that $V^H$ is equipped with its
finest convex topology.  By assumption $V$ is isomorphic to
the locally convex inductive limit $\ilim{H} V^H$, as $H$
ranges over all compact open subgroups of $\Gamma$.  Thus $V$
is also equipped with its finest convex topology.
\qed\enddemo

\proclaim{Proposition~7.2.15} 
If $V$ is an admissible locally analytic
representation of $G\times \Gamma$ that is also locally $W$-algebraic
for some object $W$ of $\Cal R$, then
$V$ is isomorphic to a representation of the form
$U \otimes_B W,$ where $B = \End_{\G}(W)$ (a finite rank semi-simple
$K$-algebra), $U$ is an admissible
smooth representation of $G \times \Gamma$ on a $B$-module (equipped
with its finest convex topology), and the $G \times \Gamma$-action
on the tensor product is defined by the diagonal action of $G$,
and the action of $\Gamma$ on the first factor.
Conversely, any such tensor product is an admissible
locally analytic representation of $G \times \Gamma$.
\endproclaim
\demo{Proof}
Proposition~4.2.10 yields a $G\times \Gamma$-equivariant topological
isomorphism $U \otimes_B W \iso V.$
Let $\check{W}$ denote the contragredient representation to $W$.
Proposition~7.2.12 shows that $\Hom(W, V)$ (which is naturally isomorphic 
to the tensor product $U \otimes_K \check{W}$) is an
admissible locally analytic representation of $G \times \Gamma$,
and hence (by proposition~7.2.8(ii)) the same is true of its closed subspace
$U = \Hom_{\lie{g}}(W,V),$ equipped with its natural $G\times \Gamma$
action.  Thus there is a topological isomorphism
$\ilim{H} U^H \iso U$ (the locally convex inductive
limit being taken over all compact open subgroups $H$ of $\Gamma$).
By corollary~6.3.3, each of the spaces $U^H$ is an
admissible smooth representation of $G$, equipped with its finest
convex topology.  Thus $U$ is an admissible smooth representation
of $G\times \Gamma,$ and is equipped with its finest convex topology,
as claimed.

Now suppose given an object $W$ of $\Cal R$, with algebra of endomorphisms
$B$, and a $B$-linear admissible smooth representation $U$
of $G \times \Gamma$ (equipped with its finest convex topology).
The remark preceding the statement of proposition~4.2.4 shows that
$U\otimes_B W$ is a locally $W$-algebraic representation of $G\times \Gamma$.
If $H$ is any compact open subgroup of $\Gamma,$ then the natural map
$U^H \otimes_B W \rightarrow (U \otimes_B W)^H$ is obviously an isomorphism.
The source of this map is the tensor product of an admissible smooth
representation of $G$ and an object of $\Cal R,$ and hence by
proposition~6.3.9 is
an admissible locally analytic representation of $G$.
Since the natural map $\ilim{H} U^H\otimes_B W \rightarrow U\otimes_B W$
is obviously an isomorphism, we see that $U\otimes_B W$ is an
admissible locally analytic isomorphism of $G \times \Gamma$, as claimed.
\qed\enddemo

\enddocument